\documentclass[12pt,spanish]{book}
\usepackage[titletoc]{appendix}
\usepackage[spanish,activeacute]{babel}
\usepackage[utf8x]{inputenc}
\usepackage{amssymb,amsfonts,amsthm,amsmath,calc,graphicx,hyperref,imakeidx,multicol,subfigure,accents}
\makeindex
\newtheorem{teo}{Teorema}[chapter]
\newtheorem{lema}[teo]{Lema}
\newtheorem{prop}[teo]{Proposición}
\newtheorem{cor}[teo]{Corolario}
\newtheorem{df}[teo]{Definición}
\newtheorem{ej}[teo]{Ejemplo}
\newtheorem{obs}[teo]{Observación}
\spanishdecimal{.}
\title{\eltitulo}
\author{\elnombre}
\date{CD. UNIVERSITARIA, CDMX., \lafecha}

\newcommand{\titulo}[1]{\def\eltitulo{#1}}
\newcommand{\carrera}[1]{\def\lacarrera{#1}}
\newcommand{\nombre}[1]{\def\elnombre{#1}}
\newcommand{\director}[1]{\def\eldirector{#1}}
\newcommand{\fecha}[1]{\def\lafecha{#1}}
\newcommand{\comite}[1]{\def\elcomite{#1}}
\titulo{\textbf{\large{COMPETENCIA Y DEPREDACIÓN: EL ARRIBO DE UNA ESPECIE EXÓTICA AL MEDIO}}}
\nombre{\uppercase{\textbf{ALVARO REYES GARCÍA}}}
\carrera{\textbf{DOCTOR EN CIENCIAS MATEMÁTICAS}}
\director{\uppercase{\textbf{DR. MANUEL JESÚS FALCONI MAGAÑA 
\\FACULTAD DE CIENCIAS, UNAM}}}
\comite{\uppercase{\textbf{DR. ABDÓN EDDY CHOQUE RIVERO \\ INSTITUTO DE FISICA Y MATEMÁTICAS, UMSNH
\\ $ $ \\ DR. RAMÓN GABRIEL PLAZA VILLEGAS \\ IIMAS, UNAM}}}
\fecha{8 DE ABRIL DE 2022}
\begin{document}
\thispagestyle{empty}
\begin{minipage}[c][15cm][s]{15cm}
\begin{center}
\includegraphics[height=2.6cm]{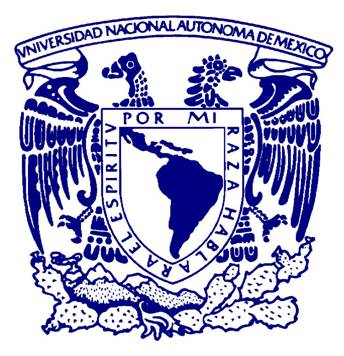}                                              \vspace{1cm}

{\large \scshape \textbf{UNIVERSIDAD NACIONAL AUTÓNOMA DE MÉXICO}}                       \vspace{.9cm}

{\scshape \textbf{POSGRADO EN CIENCIAS MATEMÁTICAS}}

{\scshape \textbf{FACULTAD DE CIENCIAS}}                                                 \vspace{.9cm}

\eltitulo                                                                               \vspace{1.4cm}

MODALIDAD DE GRADUACIÓN:

\textbf{TESIS DOCTORAL}                                                                  \vspace{.9cm}

QUE PARA OPTAR POR EL GRADO DE:

\lacarrera                                                                               \vspace{.9cm}

PRESENTA:

\elnombre                                                                              \vspace{1.4cm}

DIRECTOR DE TESIS:

\eldirector                                                                              \vspace{.9cm}

MIEMBROS DEL COMITÉ TUTOR:

\elcomite                                                                              \vspace{1.4cm}

CD. UNIVERSITARIA, CDMX, A \lafecha
\end{center}
\end{minipage}
$ $ \newpage
\thispagestyle{empty}
$ $ \newpage
\frontmatter
\maketitle
\newpage
\thispagestyle{empty}
\chapter*{}
\thispagestyle{empty}
\begin{flushright}
\emph{A Manuel Andrade,\\$ $\\el conocimiento es el pequeño\\Manuel que todos llevamos dentro.}
\end{flushright}
\tableofcontents
\chapter{Agradecimientos}
El agradecimiento más profundo es hacia el Dr. Manuel Falconi. Fue él quién hizo que todo esto fuera posible. El Dr. Falconi dirigió este proyecto de manera impecable, y es una persona que se ha ganado toda mi admiración, mi cariño y mi respeto. No cambiaría la oportunidad que tuve de trabajar con él por nada.

A Oriana Paola Neri González le extiendo todo mi cariño y toda mi gratitud por haberme acompañado en las largas horas que estuve sentado frente a la computadora.

Estoy agradecido con la Coordinación del Posgrado en Ciencias Matemáticas por haberme dado su confianza y su apoyo a lo largo de este proceso: Dra. Lourdes Esteva, Dra. Silvia Ruíz, María Inés León, María Teresa Martínez y Lucía Hernández.

A los miembros del comité tutor el Dr. Abdón Eddy Choqué Rivero y el Dr. Ramón Gabriel Plaza Villegas, así como también a mis sinodales el Dr. Pedro Miramontes, el Dr. Gamaliel Blé, el Dr. Ernesto Rosales, la Dra. Lourdes Esteva y el Dr. Marco Herrera les agradezco todos sus valiosos consejos y comentarios que me compartieron.

A mis todos mis profesores les agradezco por el gran trabajo que hicieron. En particular, a Miguel Martínez Huerta, a quien tuve el honor de conocer cuando cursaba la preparatoria, y que me enseñó la importancia del trabajo, de la ética y del profesionalismo, al Dr. Carlos Islas Moreno, a quien conocí en el verano de 2005 en el curso de Álgebra Superior I y que jamás me ha negado su apoyo, al Mtro. Ángel Manuel Carrillo Hoyo de quien fui alumno en los cursos de Análisis Matemático, y que con exposiciones muy interesantes y muy amenas, logró despertar en mí un interés genuino, y también al Dr. Francisco Marcos López García, quien me asesoró entre los años 2009 y 2011 para obtener el título de Matemático con la tesis \textit{El Problema de los Momentos}. 

A mi familia la tengo muy presente en mi pensamiento y en mi cotidianeidad: mi madre y a mi padre, mis hermanos Alejandro y Aarón, mi abuela, mis tíos Fátima, Conka, María, Guadalupe, Lourdes, Beatriz, Martín, Pablo, Carlos, Gustavo, Manuel y Héctor, mis primos Carlitos, Juan, Andrea, Natalia, Gus, Emi, Esteban, Cristobal, Pablito, Diego, Maya, Isabel y Oliver. 

También estoy muy agradecido con Jorge Alejandro, a quien conozco desde la infancia y que estuvo apoyándome en los momentos difíciles que se vivieron durante en el año 2020, con Raybel que es un gran matemático y una gran persona que trabaja dura y decididamente para cumplir todos sus objetivos, y con Pedro que me acompañaba a hacer ejercicio y siempre me hacía reír. Además existe un agradecimiento especial hacia Norma Escobar Prado, Marie José López Olivares y Tanya Arenas Reséndiz por la amistad, la camaradería y el cariño. Si le dedicara una línea a cada uno de mis amigos, demoraríamos demasiado. Sin embargo estoy seguro de que ellos saben que los quiero y les agradezco haber estado ahí para mí.

Le agradezco también el cariño y el apoyo brindado a cada uno de los integrantes de la familia Berumen Glinz.

Finalmente, cada uno de mis colegas del Equipo Joven estuvieron cerca durante la elaboración de la tesis: Manuel Falconi, Luis Osorio, Gamaliel Blé, Luis Miguel Valenzuela, Víctor Castellanos, María Fernanda Jiménez, Nashiely Juanita López, Yrina Vera, Rusby Contreras, Néstor Anaya y Erick Darinel.

Este trabajo fue apoyado parcialmente por CONACYT que me concedió la oportunidad de estar dentro de su programa de becarios con número de CVU 384637, y por la UNAM a través del proyecto PAPIIT IN116018.
\chapter{Prólogo}
Para comprender el mundo que nos rodea, los seres humanos hemos creado herramientas e implementado su uso. Se sabe que Anton van Leeuwenhoek usaba microscopios ópticos en su investigación \cite{Boutibonnes}, Galileo hizo una profunda exploración de los cuerpos celestes gracias al telescopio, por nombrar algunos ejemplos. En ambos casos, se trata de personajes cuyos trabajos fueron hechos hace más de trescientos años. En la actualidad se cuenta con submarinos, aviones, cámaras digitales, computadoras, equipo clínico, barómetros, GPS y muchos otros dispositivos que hacen que la investigación sea mucho más profunda y precisa en el día a día. Aunque podría parecer que las matemáticas no quedan dentro de este conjunto (de herramientas), la literatura reciente sugiere que a través de modelos matemáticos, uno puede describir y estudiar hechos y fenómenos que se presentan en nuestro entorno (ver por ejemplo \cite{peral}). 

En relación con la biología, la modelación matemática se ha vuelto una herramienta imprescindible en las décadas más recientes. Entre las disciplinas que han establecido conexión con la biomatemática (o biología matemática), se encuentran la ecología, la epidemiología, el crecimiento de tumores y las neurociencias. En esta tesis, se presentan modelos matemáticos que fueron creados para estudiar la depredación y la competencia en ecología, aunque no se descarta la posibilidad de incorporarlos posteriormente a otra rama.

Los humanos trabajamos con la ciencia desde tiempos remotos. A continuación se mencionan algunos aspectos que hacen notar la tendencia que muestra el ser humano a razonar científicamente, y lo cotidiano de este raciocinio que muchas veces es llevado a cabo sin siquiera darse cuenta de ello.

Se puede observar que los humanos desde muy temprana edad exploramos el medio y los objetos que nos rodean. Los niños pequeños se llevan objetos como canicas, plásticos o taparroscas de refrescos a la nariz y a la boca. ¿Por qué ocurre esto? En muchos casos, el infante ni siquiera sabe caminar, sin embargo ya está explorando el entorno en el que se encuentra inmerso. 

De forma similar, un adulto hace una exploración cuando quiere hacer una compra. Por ejemplo, si desea adquirir una nueva computadora, él se pregunta más de una vez si el equipo $A$ es el que satisface sus necesidades laborales y de ocio, o si le resulta más conveniente llevarse el equipo $B$. Entonces él hace una o varias consultas relacionadas a las computadoras $A$ y $B$ con colegas en su oficina, busca en \textit{Google} los equipos $A$ y $B$ para conocer sus características, hace una evaluación y con base en eso toma una decisión.

Existen aún más ejemplos cotidianos que muestran como las personas exploran opciones y obtienen resultados: ¿Qué ruta tomaré para llegar de la casa al trabajo?, ¿Qué platillo ordenaré en mi próxima visita a cierto restorán?, ¿Qué candidato o candidatos tendrán mi voto en la próxima elección federal?, y los ejemplos pueden seguir brotando; algunas de las preguntas anteriores pueden parecer fáciles de contestar. Sin embargo, ¿qué pasa si mi ruta favorita al trabajo está bloqueada por obras públicas? ¿O se encuentra agotado mi platillo favorito en el restorán? Entonces me veré obligado a tomar una nueva decisión, y tal vez deba preguntarme cuál es la alternativa que me genera mayor bienestar, o pedir que algún tercero me haga sugerencias al respecto.

En resumen, desde muy pequeños y de manera cotidiana, los seres humanos exploramos e investigamos, inclusive en las situaciones más simples y sencillas. En ocasiones, lo hacemos de forma natural y ni siquiera lo notamos. Al menos no se escucha a nadie en el supermercado diciendo: ``El cereal $A$ tiene contenido calórico superior al del cereal $B$, y además la caja de $A$ ofrece $650$gr. por \$$42$, mientras que $B$ cuesta \$$38$ por $600$gr. Con base en lo anterior, concluimos que para el cereal $A$ estaría pagando por cada $100$ gramos a \$$6.46$, mientras que para el cereal $B$ pagaría \$$6.33$ cada centena de gramos. Por lo tanto $A$ es más caro''. Quizás un peculiar discurso como el anterior no solamente parezca ridículo, sino que incluso se podría pensar que la persona que lo emitió pudiera llegar a tener un trastorno obsesivo. Sin embargo, una conversación común en una cantina o bar entre dos amigos muchas veces se centra en la efectividad o eficiencia de equipos de fútbol profesional. Podemos imaginar que uno de los amigos le dice al otro que el América será campeón de la Liga Mexicana de Fútbol y le recita (¡de memoria!) el nombre de todos los refuerzos que han llegado al América, junto con sus edades, estaturas, posiciones en el campo, equipos de procedencia y promedio de goles por partido. La otra persona le responde cuáles son los refuerzos de las Chivas en esta temporada y le dice a su amigo toda la información análoga concerniente a cada jugador de las Chivas. Evidentemente, aquí comienza el debate, y este puede prolongarse toda la noche y todo el día que sigue. Los contendientes (los dos amigos) justifican la razón por la cual su equipo saldrá campeón esta temporada. El primero comenta: ``Este año viene Nicolás Castillo al América. Este jugador anotó 26 goles en los 45 partidos que jugó con los Pumas. Eso quiere decir que tiene un promedio de más de un gol cada dos juegos. Por lo tanto veremos al menos un gol suyo cada dos jornadas''. Su amigo responde, ``pues las Chivas tienen en su delantera a Alexis Vega que es un hombre gol: la temporada pasada anotó 6 goles en 19 juegos. Su promedio es muy cercano al de un gol cada tres juegos, con la diferencia de que este jugador tiene apenas 21 años''. Esto provoca la réplica, ``es demasiado joven para ser factor en su equipo. En cambio, Castillo ya ganó la Copa América 2016 con la selección de Chile, y por lo tanto es un refuerzo más experimentado''. Algunos de los argumentos que los dos amigos dan, pueden parecer increíblemente parecidos a los que se mostraron con el sujeto del supermercado. ¿Por qué nos parece inusual hacer números respecto al contenido o peso de una caja de cereal, pero nos resulta perfectamente normal ver a dos hombres haciendo números para protagonizar una discusión acerca de dos equipos de fútbol? La razón es simple, es así como nos hemos acostumbrado a verlo. Pero más importante aún, es notar que en ambos casos se está haciendo uso de la capacidad cognitiva del ser humano para poder realizar una exploración sobre un objeto de estudio (en el primer caso una caja de cereal, y en el segundo un equipo de fútbol y sus refuerzos) para poder arrojar conclusiones.

De esta forma, los seres humanos constantemente estamos dando puntos de vista y sustentando con argumentos estas posturas. La discusión y el debate son practicados casi a diario, e inevitablemente se recurre a las matemáticas para decidir si algún objeto es más grande, más caro o (simple y llanamente) mejor. 

Como se acaba de mostrar, las matemáticas son invocadas en incontables ocasiones para resolver problemas. En esta tesis se muestra una aplicación de las matemáticas a la ecología. Dicha aplicación es la modelación a través de los sistemas de ecuaciones diferenciales ordinarias acopladas. 

Con respecto a la modelación matemática, un suceso previo a la redacción de este trabajo se presenta a continuación: en agosto de 2017, coincidí con el Dr. José Geiser Villavicencio Pulido (ver por ejemplo \cite{Olmos-Liceaga2017DiffusionModel, Villavicencio-Pulido2018CalculusPractices}) en la conferencia Solabima, en la ciudad de Cusco, Perú. En su presentación, el Dr. Villavicencio mencionó que un modelo matemático es como una caricatura, pues su objetivo es dar una descripción de algún fenómeno o evento (biológico, físico, económico, demográfico, etc.), usando solamente una cantidad limitada de información. En conclusión, un buen modelo (al igual que una buena caricatura) consigue acercarse suficientemente a la realidad sin tener que recurrir a todas las características del objeto que se está estudiando. 

Como se dijo anteriormente, los modelos matemáticos que aquí se trabajan son aplicables a ecología. Una obra de arte relacionada con temas de ecología, es presentada en el trabajo de Jorge Soberón Mainero \cite[Capítulo III]{SoberonMainero2000}. Aquí se habla de un gran mural del siglo XVI que se encuentra en el templo agustino en Malinalco. En este mural se describen muchas de las interacciones que existen entre flora y fauna silvestre: un bosque en el que algunos conejos están buscando alimento en el verde de las plantas, y unos colibríes se encuentran al acecho de algún néctar que extraer. Esta representación es un ejemplo muy claro de como las especies dependen unas de otras para poder coexistir. Cuando este mural es pintado hace aproximadamente cuatrocientos cincuenta años, muy probablemente ya se sabía acerca de las interacciones (como la depredación, la competencia, el mutualismo, etc.) que propician la biodiversidad.

La parte teórica de este trabajo incluye resultados clásicos y resultados modernos de la teoría de ecuaciones diferenciales. Estos son utilizados para analizar los modelos aquí propuestos. Además, se usan algunas técnicas propias basadas en geometría y en cálculo diferencial y que son parte de la investigación doctoral. También, algunos ejemplos de bifurcaciones y un sistema con dinámica caótica son estudiados.

La redacción de esta tesis ha permitido reafirmar que la matemática, más allá de ser un medio que nos permite resolver problemas y responder una gran variedad de preguntas, es también la más pura herramienta de análisis descriptivo o cualitativo. En efecto, a través de la modelación matemática, el investigador presenta conclusiones, desarrolla conjeturas y se plantea preguntas. Las personas que se dedican a estudiar modelos matemáticos, tienen al frente suyo una fuente inagotable de retos. Cada uno de estos retos tiene detrás una o varias preguntas que lo impulsan, y una vez completado el reto, este último da lugar a nuevas interrogantes que motivan un reto nuevo. Con base en lo anterior, se puede decir que el trabajo científico tiene una estructura cíclica.
\vspace{2in}
\begin{center}
Alvaro Reyes García
\end{center}
\mainmatter
\chapter{Introducción}
El mundo de los seres vivos es simplemente fascinante. Las diversas formas de vida en las que la naturaleza\index{naturaleza} se manifiesta, se pueden ver en un bosque, en una selva, en un pastizal o incluso en el fondo del océano. Es posible encontrar criaturas minúsculas como los insectos, y seres gigantescos como las ballenas. Por otra parte, están también las plantas, que proporcionan el oxígeno que nos permite a los seres aerobios respirar, y sin el cual pereceríamos. Las plantas, en muchos casos sirven de alimento a los animales herbívoros como los conejos o los venados. A su vez, estos últimos se convierten eventualmente en la presa\index{presa} de algún animal carnívoro, como por ejemplo un tigre o una hiena. 

Hechos como los anteriores dan lugar a una intrincada red de relaciones, pues las comunidades que cohabitan en un ecosistema\index{ecosistema} se mantienen en constante interacción. Un ave que desea alimentar a sus polluelos busca algunos insectos, la leona va de cacería para poder dar de comer a sus crías. Se pueden encontrar ejemplos como estos en la naturaleza\index{naturaleza}, y a partir de ello notar la importante contribución que están haciendo las presas\index{presa} a los depredadores\index{depredador}: al igual que el resto de los seres vivos, los depredadores\index{depredador} deben alimentarse y proveer alimento a sus crías (según sea el caso) para no morir. Por lo tanto, las presas\index{presa} constituyen la fuente de alimentación de los depredadores\index{depredador}. Sin embargo, la depredación elevada podría traer consecuencias en el comportamiento social de las especies, como se muestra en \cite{Ioannou2017High-predationFish}. Una tasa de captura alta se traduce en una mortalidad alta para las presas\index{presa}, y se puede suponer que esto podría llevar a las presas\index{presa} a extinguirse\index{extinción}. De manera similar, una tasa de natalidad\index{natalidad} baja por parte de las presas\index{presa} puede ser un factor que provoque su desaparición. Si la presas\index{presa} se extinguen\index{extinción} o emigran, resulta interesante preguntarse si los depredadores\index{depredador} también desaparecerán o no. Lo anterior depende de la capacidad de los depredadores\index{depredador} para sustituir su fuente alimentaria.

Para estudiar los efectos de estas interacciones, diversos modelos matemáticos han sido propuestos y analizados. La modelación matemática ha sido de gran utilidad para poder entender el comportamiento de las curvas poblacionales de las especies. En algunos casos, también ha contribuido a incrementar las interrogantes o a agudizar la curiosidad del observador.

El trasfondo histórico y la motivación del trabajo se dan a conocer a continuación.

Algunas de las primeras contribuciones que fueron hechas en relación con los modelos poblacionales datan de hace aproximadamente 800 años: en 1202 Leonardo de Pisa (Fibonacci\index{Fibonacci}) usó su famosa sucesión, que representa el crecimiento de una población de conejos. También está el trabajo del economista anglicano Thomas Robert Malthus\index{Malthus}, que en 1798 concluye que el crecimiento de la curva poblacional de la humanidad está dado por una sucesión geométrica, y que será inevitable la extinción\index{extinción} de los seres humanos dado que los alimentos tienen un crecimiento de tipo aritmético, y así los humanos enfrentaríamos una eventual escasez de víveres.  

Más de un siglo pasó desde Malthus\index{Malthus}, para que aparecieran los primeros modelos de depredación que hoy se conocen y que inspiraron muchos modelos con los que se trabaja en la actualidad. Los ecólogos propusieron dar un enfoque matemático a los problemas que se encontraban estudiando y así se desarrolló un vínculo entre los sistemas dinámicos y la ecología. Entre los años 1925 y 1926, Alfred Lotka \cite{Lotka1925} -- que nació en Ucrania --, y Vito Volterra \cite{Volterra1926} -- de nacionalidad italiana -- presentaron casi simultáneamente su trabajo, que resultó ser uno de los primeros modelos de depredación en el mundo de la biología matemática. En esos momentos, la llegada de un modelo matemático como el de Lotka y Volterra fue innovador y prometedor para futuros estudios en la disciplina, pues el modelo permitía predecir, estudiar y especular acerca del comportamiento de las poblaciones de las especies de un determinado lugar a través de un sistema de ecuaciones diferenciales acopladas. Esta gran contribución dio lugar a diversas perspectivas en la modelación matemática, como se menciona en \cite{Berryman2011}. Algunos trabajos recientes incorporan el modelo de Lotka-Volterra. Las interacciones de la fauna neozelandesa son analizadas a través de modelos de tipo Lotka-Volterra en \cite{Korobeinikov1999}, para poder así concluir que los kiwis enfrentarían una eventual extinción\index{extinción}. En \cite{Li2000} se proveen condiciones para la existencia de soluciones positivas en un sistema de Lotka-Volterra con retardo. Algoritmos computacionales de programación mixta (binaria-lineal) han sido implementados para poder encontrar los puntos de equilibrio de un sistema Lotka-Volterra, así como sus propiedades dinámicas. Otra perspectiva interesante es la incorporación de la difusión como modelo del movimiento de los individuos de una población, como se propone en \cite{Skellam1973ThePopulations} y que da lugar a publicaciones posteriores como \cite{Azevedo2012,Hastings1978}. En \cite{Llibre2014} se estudia un modelo de tipo Lotka-Volterra, con dos depredadores\index{depredador} que comparten un único recurso. Los autores utilizan capacidades de carga finita e infinita: cuando la capacidad es finita, aparecen escenarios con exclusión y escenarios con coexistencia; y cuando la capacidad es infinita, hay casos en los que aparecen soluciones fluctuantes, y otros más donde nuevamente hay exclusión por parte de una de las especies.

En 1934 Gause\index{Gause} \cite{Gause2003} presenta el Principio de Exclusión Competitiva\index{Principio de Exclusión Competitiva} (que establece que a lo más $n$ especies pueden coexistir sobre $n$ recursos). Esto llamó la atención de muchos investigadores, como se muestra en \cite{Armstrong1976}, \cite{Hardin2009}, \cite{Kuang2003BiodiversityCompetition} y \cite{McGehee1977}. En el trabajo de McGehee y Armstrong \cite{McGehee1977}, el Principio de Exclusión Competitiva\index{Principio de Exclusión Competitiva} es estudiado desde una nueva perspectiva: los espacios de funciones son la columna vertebral de un modelo matemático que sirve para comprender y analizar la migración, la invasión y el asentamiento de las especies. Se usa la ecuación diferencial
\begin{eqnarray}
\dot x_i &=& x_ig_i(x_1,x_2,\hdots,x_n),                                               \nonumber
\end{eqnarray}
para modelar el crecimiento poblacional de la $i$-ésima especie $x_i$ en el ecosistema\index{ecosistema}. Se toman en cuenta únicamente las soluciones con todas sus coordenadas positivas, es decir las soluciones en $(0,\infty)^n$, y las ecuaciones se reescriben
\begin{eqnarray}
\dot x_i &=& x_iu_i(r_1(x),r_2(x),\hdots,r_m(x)),                                      \nonumber
\end{eqnarray}
donde $r_j$ representa cada factor limitante (por ejemplo, espacio, oxígeno o luz). Tiene sentido entonces escribir a $r_j$ en términos de $x$, pues los recursos dependen de la densidad de cada especie en el medio. Así, queda definida la función $r = (r_1,r_2,\hdots,r_m)\colon (0,\infty)^n\longrightarrow \mathbb R^m$. Las comunidades ecológicas con $n$ especies y $m$ factores limitantes se escriben
\begin{eqnarray}
\Phi_m^n &=& \{g=u\circ r \colon r\in C^\infty((0,\infty)^n,\mathbb R^m),u\in C^\infty(\mathbb R^m,\mathbb R^n)\}.                                                                    \nonumber
\end{eqnarray}
Nótese que $\Phi_m^n$ es un subconjunto de $C^\infty ((0,\infty)^n,\mathbb R^n)$ para cualquier $n$ y para cualquier $m$. Las especies $x_1,x_2,\hdots,x_n$ coexisten si la función correspondiente $g\in \Phi_m^n$ da lugar a un sistema diferencial persistente (es decir, un sistema con atractor en $(0,\infty)^n$). El análisis anterior se puede consultar en \cite{Falconi2002UNAMatematica}.

Otro resultado controversial que surgió de los modelos de depredación es la Paradoja de Enriquecimiento\index{Paradoja de Enriquecimiento}, presentada por Rosenzweig\index{Rosenzweig} \cite{Rosenzweig1971} en 1971, donde se muestran seis modelos planos en los que el incremento de la capacidad del medio provoca la desaparición de los puntos de equilibrio positivos\footnote{Un punto de equilibrio positivo es aquel cuyas coordenadas son todas positivas.} atractores correspondientes a cada modelo. La forma general del sistema que se estudia en \cite{Rosenzweig1971} es
\begin{eqnarray}
\dot x &=& x G(x) - y P(x), \nonumber \\
\dot y &=& y (P(x) - D),    \nonumber
\end{eqnarray}
donde $x$ es la población de la presa, $y$ es la población del depredador, $D$ es la tasa de mortalidad, $G$ es una función de crecimiento\index{función de crecimiento} denso-dependiente con tasa intrínseca de crecimiento $R > 0$ y capacidad de carga $K > 0$ que puede estar dada por un modelo logístico\index{modelo logístico}
\begin{eqnarray}
G(x) &=& R(1 - x/K), \nonumber
\end{eqnarray}
una función de Gompertz\index{función de Gompertz}
\begin{eqnarray}
G(x) &=& R\left(\ln K - \ln x\right), \nonumber
\end{eqnarray}
o bien una función presentada por los autores
\begin{eqnarray}
G(x) &=& R\left(x^{-\beta}/K - 1\right) \quad (0 < \beta \leq 1), \nonumber
\end{eqnarray}
y la respuesta funcional $P$ con tasa de captura $Q > 0$ puede estar dada por una función de Lotka-Volterra generalizada
\begin{eqnarray}
P(x) &=& Q x^\alpha \quad (0 < \alpha \leq 1), \nonumber
\end{eqnarray}
o bien la función
\begin{eqnarray}
P(x) &=& Q\left(1 - \exp(\lambda x)\right) \quad (\lambda > 0). \nonumber
\end{eqnarray}
\begin{figure}[h]
    \centering
    \includegraphics[width=324pt]{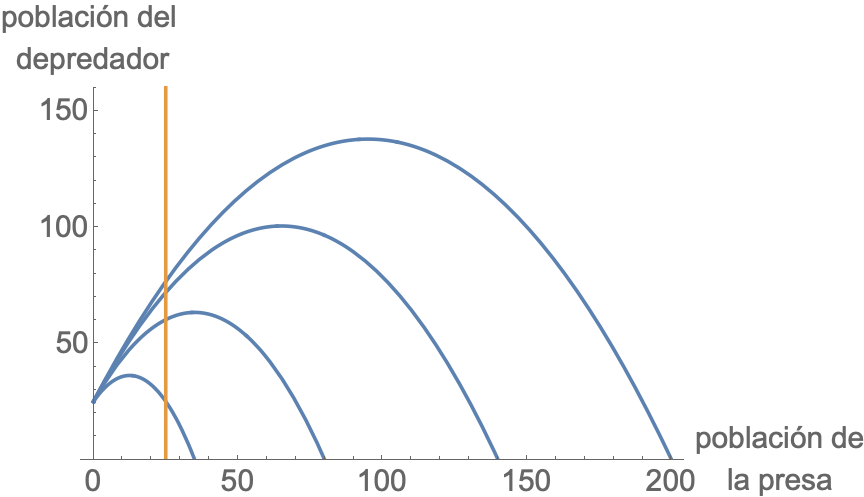}
    \caption{Se muestran algunos ejemplos de isoclinas: una isoclina de la especie depredadora\index{depredador} ($x_0 = 25$), y cuatro isoclinas para la presa\index{presa} ($K = 34$, $80$, $140$, $200$). Solamente la curva $K=34$ tiene intersección con la recta vertical en el costado de la curva que tiene pendiente negativa, y por lo tanto será la única que dé lugar a un punto de equilibrio estable.}
    \label{fig1.1}
\end{figure}

Usando las dos respuestas funcionales y las tres funciones de crecimiento mostradas, los autores obtienen los seis modelos de depredación. En todos los casos, la isoclina $\dot x = 0$ es una curva que abre hacia abajo, y la isoclina $\dot y = 0$ es una recta vertical dada por $x=x_0$ como lo muestra la Figura \ref{fig1.1}. Las isoclinas se intersectan en un único punto de equilibrio positivo, y su estabilidad asintótica depende de la ubicación del punto de equilibrio con respecto de la cúspide de la curva correspondiente a la isoclina $\dot x = 0$ (ya sea a la derecha que corresponde a un punto asintóticamente estable, o a la izquierda que corresponde a un punto inestable). Supóngase que la curva $\dot x = 0$ tiene cúspide con abscisa $x=v$. Entonces se demuestra que
\begin{eqnarray}
\frac{\partial v}{\partial K} &>& 0 \nonumber
\end{eqnarray}
en los seis modelos. Como el valor $x_0$ es fijo y el valor $v$ se desplaza hacia la derecha conforme aumenta la capacidad de carga $K$, entonces eventualmente la recta vertical intersecta a la curva $\dot x = 0$ del lado izquierdo de $v$, dando lugar a un punto de equilibrio positivo inestable.

La Paradoja de Enriquecimiento\index{Paradoja de Enriquecimiento} es la base de publicaciones subsecuentes como \cite{Abrams1996} y \cite{Vos2004}. En \cite{Abrams1996} se muestra que la densidad poblacional de las presas\index{presa} se puede subdividir en vulnerable e invulnerable para mostrar que si la transición de una categoría a la otra es adecuada, entonces el incremento de la capacidad de carga de la presa\index{presa} no necesariamente implica la inestabilidad del equilibrio positivo del modelo. En \cite{Vos2004} se muestra a través de dos modelos que la Paradoja de Enriquecimiento\index{Paradoja de Enriquecimiento} puede no cumplirse si las plantas inducen mecanismos de defensa ante la llegada de las especies herbívoras y si la mortalidad del depredador\index{depredador} en el nivel trófico más alto se encuentra ubicada en un rango determinado. 

En relación con el tema de las respuestas funcionales, en \cite{Cosner1999}, aparecen algunas de las respuestas funcionales más comunes en la literatura. Cada una de ellas refleja la forma en que se agrupan o se distribuyen los depredadores\index{depredador}. Específicamente, se discute la importancia de la respuesta funcional en un modelo matemático. Si la distribución de la especie depredadora\index{depredador} es homogénea, entonces la respuesta funcional $p$ es dependiente únicamente de la población de la presa\index{presa}:
\begin{eqnarray}
p(x,y) &=& p(x). \nonumber
\end{eqnarray}
Sin embargo el resultado es distinto cuando las colonias de cazadores se encuentran acomodadas, por ejemplo, de acuerdo al relieve o a la humedad del ecosistema\index{ecosistema}. Es así como surgen las respuestas funcionales que dependen también de la población del depredador\index{depredador}. Dichas respuestas se ven afectadas si la población del depredador\index{depredador} crece. 

Para el análisis de las respuestas funcionales, se establece la siguiente estructura
\begin{eqnarray}
p(x,y) &=& \frac{q e(x,y)/y}{1 + a e(x,y)/y}, \nonumber
\end{eqnarray}
donde $e(x,y)$ representa el número de encuentros, por depredador\index{depredador}, por unidad de tiempo. Además $q$ indica la cantidad de biomasa de la presa\index{presa} que será eliminada por el depredador\index{depredador} en cada encuentro, y la constante $a$ es el tiempo de manejo de la presa\index{presa} capturada.

El primer caso que se presenta en \cite{Cosner1999} es cuando los encuentros ocurren de acuerdo con la ley de acción de masas: $e(x,y) = e_0 x y$. Haciendo una simple modificación, se puede obtener una respuesta funcional de tipo Holling II
\begin{eqnarray}
p(x) &=& \frac{q x}{1 + a x}. \nonumber
\end{eqnarray}
Si el número de depredadores\index{depredador} en el área no afecta el número de encuentros, entonces $e(x,y) = e_0 x$, y así
\begin{eqnarray}
p(x) &=& \frac{q x}{y + a x}. \nonumber
\end{eqnarray}
Se puede observar el hecho de que cuando los depredadores\index{depredador} se organizan en manadas o enjambres para llevar a cabo la cacería, las formaciones bidimensionales dan lugar a respuestas funcionales
\begin{eqnarray}
p(x,y) &=& \frac{q x}{y^{1/2} + a x}, \nonumber
\end{eqnarray}
mientras que las tridimensionales dan lugar a
\begin{eqnarray}
p(x,y) &=& \frac{q x}{y^{1/3} + a x}. \nonumber
\end{eqnarray}
Los dos escenarios anteriores se conocen comúnmente como funciones de tipo Hassel-Varley. También en \cite{Cosner1999} se estudia el caso en el que la cantidad de depredadores\index{depredador} al acecho es importante si $y$ es pequeño (respecto a $y_0$), pero que al incrementarse $y$ deja de ser relevante. En este caso, el número de encuentros está dado por 
\begin{eqnarray}
e(x,y) &=& \frac{e_0 y_0 x y}{y_0 + y}, \nonumber
\end{eqnarray}
y así se obtiene la respuesta funcional de Beddington-DeAngelis\index{respuesta funcional de Beddington-DeAngelis}
\begin{eqnarray}
p(x,y) &=& \frac{qx}{1 + ax + by}.  \nonumber
\end{eqnarray}

Con respecto a la competencia y al papel que juega en el ecosistema\index{ecosistema}, en \cite{Chesson2008} se menciona que, junto con la depredación, la competencia promueve la conservación de la biodiversidad. En \cite{Lehman2000} se estudian tres modelos distintos para analizar la relación positiva entre la competencia y la productividad del medio, en los que se muestra la importancia que tiene la competencia. 

En otros casos, se ha optado por dar una perspectiva más realista a los modelos al incorporar retardo a las ecuaciones. En \cite{Li2015} se aborda este tópico desde el punto de vista de las ecuaciones diferenciales parciales para conocer los casos en los que existe un punto de equilibrio positivo estable, un punto de equilibrio inestable, o bien se presentan órbitas periódicas. En otros casos, el teorema de continuación de Mawhin\index{Teorema de Continuación de Mawhin} \cite{Gaines1997CoincidenceEquations} es usado para hallar órbitas periódicas en modelos de depredación, como por ejemplo en \cite{Yang2008}, donde se trabaja con un modelo de depredación con retardo estructurado con tres etapas, en \cite{Li2000} donde los autores estudian un modelo de tipo Lotka-Volterra con retardo, y en \cite{Lin2013} donde se trabaja un modelo con retardo con respuesta funcional de tipo Beddington-DeAngelis. En \cite{Cai2015} se muestra que en modelos con retardo la estabilidad asintótica de soluciones periódicas es posible y se presentan simulaciones numéricas; en \cite{Tripathi2015} se presenta un modelo con respuesta funcional de tipo Crowley-Martin\index{respuesta funcional de Crowley-Martin} y se comparan los resultados obtenidos al incorporar retardo (o gestación, como se le nombra en ese trabajo) a la ecuación; en \cite{Nindjin2006AnalysisDelay} se estudia un modelo Leslie-Gower con respuesta funcional Holling II y se discute acerca del efecto del retardo en el modelo propuesto. 

El efecto Allee es tema de interés de algunos trabajos. En \cite{Cai2015}, se muestra como se modifica la estructura del retrato fase al incorporar el efecto; en \cite{Etienne2002} un modelo con efecto Allee basado en los hábitos de la \textit{Drosophila} es construido; en \cite{GASCOIGNE2004} se discute acerca de las condiciones biológicas que se deben cumplir para que el efecto Allee pueda hacerse presente; en \cite{Hadjiavgousti2008} los umbrales que permiten la coexistencia en un modelo discreto con efecto Allee son analizados. 

El refugio o la defensa es un factor a considerarse en trabajos recientes. En \cite{Matsuda1994} se estudian las consecuencias potenciales de la efectividad de los mecanismos de defensa de las presas\index{presa} a través de algunos modelos matemáticos. En \cite{Tang2014} un modelo con respuesta funcional de tipo Holling II y refugio constante es presentado y usando el teorema de Lyapunov y el teorema de Bendixson-Dulac se concluye la estabilidad asintótica de un punto de equilibrio positivo;  en \cite{Vos2004} se estudian modelos bitróficos y tritróficos, se incorporan parámetros que representan mecanismos de defensa variables y permanentes, y se dan condiciones suficientes para que haya estabilidad en los modelos. 

Con respecto a los modelos con difusión, en \cite{ko} un modelo de tipo Hassel-Varley con difusión es presentado y se hace uso de la teoría del punto fijo para poder conocer las condiciones que garanticen la existencia de órbitas periódicas; además, los modelos matemáticos con difusión muestran como la población de roedores \cite{Azevedo2012} y de anfibios \cite{Fonseca2013} tiende a disminuir en el Amazonas en los años más recientes.

Los trabajos mencionados muestran de manera muy general la riqueza de la problemática en ecología, y constituyen el marco de la tesis. A continuación se da una descripción del problema de interés.

El objetivo general de este trabajo, que conecta de forma natural la ecología con las matemáticas (específicamente los sistemas dinámicos), es estudiar las consecuencias que acarrea consigo la invasión de una especie exótica a un ecosistema\index{ecosistema}. La llegada de un elemento externo a un entorno que se encuentra momentáneamente en equilibrio puede tener efectos devastadores en las poblaciones de las especies residentes. Ejemplo de lo anterior es el daño producido al ajolote (\textit{Ambystoma mexicanum}), en el sistema de canales de Xochimilco debido a la llegada de especies exóticas al medio (ver por ejemplo \cite{Contreras2009, Robles2009, Zambrano2007, Zambrano2010}). En \cite[Capítulo VII]{SoberonMainero2000} se menciona que alrededor del año 1840 fue introducida en Australia la especie exótica de nopal \textit{Opuntia stricta}, la cual se convirtió rápidamente en una plaga que comenzó a depredar los ecosistemas\index{ecosistema} endémicos. Por otro lado, en \cite{Ingenloff2017PredictableDecaocto} se presenta un estudio hecho a través de los nichos ecológicos para conocer las consecuencias de la invasión de la especie \textit{Streptopelia decaocto} en Norteamérica. El estudio de estos modelos presenta preguntas interesantes desde el punto de vista de la dinámica misma.

El problema se estudia desde una perspectiva general, lo que permite obtener robustez en los modelos matemáticos, en el sentido de que estos pueden ser aplicados en cualquier escenario ecológico que cumpla con las condiciones que se plantean.

El trabajo se divide en dos partes. 

La Parte \ref{partei} son los Capítulos \ref{capitulo2} y \ref{capitulo3}, donde se estudia el impacto que tiene la competencia en un entorno en el que una presa\index{presa} hace frente a varios depredadores\index{depredador}. Modelos matemáticos son propuestos y analizados, y se dan condiciones suficientes para que los puntos de equilibrio positivos correspondientes a cada modelo sean atractores. Lo anterior permite concluir que la coexistencia de las especies es factible. El aporte de la Parte \ref{partei} está en el Teorema \ref{teo2.25}, en el que se muestra la posibilidad de variar los parámetros correspondientes competencia de los consumidores sin que se vean afectadas algunas propiedades como la estabilidad asintótica del punto de equilibrio positivo en el modelo principal del Capítulo \ref{capitulo2}, y en los Teoremas \ref{teo3.9} y \ref{teo3.10} donde de forma respectiva se proveen condiciones para 
la estabilidad asintótica de un punto de equilibrio positivo en el modelo principal del Capítulo \ref{capitulo3}, y se muestra la existencia de una bifurcación de Hopf en el mismo modelo.

La Parte \ref{parteii} son los Capítulos \ref{capitulo4} y \ref{capitulo5}. Ahí se estudia el problema de la invasión a través de un modelo con dos presas\index{presa} y un depredador\index{depredador}. En el Capítulo \ref{capitulo4} se presenta un modelo matemático en su forma más general y se siguen ideas expuestas en \cite{Abrams1999} para mostrar que la tasa intrínseca de crecimiento de la especie invasora debe superar cierta cota para que esta última pueda ingresar al medio sin que esto conlleve a la desaparición o a la migración de los residentes. La coexistencia de las especies (vista a través de los sistemas dinámicos) se puede obtener a través de un punto de equilibrio atractor, de un ciclo límite atractor, o incluso a través de dinámica caótica; también en el Capítulo \ref{capitulo4} se realiza una construcción geométrica que permite hallar puntos de equilibrio positivos para el modelo propuesto y se proveen condiciones suficientes para que un punto de equilibrio sea un atractor, además de que se trabaja con las bifurcaciones de dicho modelo matemático. Las contribuciones anteriores aparecen a los Teoremas \ref{teo4.18} y \ref{teo4.20}; en el Capítulo \ref{capitulo5} se estudia la existencia de oscilaciones caóticas en el modelo y se analizan las propiedades correspondientes.
\part{Competencia: una presa y varios depredadores}
\label{partei}
\chapter*{Parte \ref{partei} \\ $ $ \\ Competencia: una presa \\ y varios depredadores}
Los modelos matemáticos son una herramienta notable en ecología, pues a través de ellos, es posible estudiar las interacciones de las especies y obtener resultados aplicables en el laboratorio o en el campo. Con respecto a estas interacciones, la depredación y la competencia son incorporadas en los modelos que aquí se proponen. Los depredadores\index{depredador} que aparecen en el modelo matemático pueden ser carnívoros que se alimentan de otra especie animal, herbívoros que consumen un pastizal, incluso bacterias o protozoarios.

 En esta parte del trabajo se presenta el caso en el que una especie hace frente a varios depredadores\index{depredador} comunes, y se analiza el impacto que tiene la competencia a través de un sistema dinámico que modela las interacciones entre las especies que comparten un recurso único y no remplazable.
 
 En el Capítulo \ref{capitulo2} se presenta el modelo \eqref{2.1} que usa la respuesta funcional de tipo Lotka-Volterra (que refleja tasa de depredación constante). Una característica de este modelo es que los resultados obtenidos son válidos para cualquier función $g$ de crecimiento\index{función de crecimiento} del recurso decreciente y diferenciable en el intervalo $[0,K]$. La Proposición \ref{prop2.17} muestra que el incremento del coeficiente de la competencia entre los consumidores tiende a hacer crecer el valor de equilibrio de la población de la especie recurso. En el Teorema \ref{teo2.19} se dan condiciones suficientes para la estabilidad asintótica global de un punto de equilibrio positivo del sistema; lo anterior permite dar respuesta al problema de la inserción de especies exóticas en un ecosistema\index{ecosistema}, pues muestra que nuevos consumidores pueden ser introducidos exitosamente si la competencia intraespecífica que los nuevos miembros de la comunidad ejecutan es suficientemente alta en comparación con la competencia intra- e interespecífica de la especie nativa.

En el Capítulo \ref{capitulo3} se analiza el modelo \eqref{3.1} que contempla dos especies depredadoras\index{depredador} y una especie recurso. La función de crecimiento\index{función de crecimiento} del recurso en este caso, es el modelo logístico\index{modelo logístico}, mientras que la tasa de depredación está dada por la respuesta funcional de tipo Holling II. La competencia intraespecífica está dada por valores constantes, y se desprecia la competencia interespecífica. A partir de este análisis, se muestra la existencia de un umbral de bifurcación de Hopf. En el Apéndice \ref{apendicea} aparece el submodelo plano que se obtiene al retirar uno de los consumidores del sistema original (que en la literatura se conoce comúnmente como el modelo de Bazykin). Este submodelo muestra como es el ecosistema\index{ecosistema} antes de la llegada de una especie exótica que se alimenta del mismo recurso que la especie nativa. 
\chapter{Tasa de captura constante}
\label{capitulo2}
\begin{flushright}
\emph{“Solo después de haber conocido la\\superficie de las cosas, se puede\\uno animar a buscar lo que hay\\debajo. Pero la superficie de las\\cosas es inagotable.”\\$ $\\ Italo Calvino \cite{Calvino2017Palomar}}
\end{flushright}

Los modelos de depredación han sido usados para analizar las interacciones entre una especie carnívora y una herbívora, o una especie herbívora y el pastizal que representa su fuente alimentaria, o incluso para analizar el comportamiento de bacterias y protozoarios. Estos modelos han presentado muchas variantes, como por ejemplo la estructura de edades. En \cite{Ma2015} un modelo con dos presas\index{presa} y un solo depredador\index{depredador} en dos diferentes etapas es introducido, y se muestra que en algunos casos reales, la inserción de pesticidas y la implementación de medidas de prevención de enfermedades contribuye a estabilizar un punto de equilibrio positivo. También en \cite{Song2006}, los autores muestran que un sistema con estructura de edades presenta una bifurcación de ruptura de simetría y caos; en \cite{Shi2011} un sistema con estructura de edades para el depredador\index{depredador} y respuesta funcional de Crowley-Martin\index{respuesta funcional de Crowley-Martin} es estudiado, y se dan condiciones suficientes para la existencia de puntos de equilibrio positivo que sean atractores globales para el sistema. En \cite{Khajanchi2014} un modelo de depredación con estructura de edades con respuesta funcional de tipo Beddington-DeAngelis\index{respuesta funcional de Beddington-DeAngelis} es analizado, y se muestra que la estructura de edades para especies depredadoras\index{depredador} implica la existencia de órbitas periódicas. Además de lo anterior, modelos con depredación generalista aparecen en \cite{Ali2013b} y en \cite{Yu2014}. En el primero se muestra un modelo con respuesta funcional de tipo Beddington-DeAngelis\index{respuesta funcional de Beddington-DeAngelis}, y en el segundo se muestra un modelo con respuesta funcional de tipo Crowley-Martin\index{respuesta funcional de Crowley-Martin}.

La competencia también ha sido estudiada desde una perspectiva teórica global en \cite{Armstrong1976} y \cite{McGehee1977}. El Principio de Exclusión Competitiva\index{Principio de Exclusión Competitiva} (que establece que a lo más $n$ especies pueden coexistir en $n$ recursos) fue enunciado por Gause\index{Gause} \cite{Gause2003} en 1934 y consiguió atrapar la atención de algunos investigadores (ver por ejemplo \cite{Armstrong1976}, \cite{Hardin2009}, \cite{Kuang2003BiodiversityCompetition} y \cite{McGehee1977}).En \cite{Kuang2003BiodiversityCompetition} un modelo de tipo Lotka-Volterra es analizado y se muestra que la competencia intraespecífica en los consumidores favorece la coexistencia. En \cite{Loladze2004} se introduce un sistema con dos consumidores en un solo recurso, y a partir de la composición química del recurso, se concluye que la heterogeneidad química de las especies propicia la coexistencia a través de un punto de equilibrio positivo estable en el sistema. En \cite{Lehman2000} se realiza un análisis comparativo de tres modelos multiespecíficos para poder discutir el vínculo entre biodiversidad y estabilidad a nivel población y a nivel comunidad. En \cite{Lin2013} el efecto de retardo en coeficientes intra- e interespecíficos es estudiado, y se obtienen órbitas periódicas que numéricamente aparentan ser estables. La manera en la que la depredación y la competencia actúan para favorecer la biodiversidad es un tema de interés; en \cite{Chesson2008} se establece que para preservar la diversidad de las especies, los lazos intraespecíficos deben ser más fuertes que los lazos interespecíficos. Para mostrar esto, la depredación y la competencia son cuantificadas y analizadas a través de un sistema tritrófico. En \cite{Abrams1989}, se estudian los nutrientes en un recurso y los recursos son clasificados en no-remplazables y complementarios.

Tomando como base las publicaciones anteriores, se propone analizar un entorno en el que varios consumidores compiten por un solo recurso. Con esa finalidad, se considera un modelo matemático con $n+1$ ecuaciones diferenciales que representan las poblaciones de un recurso y de $n$ consumidores. El crecimiento del recurso está dado por una función $g$, y a su vez tiene un decremento causado por el $i$-ésimo consumidor a razón de la tasa de captura $q_i$. Adicionalmente, la tasa de natalidad\index{natalidad} del $i$-ésimo consumidor es proporcional a la biomasa capturada, la tasa de mortalidad se denota con $\mu_i$ y los coeficientes de competencia con $m_{ij}$. Concretamente, el modelo está dado por:
\begin{eqnarray}
\label{2.1}
\dot x   &=& x  \left(r\, g(x) - q_1y_1 - q_2y_2 - \hdots \hdots - q_ny_n\right),                      \nonumber \\
\dot y_1 &=& y_1\left(c_1q_1x - \mu_1 - m_{11}y_1 - m_{12}y_2 - \hdots \hdots - m_{1n}y_n\right),  \nonumber \\
\dot y_2 &=& y_2\left(c_2q_2x - \mu_2 - m_{21}y_1 - m_{22}y_2 - \hdots \hdots - m_{2n}y_n\right),            \\
&\vdots&                                                                                           \nonumber \\
\dot y_n &=& y_n\left(c_nq_nx - \mu_n - m_{n1}y_1 - m_{n2}y_2 - \hdots \hdots - m_{nn}y_n\right),  \nonumber \\
                                                                                                    \nonumber \\
& & n>0\text{, }x(0)>0\text{, }y_1(0)>0\text{, }y_2(0)>0\text{, }\hdots \hdots \text{, }y_n(0)>0,  \nonumber
\end{eqnarray}
donde $x$ representa la población del recurso, mientras que $y_i$ es la población del $i$-ésimo consumidor para $i=1,\hdots,n$. La función de crecimiento\index{función de crecimiento} $g$ es diferenciable, $g'(x)<0$ para toda $x>0$, y $g(K) = 0$. Los valores $r$ y $K$ son la tasa intrínseca de crecimiento y la capacidad de carga respectivamente. La tasa de captura de cada consumidor es constante (y está dada por $q_i$); $c_i$ es la eficiencia en la conversión de biomasa de recurso capturada a biomasa del consumidor. Como ya se dijo, $\mu_i$ es la mortalidad del $i$-ésimo depredador y $m_{ij}$ es la competencia del $i$-ésimo depredador sobre el $j$-ésimo depredador. En el resultado que aparece en \cite[Teorema 3.1]{Kuang2003BiodiversityCompetition} se presenta un modelo similar, pero que solamente considera competencia intraespecífica. La competencia interespecífica ha sido estudiada en \cite{Connell1961}, \cite{Hairston1960}, \cite{Holt1977}, \cite{Schoener1978} y \cite{Sinclair1985}. En \cite{Holt1977}, se menciona que dos especies están en \textit{competencia aparente} si la presencia de una de las especies induce una disminución en la densidad de equilibrio de las otras especies, y se muestran dos casos concretos (en el primero se estudian liebres árticas, y en el segundo se estudia zooplancton), en los que el arribo de competidores ocasiona la migración de las especies residentes.

Se asume que $r$, $K$, $c_i$, $q_i$, $\mu_i$, $m_{ii} > 0$, y $m_{ij} \geq 0$, (para $i\neq j$). A pesar de que \eqref{2.1} puede tener hasta $2^n + 1$ diferentes puntos de equilibrio (esto se verifica haciendo cero el lado izquierdo de \eqref{2.1} y resolviendo el sistema correspondiente), \eqref{2.1} tiene solamente un punto de equilibrio $E^*$ con coordenadas positivas. La existencia de $E^*$ resulta ser de particular interés, pues representa un estado en el que todas las especies pueden coexistir en el medio.

A partir del modelo propuesto, surgen algunas preguntas: ¿Pueden los competidores (que se alimentan de un único recurso) alcanzar la coexistencia? En caso de que lo anterior ocurra, ¿cuáles son las estrategias que les permitirán llegar a ese estado? ¿Cómo se refleja la competencia interespecífica en la coexistencia de las especies? Para responder estas preguntas, se hace un análisis detallado del punto de equilibrio positivo $E^*$.
\section{Propiedades generales del modelo}
\label{seccion2.1}
 Se analiza la disipatividad del sistema \eqref{2.1} en el Teorema \ref{teo2.1}, aparece la definición de umbral de eficiencia energética (un concepto que permite profundizar en el análisis del modelo), se calcula la matriz jacobiana de \eqref{2.1} para conocer las propiedades de estabilidad de dos de sus puntos de equilibrio en la frontera y se estudia la persistencia uniforme de \eqref{2.1} en el Teorema \ref{teo2.7}.
 
La demostración del siguiente teorema es similar a la que se hace en \cite{Sarwardi2013DynamicalRefuge}.
\begin{teo}
\label{teo2.1}
Existe un subconjunto compacto $C$ de $[0,\infty)^{n+1}$ de tal manera que para toda solución no-negativa $(x,y_1,\hdots,y_n)$ de \eqref{2.1}, existe $T > 0$ tal que $(x(t),y_1(t),\hdots,y_n(t)) \in C$ para toda $t\geq T$.
\end{teo}
\begin{proof}
Se considera $\displaystyle w = x + \sum_{i=1}^n c_i^{-1} y_i$, se toma $\varphi$ en el intervalo $\left(0,\min\left\{\mu_i\right\}_{i=1}^n\right]$, y se define \newline $\displaystyle \rho = \frac{K(r + \varphi)^2}{4 r}$. Entonces se tiene
\begin{eqnarray}
w' + \varphi w &=& x(r(1 - x/K)) - \sum_{i=1}^nc_i^{-1}\mu_i y_i - \sum_{i=1}^n\left(\sum_{j=1}^n m_{ij}y_iy_j\right) + \varphi x + \sum_{i=1}^nc_i^{-1}\varphi y_i \nonumber \\
&\leq& x\left(r + \varphi - \left(\frac{r}{K}\right)x\right) - \sum_{i=1}^n c_i^{-1}(\mu_i - \varphi)y_i \nonumber \\
&\leq& x\left(r + \varphi - \left(\frac{r}{K}\right)x\right) \nonumber \\
&\leq& \frac{K(r + \varphi)^2}{4 r} \nonumber \\
&=& \rho. \nonumber
\end{eqnarray}

Del análisis previo, se sigue que $\displaystyle w \leq e^{-\varphi t}w_0 + \frac{\rho}{\varphi}(1 - e^{-\varphi t})$. 

Así, se tiene que $\displaystyle \limsup_{t\to \infty} w(t) \leq \frac{\rho}{\varphi}$. 

Se concluye entonces que para cualquier solución no-negativa $(x,y_1,\hdots,y_n)$ de \eqref{2.1}, existe $T>0$ tal que
\begin{eqnarray}
(x(t),y_1(t),\hdots,y_n(t)) &\in& \left\{(x,y_1,\hdots,y_n)\in [0,\infty)^{n+1}\colon x + \sum_{i=1}^n c_i^{-1}y_i\leq \frac{\rho}{\varphi}\right\}\text{, si } t \geq T. \nonumber
\end{eqnarray}
\end{proof}
Como las curvas integrales de \eqref{2.1} están eventualmente acotadas, se puede realizar el estudio del espacio de estados en una celda compacta de $[0,\infty)^{n+1}$. Además es posible verificar que el intervalo maximal de existencia de una solución $\boldsymbol x(t)$ de \eqref{2.1} es el intervalo no-negativo $[0,\infty)$, (ver \cite[Sección 2.4, Corolario 2]{Perko2001}).

\begin{df}\cite{Kuang2003BiodiversityCompetition}
Dado el sistema \eqref{2.1}, el \textbf{umbral de eficiencia energética (u.e.e.) de la población del recurso con respecto al consumidor} $\boldsymbol{y_i}$ es la cantidad de recurso requerido en el medio para que la tasa de crecimiento y la tasa de mortalidad de la $i$-ésima especie depredadora\index{depredador} coincidan. En otras palabras, el u.e.e. $x_i^b$ es el valor positivo que satisface
\begin{eqnarray}
x_i^b &=& \frac{\mu_i}{c_iq_i}.                                                                                    \nonumber
\end{eqnarray}
\end{df}
El u.e.e. provee información acerca de la adaptabilidad o de las posibilidades de la especie $y_i$ para permanecer en el hábitat. De acuerdo con Kuang et al. \cite[Proposición 2.1]{Kuang2003BiodiversityCompetition}, si el u.e.e. $x_i^b$ (de la especie depredadora\index{depredador} $y_i$) es mayor que la capacidad de carga $K$, entonces la especie morirá. Se asume entonces que $x_i^b < K$ para $i=1,2,\hdots,n$.

A continuación se muestra que los puntos de equilibrio $E_{00} = (0,\hdots,0)$, $E_0 = (K,0,\hdots,0)$ son puntos silla. Además se muestran las variedades estable\index{variedad estable} e inestable\index{variedad inestable} de $E_{00}$, y la variedad estable\index{variedad estable} de $E_0$. El punto de equilibrio en la frontera en \eqref{2.1} $E_{00}$ representa la exclusión de todas las especies, mientras que $E_0$ representa la exclusión de todos los consumidores y la saturación del recurso.

Si $\boldsymbol x = (x,y_1,y_2,\hdots,y_n)$, entonces la matriz jacobiana del sistema $\eqref{2.1}$ en el punto $\boldsymbol x$ es
\begin{eqnarray}
J(\boldsymbol x) &=& \left(\begin{array}{ccccc} 
r\, (xg'(x) + g(x))  & -q_1x                    & -q_2x                    
                     & \hdots                   & -q_nx                  \\
- q_1y_1 - q_2y_2    &                          &
                     &                          &                        \\
- \hdots - q_ny_n    &                          &
                     &                          &                       \\\\
c_1q_1y_1            & c_1q_1x - \mu_1          & -m_{12}y_1    
                     & \hdots                   & -m_{1n}y_1             \\
                     & - 2m_{11}y_1 - m_{12}y_2 &
                     &                          &                        \\
                     & - \hdots - m_{1n}y_n     &
                     &                          &                      \\\\
c_2q_2y_2            & -m_{21}y_2               & c_2q_2x - \mu_2 
                     & \hdots                   & -m_{2n}                \\
                     &                          & -m_{21}y_1-2m_{22}y_2 
                     &                          &                        \\
                     &                          & -\hdots-m_{2n}y_n 
                     &                          &                       \\\\
\vdots               & \vdots                   & \vdots        
                     & \ddots                   & \vdots                \\\\
c_nq_ny_n            & -m_{n1}y_n               & -m_{n2}y_n  
                     & \hdots                   & c_nq_nx - \mu_n         \\
                     &                          &                          
                     &                          & - m_{n1}y_1 - m_{n2}y_2 \\
                     &                          &                          
                     &                          & -\hdots - 2m_{nn}y_n 
\end{array}\right).\nonumber
\end{eqnarray}

Si $\displaystyle \lim_{x\to 0^+} xg'(x) + g(x) = \beta > 0$, entonces $J$ evaluada en $E_{00}$ está dada por
\begin{eqnarray}
J(E_{00}) &=& \left(\begin{array}{ccccc} 
r \beta & 0        & 0        & \hdots & 0      \\ 
0       & -\mu_1   & 0        & \hdots & 0      \\
0       & 0        & -\mu_2   & \hdots & 0      \\
\vdots  & \vdots   & \vdots   & \ddots & \vdots \\
0       & 0        & 0        & \hdots & -\mu_n   \end{array}\right),\nonumber
\end{eqnarray}
por lo que $E_{00}$ es un punto silla en \eqref{2.1} con variedad estable\index{variedad estable}
\begin{eqnarray}
W^s(E_{00}) &=& \{0\} \times [0,\infty)^n \nonumber,
\end{eqnarray}
y variedad inestable\index{variedad inestable}
\begin{eqnarray}
W^u(E_{00}) &=& [0,\infty) \times \{0\}^n \nonumber.
\end{eqnarray}
Análogamente se obtiene
\begin{eqnarray}
J(E_{0}) &=& \left(\begin{array}{ccccc} 
rKg'(K)      & -q_1K         & -q_2K         & \hdots & -q_nK         \\ 
 0          & c_1q_1K-\mu_1 & 0             & \hdots & 0             \\
 0          & 0             & c_2q_2K-\mu_2 & \hdots & 0             \\
 \vdots     & \vdots        & \vdots        & \ddots & \vdots        \\
 0          & 0             & 0             & \hdots & c_nq_nK-\mu_n 
\end{array}\right).\nonumber
\end{eqnarray}
Dado que los valores u.e.e. $x_i^b$ son menores que $K$, entonces $E_0$ es también un punto silla y
\begin{eqnarray}
W^s(E_{0}) &=& [0,\infty) \times \{0\}^n. \nonumber \end{eqnarray}
\begin{teo}
\label{teo2.3}
Si $x_i^b < K$ para toda $i = 1,2,\hdots,n$, entonces $E_{00}$ y $E_0$ son puntos silla en \eqref{2.1}.
\end{teo}
Otros equilibrios en la frontera distintos a los mostrados anteriormente son de la forma $\widetilde E = (\widetilde x,\widetilde y_1,\hdots,\widetilde y_k,0,\hdots,0)$, con $\widetilde x, \widetilde y_1, \hdots, \widetilde y_n > 0$, y representan la exclusión de algunos de los consumidores. Estos casos no son estudiados aquí.

Para analizar la persistencia de \eqref{2.1}, se usan algunas definiciones y resultados que aparecen en \cite{Freedman1994UniformSet}. 
\begin{df} 
Dado un sistema dinámico con función de flujo $\boldsymbol \varphi$ y un conjunto invariante $\Gamma$ para $\boldsymbol \varphi$, se dice que el sistema es
\begin{enumerate}
    \item \textbf{Débilmente persistente en} $\boldsymbol \Gamma$, si para toda $\overline x \in int(\Gamma)$ se cumple
    \begin{eqnarray}
    \limsup_{t\to \infty} d(\boldsymbol \varphi(\overline x,t),\partial \Gamma) &>& 0. \nonumber
    \end{eqnarray}
    \item \textbf{Persistente en} $\boldsymbol \Gamma$ si para toda $\overline x \in int(\Gamma)$ se tiene
    \begin{eqnarray}
    \liminf_{t\to \infty}d(\boldsymbol \varphi(\overline x,t),\partial \Gamma) &>& 0. \nonumber
    \end{eqnarray}
    \item \textbf{Uniformemente persistente en} $\boldsymbol \Gamma$ si existe $\varepsilon > 0$ tal que para toda $\overline x\in int(\Gamma)$, se tiene
    \begin{eqnarray}
    \liminf_{t\to \infty}d(\boldsymbol \varphi(\overline x,t),\partial \Gamma) &>& \varepsilon.     \nonumber
    \end{eqnarray}
\end{enumerate}
\end{df}
\begin{teo}
\label{teo2.5}
Sea $\boldsymbol \varphi$ una función de flujo y sea $\Gamma$ un conjunto invariante. Si adicionalmente existe $r > 0 $ tal que el flujo es disipativo en $N_r (\Gamma) \setminus \Gamma$, entonces se cumple exactamente una de las siguientes afirmaciones:
\begin{enumerate}
    \item Para todo $\varepsilon > 0$, existe un conjunto invariante $\mathcal K \subset N_\varepsilon (\Gamma) \setminus \Gamma$.
    \item Existe $\overline x \in N_r(\Gamma) \setminus \Gamma$ de tal manera que $\omega(\overline x) \subset \Gamma$.
    \item Existe $\varepsilon > 0$ tal que para toda $\overline x \in N_r(\Gamma)\setminus \Gamma$, se cumple
    \begin{eqnarray}
    \lim_{t\to \infty}d(\boldsymbol \varphi(\overline x,t),\partial \Gamma) &\geq& \varepsilon \nonumber
    \end{eqnarray}
\end{enumerate}
\end{teo}
La demostración del Teorema \ref{teo2.5} se encuentra en \cite{Freedman1994UniformSet}. A continuación se usa el Teorema \ref{teo2.5} para verificar que el sistema \eqref{2.1} es uniformemente persistente.
\begin{lema}
\label{lema2.6}
Si $n=1$ y $x_1^b < K$, entonces el sistema \eqref{2.1} es débilmente persistente en $[0,\infty)^2$.
\end{lema}
\begin{proof}
Se procede por contradicción. Supóngase que existe $\overline x \in [0,\infty)^2$ tal que 
\begin{eqnarray}
\limsup_{t\to \infty} d(\boldsymbol \varphi(\overline x, t),\partial ([0,\infty]^2)) &=& 0.    \nonumber
\end{eqnarray}
Se analiza una a una las entradas de $\boldsymbol \varphi (\cdot) = (x(\cdot),y_1(\cdot))$.
\begin{enumerate}
    \item Supóngase que $\displaystyle \limsup_{t\to \infty} x(t) = 0$. Considérese $\varepsilon > 0$ arbitrario tal que $\varepsilon < x_1^b$. Existe $t_0 > 0$ de tal manera que si $t \geq t_0$, entonces $x(t) < \varepsilon$. Así $\dot y_1 < y_1(c_1q_1 \varepsilon - \mu_1) < 0$. De lo anterior se sigue que
    \begin{eqnarray}
    \lim_{t\to \infty}\boldsymbol \varphi(t) = \lim_{t\to \infty}(x(t),y_1(t)) &=& (0,0),                \nonumber
    \end{eqnarray}
    lo cual es imposible, pues $E_{00}=(0,0)$ es un punto silla con variedad estable\index{variedad estable} $\{0\}\times [0,\infty)$ (ver Teorema \ref{teo2.3}).
    \item Supóngase que
    \begin{eqnarray}
    \limsup_{t\to \infty} y_1(t) = 0. \nonumber
    \end{eqnarray}
     Cuando $t$ tiende a infinito, la solución $\boldsymbol \varphi (t) = (x(t),y_1(t))$ no puede converger a $E_0 = (K,0)$ ya que la variedad estable\index{variedad estable} de $E_0$ es $[0,\infty) \times \{0\}$ ni tampoco puede converger a $E_{00}$.
     Al no haber más puntos de equilibrio en el eje $X$ (con coordenada $y_1 = 0$), se concluye que $y_1(t)$ no tiende a cero.
\end{enumerate}
\end{proof}
En el siguiente resultado, se demuestra la persistencia uniforme usando el Lema \ref{lema2.6}
\begin{teo}
\label{teo2.7}
Si $n=1$ y $x_1^b<K$, entonces el sistema \eqref{2.1} es uniformemente persistente.
\end{teo}
\begin{proof}
Primero se analiza el eje $X$. Sea entonces $\Gamma_1 = [0,\infty) \times \{0\}$. Ya se demostró (Teorema \ref{teo2.1}) que el sistema \eqref{2.1} es disipativo, entonces se cumple exactamente una de las afirmaciones del Teorema \ref{teo2.5}. La afirmación 1. falla porque $\Gamma_1$ es una variedad invariante maximal de \eqref{2.1} (ver demostración del Lema \ref{lema2.6}). La afirmación 2. tampoco se cumple ya que de acuerdo al Lema \ref{lema2.6}, el sistema \eqref{2.1} es débilmente persistente. De lo anterior, se sigue que
\begin{eqnarray}
\liminf_{t \to \infty} d(\boldsymbol \varphi(x,t),\partial \Gamma_1) &>& \varepsilon_0\text{, para algún } \varepsilon_0 >0.     \nonumber 
\end{eqnarray}
De forma similar, se verifica que lo anterior se cumple para $\Gamma_2 = \{0\} \times [0,\infty)$. Se concluye entonces que el sistema \eqref{2.1} es uniformemente persistente.

\end{proof}
Si $n > 1$, entonces también es posible mostrar que el sistema \eqref{2.1} es persistente a través de un procedimiento análogo. 
\section{Dinámica alrededor del punto de equilibrio positivo}
\label{seccion2.2}
En esta sección se presentan las propiedades principales del punto de equilibrio positivo $E^* = (x^*,y_1^*,y_2^*,\hdots,y_n^*)$. Además, se analiza como varía la coordenada que representa la población del recurso con respecto a la intensidad de la competencia.

En lo que sigue, los vectores positivos se denotan con $\boldsymbol y = (y_1,y_2,\hdots,y_n)$, 
$\boldsymbol q = (q_1,q_2,\hdots,q_n)$, $\boldsymbol \mu = (\mu_1,\mu_2,\hdots,\mu_n)$, $\boldsymbol c = (c_1,c_2,\hdots,c_n)$, $\boldsymbol \kappa = (c_1q_1,c_2q_2,\hdots,c_nq_n)$. La matriz de competencia se denota con
\begin{eqnarray}
M &\colon =& \left( \begin{array}{cccc}
m_{11} & m_{12} & \hdots & m_{1n} \\
m_{21} & m_{22} & \hdots & m_{2n} \\
\vdots & \vdots & \ddots & \vdots \\
m_{n1} & m_{n2} & \hdots & m_{nn} 
\end{array} \right), \nonumber
\end{eqnarray}
y se asume que
\begin{eqnarray}
\det(M) &\neq& 0. \nonumber
\end{eqnarray}
\begin{df}
El \textbf{vector de balance competitivo} denotado con $\boldsymbol \alpha = (\alpha_1,\alpha_2,\hdots,\alpha_n)$ está dado por
\begin{eqnarray}
\boldsymbol \alpha^T &=& M^{-1} \boldsymbol q^T. \nonumber
\end{eqnarray}
\end{df}
Es claro que $\boldsymbol \alpha$ es la solución al sistema lineal
\begin{eqnarray}
m_{11}\alpha_1 + m_{12}\alpha_2 + \hdots + m_{1n}\alpha_n &=& q_1, \nonumber \\
m_{21}\alpha_1 + m_{22}\alpha_2 + \hdots + m_{2n}\alpha_n &=& q_2, \nonumber \\
                                                      &\vdots&     \nonumber \\
m_{n1}\alpha_1 + m_{n2}\alpha_2 + \hdots + m_{nn}\alpha_n &=& q_n. \nonumber 
\end{eqnarray}
Obsérvese que puede ocurrir $\alpha_j < 0$ para alguna $j\in \{1,\hdots,n\}$.

A continuación se introduce la función
\begin{eqnarray}
h(x) &=& \boldsymbol q M^{-1}\left(x  \boldsymbol \kappa^T -  \boldsymbol \mu^T\right) - g(x), \quad \text{con }x\in [0,K].  \nonumber
\end{eqnarray}
\begin{obs}
\label{obs2.9}
$E^* = (x^*,y_1^*,y_2^*,\hdots,y_n^*)$ es un punto de equilibrio positivo del sistema \eqref{2.1} si y solo si las igualdades
\begin{eqnarray}
h(x^*) &=& 0,                                                                     \nonumber   \\
(\boldsymbol y^*)^T &=& M^{-1}\left(x^* \boldsymbol \kappa^T - \boldsymbol \mu^T\right)\in (0,\infty)^n.  \nonumber
\end{eqnarray}
se cumplen.
\end{obs}
En el Teorema \ref{teo3.7} se muestran condiciones para la existencia de un punto de equilibrio positivo para el modelo \eqref{3.1} estudiado en el Capítulo \ref{capitulo3}. De hecho, el modelo \eqref{2.1} es el caso particular de \eqref{3.1} con $a_1 = a_2 = 0$. 
\begin{df} 
El \textbf{conjunto diagonal} $\Delta$ del espacio de matrices está dado por
\begin{eqnarray}
\Delta &=& \{M\in \mathcal M_{n\times n} \colon m_{ii}>0\text{, and }m_{ij}=0 \text{ si }i\neq j\}. \nonumber
\end{eqnarray}
\end{df}
\begin{obs}
\label{obs2.11}
Si $M\in \Delta$, entonces $\boldsymbol \alpha$ pertenece a $(0,\infty)^n$ y
\begin{eqnarray}
h'(x) &=& \left(\sum_{i=1}^n m_{ii}^{-1}c_iq_i^2\right) - g'(x) > 0\text{, para toda } x\in[0,K]. \nonumber
\end{eqnarray}
\end{obs}
\begin{obs}
\label{obs2.12}
Si $M\in \Delta$, entonces $M^{-1}\in \Delta$ y por lo tanto la correspondiente función $h$ satisface
\begin{eqnarray}
h(0) \quad < \quad 0 \quad < \quad h(K). \nonumber
\end{eqnarray} 
Entonces se tiene por el teorema del valor intermedio que la primera igualdad de la Observación \ref{obs2.9} se cumple para algún valor en $(0,K)$.
\end{obs}
\begin{obs}
\label{obs2.13}
Para toda $\widetilde M\in \Delta$, existe $\varepsilon_1(\widetilde M)>0$ tal que la función correspondiente $h$ satisface la cadena de desigualdades en la Observación \ref{obs2.12} para toda $M\in \mathcal M_{n\times n}$ con $\|M-\widetilde M\| < \varepsilon_1(\widetilde M)$.
\end{obs}
\subsection{Efectos de la competencia en la población del recurso}
Dada una matriz cuadrada de $n \times n$
\begin{eqnarray}
A &=& (a_{ij}), \nonumber
\end{eqnarray}
se denota con $A^{ij}$ a la matriz de $(n-1)\times (n-1)$ que resulta de la eliminación de la $i$-ésima fila y la $j$-ésima columna de $A$. Además, el menor de $A$ correspondiente a la $k$-ésima fila y la $\ell$-ésima columna será denotado con $|A|_{k\ell}$. Obsérvese que
\begin{eqnarray}
\label{2.2}
\det(A) &=& \sum_{k=1}^n (-1)^{k+j}a_{kj}|A|_{kj} = \sum_{\ell =1}^n (-1)^{i+\ell}a_{i\ell}|A|_{i\ell}, \quad \text{para toda }i,j\in\{1,\hdots,n\}.
\end{eqnarray}
Las siguientes matrices son estudiadas en los Lemas \ref{lema2.14} y \ref{lema2.15}.

\begin{scriptsize}
\begin{eqnarray}
\overline A_{ij} = & & \nonumber
\end{eqnarray}
\begin{eqnarray}
\left( \begin{array}{cccccccc}
|A^{ij}|_{11}           & -|A^{ij}|_{21}                  & \hdots 
                        & (-1)^i|A^{ij}|_{(i-1)1}         & 0
                        & (-1)^i|A^{ij}|_{(i+1)1}         & \hdots
                        & (-1)^{n+1}|A^{ij}|_{n1}                  \\
-|A^{ij}|_{12}          & |A^{ij}|_{22}                   & \hdots 
                        & (-1)^{i+1}|A^{ij}|_{(i-1)2}     & 0
                        & (-1)^{i+1}|A^{ij}|_{(i+1)2}     & \hdots
                        & (-1)^n|A^{ij}|_{n2}                      \\
\vdots                  & \vdots                          & \ddots 
                        & \vdots                          & \vdots 
                        & \vdots                          & \ddots 
                        & \vdots                                   \\
(-1)^j|A^{ij}|_{1(j-1)} & (-1)^{j+1}|A^{ij}|_{2(j-1)}     & \hdots 
                        & (-1)^{i+j}|A^{ij}|_{(i-1)(j-1)} & 0      
                        & (-1)^{i+j}|A^{ij}|_{(i+1)(j-1)} & \hdots 
                        & (-1)^{n+j+1}|A^{ij}|_{n(j-1)}            \\
0                       & 0                               & \ddots 
                        & 0                               & 0
                        & 0                               & \ddots 
                        & 0                                      \\\\
(-1)^j|A^{ij}|_{1(j+1)} & (-1)^{j+1}|A^{ij}|_{2(j+1)}     & \hdots 
                        & (-1)^{i+j}|A^{ij}|_{(i-1)(j+1)} & 0      
                        & (-1)^{i+j}|A^{ij}|_{(i+1)(j+1)} & \hdots 
                        & (-1)^{n+j+1}|A^{ij}|_{n(j+1)}            \\
\vdots                  & \vdots                          & \ddots 
                        & \vdots                          & \vdots 
                        & \vdots                          & \ddots 
                        & \vdots                                   \\
(-1)^{n+1}|A^{ij}|_{1n} & (-1)^{n+2}|A^{ij}|_{2n}         & \hdots 
                        & (-1)^{i+n+1}|A^{ij}|_{(i-1)n}   & 0      
                        & (-1)^{i+n+1}|A^{ij}|_{(i+1)n}   & \hdots 
                        & |A^{ij}|_{nn}
\end{array} \right), \nonumber
\end{eqnarray}
\end{scriptsize}

\begin{footnotesize}
\begin{eqnarray}
\overline{\overline A}_{ij} = & & \nonumber
\end{eqnarray}
\begin{eqnarray}
\underset{j\text{-ésima fila}\longrightarrow}{}
\left( \begin{array}{cccccccc}
(-1)^{i+j}|A|_{ij} & 0                  & \hdots 
                   & 0                  & \overbrace{(-1)^i|A|_{i1}}^{j\text{-ésima columna}} 
                   & 0                  & \hdots   
                   & 0                                                \\
0                  & (-1)^{i+j}|A|_{ij} & \hdots 
                   & 0                  & (-1)^{i+1}|A|_{i2}  
                   & 0                  & \hdots   
                   & 0                                                \\
\vdots             & \vdots             & \ddots 
                   & \vdots             & \vdots                                   
                   & \vdots             & \ddots
                   & \vdots                                           \\
0                  & 0                  & \hdots 
                   & (-1)^{i+j}|A|_{ij} & (-1)^{i+j}|A|_{i(j-1)}           
                   & 0                  & \hdots   
                   & 0                                                \\
0                  & 0                  & \ddots 
                   & 0                  & 0
                   & 0                  & \ddots   
                   & 0                                              \\\\
0                  & 0                  & \hdots 
                   & 0                  & (-1)^{i+j} |A|_{i(j+1)}
                   & (-1)^{i+j}|A|_{ij} & \hdots   
                   & 0                                                 \\
\vdots             & \vdots             & \ddots 
                   & \vdots             & \vdots
                   & \vdots             & \ddots   
                   & \vdots                                            \\
0                  & 0                  & \hdots 
                   & 0                  & (-1)^{i+n+1}|A|_{in}
                   & 0                  & \hdots   
                   & (-1)^{i+j}|A|_{ij} 
\end{array} \right). \nonumber
\end{eqnarray}
\end{footnotesize}

\begin{lema}
\label{lema2.14}
 $\overline A_{ij} A = \overline{\overline A}_{ij}$.
\end{lema}
La demostración se sigue de la ecuación \eqref{2.2} y del hecho de que intercambiar filas en una matriz cambia el signo del determinante.
\begin{lema}
\label{lema2.15}
Si $\boldsymbol u$, $\boldsymbol v \in \mathbb R^n$ y $A^*$ representa a la matriz adjunta de $A$, entonces $\displaystyle \frac{\partial}{\partial a_{ij}} \boldsymbol u A^* \boldsymbol v^T = \boldsymbol u \overline A_{ij} \boldsymbol v^T$.
\end{lema}
\begin{proof}
Es fácil ver que
\begin{eqnarray}
A^* \quad = \quad \det(A)(A^{-1}) &=& \left( \begin{array}{cccc}
|A|_{11}           & -|A|_{21}      & \hdots & (-1)^{n+1}|A|_{n1} \\
-|A|_{12}          & |A|_{22}       & \hdots & (-1)^n|A|_{n2}     \\
\vdots             & \vdots         & \ddots & \vdots             \\
(-1)^{n+1}|A|_{1n} & (-1)^n|A|_{2n} & \hdots & |A|_{nn}
\end{array} \right). \nonumber
\end{eqnarray}
Por lo tanto $\displaystyle \frac{\partial}{\partial a_{ij}} A^* = \overline A_{ij}$ y el resultado se sigue.
\end{proof}
El Teorema \ref{teo2.16} establece que la razón de cambio de la población de equilibrio del recurso con respecto al coeficiente de competencia es directamente proporcional a la población de equilibrio de la especie competidora, así como al balance competitivo de la especie que recibe esta competencia. Obsérvese que la población de equilibrio del recurso $x^*(M)$ está definida implícitamente (por la primera ecuación de la Observación \ref{obs2.9} y el teorema de la función implícita) en un subconjunto del espacio vectorial $\mathcal M_{n\times n}$. Nótese que a pesar de que la función $h$ está dada en términos de $x$, depende también de los parámetros $m_{ij}$. Para evitar ambigüedad, se usa la notación $h'$ para referirse a $\partial h/(\partial x)$.
\begin{teo}
\label{teo2.16}
Si $E^*=(x^*,y_1^*,y_2^*,\hdots,y_n^*)$ es un punto de equilibrio positivo de \eqref{2.1}, 
entonces
\begin{eqnarray}
\frac{\partial x^*}{\partial m_{ij}} &=& \frac{\alpha_i y_j^*}{\det(M)h'(x^*)}. \nonumber
\end{eqnarray}
\end{teo}
\begin{proof}
De los Lemas \ref{lema2.14} y \ref{lema2.15} se tiene que
\begin{eqnarray}
\frac{\partial h}{\partial m_{ij}} &=& \frac{\partial}{\partial m_{ij}}
                                         \left(\boldsymbol q M^{-1}\left(x^*\boldsymbol \kappa^T -    
                                         \boldsymbol \mu^T\right)\right)              \nonumber \\
                                     &=& \frac{\partial}{\partial m_{ij}}
                                         \left(\boldsymbol q\frac{M^*}{\det(M)}\left(x^*\boldsymbol \kappa^T -
                                         \boldsymbol \mu^T\right)\right)              \nonumber \\
                                     &=& \frac{1}{\det(M)}\left(\boldsymbol q \overline M_{ij}(x^*\boldsymbol \kappa^T -
                                        \boldsymbol \mu^T)-|M|_{ij}\boldsymbol q M^{-1}(x^*\boldsymbol \kappa^T-\boldsymbol \mu^T)\right)
                                                                        \nonumber \\
                                     &=& \frac{\left(\boldsymbol q\overline{\overline M}_{ij}-|M|_{ij}\boldsymbol q
                                         \right)\boldsymbol y^*}{\det(M)}               \nonumber \\
                                     &=& -\frac{\alpha_iy_j^*}{\det(M)}.    \nonumber
\end{eqnarray}
Del teorema de la función implícita se sigue que
\[\frac{\partial x^*}{\partial m_{ij}} = -\frac{\partial h}{\partial m_{ij}}
/\left(\frac{\partial h}{\partial x^*}\right) = \frac{\alpha_i y_j^*}{\det(M)h'(x^*)}.\]
\end{proof}
La Observación \ref{obs2.11} y el Teorema \ref{teo2.16} dan lugar a lo siguiente. 
\begin{prop}
\label{prop2.17}
Si $M\in \Delta$, entonces
\begin{eqnarray}
\frac{\partial x^*}{\partial m_{ij}}(M) &>& 0. \nonumber
\end{eqnarray}
\end{prop}
En este caso, la tasa de crecimiento de la población de equilibrio del recurso respecto a la competencia es positiva. Se concluye que si los coeficientes de competencia interespecífica son igual a cero, entonces el incremento de la competencia (intra o interespecífica) de los consumidores genera un beneficio indirecto para el recurso  en \eqref{2.1}.
\begin{obs}
\label{obs2.18}
Para toda $\widetilde M\in \Delta$, existe $\varepsilon_2(\widetilde M)>0$ tal que para $M\in \mathcal 
M_{n\times n}$, $\|M-\widetilde M\| < \varepsilon_2(\widetilde M)$ implica que  $x^*(M)$ satisface la desigualdad de la Proposición \ref{prop2.17}.
\end{obs}
\subsection{Estabilidad asintótica}
Para analizar las propiedades dinámicas del sistema \eqref{2.1}, se introduce la transformación lineal $L_{\boldsymbol c}\colon \mathcal M_{n\times n}\longrightarrow \mathcal M_{n\times n}$ dada por
\begin{eqnarray}
L_{\boldsymbol c}(M)  &=&  \left( \begin{array}{cccc}
c_1^{-1}m_{11}                     & (c_1^{-1}m_{12} + c_2^{-1}m_{21})/2  & \hdots   
& (c_1^{-1}m_{1n} + c_n^{-1}m_{n1})/2 \\
(c_1^{-1}m_{12} + c_2^{-1}m_{21})/2 & c_2^{-1}m_{22}                      & \hdots   
& (c_2^{-1}m_{2n} + c_n^{-1}m_{n2})/2 \\
\vdots                             & \vdots                               & \ddots  
& \vdots                             \\
(c_1^{-1}m_{1n} + c_n^{-1}m_{n1})/2 & (c_2^{-1}m_{2n} + c_n^{-1}m_{n2})/2 & \hdots   
& c_n^{-1}m_{nn}
\end{array} \right). \nonumber
\end{eqnarray}
\begin{teo}
\label{teo2.19}
Si $L_{\boldsymbol c}(M)$ es positiva definida, entonces $E^*$ es atractor global en \eqref{2.1}.
\end{teo}
\begin{proof}
Se reescribe \eqref{2.1} de la siguiente manera
\begin{eqnarray}
\dot x   &=&   x\left(g(x)-g(x^*)      - q_1   (y_1-y_1^*) - q_2   (y_2-y_2^*) - 
\hdots - q_n   (y_n-y_n^*)\right), \nonumber \\
\dot y_1 &=& y_1\left(c_1q_1(x-x^*)    - m_{11}(y_1-y_1^*) - m_{12}(y_2-y_2^*) - 
\hdots - m_{1n}(y_n-y_n^*)\right), \nonumber \\
\dot y_2 &=& y_2\left(c_2q_2(x-x^*)    - m_{21}(y_1-y_1^*) - m_{22}(y_2-y_2^*) - 
\hdots - m_{2n}(y_n-y_n^*)\right),         \nonumber  \\
       &\vdots&                   \nonumber \\
\dot y_n &=& y_n\left(c_nq_n(x-x^*)    - m_{n1}(y_1-y_1^*) - m_{n2}(y_2-y_2^*) - 
\hdots - m_{nn}(y_n-y_n^*)\right),\nonumber
\end{eqnarray}
y se considera la función de Lyapunov $V\colon (0,\infty)^{n+1} \longrightarrow [0,\infty)$ dada por
\begin{eqnarray}
V(x,y_1,y_2,\hdots,y_n) &=& \lambda_0 \left(x  -x^*  -x^*  \ln\left(\frac{x  }{x^*  }
                            \right)\right) + \lambda_1 \left(y_1-y_1^* - 
                            y_1^*\ln\left(\frac{y_1}{y_1^*}\right)\right) + \hdots 
                            \nonumber \\
                        & & \hdots + \lambda_n\left(y_n-y_n^* -
                            y_n^*\ln\left(\frac{y_n}{y_n^*}\right)\right), 
                            \hspace{73pt} \text{con }\lambda_0, \lambda_1,\hdots, 
                            \lambda_n > 0.                           \nonumber
\end{eqnarray}
La derivada de $V$ está dada por
\begin{eqnarray}
\frac{dV}{dt} &=& \left(x - x^*\right)\left[\left(g(x) - g(x^*)\right) + 
                  \left(\lambda_0 - \lambda_1c_1\right)q_1\left(y_1 - y_1^*\right) + 
                  \hdots +
                  \left(\lambda_0 - \lambda_nc_n\right)q_n\left(y_n - y_n^*\right)
                  \right] \nonumber \\ \nonumber \\
              && - \left(\boldsymbol y - \boldsymbol y^*\right)
\left( \begin{array}{cccc}
\lambda_1m_{11}                       &  
                               \hdots & (\lambda_1m_{1n} + \lambda_nm_{n1})/2 \\
(\lambda_1m_{12} + \lambda_2m_{21})/2 &  
                               \hdots & (\lambda_2m_{2n} + \lambda_nm_{n2})/2 \\
\vdots                                & 
                               \ddots & \vdots                               \\
(\lambda_1m_{1n} + \lambda_nm_{n1})/2 & 
                               \hdots & \lambda_nm_{nn}
\end{array} \right)
\left(\boldsymbol y - \boldsymbol y^*\right)^T. \nonumber
\end{eqnarray}
Si se escoge $\lambda_0 = 1$ y $\lambda_i = 1/c_i$, la ecuación anterior queda
\begin{eqnarray}
\frac{dV}{dt} &=& \left(x - x^*\right)\left(g(x) - g(x^*)\right) -\left(\boldsymbol y - \boldsymbol y^*\right) L_{\boldsymbol c}(M) \left(\boldsymbol y - \boldsymbol y^*\right)^T. \nonumber 
\end{eqnarray}
Del hecho de que $L_{\boldsymbol c}(M)$ tiene espectro positivo y que $g$ es decreciente, se concluye que $E^*$ es un atractor global en \eqref{2.1}.

\end{proof}

En el caso $n=2$, la desigualdad de las medias aritmética y geométrica da lugar a $((m_{12}/c_1 + m_{21}/c_2)/2)^2 \geq (m_{12}m_{21})/(c_1c_2)$. Así, se tiene que si $L_{\boldsymbol c}(M)$ es una matriz positiva definida, entonces $\det(M) > 0$. 

En general, parece complicado verificar si una matriz es positiva definida. Sin embargo, en el caso particular $n=2$, las condiciones de estabilidad para el punto de equilibrio positivo $E^*$ están dadas en términos de las medias aritmética y geométrica.
\begin{cor}
\label{cor2.20}
Considérese el sistema \eqref{2.1} con $n=2$. Si se cumple
\begin{eqnarray}
\frac{c_2m_{12} + c_1m_{21}}{2} &<& \sqrt{c_1m_{11}c_2m_{22}}, \nonumber
\end{eqnarray}
entonces $E^*=(x^*,y_1^*,y_2^*)$ es un atractor global.
\end{cor}
La desigualdad del Corolario \ref{cor2.20} a su vez implica las desigualdades
\begin{eqnarray}
\frac{c_2m_{12} + c_1m_{21}}{2} &<& \frac{c_1m_{11} + c_2m_{22}}{2}, \nonumber \\ \nonumber \\
\sqrt{c_2m_{12}c_1m_{21}} &<& \sqrt{c_1m_{11}c_2m_{22}}. \nonumber
\end{eqnarray}

Un escenario particular en el que el Teorema \ref{teo2.19} y el Corolario \ref{cor2.20} se pueden aplicar es el siguiente.
\begin{ej}
\label{ej2.21}
Considérese una comunidad donde dos consumidores ($y_1$, $y_2$) comparten un único recurso $x$. Supóngase además que las poblaciones y las interacciones de estas especies están modeladas por el sistema \eqref{2.1} (con $n=2$) y que la desigualdad del Corolario \ref{cor2.20} se cumple. Entonces la invasión de un nuevo competidor $y_3$ será exitosa si se cumple
\begin{eqnarray}
\label{2.3}
\frac{4m_{11}m_{22}}{c_1c_2} - \left(\frac{m_{12}}{c_1} 
+ \frac{m_{21}}{c_2}\right)^2                                                   &>& 
\frac{c_3}{m_{33}}\left[\frac{m_{11}}{c_1}\left(\frac{m_{23}}{c_2} 
+ \frac{m_{32}}{c_3}\right)^2 + \frac{m_{22}}{c_2}\left(\frac{m_{13}}{c_1} 
+ \frac{m_{31}}{c_3}\right)^2 \right.                                              \\
& & 
\left.-\left(\frac{m_{12}}{c_1} + \frac{m_{21}}{c_2}\right)\left(\frac{m_{13}}{c_1} 
+\frac{m_{31}}{c_3}\right)\left(\frac{m_{23}}{c_2} 
+ \frac{m_{32}}{c_3}\right)\right].                                    \nonumber
\end{eqnarray}
La desigualdad \eqref{2.3} se puede deducir del determinante de la matriz correspondiente de $3\times 3$ $L_{\boldsymbol c}(M)$ en el Teorema \ref{teo2.19}. Además, el Corolario \ref{cor2.20} implica que el lado izquierdo de \eqref{2.3} es un número positivo. De acuerdo con eso, para una invasión exitosa de $y_3$, es suficiente que la especie mantenga un perfil bajo con los otros competidores (de tal forma que $(m_{23}/c_2 + m_{32}/c_3)$ y $(m_{13}/c_1 + m_{31}/c_3)$ sean suficientemente pequeños en comparación con los otros parámetros de \eqref{2.3}), o que tenga un nivel de competencia intraespecífica suficientemente grande (es decir, $c_3/m_{33}$ suficientemente pequeña).
\end{ej}

El Ejemplo \ref{ej2.21} puede adaptarse a dimensiones más altas, sin embargo los cálculos se vuelven más largos.

De la Observación \ref{obs2.12} se sigue que existe un valor de equilibrio $x^* \in (0,K)$ si $M\in \Delta$. En este caso, existe un punto de equilibrio positivo $(x^*,y_1^*,\hdots,y_n^*)$ si y solo si
\[M^{-1}(x^*\boldsymbol \kappa^T - \boldsymbol \mu^T) \in (0,\infty)^n.\]
\begin{obs} 
\label{obs2.22}
Si $M\in \Delta$, entonces las siguientes tres condiciones son equivalentes
\begin{enumerate}
\item $\boldsymbol y^* \in \mathbb (0,\infty)^n$,
\item $m_{ii}^{-1}\left(c_iq_ix^* - \mu_i\right) > 0\text{ para toda }i\in\{1,\hdots,n\}$,
\item $x^* > x_i^b\text{ para toda }i\in\{1,\hdots,n\}$.
\end{enumerate}
\end{obs}

Se denota
\begin{eqnarray}
\widetilde \Delta &=& \{M\in \Delta \colon M \text{ cumple alguna de las condiciones de la Observación \ref{obs2.22}}\}    \nonumber
\end{eqnarray}
\begin{prop}
$\widetilde \Delta \neq \emptyset$
\end{prop}
\begin{proof}
Si $n=1$, entonces
\begin{eqnarray}
h(x_1^b) \quad < \quad 0 \quad < \quad h(K), \nonumber
\end{eqnarray}
donde $h$ está dada por 
\begin{eqnarray}
h(x) &=& \boldsymbol q M^{-1}\left(x  \boldsymbol \kappa^T -  \boldsymbol \mu^T\right) - g(x), \quad \text{con }x\in [0,K],  \nonumber
\end{eqnarray}
y el resultado se sigue del teorema del valor intermedio.

Para el caso $n>1$, se escribe
\begin{eqnarray}
x_N^b &=& \max\left\{x_i^b\right\}_{i=1}^n,          \nonumber \\
\gamma     &=& \max\left\{c_iq_i(x_N^b - x_i^b)\right\}_{i=1}^{n-1}. \nonumber
\end{eqnarray}
Entonces, se escoge $M\in \Delta$, con
\[m_{ii}>\frac{(n-1)\gamma}{g(x_N^b)}\text{, para toda }i\in\{1,\hdots,n-1\}.\]
Por lo tanto $h(x_N^b)<0$, y se procede como en el primer caso.

\end{proof}
\begin{obs}
\label{obs2.24}
Para toda $\widetilde M \in \widetilde \Delta$, existe $\varepsilon_3(\widetilde M) > 0$, tal que $M\in \mathcal M_{n\times n}$, $\|M - \widetilde M\| < \varepsilon_3(\widetilde M)$ implica que $M^{-1}\left( x^* \boldsymbol \kappa^T - \boldsymbol \mu^T \right)$ es un vector positivo en $\mathbb R^n$ y $\min \sigma\left(L_{\boldsymbol c}(M)\right) > 0$, donde $\sigma$ denota al espectro de $M$. 
\end{obs}
Los resultados principales del Capítulo \ref{capitulo2} (la Proposición \ref{prop2.17} y el Teorema \ref{teo2.19}) se pueden generalizar para las matrices de competencia $M$ de un subconjunto abierto y conexo de $\mathcal M_{n\times n}$.
\begin{teo}[Reyes-García]
\label{teo2.25}
Existe un conjunto abierto y conexo $U\subset \mathcal M_{n\times n}$ de tal forma que $\widetilde \Delta\subset U$ y si $M \in U$, entonces un sistema \eqref{2.1} asociado a $M$ tiene un punto de equilibrio positivo $E^*$ que es un atractor global y además $\partial x^*/ (\partial m_{ij})(M) > 0$.
\end{teo}
\begin{proof}
Para toda $M\in \widetilde \Delta$, se escoge $\varepsilon(M) = \min\left\{\varepsilon_1(M), \varepsilon_2(M), 
\varepsilon_3(M)\right\}$, donde $\varepsilon_i(M)$ denota al valor positivo de las Observaciones \ref{obs2.13}, \ref{obs2.18} y \ref{obs2.24} respectivamente. Entonces se introduce el conjunto
\begin{eqnarray}
U &=& \bigcup \left\{B\left(M,\varepsilon(M)\right)\colon M\in \widetilde 
\Delta\right\} \nonumber
\end{eqnarray}
De las Observaciones \ref{obs2.9}, \ref{obs2.12}, \ref{obs2.13} y \ref{obs2.22}, se sigue que existe un punto de equilibrio positvo $E^*$. Posteriormente, $\partial x^*/(\partial m_{ij})(M) > 0$ se obtiene por la Observación \ref{obs2.18}. Finalmente, por la Observación \ref{obs2.24} se tiene que $E^*$ es un atractor global.

\end{proof}

Para dar un sistema de ecuaciones diferenciales \eqref{2.1} con un estado estacionario de coexistencia $E^* = (x^*,y_1^*,\hdots,y_n^*)$ que sea atractor global y tal que el valor estacionario de la población del recurso $x^*$ sea creciente respecto a la competencia $m_{ij}$, conviene considerar las matrices de competencia que pertenecen al conjunto $U \subset \widetilde \Delta$.
\section{Simulaciones numéricas}
\label{seccion2.3}
A continuación se realizan simulaciones que ilustran algunas de las características del sistema \eqref{2.1}. Considérese el sistema \eqref{2.1} con $n=2$ (dos consumidores alimentándose de un recurso), y supóngase que todos los parámetros con excepción de $m_{12}$ y $m_{21}$ están fijos: los valores que se escogen para los parámetros son $r=3$, $K=8$, $c_1=1$, $q_1=4$, $\mu_1=3$, $m_{11}=1.5$, $c_2=1.2$, $q_2=2$, $\mu_2=1$ y $m_{22}=2$; las condiciones iniciales están dadas por $x(0)=1$, $y_1(0)=0.8$, $y_2(0)=0.6$. Se mostrará (para $g$ dada por el modelo logístico) que existe un subconjunto $\Omega\subset [0,\infty)^2$ tal que para $(m_{12},m_{21})\in \Omega$, la matriz correspondiente $M$ pertenece a $U$ (ver el Teorema \ref{teo2.25}) y por lo tanto el sistema \eqref{2.1} tiene un punto de equilibrio $E^*$ que es un atractor global. Además se muestra que se cumple la desigualdad
\begin{eqnarray}
\frac{\partial x^*}{\partial m_{ij}}(M) &>& 0. \nonumber
\end{eqnarray}

Para hallar un punto de equilibrio positivo, se resuelve el siguiente sistema lineal
\begin{eqnarray}
3(1 - x^*/8) - 4y_1^* - 2y_2^*    &=& 0, \nonumber \\
4x^* - 1.5y_1^* - m_{12}y_2^* - 3 &=& 0, \nonumber \\
2.4x^* - m_{21}y_1^* - 2y_2^* - 1 &=& 0, \nonumber
\end{eqnarray}
para $x^*$, $y_1^*$, $y_2^*$.

Se sigue que
\begin{eqnarray}
x^*   &=& \frac{8(36 - 4m_{12} - 6m_{21} - 3m_{12}m_{21})}{322.6 - 76.8m_{12} - 64m_{21} - 3m_{12}m_{21}} \nonumber \\
y_1^* &=& \frac{122.8 - 54.6m_{12}}{322.6 - 76.8m_{12} - 64m_{21} - 3m_{12}m_{21}},                     \nonumber \\
y_2^* &=& \frac{184.3 - 87m_{21}}{322.6 - 76.8m_{12} - 64m_{21} - 3m_{12}m_{21}}.                       \nonumber
\end{eqnarray}
Es claro que $y_1^*$ tiende a cero si
\begin{eqnarray}
m_{12} &\to& 2.249, \nonumber
\end{eqnarray}
y que $y_2^*$ tiende a cero si
\begin{eqnarray}
m_{21} &\to& 2.1183. \nonumber
\end{eqnarray}
Por lo tanto $(y_1^*,y_2^*)$ es una pareja ordenada positiva si 
\begin{eqnarray}
(m_{12},m_{21}) \in \Omega_1 &\colon =& [0,2.249)\times [0,2.1183). \nonumber
\end{eqnarray}
Por otro lado, el Corolario \ref{cor2.20} implica que $L_{\boldsymbol c}(M)$ es una matriz positiva definida si $(m_{12},m_{21})$ pertenece al conjunto
\begin{eqnarray}
\Omega_2 &\colon =& \left\{(m_{12},m_{21})\in[0,\infty)^2\colon \frac{1.2m_{12} + m_{21}}{2} < \sqrt3.6\right\}. \nonumber
\end{eqnarray}
Finalmente, las derivadas parciales $\displaystyle \frac{\partial x^*}{\partial m_{12}}$, $\displaystyle \frac{\partial x^*}{\partial m_{21}}$ son positivas si 
\begin{eqnarray}
(m_{12},m_{21})\in \Omega_3 &\colon =& [0,0.75)\times [0,2.1183). \nonumber
\end{eqnarray}
Nótese que $\Omega_3 \subset \Omega_1$. Del análisis anterior, se concluye que si $(m_{12},m_{21})$ es un elemento de $\Omega\colon = \Omega_1\cap \Omega_2\cap \Omega_3 = \Omega_1\cap \Omega_2$, entonces la matriz de competencia correspondiente $M$ pertenece al conjunto $U$ del Teorema \ref{teo2.25}.

\begin{figure}
\centering
\vspace{.5cm}
\includegraphics[width=242pt]{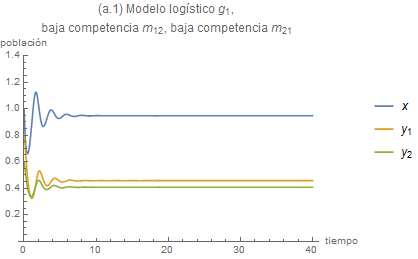}
\includegraphics[width=242pt]{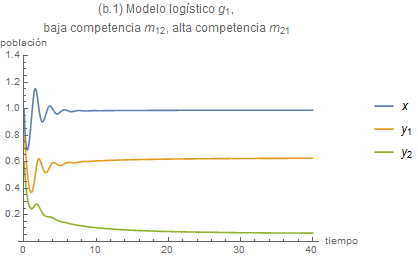}
\vspace{.5cm}
\includegraphics[width=242pt]{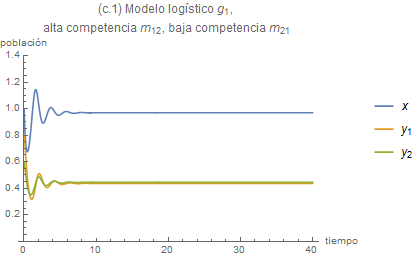}
\includegraphics[width=242pt]{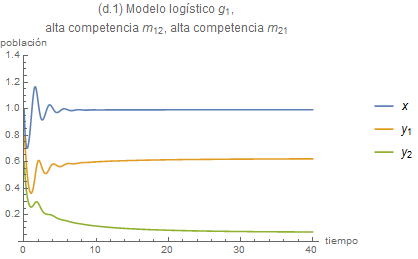}
\caption{Soluciones del sistema \eqref{2.1} con el modelo logístico como función de crecimiento y parámetros de competencia interespecífica $m_{12} = 0.25$, $0.5$; $m_{21} = 1$, $2$.} 
\label{fig2.1}
\end{figure}

\begin{figure}
\centering
\vspace{.5cm}
\includegraphics[width=242pt]{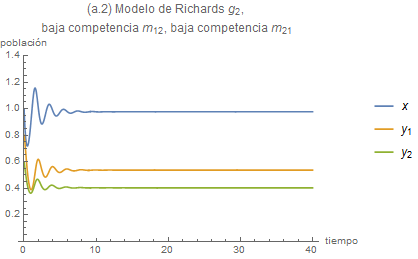}
\includegraphics[width=242pt]{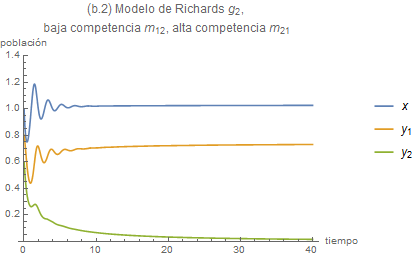}
\vspace{.5cm}
\includegraphics[width=242pt]{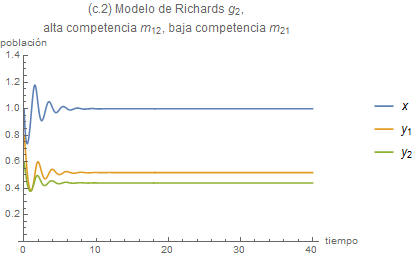}
\includegraphics[width=242pt]{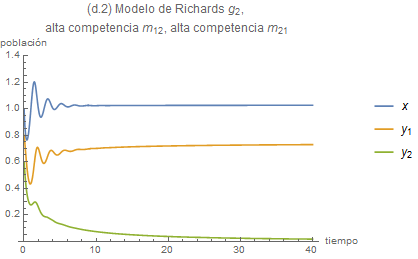}
\caption{Soluciones del sistema \eqref{2.1} con el modelo de Richards como función de crecimiento y parámetros de competencia interespecífica $m_{12} = 0.25$, $0.5$; $m_{21} = 1$, $2$.}
\label{fig2.2}
\end{figure}

\begin{figure}
\centering
\vspace{.5cm}
\includegraphics[width=242pt]{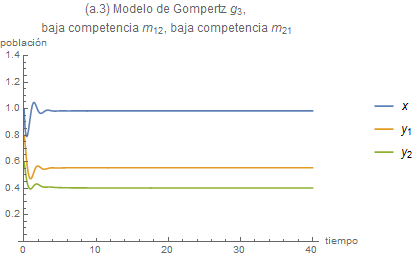}
\includegraphics[width=242pt]{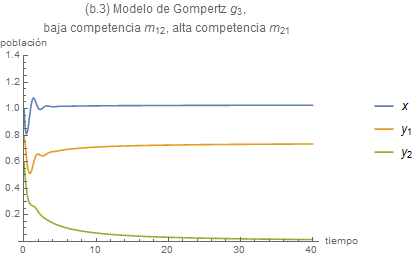}
\vspace{.5cm}
\includegraphics[width=242pt]{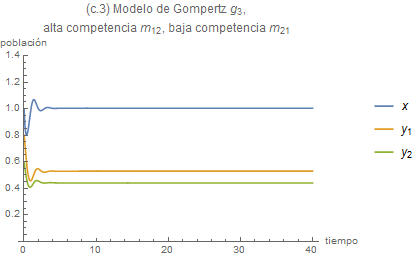}
\includegraphics[width=242pt]{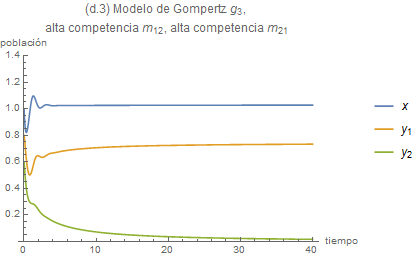}
\caption{Soluciones del sistema \eqref{2.1} con el modelo de Gompertz como función de crecimiento y parámetros de competencia interespecífica $m_{12} = 0.25$, $0.5$; $m_{21} = 1$, $2$.}
\label{fig2.3}
\end{figure}

\begin{figure}
\centering
\vspace{.5cm}
\includegraphics[width=242pt]{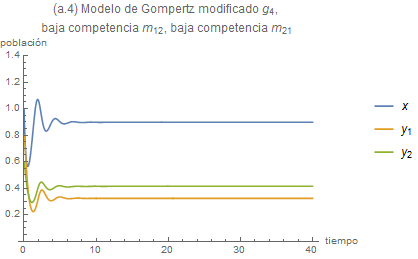}
\includegraphics[width=242pt]{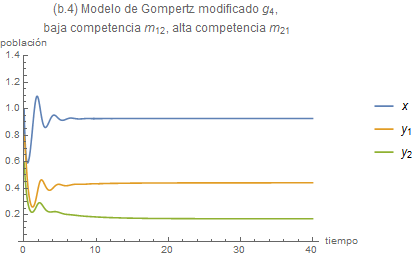}
\vspace{.5cm}
\includegraphics[width=242pt]{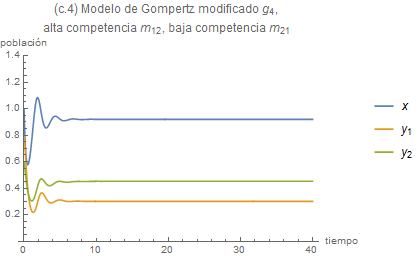}
\includegraphics[width=242pt]{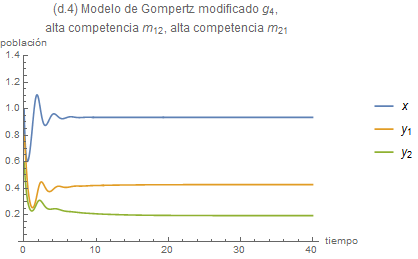}
\caption{Soluciones del sistema \eqref{2.1} con el modelo de Gompertz modificado como función de crecimiento y parámetros de competencia interespecífica $m_{12} = 0.25$, $0.5$; $m_{21} = 1$, $2$.}
\label{fig2.4}
\end{figure}

\begin{figure}
\centering
\vspace{.5cm}
\includegraphics[width=242pt]{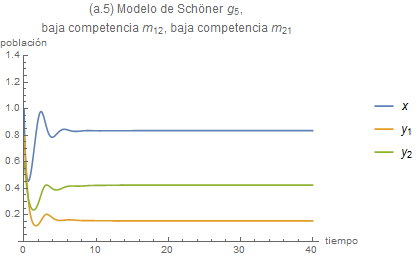}
\includegraphics[width=242pt]{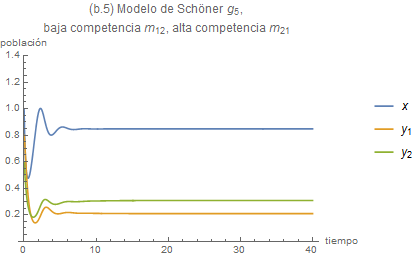}
\vspace{.5cm}
\includegraphics[width=242pt]{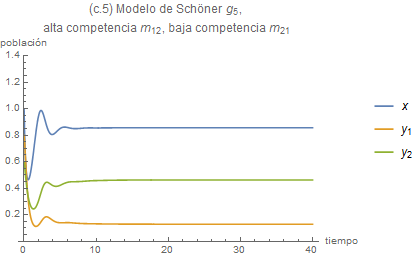}
\includegraphics[width=242pt]{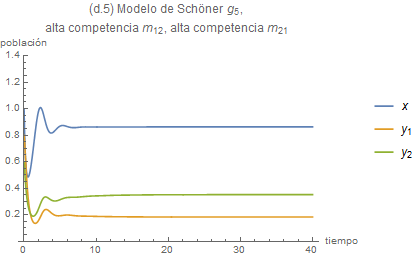}
\caption{Soluciones del sistema \eqref{2.1} con el modelo de Schöner como función de crecimiento y parámetros de competencia interespecífica $m_{12} = 0.25$, $0.5$; $m_{21} = 1$, $2$.}
\label{fig2.5}
\end{figure}

Los resultados presentados para el sistema \eqref{2.1} en la Sección \ref{seccion2.2} son robustos en el sentido de que se cumplen para cualquier función decreciente que satisfaga $g(K) = 0$. En las Figuras \ref{fig2.1}, \ref{fig2.2}, \ref{fig2.3}, \ref{fig2.4} y \ref{fig2.5} se comparan las soluciones de \eqref{2.1} para cuatro conjuntos de valores para los parámetros de competencia interespecífica en $\Omega$, y cinco diferentes funciones de crecimiento. Nótese que para algunas funciones $g$, no existe un punto de equilibrio $E^*$ en \eqref{2.1}, a pesar de que $(m_{12},m_{21})$ pertenece a $\Omega$.

En las gráficas de las Figuras \ref{fig2.1}-\ref{fig2.5}, los coeficientes de competencia interespecífica toman los valores $m_{12}\in\left\{0.25,0.5\right\}$; $m_{21}\in\left\{1,2\right\}$; las funciones de crecimiento están dadas por:
\begin{enumerate}
    \item Modelo logístico: $g_1(x) = 1 - x/K$; 
    \item Modelo de Richards (ver \cite{Mueller1981}, \cite{Nelder1961}): $g_2(x) = 1 - (x/K)^2$; 
    \item Modelo de Gompertz (ver \cite{Jiang2014OnModels}, \cite{Liu2017OnSystem}): $g_3(x) = (\ln K)^{-1} \ln (K/x)$; 
    \item Modelo de Gompertz modificado: $g_4(x) = \left(\ln(K + 1)\right)^{-1}\ln\left((K + 1)/(x + 1)\right)$;
    \item Modelo de Sch\"oner (ver \cite{Gilpin1976SchoenersCompetition}): $g_5(x) = K^{-1}\left((K + 1)/(x + 1) - 1\right)$.
\end{enumerate}

\begin{table}[h]
\centering
\begin{tabular}{||c||c|c|c|c|c||}\hline 
                    & (1) Logístico & (2) Richards & (3) Gompertz & (4) Gompertz & (5) Sch\"oner \\ 
                    &               &              &              & modificado   &               \\ \hline\hline
(a) $m_{12}=0.25$,  & $(0.947,$     & $(0.9766,$   & $(0.9832,$   & $(0.8974,$   & $(0.8343,$    \\ 
$m_{21}=1$          & $0.4574,$     & $0.5371,$    & $0.5548,$    & $0.3239,$    & $0.1541,$     \\
                    & $0.4077)$     & $0.4033)$    & $0.4024)$    & $0.4149)$    & $0.424)$      \\ \hline
(b) $m_{12}=0.25$,  & $(0.9891,$    & no           & no           & $(0.926,$    & $(0.8475,$    \\
$m_{21}=2$          & $0.6277,$     & existe       & existe       & $0.4412,$    & $0.2087,$     \\
                    & $0.0592)$     &              &              & $0.17)$      & $0.3082)$     \\ \hline
(c) $m_{12}=0.5$,   & $(0.9694,$    & $(0.9993,$   & $(1.0033,$   & $(0.9194,$   & $(0.8563,$    \\
$m_{21}=1$          & $0.4367,$     & $0.5182,$    & $0.529,$     & $0.3009,$    & $0.1292,$     \\
                    & $0.4449)$     & $0.4401)$    & $0.4394)$    & $0.4528)$    & $0.4629)$     \\ \hline
(d) $m_{12}=0.5$,   & $(0.9922,$    & no           & no           & $(0.9346,$   & $(0.8627,$    \\
$m_{21}=2$          & $0.6234,$     & existe       & existe       & $0.4278,$    & $0.1831,$     \\
                    & $0.0672)$     &              &              & $0.1937)$    & $0.352)$      \\ \hline
\end{tabular}
\caption{Coordenadas del punto de equilibrio positivo $E^* = (x^*,y_1^*,y_2^*)$.} 
\label{cuadro2.1}
\end{table}
Nótese que las poblaciones de equilibrio $x^*$, $y_1^*$ son llevadas a una ubicación más baja para las funciones de crecimiento $g_4$ y $g_5$, en comparación con los casos $g_1$, $g_2$ y $g_3$. Además, el sistema \eqref{2.1} con los modelos de crecimiento $g_2$ y $g_3$, y parámetro de competencia interespecífica $m_{21} = 2$ no tiene un punto de equilibrio $E^*$ en el retrato fase. En estos casos, el punto de equilibrio $E_1$ es un atractor para el sistema \eqref{2.1}.

Los puntos de equilibrio en la frontera $E_{00}$ y $E_0$ no dependen de la función de crecimiento ni de los parámetros $m_{12}$ y $m_{21}$, mientras que los puntos de equilibrio $E_1$ y $E_2$ cambian si la función de crecimiento cambia, pero permanecen estáticos si los parámetros $m_{12}$ y $m_{21}$ son modificados. De acuerdo con el Teorema \ref{teo2.25} el valor estacionario $x^*$ crece si $m_{12}$ y $m_{21}$ aumentan, y entonces el punto de equilibrio positivo $E^*$ varía. A través de un cálculo directo se puede conocer la posición de los puntos de equilibrio positivo para cada caso de las Figuras \ref{fig2.1}-\ref{fig2.5}. Esto se muestra en el Cuadro \ref{cuadro2.1}.

\section{Conclusiones}
En este capítulo, el modelo de competencia y depredación \eqref{2.1} es introducido y analizado para profundizar en la pregunta sobre el papel de la competencia intra- e interespecífica en la estabilidad de la comunidad (ver por ejemplo \cite{Kuang2003BiodiversityCompetition}). La estabilidad asintótica para el punto de equilibrio positivo en \eqref{2.1} fue demostrada en el Teorema \ref{teo2.19} a través de una función clásica de Lyapunov introducida por Goh en \cite{Goh1978}. Hasta donde se sabe, una demostración para un modelo con coeficientes de competencia interespecífica no nulos y una función de crecimiento general no ha sido dada antes. Los resultados del capítulo muestran que la coexistencia de las especies es factible siempre que un equilibrio competitivo entre los parámetros intra- e interespecífico es alcanzado; en particular véase el Corolario \ref{cor2.20}, donde las condiciones para estabilidad global han sido establecidas en términos de medias geométricas y aritméticas que involucran los parámetros de competencia. En el Ejemplo \ref{ej2.21}, la invasión de un depredador exótico hacia un entorno que está ocupado por una presa nativa y dos depredadores es estudiado. En este caso, la condición para una invasión exitosa (desigualdad \eqref{2.3}) está expresada en términos de medias aritméticas y geométricas de los parámetros de competencia.

La competencia entre consumidores puede propiciar un mayor crecimiento para el recurso, como se muestra en el Teorema \ref{teo2.16} y en la Proposición \ref{prop2.17}. Los ejemplos de la Sección \ref{seccion2.3} muestran que lo anterior es válido para diferentes funciones de crecimiento. De hecho, se puede verificar que
\begin{eqnarray}
\frac{\partial x^*}{\partial m_{ij}} &=& \frac{\det(M_i) y_j^*}{(\det(M))^2 h'(x^*)}. \nonumber
\end{eqnarray}
Esta derivada parcial es positiva si y solo si $\det(M_i)$ y $h'(x^*)$ tienen signos iguales. Para el caso $n=2$, es fácil ver que si $\det(M) > 0$ y $m_{21}/m_{22} \leq q_1/q_2 \leq m_{11}/m_{12}$, entonces $h'(x^*) > 0$. Análogamente, si $\det(M) < 0$ y $m_{11}/m_{12} \leq q_1/q_2 \leq m_{21}/m_{22}$, entonces $h'(x^*) > 0$. Así, la razón de los coeficientes de depredación y competencia se puede usar para determinar el signo de la derivada parcial de $x^*$.

En \cite{Chesson2008}, se establece que una condición suficiente para alcanzar la coexistencia, es que los parámetros intraespecíficos sean mayores que los parámetros interespecíficos. Sin embargo, se puede demostrar que no es una condición necesaria; una selección de parámetros conveniente puede no cumplir esa condición pero sí cumplir las hipótesis del Corolario \ref{cor2.20}. Por lo tanto el Teorema \ref{teo2.19} y el Corolario \ref{cor2.20} confirman y complementan ese trabajo.

La región de coexistencia $U \subset \mathcal M_{n\times n}$ es introducida en el Teorema \ref{teo2.25}. A pesar de que se sabe que $U$ es un conjunto abierto y no vacío, poco se sabe acerca del diámetro de la región. En la Sección \ref{seccion2.3} se muestra que la función de crecimiento afecta la población de equilibrio de las especies, a pesar de que los parámetros permanecen fijos. Entonces es una tarea futura investigar acerca del conjunto $U$ para diferentes funciones de crecimiento en el modelo.

Con el fin de hacer más realista la modelación, se considera una respuesta funcional de tipo Holling II que es denso-dependiente en el Capítulo \ref{capitulo3}.
\chapter{Tasa de captura denso-dependiente}
\label{capitulo3}
\begin{flushright}
\emph{``Las matemáticas son el arte de dar\\el mismo nombre a cosas diferentes.''\\$ $\\Henri Poincaré}
\end{flushright}

Para continuar con el estudio del problema de la coexistencia de las especies, se considera ahora un modelo matemático con dos depredadores\index{depredador}, una función de crecimiento\index{función de crecimiento} dada por el modelo logístico\index{modelo logístico} para la presa\index{presa}, valores de competencia interespecífica $m_{12} = m_{21} = 0$, y una respuesta funcional de tipo Holling II
\begin{eqnarray}
p(x) &=& \frac{q x}{1 + a x}.                                                                         \nonumber
\end{eqnarray}
El parámetro $a$ representa el coeficiente de saturación o el tiempo de manejo de la presa\index{presa} capturada por parte del depredador\index{depredador}. Conjuntando las características anteriores, se obtiene el modelo que aparece a continuación

\begin{eqnarray}
\label{3.1}
\dot x &=& x\left(r(1 - x/K) -\frac{q_1 y}{1 + a_1 x} - \frac{q_2 z}{1 + a_2 x}\right) \quad \colon = \quad xF(x,y,z), \nonumber \\
\dot y &=& y\left(\frac{c_1q_1x}{1 + a_1x} - \mu_1 - m_1y\right) \quad \quad \quad \quad \quad \quad \quad \colon = \quad yG(x,y),                                 \\
\dot z &=& z\left(\frac{c_2q_2x}{1 + a_2x} - \mu_2 - m_2z\right) \quad \quad \quad \quad \quad \quad \quad \colon = \quad zH(x,z),                       \nonumber \\ \nonumber \\
& &\text{con }x(0)>0\text{, }y(0)>0\text{, }z(0)>0\text{. }                                                    \nonumber
\end{eqnarray}
Es claro que si se considera $a_1 = a_2 = 0$, entonces se obtiene el escenario con un recurso y dos consumidores que fue estudiado en el Capítulo \ref{capitulo2}. En lo que sigue, los parámetros $r$, $K$, $q_1$, $q_2$, $a_1$, $a_2$, $c_1$, $c_2$, $\mu_1$, $\mu_2$, $m_1$, $m_2$ son números positivos.

Algunas de las contribuciones hechas a los modelos con una presa\index{presa} y dos depredadores\index{depredador} se enlistan a continuación: en \cite{Loladze2004} se analiza el impacto de la constitución química en un sistema de dos consumidores de un mismo recurso. Esta constitución química es modelada a través de las concentraciones de dos elementos químicos (fósforo y carbono). A través de un análisis numérico se encuentra que las restricciones estequiométricas limitan de forma natural las densidades de población, y eso hace posible la coexistencia; en \cite{Llibre2014} se analiza la dinámica global de un modelo de dos consumidores de un mismo recurso con respuesta funcional de tipo Lotka-Volterra; en \cite{Agrawal2015} se propone un modelo dependiente de la razón de tipo Holling-Tanner, se analiza la dinámica de los submodelos planos y posteriormente se muestra que en el modelo hay presencia de caos; en \cite{Lehman2000} se consideran modelos estocásticos para analizar la estabilidad temporal de las comunidades. Se lleva a cabo un análisis comparativo de los resultados que se obtienen a través de cada modelo; en \cite{Savitri2019APredators} se construye un modelo a partir de otros dos ya conocidos, y se realiza un análisis local de la dinámica del punto de equilibrio interior y de algunos puntos de equilibrio en la frontera; en \cite{Sarwardi2013DynamicalRefuge} se presenta un modelo que es muy similar al que se trabaja en este capítulo, y se muestra una posible bifurcación de Hopf, junto con varias simulaciones numéricas; en \cite{Dubey2004PersistenceSystem} se lleva a cabo un modelo cuya respuesta funcional es dependiente de razón, y se ilustran los resultados a través de simulaciones numéricas. Además se compara un escenario en el que los dos consumidores son excluidos con un resultado análogo obtenido en \cite{Hsu1981OnInterference}; en \cite{Hsu2001RichModel} se muestra que para un modelo con respuesta funcional dependiente de la razón (similar a una función de tipo Beddington-DeAngelis\index{respuesta funcional de Beddington-DeAngelis}), el Principio de Exclusión Competitiva\index{Principio de Exclusión Competitiva} parece cumplirse en muchos escenarios, e incluso pudiendo ocurrir la desaparición de los dos depredadores\index{depredador} o de la presa\index{presa} misma, lo que orillaría automáticamente a la muerte a los consumidores; en \cite{Ton2011DynamicsPrey} se usan los espacios de Banach para poder dar a conocer las propiedades dinámicas de un modelo con respuesta funcional de tipo Beddington-DeAngelis\index{respuesta funcional de Beddington-DeAngelis}, y se ilustran los resultados con algunas simulaciones.

\section{Propiedades generales del modelo}
En esta sección se definen los umbrales de eficiencia energética (u.e.e.), se dan a conocer dos puntos de equilibrio en la frontera y se demuestra que el sistema \eqref{3.1} es disipativo.

Los valores de $\dot y, \dot z$ en \eqref{3.1} tienen el mismo signo que $G(x,y)$ y $H(y,z)$ respectivamente, ya que $y$, $z$ siempre son no-negativos. De lo anterior se sigue que
\begin{eqnarray}
\dot y \quad > \quad 0 &\Leftrightarrow& \frac{c_1 q_1 x}{1 + a_1 x} - \mu_1 - m_1 y \quad > \quad 0.                       \nonumber
\end{eqnarray}
Consecuentemente
\begin{eqnarray}
\dot y \quad > \quad 0 &\Rightarrow& \frac{c_1 q_1 x}{1 + a_1 x} - \mu_1 \quad > \quad 0,                                   \nonumber 
\end{eqnarray}
de donde se sigue que
\begin{eqnarray}
\dot y \quad > \quad 0 &\Rightarrow& x \quad > \quad \frac{\mu_1}{c_1 q_1 - a_1 \mu_1} .                        \nonumber
\end{eqnarray}

\begin{df}\cite{Kuang2003BiodiversityCompetition}
Los \textbf{umbrales de eficiencia energética (u.e.e.)} en el sistema \eqref{3.1}, que se denotan con $y^b$, $z^b$ están dados por 
\begin{eqnarray}
y^b \quad = \quad \frac{\mu_1}{c_1 q_1 - a_1 \mu_1}, & & z^b \quad = \quad \frac{\mu_2}{c_2 q_2 - a_2 \mu_2}.             \nonumber
\end{eqnarray}
\end{df}
Recuérdese (de la Sección \ref{seccion2.1}) que los u.e.e. indican la cantidad de población de presa\index{presa} requerida para que la tasa de crecimiento sea exactamente igual a la tasa de mortalidad del respectivo depredador\index{depredador}. Los u.e.e. $y^b$, $z^b$ son elementos decisivos en la coexistencia o la exclusión de las especies depredadoras\index{depredador}.

Dos puntos de equilibrio en la frontera que siempre aparecen en el retrato fase del modelo \eqref{3.1} son $E_{00} = (0,0,0)$, $E_0 = (K,0,0)$. El primero de ellos representa la exclusión de las tres especies, mientras que el segundo representa la permanencia de la especie recurso junto con la exclusión de las dos especies consumidoras. También es posible hallar puntos de la forma $E_1 = (x_1, y_1, 0)$, $E_2 = (x_2,0,z_2)$, con $x_1$, $y_1$, $x_2$, $z_2 > 0$, en los cuales se representa la permanencia de la especie recurso y una especie consumidora junto con la exclusión de la otra especie consumidora.
\begin{teo}
\label{teo3.2}
Los puntos de equilibrio en la frontera del sistema \eqref{3.1} dados por $E_{00} = (0,0,0)$, $E_0 = (K,0,0)$ tienen las siguientes propiedades.
\begin{enumerate}
    \item $E_{00}$ es un punto silla.
    \item Si se cumple $0 < y^b < K$, o bien $0 < z^b < K$, entonces $E_0$ es un punto silla.
    \item Si $y^b > K$, y $z^b > K$, entonces $E_0$ es un nodo localmente asintóticamente  estable.
\end{enumerate}
\end{teo}
\begin{proof}
La matriz jacobiana de \eqref{3.1} está dada por
\begin{eqnarray}
J(x,y,z) &=&
\begin{pmatrix}
r\left(1 - 2x/K\right) - q_1 y(1 + a_1 x)^{-2} & - q_1 x(1 + a_1 x)^{-1} & - q_2 x(1 + a_2 x)^{-1}             \\
 - q_2 z(1 + a_2 x)^{-2} \\\\\\
c_1 q_1 y(1 + a_1 x)^{-2}  & c_1 q_1 x(1 + a_1 x)^{-1} & 0   \\
   & - \mu_1 - 2 m_{11} y  \\\\\\
c_2 q_2 z (1 + a_2 x)^{-2} & 0 & c_2 q_2 x(1 + a_2 x)^{-1}    \\
 & & - \mu_2 - 2 m_{22} z 
\end{pmatrix}. \nonumber
\end{eqnarray}
Entonces se tiene que
\begin{eqnarray}
J(E_{00}) &=&
\begin{pmatrix}
r & 0      & 0       \\
0 & -\mu_1 & 0       \\
0 & 0      & - \mu_2
\end{pmatrix}, \nonumber
\end{eqnarray}
\begin{eqnarray}
J(E_0) &=&
\begin{pmatrix}
- r & q_1 K(1 + a_1 K)^{-1}             &  q_2 K(1 + a_2 K)^{-1}              \\\\
0   & c_1 q_1 K(1 + a_1 K)^{-1} - \mu_1 &                  0                  \\\\
0   & 0                                   & c_2 q_2 K(1 + a_2 K)^{-1} - \mu_2
\end{pmatrix}. \nonumber
\end{eqnarray}
Por lo tanto, los resultados se siguen para $E_{00}$ y $E_0$ en \eqref{3.1}.

\end{proof}
Nótese que si se tiene $y^b > K$, $z^b > K$, entonces los puntos de equilibrio de la forma $E_1 = (x_1,y_1,0)$, $E_2 = (x_2,0,z_2)$ desaparecen, y lo mismo ocurre con los puntos de equilibrio positivo positivos $E^* = (x^*,y^*,z^*)$.

A continuación se realiza el análisis de los puntos de equilibrio en la frontera restantes.
\begin{prop}
\label{prop3.3}
Si se cumple $0 < y^b < K$ ( $0 < z^b < K$), entonces \eqref{3.1} tiene al menos un punto de equilibrio de la forma $E_1 = (x_1,y_1,0)$ ($E_2 = (x_2,0,z_2)$) con $x_1$, $y_1 > 0$ (con $x_2$, $z_2 > 0$). Más aún, si se satisface
\begin{eqnarray}
\frac{K - 1/a_1}{2} \quad \leq \quad y^b \quad < \quad K \nonumber \\
\left(\frac{K - 1/a_2}{2} \quad \leq \quad z^b \quad < \quad K\right), \nonumber
\end{eqnarray}
entonces \eqref{3.1} tiene un único punto de equilibrio de la forma $E_1 = (x_1,y_1,0)$ ($E_2 = (x_2,0,z_2)$).
\end{prop}
\begin{proof}
Si $z = 0$, entonces se tiene un subsistema plano de \eqref{3.1} cuyas curvas isoclinas están dadas por
\begin{eqnarray}
f(x) &=& \frac{r}{q_1}(1 - x/K)(1 + a_1x), \nonumber \\
g(x) &=& \frac{1}{m_1}\left(\frac{c_1 q_1 x}{1 + a_1 x} - \mu_1\right), \nonumber
\end{eqnarray}
Se verifica fácilmente que $f$ y $g$ son simultáneamente no-negativas precisamente en el intervalo $[y^b,K]$. Por otra parte, $f$ es una función decreciente en $[(K - 1/a)/2, K]$, y $g$ es creciente en $[y^b,K]$. Finalmente, se tiene $f(y^b) > 0 = g(y^b)$, $f(K) = 0 < g(y^b)$. Entonces existe exactamente un punto $x_1\in (y^b,K)$ que satisface
\begin{eqnarray}
f(x_1) &=& g(x_1). \nonumber
\end{eqnarray}
\end{proof}
A continuación se estudian las propiedades dinámicas de estos puntos de equilibrio.
\begin{prop}
\label{prop3.4}
Si $0 < y^b < z^b < K$, entonces cualquier punto de equilibrio de la forma $E_2 = (x_2,0,z_2)$ es asintóticamente inestable en \eqref{3.1}. 
\end{prop}
\begin{proof}
Como el renglón intermedio de la matriz jacobiana $D(E_2)$ está dado por \newline $(0,c_1q_1x_2(1 + a_1x_2)^{-1} - \mu_1,0)$, y se tiene que $y^b < z^b < x_2$, entonces $c_1q_1x_2(1 + a_1x_2)^{-1} - \mu_1> 0$ y así $E_2$ es inestable.
\end{proof}

A través de un procedimiento similar al que se aplica en \cite{Sarwardi2013DynamicalRefuge}, se muestra que el sistema \eqref{3.1} es disipativo.
\begin{teo}
\label{teo3.5}
Existe un subconjunto compacto $C$ de $[0,\infty)^3$ de tal manera que para toda solución no-negativa $(x,y,z)$ de \eqref{3.1}, existe $T > 0$ de tal forma que $(x(t),y(t),z(t)) \in C$ para toda $t\geq T$.
\end{teo}
\begin{proof}
De forma similar a lo que se hizo en la demostración del Teorema \ref{teo2.1}, se considera $w = x + c_1^{-1}y + c_2^{-1}z$, y posteriormente se definen $\varphi \in (0,\min\{\mu_1,\mu_2\}]$, $\displaystyle \rho = \frac{K(r + \varphi)^2}{4 r}$. Así
\begin{eqnarray}
w' + \varphi w &=& x(r(1 - x/K)) - c_1^{-1}(\mu_1 y + m_1 y^2)  \nonumber \\
& & - c_2^{-1}(\mu_2 z + m_2 z^2) + \varphi(x + c_1^{-1}y + c_2^{-1}z) \nonumber \\
&\leq& x\left(r + \varphi - \left(\frac{r}{K}\right)x\right) - c_1^{-1}(\mu_1 - \varphi)y - c_2^{-1}(\mu_2 - \varphi)z \nonumber \\
&\leq& x\left(r + \varphi - \left(\frac{r}{K}\right)x\right) \nonumber \\
&\leq& \frac{K(r + \varphi)^2}{4 r} \nonumber \\
&=& \rho. \nonumber
\end{eqnarray}

De lo anterior se tiene $\displaystyle w \leq e^{-\varphi t}w_0 + \frac{\rho}{\varphi}(1 - e^{-\varphi t})$. 

Y por lo tanto $\displaystyle \limsup_{t\to \infty} w(t) \leq \frac{\rho}{\varphi}$. 

Se concluye entonces que para cualquier solución no-negativa $(x,y,z)$ de \eqref{3.1}, existe $T>0$ tal que
\begin{eqnarray}
(x(t),y(t),z(t)) &\in& \left\{(x,y,z)\in \mathbb [0,\infty)\colon x + c_1^{-1}y + c_2^{-1}z\leq \frac{\rho}{\varphi}\right\}\text{, si } t \geq T. \nonumber
\end{eqnarray}
\end{proof}
\section{Dinámica alrededor de los puntos de equilibrio positivos}
Recuérdese que los puntos de equilibrio positivos tienen un significado ecológico relevante, pues representan estados de coexistencia de todas especies. En el Teorema \ref{teo3.7} se muestran algunas condiciones que deben cumplir los parámetros para que exista un punto de equilibrio positivo $E^* = (x^*,y^*,z^*)$ en \eqref{3.1}.

El sistema \eqref{3.1} es simétrico en el sentido de que los papeles que juegan $y$, $z$ son conmutables. Por lo tanto se puede suponer sin pérdida de generalidad que siempre se cumple $0 < y^b < z^b < K$. 
\begin{lema}
\label{lema3.6}
Si $0 < y^b < z^b < K$, entonces se cumplen las desigualdades 
\begin{eqnarray}
c_2q_2 - a_2\mu_2 &>& 0, \nonumber \\
\mu_2(c_1q_1 - a_1\mu_1) - \mu_1(c_2q_2 - a_2\mu_2) &>& 0, \nonumber \\ 
K(c_2q_2 - a_2\mu_2) - \mu _2 &>& 0. \nonumber
\end{eqnarray}
\end{lema}
\begin{proof}
La primera desigualdad se verifica de manera inmediata. En la segunda desigualdad se observa que
\begin{eqnarray}
z^b - y^b > 0 &\Rightarrow&
\frac{\mu_2}{c_2q_2 - a_2\mu_2} - \frac{\mu_1}{c_1q_1 - a_1\mu_1} > 0 \nonumber \\ \nonumber \\ &\Rightarrow& \mu_2(c_1q_1 - a_1\mu_1) - \mu_1(c_2q_2 - a_2\mu_2) > 0. \nonumber
\end{eqnarray}
Finalmente, la tercera desigualdad se comprueba así
\begin{eqnarray}
K - z^b > 0 &\Rightarrow& K - \frac{\mu_2}{c_2q_2 - a_2\mu_2} >0 \nonumber \\ \nonumber \\
&\Rightarrow& K(c_2q_2 - a_2\mu_2) - \mu _2 > 0. \nonumber
\end{eqnarray}
\end{proof}
 Es importante tener en cuenta la positividad de los valores que se analizan en la proposición anterior para la demostración del siguiente resultado.
\begin{teo}
\label{teo3.7}
Supóngase que $0 < y^b <  z^b < K$. Entonces 
\begin{eqnarray}
r &>& \left(\frac{c_2q_2 - a_2\mu _2}{a_1\mu_2 - a_2 \mu_2 + c_2q_2}\right)^2 \frac{K q_1 (\mu_2(c_1q_1 - a_1\mu_1) - \mu_1(c_2q_2 - a_2\mu_2))}{m_1(K(c_2q_2 - a_2\mu_2) - \mu _2)}, \nonumber
\end{eqnarray}
si y solo si \eqref{3.1} tiene al menos un punto de equilibrio positivo.
\end{teo}
\begin{proof}
Del Lema \ref{lema3.6}, se sigue que el lado derecho de la desigualdad anterior es mayor que cero. Se definen $f$, $g$, $h$ como las funciones isoclinas de $\dot x$, $\dot y$, $\dot z$ (respectivamente). Estas últimas están dadas por 
\begin{eqnarray}
f(x,y) &\colon =& \frac{1}{q_2}\left(r(1 - x/K) - \frac{q_1 y}{1 + a_1x}\right)(1 + a_2x),            \nonumber \\
g(x) &\colon =& \frac{1}{m_1}\left(\frac{c_1q_1 x}{1 + a_1x} - \mu_1\right),                          \nonumber \\
h(x) &\colon =& \frac{1}{m_2}\left(\frac{c_2q_2 x}{1 + a_2x} - \mu_2\right).                          \nonumber
\end{eqnarray}
A través de un cálculo se verifica que la hipótesis se satisface si y solo si
\begin{eqnarray}
r\left(\frac{a_2\mu _2}{q_2(c_2q_2 - a_2\mu_2)} + \frac{1}{q_2} - \frac{a_2\mu _2^2}{Kq_2(c_2q_2 - a_2\mu_2)^2} - \frac{\mu_2}{Kq_2(c_2q_2 - a_2\mu _2)}\right) \quad > \nonumber
\end{eqnarray}
\begin{eqnarray}
\frac{c_1q_1^2\mu_2}{m_1q_2(c_2q_2 - a_2\mu_2) \left(\frac{a_1\mu_2}{c_2q_2 - a_2\mu_2} + 1\right)^2} &+& \frac{a_2c_1q_1^2\mu_2^2}{m_1q_2(c_2q_2 - a_2\mu_2)^2\left(\frac{a_1\mu_2}{c_2q_2 - a_2\mu_2} + 1\right)^2} \nonumber \\ \nonumber \\
- \frac{q_1\mu_1}{m_1q_2\left(\frac{a_1\mu_2}{c_2q_2 - a_2\mu_2} + 1\right)} &-& \frac{a_2q_1\mu_1\mu_2}{m_1q_2(c_2q_2 - a_2\mu_2) \left(\frac{a_1\mu_2}{c_2q_2 - a_2\mu_2} + 1\right)}, \nonumber
\end{eqnarray}
que es equivalente a  
\begin{eqnarray}
f(z^b,g(z^b)) &>& 0. \nonumber
\end{eqnarray}
Es fácil ver que 
\begin{eqnarray}
\quad h(z^b) &=& 0. \nonumber
\end{eqnarray}
y que
\begin{eqnarray}
h(K) \quad > \quad 0 &>& f(K,g(K)). \nonumber
\end{eqnarray}
Del teorema del valor intermedio, se sigue que existe $x^* \in (0,K)$ tal que a $h(x^*) = f(x^*,g(x^*))$. Así, $E^* = (x^*,y^*,z^*)$ es un punto de equilibrio positivo en \eqref{3.1}, con $y^* = g(x^*)$, $z^* = h(x^*)$.

\end{proof}
Se puede observar que el lado derecho de la desigualdad del Teorema \eqref{teo3.7} se puede reescribir
\begin{eqnarray}
\left(\frac{1}{a_1 z^b + 1}\right)^2\left(\frac{K q_1 \mu_1 (z^b - y^b)}{m_1 y^b(K - z^b)}\right). \nonumber
\end{eqnarray}
El único parámetro que no aparece dicha expresión es precisamente $m_2$. Entonces la competencia que ejerce el depredador\index{depredador} débil es independiente para existencia de puntos de equilibrio positivos. Por otro lado, se puede apreciar que al incrementar los valores de $K$, $q_1$, crece también el lado derecho de la desigualdad. Lo anterior dificulta la existencia de dichos estados estacionarios, mientras que el incremento en $m_1$ aparentemente propicia la existencia de estos últimos.

Es claro que el Teorema \ref{teo3.7} se puede aplicar al modelo propuesto en el Capítulo \ref{capitulo2} si consideramos el caso $n=2$ con el modelo logístico como función de crecimiento y $m_{12} = m_{21} = 0$. Solo hace falta hacer cero los parámetros $a_1$, $a_2$.

\begin{cor}
\label{cor3.8}
Si $(K - 1/a_1)/2 < y^b <z^b < K$ y además el punto de equilibrio de la forma $E_1 = (x_1,y_1,0)$ es asintóticamente estable en \eqref{3.1}, entonces el sistema \eqref{3.1} no cuenta con puntos de equilibrio positivo $E^* = (x^*,y^*,z^*)$.
\end{cor}
\begin{proof}
De la Proposición \ref{prop3.3}, se tiene que $E_1 = (x_1,y_1,0)$ es el único punto de equilibrio con esta estructura ($z = 0$). El último renglón de la matriz jacobiana $D(E_1)$ está dado por $(0,0,c_2q_2x_1(1 + a_2x_1)^{-1} - \mu_2)$. Como $E_1$ es estable, entonces $c_2q_2x_1(1 + a_2x_1)^{-1} - \mu_2 < 0$, de donde se sigue que $x_1 < z^b$. Por otro lado, se tiene que $x_1 > 2^{-1}(K - 1/a_1)$. Obsérvese que las funciones definidas en la demostración del Teorema \ref{teo3.7} cumplen:
\begin{enumerate}
    \item $g$ es creciente respecto a $x$.
    \item $f$ es decreciente respecto a $y$.
    \item $f$ es decreciente respecto a $x$ en $[(K - 1/a_1)/2,\infty)$.
\end{enumerate}
Como $f(x_1,g(x_1)) = 0$, entonces
\begin{eqnarray}
f(z^b,g(z^b)) &<& 0, \nonumber
\end{eqnarray}
y esto ocurre si y solo si
\begin{eqnarray}
r &<& \left(\frac{c_2q_2 - a_2\mu _2}{a_1\mu_2 - a_2 \mu_2 + c_2q_2}\right)^2 \frac{K q_1 (\mu_2(c_1q_1 - a_1\mu_1) - \mu_1(c_2q_2 - a_2\mu_2))}{m_1(K(c_2q_2 - a_2\mu_2) - \mu _2)}. \nonumber
\end{eqnarray}
Del Teorema \ref{teo3.7} se sigue que \eqref{3.1} no tiene puntos de equilibrio positivo.

\end{proof}

A continuación se presenta un resultado que trata la estabilidad asintótica local de los puntos de equilibrio positivos. En la demostración se aplica el teorema de Hartman-Grobman a la matriz jacobiana de \eqref{3.1}, y también se usa el criterio de Routh-Hurwitz para estudiar las raíces del polinomio característico de la matriz jacobiana.
\begin{teo}[Reyes-García]
\label{teo3.9}
Si \eqref{3.1} tiene un punto de equilibrio positivo en $E^* = (x^*,y^*,z^*)$, y además \newline $y^b,z^b \in (K/2,K)$, entonces $E^*$ es el único punto de equilibrio positivo, y $E^*$ es un atractor local en \eqref{3.1}.
\end{teo}
\begin{proof}
Siguiendo la notación del Teorema \ref{teo3.7}, se tiene $f(x^*,g(x^*)) = h(x^*)$, donde $K/2 < x^* < K$. Además, $f(\cdot,g(\cdot))$ es decreciente y $h$ es creciente en el intervalo $(K/2,K)$, de donde se sigue que $E^*$ es único. 

Considérese la matriz jacobiana del sistema \eqref{3.1} aplicada al punto de equilibrio positivo $E^* = (x^*,y^*,z^*)$,
\begin{eqnarray}
D(E^*) &=& \left( \begin{array}{ccc}
                A         & - B_1    & - B_2  \nonumber \\
                C_1       & - M_1    & 0      \nonumber \\
                C_2       & 0        & - M_2   \nonumber
                \end{array} \right),                                  \nonumber
\end{eqnarray}
con 
\begin{eqnarray}
A = r\left(1 - \frac{2x^*}{K}\right) - \frac{q_1y^*}{(1 + a_1x^*)^2} &-& \frac{q_2z^*}{(1 + a_2x^*)^2},                                                            \nonumber \\
B_1 = \frac{q_1x^*}{1 + a_1x^*}, &&
B_2 = \frac{q_2x^*}{1 + a_2x^*},                                    \nonumber \\
C_1 = \frac{c_1q_1y^*}{(1 + a_1x^*)^2}, &&
C_2 = \frac{c_2q_2z^*}{(1 + a_2x^*)^2},                                \nonumber \\ \nonumber \\
M_1 = m_1 y^*, && M_2 = m_2 z^*. \nonumber
\end{eqnarray}
Se puede ver que los valores de $B_1$, $B_2$, $C_1$, $C_2$, $M_1$, $M_2$ son positivos. El polinomio característico de la matriz $D$ está dado por
\begin{eqnarray}
P(\lambda) &=& \lambda^3 + \Omega_2 \lambda^2 + \Omega_1 \lambda + \Omega_0,    \nonumber \\
\text{donde}                                                                    \nonumber \\
\Omega_2 &=& M_1 + M_2 - A,               \nonumber \\
\Omega_1 &=& M_1M_2 + B_1C_1 + B_2C_2 - A(M_1 + M_2),      \nonumber \\
\Omega_0 &=& B_1C_1M_2 + B_2C_2M_1 - AM_1M_2. \nonumber
\end{eqnarray}
El Criterio de Routh-Hurwitz dice que si los valores $\Omega_0$, $\Omega_2$ y $\Omega_2\Omega_1 - \Omega_0$ son todos positivos, entonces las raíces de $P$ tienen parte real negativa. Haciendo un cálculo se obtiene
\begin{eqnarray}
\Omega_2\Omega_1 - \Omega_0 &=& B_1C_1(M_1 - A) + B_2C_2(M_2 - A) + (M_1 + M_2)(M_1 - A)(M_2 - A).    \nonumber
\end{eqnarray}
Como $\min \{y^b,z^b\} > K/2$ entonces $x^* > K/2$, y por lo tanto $A < 0$. De lo anterior, se concluye que $E^*$ es un atractor local para \eqref{3.1}.

\end{proof}
Se puede observar que la relación de los parámetros de competencia que se establece en el Corolario \ref{cor2.20} aparece también en el Teorema \ref{teo3.9}, pues se tiene $m_{12} = m_{21} = 0$.
\section{Bifurcación de Hopf}
\label{seccion3.3}
Se realiza un procedimiento para verificar la existencia de una bifurcación de Hopf en el sistema \eqref{3.1}. Supóngase que \eqref{3.1} tiene un punto de equilibrio con coordenadas positivas en $E^* = (x^*,y^*,z^*)$. Entonces el primer paso consiste en efectuar un cambio de variables como se muestra a continuación.
\begin{eqnarray}
\overline x = x/x^*,\quad  
\overline y = y/y^*, \quad
\overline z = z/z^*, \quad
\overline r = r, \nonumber \\ 
\overline K = K/x^*, \quad
\overline q_1 = q_1 y^*, \quad
\overline q_2 = q_2 z^*, \quad 
\overline a_1 = a_1 x^*, \nonumber \\
\overline a_2 = a_2 x^*, \quad    
\overline c_1 = c_1 x^*/y^*, \quad
\overline c_2 = c_2 x^*/z^*, \nonumber \\  
\overline \mu_1 = \mu_1, \quad         
\overline \mu_2 = \mu_2, \quad         
\overline m_1 = m_1 y^*,        \quad
\overline m_2 = m_2 z^*,        \nonumber
\end{eqnarray}
Después de realizar el cambio, se retiran las barras para simplificar la notación. Se puede verificar que \eqref{3.1} tiene ahora un punto de equilibrio en $(1,1,1)$. Se propone hacer igual a uno los siguientes parámetros 
\begin{eqnarray}
q_2 \quad = \quad a_1 \quad = \quad a_2 \quad = \quad c_1 \quad = \quad c_2 \quad = \quad m_1 \quad = \quad 1. \nonumber
\end{eqnarray}
A partir de lo anterior, \eqref{3.1} queda transformado en
\begin{eqnarray}
\label{3.2}
\dot x &=& x\left(r(1 - x/K) -\frac{q_1 y}{1 + x} - \frac{z}{1 + x}\right), \nonumber \\
\dot y &=& y\left(\frac{q_1x}{1 + x} - \mu_1 - y\right),                                 \\
\dot z &=& z\left(\frac{x}{1 + x} - \mu_2 - m_2z\right).                       \nonumber 
\end{eqnarray}
Como $(1,1,1)$ es punto de equilibrio para \eqref{3.2}, es posible despejar algunos parámetros
\begin{eqnarray}
r &=& \frac{K(q_1 + 1)}{2(K - 1)}, \nonumber \\
\mu_1 &=& \frac{q_1 - 2}{2}, \nonumber \\
\mu_2 &=& \frac{1 }{2}-m_2. \nonumber
\end{eqnarray}
Así, se obtiene un nuevo sistema de ecuaciones diferenciales
\begin{eqnarray}
\label{3.3}
\dot x &=& x\left(\frac{K(q_1 + 1)}{2(K - 1)}(1 - x/K) -\frac{q_1 y}{1 + x} - \frac{z}{1 + x}\right), \nonumber \\
\dot y &=& y\left(\frac{q_1x}{1 + x} - \left(\frac{q_1 - 2}{2}\right) - y\right),                                 \\
\dot z &=& z\left(\frac{x}{1 + x} - \left(\frac{1}{2} - m_2\right) - m_2z\right).                       \nonumber 
\end{eqnarray}
\begin{teo}[Reyes-García]
\label{teo3.10}
Supóngase que $K > 1$, $q_1 >2$, y adicionalmente se cumple
\begin{eqnarray}
q_1(q_1((K - 3)q_1(4(K(K(13K - 37) + 47)  & & \nonumber \\
- 31) - (K - 3)(K(7K - 10) - 1)q_1)   & & \nonumber \\
+ 2K(K(- 51(K - 4)K - 322) + 140) + 154)  & & \nonumber \\
- 4(K - 3)^2(K(17K - 30) + 9)) + (K (17 K - 30) + 9)^2 &>& 0. \nonumber
\end{eqnarray}
Entonces el sistema \eqref{3.3} admite configuraciones de parámetros para las cuales el espectro de la matriz jacobiana del punto de equilibrio positivo $(1,1,1)$ tiene la forma $\{\nu,i\omega,-i\omega\}$.
\end{teo}
\begin{proof}
La matriz jacobiana de \eqref{3.3} aplicada al punto $(1,1,1)$ está dada por
\begin{eqnarray}
J(1,1,1) &=& \begin{pmatrix}
(K - 3)(q_1 + 1)(4(K - 1))^{-1} & - q_1/2 & - 1/2 \\\\
q_1/4                             & - 1     & 0     \\\\
1/4                               & 0       & -m_2
\end{pmatrix},
\nonumber
\end{eqnarray}
y el polinomio característico de $D$ es
\begin{eqnarray}
P(\lambda) &=& - \lambda^3 + A \lambda^2 + B \lambda + C\text{, donde} \nonumber
\end{eqnarray}
\begin{eqnarray}
A &=& \frac{(K - 3)\left(q_1 + 1\right)}{4(K - 1)} - m_2 - 1, \nonumber \\ \nonumber \\
B &=& \frac{2 (K - 3) \left(m_2 + 1\right) q_1 - 6 K m_2 - (K - 1) q_1^2 + K + 2 m_2 - 5}{8 (K - 1)}, \nonumber \\ \nonumber \\
C &=& \frac{(K-3) m_2 \left(q_1+1\right)}{4 (K-1)}-\frac{1}{8} m_2 q_1^2-\frac{1}{8}. \nonumber
\end{eqnarray}
Se desea que el polinomio $P$ tenga raíces $\{\nu,i\omega,i\omega\}$. Entonces $P$ debe tener la forma siguiente
\begin{eqnarray}
P(\lambda) &=& (\nu - \lambda)(i\omega - \lambda)(- i\omega - \lambda) \nonumber \\
&=& -\lambda^3 + \nu \lambda^2 - \omega^2 \lambda + \nu \omega^2. \nonumber
\end{eqnarray}
Si $B$ es un número negativo, y $A$ tiene el mismo signo que $C$, entonces se escribe
\begin{eqnarray}
A = \nu, \quad B = - \omega^2, \quad C = \nu \omega^2, \nonumber
\end{eqnarray}
de donde se sigue
\begin{eqnarray}
- AB = C. \nonumber
\end{eqnarray}
Al hacer la sustitución de los valores $A$, $B$, $C$ se obtiene
\begin{eqnarray}
-\left(\frac{(K - 3)\left(q_1 + 1\right)}{4(K - 1)} - m_2 - 1\right) \times & & \nonumber \\
\left(\frac{2 (K - 3) \left(m_2 + 1\right) q_1 - 6 K m_2 - (K - 1) q_1^2 + K + 2 m_2 - 5}{8 (K - 1)}\right) & & \nonumber 
\end{eqnarray}
\begin{eqnarray}
&=& \frac{(K-3) m_2 \left(q_1+1\right)}{4 (K-1)}-\frac{1}{8} m_2 q_1^2-\frac{1}{8}. \nonumber
\end{eqnarray}
Despejando $m_2$ de la ecuación anterior, se obtiene
\begin{eqnarray}
m_2 &=& 8^{-1}(K - 1)^{-1} \left((K - 3)q_1 - 3K + 1\right)^{-1}((K - 3)^2q_1^2 - 2(K - 3)(3K - 1)q_1 \nonumber \\
& & + K(11K - 10)+ 3 \pm (q_1(q_1((K - 3)q_1(4(K(K (13K - 37) + 47) - 31)  \nonumber \\
& & - (K - 3)(K(7K - 10) - 1)q_1) + 2K(K(- 51 (K - 4)K - 322) + 140) + 154) \nonumber \\
& &  - 4(K - 3)^2(K(17 K - 30) + 9)) + (K(17 K - 30) + 9)^2)^{1/2}). \nonumber
\end{eqnarray}
Es precisamente en ese valor en el que el espectro tiene la forma deseada.

\end{proof}
Para poder concluir la existencia de una bifurcación de Hopf en el retrato fase de \eqref{3.1}, es necesario que se cumplan las condiciones de transversalidad y de genericidad. Con respecto a la condición de transversalidad en una bifurcación de Hopf, se puede mostrar que para determinadas selecciones de parámetros, la derivada respecto a $m_2$ de la parte real de los valores propios de la matriz jacobiana $J(1,1,1)$ que pertenecen al conjunto $\mathbb C \setminus \mathbb R$ es distinta de cero. 

La bifurcación resulta ser genérica para esas selecciones de parámetros. Esto se observa al calcular los valores del primer coeficiente de Lyapunov $\ell_1$. En todos los experimentos hechos se obtuvo $\ell_1 \neq 0$. En los primeros seis casos se tiene un valor negativo para $\ell_1$ (lo que da lugar a una bifurcación de Hopf supercrítica y la existencia de ciclos límite estables), y en el último de ellos $\ell_1$ es positivo (lo que implica una bifurcación de Hopf subcrítica y la existencia de un ciclo límite inestable).
$ $ \newpage
\[\]\[\]\[\]\[\]
\begin{figure}[h]
    \centering
    \includegraphics[width=488pt]{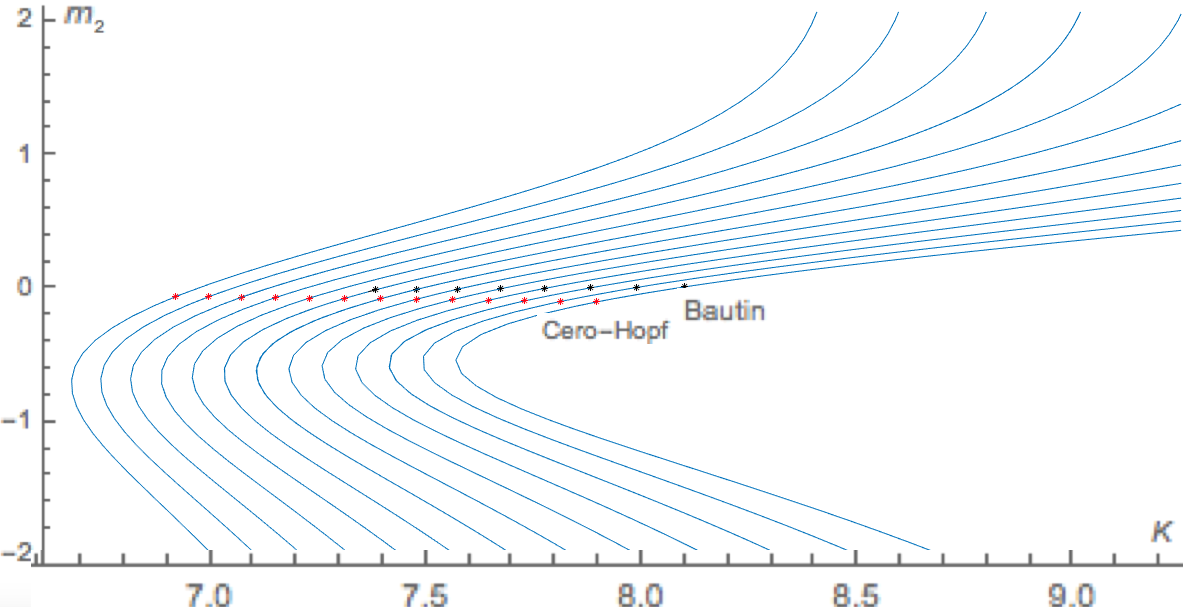}
    \caption{Para los niveles $q_1 = 4$, $4.05$, $4.1$, $4.15$, $4.2$, $4.25$, $4.3$, $4.35$, $4.4$, $4.45$, $4.5$, $4.55$, $4.6$ se grafican de derecha a izquierda las curvas formadas por las parejas ordenadas $(K,m_2)$ en las que se presenta un umbral de bifurcación de Hopf (curvas de Hopf) para el sistema \eqref{3.3}. Los puntos rojos corresponden a umbrales de una bifurcación de cero-Hopf, mientras que los puntos negros son umbrales de bifurcaciones de Bautin. Un aspecto que conviene mencionar, es que los umbrales de bifurcación de cero-Hopf y de Bautin exhibidos en este diagrama, se encuentran por debajo del eje $m_2=0$, y por lo tanto están fuera de la región de interés ecológico.}
    \label{fig3.1}
\end{figure}
\newpage
\begin{eqnarray}
K = 8.0,   \quad q_1 = 4.00     &\Rightarrow& \frac{\partial Re(\lambda)}{\partial m_2} = -0.383873, \quad \ell_1 = -0.383873. \nonumber \\
K = 7.9, \quad q_1 = 4.15  &\Rightarrow& \frac{\partial Re(\lambda)}{\partial m_2} = -0.102561, \quad \ell_1 = -3.59672 \nonumber \\
K = 7.8, \quad q_1 = 4.30   &\Rightarrow& \frac{\partial Re(\lambda)}{\partial m_2} = -0.0977396, \quad \ell_1 = -4.39827. \nonumber \\
K = 7.7, \quad q_1 = 4.45  &\Rightarrow& \frac{\partial Re(\lambda)}{\partial m_2} = -0.0871043, \quad \ell_1 = -4.74064. \nonumber \\
K = 7.6, \quad q_1 = 4.60   &\Rightarrow& \frac{\partial Re(\lambda)}{\partial m_2} = -0.0706634, \quad \ell_1 = -4.92111. \nonumber \\
K = 7.5, \quad q_1 = 4.45   &\Rightarrow& \frac{\partial Re(\lambda)}{\partial m_2} = -0.0877358, \quad \ell_1 = -4.58241. \nonumber \\
K = 7.4, \quad q_1 = 4.30   &\Rightarrow& \frac{\partial Re(\lambda)}{\partial m_2} = -0.0774249, \quad \ell_1 = 52.8238. \nonumber
\end{eqnarray}
Con base en los puntos anteriores, se concluye la existencia de bifurcaciones de tipo Hopf en el modelo de interés. El diagrama de la Figura \ref{fig3.1} (que fue realizado con \textit{Matcont}, una interfaz creada en \textit{Matlab}) permite ubicar los puntos en los que aparecen los umbrales de bifurcaciones de cero-Hopf así como de bifurcaciones de Bautin.
\section{Simulaciones numéricas}
\label{seccion3.4}
Ahora se muestran algunos ejemplos de la dinámica en el sistema \eqref{3.1}. En el primero de ellos se observa la existencia de más de un punto de equilibrio positivo, mientras que en el segundo aparentemente hay un ciclo límite atractor.
\begin{ej}
\label{ej3.11}
Uno de los escenarios en los que hay más de un punto de equilibrio positivo en el modelo \eqref{3.1} es presentado a continuación. Para la configuración de parámetros que abajo se muestra
\begin{eqnarray}
r = 1.3, \quad 
K = 4, && 
q_1 = 1.7, \quad 
a_1 = 3, \nonumber \\
q_2 = 2, \quad
a_2 = 5, &&    
c_1 = 2, \quad
c_2 = 1, \nonumber \\
\mu_1 = 0.3, \quad
\mu_2 = 0.2, &&
m_1 = 0.98, \quad         
m_2 = 0.06, \nonumber 
\end{eqnarray}
se tienen tres puntos de equilibrio positivo en \eqref{3.1}. Estos últimos se encuentran ubicados en 
\begin{eqnarray}
E_1^* = (0.622294,0.446953,1.71185) & & \text{ que es un punto silla, } \nonumber \\
E_2^* = (1.06112,0.573898,2.27608) & & \text{ que también es un punto silla, } \nonumber \\
E_3^* = (1.85518,0.674199,2.68457) & & \text{que es un foco atractor.} \nonumber
\end{eqnarray}

Como solo uno de los puntos de equilibrio mostrados es atractor, se puede concluir que en este caso particular, no hay multiestabilidad en las soluciones positivas.

    \begin{figure}[h]
        \centering
        \includegraphics[width=160pt]{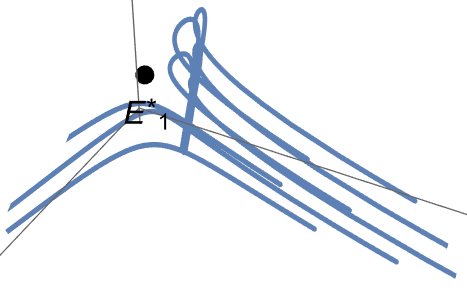}
        \includegraphics[width=160pt]{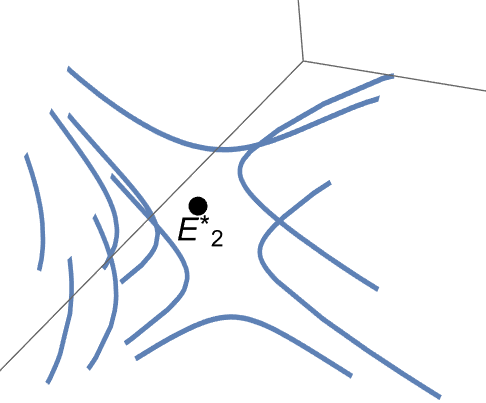}
        \includegraphics[width=160pt]{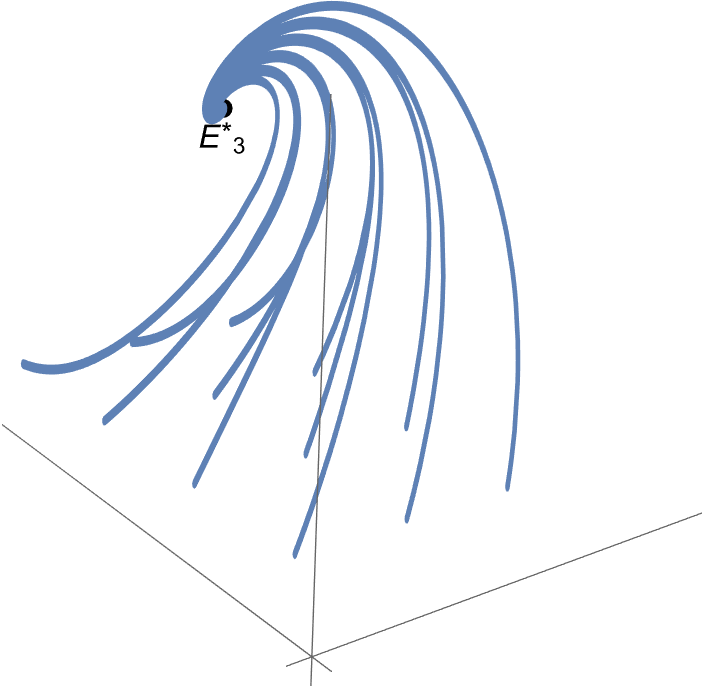}
        \caption{Dinámica local de los puntos de equilibrio del sistema \eqref{3.1} con la selección de parámetros especificada en el Ejemplo \ref{ej3.11}. Se bosquejan algunas de las trayectorias, y se observa su comportamiento asintótico. $E_1^* = (0.622294,0.446953,1.71185)$, $E_2^* = (1.06112,0.573898,2.27608)$, $E_3^* = (1.85518,0.674199,2.68457)$}
    \end{figure}
    
Se puede observar que los parámetros seleccionados en el Ejemplo \ref{ej3.11} satisfacen las hipótesis del Teorema \ref{teo3.7}, pues
\begin{eqnarray}
r &=& 1.3, \nonumber
\end{eqnarray}    

\begin{eqnarray}
\left(\frac{c_2q_2 - a_2\mu _2}{a_1\mu_2 - a_2 \mu_2 + c_2q_2}\right)^2 \frac{K q_1 (\mu_2(c_1q_1 - a_1\mu_1) - \mu_1(c_2q_2 - a_2\mu_2))}{m_1(K(c_2q_2 - a_2\mu_2) - \mu _2)} &=& 0.142656. \nonumber
\end{eqnarray}
    
Los u.e.e. están dados por 
 \begin{eqnarray}
 y^b \quad = \quad 0.12, && z^b \quad = \quad 0.2, \nonumber
 \end{eqnarray}
 mientras que la capacidad de carga está dada por 
 \begin{eqnarray}
 K = 4. \nonumber
 \end{eqnarray}
 Por lo tanto las hipótesis del Teorema \ref{teo3.9} no se cumplen, sin embargo el sistema del Ejemplo \eqref{ej3.11} sí cuenta con un equilibrio positivo que es localmente estable.
\end{ej}

\begin{ej}
\label{ej3.12}
Se exhibe un caso que da lugar a un punto de equilibrio positivo en \eqref{3.1} que es asintóticamente inestable. Los parámetros que se consideran son
\begin{eqnarray}
r = 3, \quad 
K = 4, && 
q_1 = 1.8, \quad 
a_1 = 1, \nonumber \\
q_2 = 2, \quad
a_2 = 1, &&    
c_1 = 1, \quad
c_2 = 1, \nonumber \\
\mu_1 = 0.6, \quad
\mu_2 = 0.8, &&
m_1 = 1, \quad         
m_2 = 0.05. \nonumber 
\end{eqnarray}

El único punto de equilibrio positivo se localiza en $E^* = (0.816562,0.209117,1.98039)$. Al ser inestable $E^*$, se tiene que las trayectorias se alejan de $E^*$, y se aproximan a un ciclo límite. Lo anterior se ilustra en las Figuras \ref{fig3.3} y \ref{fig3.4}.

Un aspecto que muestra el Ejemplo \ref{ej3.12} es la existencia de un depredador marginal\index{depredador!marginal}, que en este caso está representado con la variable $y$.
\begin{figure}[h]
    \centering
    \includegraphics[width=380pt]{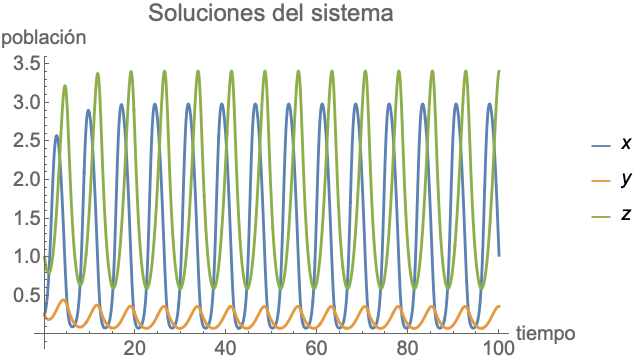}
    \caption{Soluciones del sistema \eqref{3.1} con los parámetros indicados en el Ejemplo \ref{ej3.12} y condiciones iniciales $x(0) = 0.25$, $y_0 = 0.25$, $z_0 =1$.}
    \label{fig3.3}
\end{figure}

Obsérvese que la hipótesis del Teorema \ref{teo3.7} se cumple, ya que
\begin{eqnarray}
r &=& 3, \nonumber
\end{eqnarray}
\begin{eqnarray}
\left(\frac{c_2q_2 - a_2\mu _2}{a_1\mu_2 - a_2 \mu_2 + c_2q_2}\right)^2 \frac{K q_1 (\mu_2(c_1q_1 - a_1\mu_1) - \mu_1(c_2q_2 - a_2\mu_2))}{m_1(K(c_2q_2 - a_2\mu_2) - \mu _2)} &=& 0.15552. \nonumber \nonumber
\end{eqnarray}
También se puede verificar que los u.e.e. están dados por
\begin{eqnarray}
y^b \quad = \quad 0.5, & & z^b \quad = \quad 0.6666, \nonumber
\end{eqnarray}
y la capacidad de carga es
\begin{eqnarray}
K &=& 4, \nonumber
\end{eqnarray}
de tal forma que la hipótesis del Teorema \ref{teo3.9} no se cumple. Esta hipótesis evidentemente fallaría en cumplirse al ser inestable el punto de equilibrio positivo $E^*$, y de hecho se puede verificar que en este caso \eqref{3.1} no tiene puntos de equilibrio estables.

\begin{figure}[h]
    \centering
    \includegraphics[width=180pt]{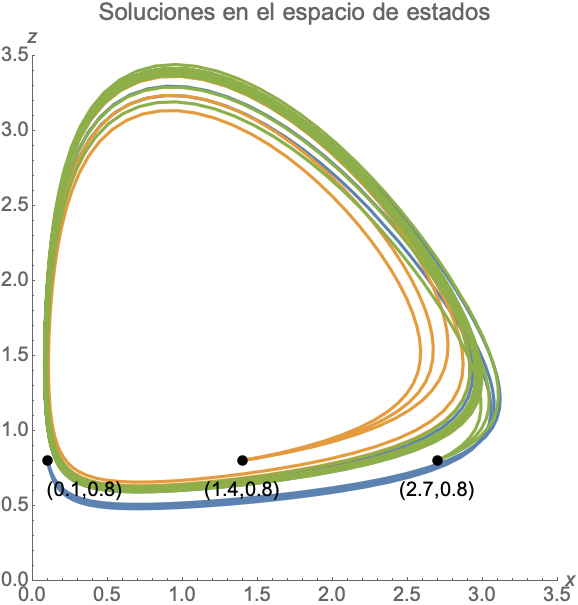}
    \includegraphics[width=180pt]{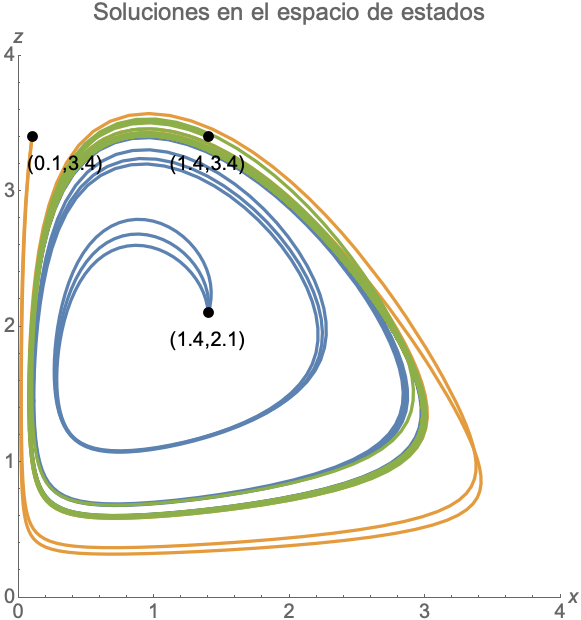}
    \caption{Se grafican las coordenadas $x$, $z$ de algunas soluciones particulares de \eqref{3.1} con los parámetros del Ejemplo \ref{ej3.12}. Las condiciones $x_0$, $z_0$ se muestran en cada caso, mientras que $y_0$ toma alguno de los valores $0.1$, $0.2$, $0.3$. Se puede apreciar la convergencia hacia una órbita periódica no trivial en ambos casos.}
    \label{fig3.4}
\end{figure}
\end{ej}

\begin{ej} 
\label{ej3.13}
Una tercera simulación muestra que existen casos en los que aparecen tres o más puntos de equilibrio en la frontera. En relación con lo que se establece en el segundo inciso del Lema \ref{lema3.6}, la existencia de puntos de equilibrio positivo no implica necesariamente la inestabilidad de los puntos de equilibrio en la frontera. Sin embargo, se puede observar que el único escenario en el que hay puntos de equilibrio estables en la frontera, y puntos de equilibrio positivo, es cuando hay tres puntos de equilibrio en la frontera. Los parámetros están dados por
\begin{eqnarray}
r=0.31, \quad K=4.5, && q_1=1.7, \quad a_1=3, \nonumber \\
q_2=2, \quad a_2=1, && c_1=2, \quad c_2=1, \nonumber \\
\mu _1=0.3, \quad \mu _2=0.95, && m_1=0.98, \quad m_2=0.8. \nonumber
\end{eqnarray}

Aquí aparecen siete puntos de equilibrio en total. Seis de ellos en la frontera, y únicamente uno en el interior. El punto de equilibrio positivo $E^* = (0.906473,0.539415,0.000961756)$ es un punto silla. A pesar de la existencia de un punto de equilibrio positivo, el punto frontera $E_1 = (0.702991,0.478364,0)$ es un atractor en el sistema \eqref{3.1}. Nótese que esto pudo ocurrir debido a que \eqref{3.1} tiene más de un punto de equilibrio con $x,y > 0$, $z=0$.
\begin{figure}[h]
    \centering
    \includegraphics[width=171pt]{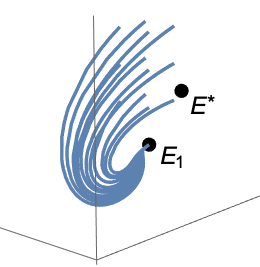}
    \caption{Dinámica local del nodo atractor de \eqref{3.1} ubicado en la frontera. A pesar de la existencia de puntos de equilibrio positivos, es posible hallar un atractor no interior. En este caso, los puntos de equilibrio tienen ubicación
    $E_1 = (0.702991,0.478364,0)$, $E^* = (0.906473,0.539415,0.000961756)$, y valores propios $\sigma(E_1) = \{-0.32693,-0.124404,-0.0128568\}$, $\sigma(E^*) = \{-0.428367,0.00871803 + 0.0231181i,0.00871803 - 0.0231181i\}$.}
    \label{}
\end{figure}
\end{ej}
\section{Conclusiones}
A través de un análisis del sistema \eqref{3.1}, fue posible conocer las condiciones para la existencia de puntos de equilibrio positivo y para la estabilidad de estos últimos, y se dieron condiciones para la existencia de un umbral de bifurcación. El Teorema \ref{teo3.9} es una generalización de la Proposición \ref{propa.1} (del Apéndice \ref{apendicea}) en la que aparece el modelo de Bazykin \eqref{a.1} (el modelo de Bazykin representa el entorno previo a la llegada del competidor foráneo). Esto permite notar la compatibilidad que existe entre los modelos \eqref{3.1} y \eqref{a.1}. En el Ejemplo \ref{ej3.12} se puede apreciar que el equilibrio positivo $E^*$ es inestable, pero el submodelo \eqref{a.1} correspondiente (con $x$, $y$) tiene un punto de equilibrio positivo estable. ¿Será posible hallar un caso en el que el punto de equilibrio positivo sea inestable, pero que en los dos submodelos inducidos $(x,y)$, $(x,z)$ el punto de equilibrio positivo sea estable? ¿O viceversa (un modelo con un equilibrio positivo estable, cuyos submodelos inducidos tengan ambos un equilibrio positivo inestable)?  

En contraste con lo hecho en \cite{Sarwardi2013DynamicalRefuge} (donde se trabaja con competencia interespecífica y sin competencia intraespecífica, y se obtiene un punto de equilibrio positivo que no es atractor en ninguno de los casos analizados), el modelo \eqref{3.1} presenta coeficientes no-nulos de competencia intraespecífica, y se proveen condiciones para que un punto de equilibrio positivo sea atractor. El análisis hecho muestra que existe una gran diferencia entre la competencia intraespecífica y la competencia interespecífica. Con base en lo anterior, es válido preguntarse cuál es la proporción entre estas competencias para poder asegurar la estabilidad asintótica del punto de equilibrio positivo del sistema. Esto podría ser analizado, por ejemplo, desde la perspectiva de la estabilización dispersiva (véase \cite{Metivier2010DispersiveStabilization}).

\part{Invasión: dos presas y un depredador}
\label{parteii}
\chapter*{Parte \ref{parteii} \\ $ $ \\ Invasión: dos presas \\ y un depredador}
Algunas de las preguntas que surgen cuando se estudian un ecosistema\index{ecosistema} y las especies que lo habitan son las siguientes: ¿Cuántas especies en el mismo nivel trófico pueden coexistir en el hábitat? ¿Pueden compartir dos o más especies el entorno aún si estas dependen del mismo recurso alimentario (ver por ejemplo \cite{Armstrong1976})? ¿Cuáles son las facilidades o las dificultades que encontrará una especie al arribar a un medio en el que ya se encuentran otras especies previamente establecidas? Si dos especies mantienen un conjunto de interacciones que les permiten coexistir, ¿cuáles son las características que debe cumplir una especie migrante para poder mezclarse con las nativas sin que esto implique la desaparición de alguna? En particular, si dos competidores comparten un ecosistema\index{ecosistema}, y hacen frente a un depredador\index{depredador} común, ¿cuáles son los rasgos que definirán a la presa\index{presa} mejor adaptada? ¿Cuál de las dos presas\index{presa} tiene mejores posibilidades de sobrevivir? ¿Qué condiciones se deben cumplir para la coexistencia de las tres especies ya sea en equilibrio, o bien con fluctuaciones en las densidades? Diversos modelos matemáticos han sido propuestos y analizados, lo que ha permitido obtener conclusiones, pero también ha dado pie a nuevas interrogantes y a una nutrida discusión. 

En esta parte del trabajo se analizan las consecuencias que trae consigo la llegada de un nuevo competidor al ecosistema\index{ecosistema}. Para ello, se plantea un modelo matemático en el que se describen las poblaciones de tres especies. Una de las especies que aparece en el modelo, es depredadora\index{depredador} de las otras dos (problema similar al que se estudia en \cite{Elettreby2009b}). Además las presas\index{presa} se encuentran en constante competencia intra e interespecífica (ver \cite{Ayala1973} y \cite{Leon}), y se asume que su crecimiento en ausencia de las demás especies está descrito por un modelo denso-dependiente que se conoce comúnmente como modelo logístico\index{modelo logístico}. Para la actividad del depredador\index{depredador} se usa una respuesta funcional de tipo Holling II, que refleja un nivel de saturación en la tasa de captura del depredador\index{depredador}.

En el Capítulo \ref{capitulo4}, se muestra que el modelo matemático es disipativo y persistente, y se lleva a cabo un análisis del modelo similar al que realiza Peter A. Abrams en \cite{Abrams1999} para conocer las características que debe reunir la especie invasora que le permita coexistir con las especies residentes. También en el Capítulo \ref{capitulo4} se analiza la dinámica local del modelo y la existencia de bifurcaciones; en el Capítulo \ref{capitulo5} se aborda el tema de la dinámica caótica.
\chapter{Características e interpretación del modelo}
\label{capitulo4}
\begin{flushright}
\emph{``La interacción de la competencia\\con la depredación puede influir\\enormemente en la coexistencia\\o en la exclusión de las especies,\\y por lo tanto en la composición\\de una comunidad.''\\$ $\\Robert Holt, Gary Polis \cite{Holt1997b}}
\end{flushright}

A continuación se hace uso de los sistemas dinámicos y de la modelación matemática para estudiar la invasión de un ecosistema\index{ecosistema}. En la literatura hay algunas ideas que están en conexión con esto. Por ejemplo, en \cite{ShivaReddy2016DynamicsFluctuations} se analiza un caso en el que el aprovechamiento de un depredador\index{depredador} (que se alimenta de dos presas\index{presa}) viene acompañado de un retardo, y se muestra la existencia de puntos de equilibrio, de ciclos límite y de ruido blanco; otro modelo en el que se trabaja con dos presas\index{presa} y un depredador\index{depredador} aparece en \cite{Ali2017DynamicsHelp}, donde se incorporan coeficientes que representan cooperación entre las especies que enfrentan al depredador\index{depredador}, y algunos umbrales de bifurcación de Hopf son mostrados, así como condiciones para la estabilidad de los puntos de equilibrio del sistema; también se trabaja un modelo con coeficientes de cooperación en \cite{Elettreby2009b}; en \cite{Lee2008THEMODEL} aparece un modelo con dos presas\index{presa} en el que una de ellas se ve beneficiada al convertir en alimento los cadáveres de la otra; en \cite{WEI2010} se bosquejan diagramas de bifurcación y atractores extraños para un modelo de dos presas\index{presa} y un depredador\index{depredador} similar al modelo que se presenta aquí.

Los trabajos anteriores motivan el planteamiento un nuevo problema: una especie exótica llega a un ecosistema\index{ecosistema}, con una población baja en comparación con las de las especies nativas. Se propone usar los sistemas dinámicos para analizar las curvas poblacionales de las especies (endémicas y foráneas). Concretamente, se presenta un modelo determinista (triespecífico) de depredación y competencia. Las tres ecuaciones diferenciales ordinarias que conforman este modelo matemático representan la tasa de crecimiento poblacional de las especies. Las dos primeras corresponden a los dos competidores, y la última corresponde a un depredador\index{depredador}. El crecimiento de las poblaciones de las dos presas\index{presa} (denotadas con $x$, $y$) en ausencia de las otras especies en el medio está dado por el modelo logístico\index{modelo logístico} (también conocido como modelo de Verhulst). Las poblaciones $x$, $y$ disminuyen a causa de la depredación, y esto se representa con la respuesta funcional Holling II (también  conocida como Michaelis-Menten). Los factores que intervienen en la variación de la población del depredador\index{depredador} (denotada con $z$) son la natalidad\index{natalidad} (que es proporcional a la depredación ejercida), la tasa de mortalidad y la competencia intraespecífica. Si las dos presas\index{presa} $x$, $y$ se ausentan en el medio, entonces $z$ enfrentará la extinción\index{extinción}. El modelo propuesto es
\begin{eqnarray}
\label{4.1}
 \dot x &=& x\left(r_1 \left(1 - x/K_1\right) - \alpha_{12}y - \frac{q_1 z}{1 + a_1x} \right),                                                                                                  \nonumber \\
 \dot y &=& y\left(r_2 \left(1 - y/K_2\right) - \alpha_{21}x - \frac{q_2z}{1 + a_2y} \right),                                                                                                                     \\
 \dot z &=& z\left(\frac{c_1 q_1 x}{1 + a_1 x} + \frac{c_2q_2y}{1 + a_2y} - \mu -mz\right),                                                                                              \nonumber \\ \nonumber \\
& & \text{con }x(0)>0\text{, }y(0)>0\text{, }z(0)>0.  \nonumber 
\end{eqnarray}

Se asume que los parámetros de \eqref{4.1} dados por $r_i$, $K_i$, $a_i$, $c_i$, $\mu$, $m$ son siempre positivos, y que los parámetros $\alpha_{ij}$ son no negativos.

Los parámetros $r_i$ representan la tasa intrínseca de crecimiento de las presas\index{presa}; los parámetros $K_i$, se conocen comúnmente como capacidades de carga del medio; la competencia interespecífica de $x$, $y$ está dada por los parámetros $\alpha_{ij}$; los parámetros $a_i$ son conocidos como tiempo de manejo; los parámetros $c_i$ indican la eficiencia en la conversión de biomasa capturada; finalmente $\mu$ es la tasa de mortalidad de la especie consumidora y $m$ representa la competencia intraespecífica por parte del consumidor.
\section{Propiedades generales del modelo}
\label{seccion4.1}

Los puntos de equilibrio en la frontera del sistema \eqref{4.1} están dados por: $E_0 = (0,0,0)$, $E_1 = (K_1,0,0)$, $E_2 = (0,K_2,0)$, $E_3 = (x_3,y_3,0)$, $E_4 = (0,y_4,z_4)$, $E_5 = (x_5,0,z_5)$, con $x_3$, $y_3$, $y_4$, $z_4$, $x_5$, $z_5 > 0$. Obsérvese que puede haber más de un punto de la forma $E_4$ o de la forma $E_5$ para \eqref{4.1}. En esta sección solo se analizan los puntos de equilibrio de frontera. En la Sección \ref{seccion4.4} se analiza $E^* = (x^*,y^*,z^*)$ que es el punto de equilibrio positivo de \eqref{4.1}. La relevancia ecológica de $E^*$ se debe a que representa un estado de coexistencia de las especies, y por lo tanto se buscarán condiciones que propicien la estabilidad asintótica de ese estado.

Al igual que en los Capítulos \ref{capitulo2} y \ref{capitulo3}, se introducen los umbrales de eficiencia energética (u.e.e.) para \eqref{4.1}. Recuérdese que los u.e.e. muestran la mínima población de recurso que debe existir en el medio (en caso de que el otro recurso esté ausente) para que el consumidor pueda sobrevivir.

\begin{df} \label{def4.1} \cite{Kuang2003BiodiversityCompetition}
Los \textbf{umbrales de eficiencia energética (u.e.e.)} en el sistema \eqref{4.1}, que se denotan con $x^b$, $y^b$ están dados por
\begin{eqnarray}
x^b \quad = \quad \frac{\mu}{c_1 q_1 - a_1 \mu}, & & y^b \quad = \quad \frac{\mu}{c_2 q_2 - a_2 \mu}.             \nonumber
\end{eqnarray}
\end{df}
\begin{teo}
\label{teo4.2}
En el sistema \eqref{4.1}:
\begin{enumerate}
    \item El punto de equilibrio $E_0 = (0,0,0)$ es un punto silla.
    \item Si $\alpha_{21} < r_2/K_1$ o bien $0 < x^b < K_1$, entonces $E_1 = (K_1,0,0)$ es un punto silla. 
    \item Si $\alpha_{12} < r_1/K_2$ o bien $0 < y^b < K_2$, entonces $E_2 = (0,K_2,0)$ es un punto silla.
\end{enumerate}
Si las condiciones en 2. o en 3. no se cumplen, entonces el punto correspondiente es un atractor local.
\end{teo}
\begin{proof}
La matriz jacobiana de \eqref{4.1} está dada por

\begin{eqnarray}
J(x,y,z) &=& \left( \begin{array}{ccc}
                r_1\left(1 - 2x/K_1\right) - \alpha_{12}y & &   \\
                - q_1 z(1 + a_1 x)^{-2} & - \alpha_{12} x    & - q_1 x(1 + a_1 x)^{-1} \\\\
                & r_2\left(1 - 2y/K_2\right) - \alpha_{21}x     & \\
                - \alpha_{21} y & - q_2 z(1 + a_2 y)^{-2}  & - q_2 y(1 + a_2 y)^{-1} \\\\
                 & & c_1 q_1 x(1 + a_1 x)^{-1} \\
                 c_1 q_1 z(1 + a_1 x)^{-2}    & c_2 q_2 z(1 + a_2 y)^{-2} & + c_2 q_2 y(1 + a_2 y)^{-1}  \\
                & & - \mu - 2 m z  \end{array} \right).                                  \nonumber
\end{eqnarray}
Si se evalúa la matriz jacobiana de $E_0 = (0,0,0)$ se obtiene
\begin{eqnarray}
J(E_0) &=& \left( \begin{array}{ccc}
                r_1 & 0   & 0    \\
                0   & r_2 & 0    \\
                0   & 0   & -\mu 
                \end{array} \right),                                  \nonumber
\end{eqnarray}
y siguiendo el mismo procedimiento con $E_1 = (K_1,0,0)$ se llega a

\begin{eqnarray}
J(E_1) &=& \left( \begin{array}{ccc}
                - r_1 & - K_1 r_1             & - K_1 q_1(1 + a_1 K_1)^{-1}         \\\\
                0     & r_2 - K_1 \alpha_{21} & 0                                     \\\\
                0     & 0                     & c_1 K_1 q_1(1 + a_1 K_1)^{-1} - \mu \\
                \end{array} \right).                                  \nonumber
\end{eqnarray}
De los cálculos previos, se sigue el resultado.

\end{proof}
Además de las posibilidades mostradas en el Teorema \ref{teo4.2}, aún hace falta estudiar otras tres configuraciones para los puntos de equilibrio en la frontera  $E_3 = (x_3,y_3,0)$, $E_4 = (0,y_4,z_4)$ y $E_5 = (x_5,0,z_5)$, con $x_3$, $y_3$, $y_4$, $z_4$, $x_5$, $z_5 > 0$. 

\begin{teo}
\label{teo4.3}
Sean
\begin{eqnarray}
p_1(x) \quad = \quad \frac{c_1 q_1 x}{1 + a_1 x}, && p_2(y) \quad = \quad \frac{c_2 q_2 y}{1 + a_2 y}. \nonumber
\end{eqnarray}
En el sistema \eqref{4.1} se cumple:
\begin{enumerate}
    \item Existe un único punto de equilibrio $E_3=(x_3,y_3,0)$, con $x_3$, $y_3 > 0$ si y solo si $r_1/K_2 > \alpha_{12}$ y $r_2/K_1 > \alpha_{21}$. Si $p_1(x_3) + p_2(y_3) > \mu$, entonces $E_3$ es inestable.
    \item Si $0 < x^b < K_1$, entonces existe al menos un punto de equilibrio de la forma $E_4 = (0,y_4,z_4)$, con $y_4$, $z_4 > 0$. Si $E_4 = (0,y_4,z_4)$, con $y_4$, $z_4 > 0$, es un punto de equilibrio, y $r_1 > \alpha_{12}y_4 + q_1z_4$, entonces $E_4$ es inestable.
    \item Si $0 < y^b < K_2$, entonces existe al menos un punto de equilibrio de la forma $E_5 = (x_5,0,z_5)$, con $x_5$, $z_5 > 0$. Si $E_5 = (x_5,0,z_5)$, con $x_5$, $z_5 > 0$, es un punto de equilibrio, y $r_2 > \alpha_{21}x_5 + q_2z_5$, entonces $E_5$ es inestable.
\end{enumerate}
\end{teo}
\begin{proof}
Es fácil ver que los respectivos puntos de equilibrio existen si se cumplen las condiciones indicadas. 

El último renglón de la matriz $D(E_3)$ es $(0,0,p_1(x_3)+p_2(y_3)-\mu)$, por lo tanto $D(E_3)$ tiene un valor propio positivo y por lo tanto $E_3$ es asintóticamente inestable.

Por otra parte, se tiene que el primer renglón de la matriz $D(E_4)$ es $(r_1 - \alpha_{12}y_4 - q_1z_4,0,0)$, y entonces $E_4$ es asintóticamente inestable. Un argumento similar funciona para $E_5$.
\end{proof}
En el Teorema \ref{teo4.3} se encuentran condiciones suficientes para la inestabilidad de algunos de los puntos de equilibrio en la frontera. Al final del capítulo se muestra que en los puntos de equilibrio $E_3$, $E_4$ y $E_5$ se puede presentar también estabilidad local.

A continuación se demuestra que el sistema \eqref{4.1} es disipativo.
\begin{teo} Existe un subconjunto compacto $C$ de $[0,\infty)^3$ de tal manera que para toda solución no-negativa $(x,y,z)$ de \eqref{4.1}, existe $T > 0$ de tal forma que $(x(t),y(t),z(t)) \in C$ para toda $t \geq T$.
\end{teo}
\begin{proof}
De forma similar a como ya se hizo en las demostraciones de los Teoremas \ref{teo2.1} y \ref{teo3.5}, se toma 
\begin{eqnarray}
w &=& c_1 x + c_2 y + z, \quad \varphi \quad \in \quad (0,\mu], \nonumber \\ \displaystyle \rho &=& \frac{K_1(r_1 + \varphi)^2}{4 r_1} + \frac{K_2(r_2 + \varphi)^2}{4 r_2}. \nonumber
\end{eqnarray}
Entonces se tiene
\begin{eqnarray}
w' + \varphi w &=& c_1(x (r_1 (1 - x/K_1) - \alpha_{12}y)) + c_2(y (r_2 (1 - y/K_2) - \alpha_{21}x))  \nonumber \\
& & - (\mu z + m z^2) + \varphi(c_1 x + c_2 y + z) \nonumber \\
&\leq& c_1(x \left(r_1 + \varphi - \left(\frac{r_1}{K_1}\right)x\right) + c_2(y \left(r_2 + \varphi - \left(\frac{r_2}{K_2}\right)y\right) - (\mu - \varphi)z \nonumber \\
&\leq& c_1(x \left(r_1 + \varphi - \left(\frac{r_1}{K_1}\right)x\right) + c_2(y \left(r_2 + \varphi - \left(\frac{r_2}{K_2}\right)y\right) \nonumber \\
&\leq& \frac{c_1K_1(r_1 + \varphi)^2}{4 r_1} + \frac{c_2K_2(r_2 + \varphi)^2}{4 r_2}\nonumber \\
&=& \rho. \nonumber
\end{eqnarray}

Lo anterior da lugar a  $\displaystyle w \leq e^{-\varphi t}w_0 + \frac{\rho}{\varphi}(1 - e^{-\varphi t})$. 

Y así $\displaystyle \limsup_{t\to \infty} w(t) \leq \frac{\rho}{\varphi}$. 

Por lo tanto, para cualquier solución no-negativa $(x,y,z)$ de \eqref{4.1}, existe $T>0$ tal que
\begin{eqnarray}
(x(t),y(t),z(t)) &\in& \left\{(x,y,z)\in [0,\infty)^3\colon c_1 x + c_2 y + z \leq \frac{\rho}{\varphi}\right\}\text{, si } t \geq T. \nonumber
\end{eqnarray}
\end{proof}

Obsérvese que al escribir $y = 0$ o $z=0$ en \eqref{4.1}, se obtiene el modelo de Bazykin, (véase el Apéndice \ref{apendicea}). Por otro lado, un modelo de competencia de tipo Lotka-Volterra (véase por ejemplo \cite{buzzi}, \cite{Goh1978}, \cite{Korobeinikov1999}, \cite{Llibre2014}, \cite{soliman}) se obtiene de hacer $z = 0$ en \eqref{4.1}.

En la siguiente sección se aborda el problema de la invasión que inicialmente fue planteado. Este análisis se hace a través del modelo matemático \eqref{4.1}.
\section{Interpretación: Aparición de un nuevo competidor}
\label{seccion4.2}
Como ya se sabe, la biodiversidad y la riqueza de un hábitat es un tema que ha sido estudiado por reconocidos autores (ver por ejemplo \cite{Kuang2003BiodiversityCompetition, Lehman2000, Tilman2004}). Esta diversidad se puede ver mermada cuando se introduce una especie ajena a un ecosistema\index{ecosistema}, pues la competencia de las especies foráneas con las endémicas puede dar lugar a la disminución o a la desaparición de alguna especie en el ecosistema\index{ecosistema}. Ejemplo de lo anterior es el impacto que tuvo en el ajolote (\textit{Ambystoma mexicanum}) la introducción de la carpa \textit{Cyprinus carpio} y la tilapia \textit{Oreochromis niloticus} al lago de Xochimilco \cite{Contreras2009, Robles2009, Zambrano2007, Zambrano2010}. La introducción de especies exóticas ha colocado al ajolote al borde de la extinción\index{extinción}. Otro ejemplo es la tasa de mortalidad del crustáceo \textit{Chthamalus stellatus} que se ve afectada severamente en la temporada de alta reproductividad de su competidor \textit{Balanus balanoides}; este hecho se muestra a través de experimentos en \cite{Connell1961}. Un caso similar se presenta en la planta \textit{Centaurea maculosa}, la cual afecta a su competidora \textit{Arabis fecunda} a través de la secreción de una sustancia \cite{Lesica2016}.

Para estudiar la introducción de nuevas especies competidoras, se usa el sistema \eqref{4.1}, y se asume que la presa\index{presa} (cuya población se denota con $x$) se encuentra en coexistencia con su depredador\index{depredador} (con población denotada con $z$). Bajo estas condiciones, la especie colonizadora aparece en el nicho ecológico con una población $y$ reducida en comparación con las poblaciones $x$ y $z$. De acuerdo con las ideas de Abrams \cite{Abrams1999}, el problema se enfoca partiendo de la tasa intrínseca de crecimiento de $y$. 

\begin{obs}
En el sistema \eqref{4.1}, el crecimiento de la población $y$ está dado por
\begin{eqnarray}
\label{4.2}
\dot y &=& y\left(r_2\left(1 - \frac{y}{K_2}\right) - \alpha_{12} x - \frac{q_2 z}{1 + a_2 y}\right).
\end{eqnarray}
\end{obs}
En \eqref{4.2}, $y$ es un número positivo. Entonces el signo de $\dot y$ es el mismo signo de
\begin{eqnarray}
r_2\left(1 - \frac{y}{K_2}\right) - \alpha_{12} x - \frac{q_2 z}{1 + a_2 y}.               \nonumber
\end{eqnarray}
Por lo tanto, se tiene lo siguiente
\begin{obs}
La primera derivada $\dot y$ satisface $\dot y >0$ si y solamente si 
\begin{eqnarray}
\label{4.3}
r_2 &>&  \frac{\alpha_{21} x + q_2 z (1 + a_2 y)^{-1}}{1 - y/K_2}.
\end{eqnarray}
\end{obs}
Como ya se mencionó anteriormente, las especies $x$, $z$ se encuentran en coexistencia. Se puede suponer que en el submodelo correspondiente (que es el que aparece en el Apéndice \ref{apendicea})
\begin{eqnarray}
\label{4.4}
\dot x &=& x\left(r_1(1 - x/K_1) - \frac{q_1 z}{1 + a_1 x}\right), \nonumber \\
\dot z &=& z\left(\frac{c_1 q_1 x}{1 + a_1 x} - \mu - m z\right), 
\end{eqnarray}
existe un estado estacionario
\begin{eqnarray}
\widetilde E &=& (\widetilde x,\widetilde z). \nonumber
\end{eqnarray}
Supóngase adicionalmente que $\widetilde E$ es un atractor global del sistema en \eqref{4.4} (véase la Proposición \ref{propa.1}). Si la especie $y$ ingresa al entorno con una población escasa, entonces los valores $y/K_2$, $a_2 y$ se pueden despreciar en el lado derecho de la desigualdad \eqref{4.3}. De lo anterior, se sigue que la tasa de crecimiento $\dot y$ es positiva si y solo si
\begin{eqnarray}
r_2 &>& \alpha_{21} \widetilde x + q_2 \widetilde z.                                                         \nonumber
\end{eqnarray}
Nótese la similitud que hay entre la desigualdad anterior y la desigualdad obtenida en el tercer punto del Teorema \ref{teo4.3}. 

En el submodelo \eqref{4.4} se puede verificar que $\widetilde z = r_1q_1^{-1} (1 - \widetilde x/K_1)(1 + a_1\widetilde x)$. Lo anterior da lugar a lo siguiente.
\begin{obs}
\label{obs4.7}
Si $\widetilde E = (\widetilde x, \widetilde z)$ es un punto de equilibrio positivo del submodelo \eqref{4.4} que es atractor global, entonces se cumple $\dot y > 0$ si y solo si
\begin{eqnarray}
\label{4.5}
r_2 &>& \alpha_{21} \widetilde x + \frac{r_1 q_2}{q_1}\left(1 - \frac{\widetilde x}{K_1}\right)\left(1 + a_1\widetilde x\right).
\end{eqnarray}
\end{obs}
Si se cumple la condición \eqref{4.5}, entonces se concluye que la introducción de $y$ al medio será exitosa. La especie $y$ podrá mezclarse con las especies residentes si las características del medio y de las especies se mantienen similares a como eran al momento de la introducción. En \eqref{4.5}, el valor de $\widetilde x$ está dado por una raíz del polinomio de tercer grado (que es el mismo que aparece en la demostración del Teorema \ref{teoa.2})

\begin{eqnarray}
P(x) &=& -a_1^2mr_1x^3 + (-2a_1mr_1 + a_1^2K_1mr_1) x^2                                 \nonumber \\
     & &+ (a_1\mu K_1q_1 - c_1K_1q_1^2 - m r_1 + 2a_1K_1mr_1) x + (\mu K_1 q_1 + K_1 m r_1).                                                                                             \nonumber 
\end{eqnarray}
La fórmula de Cardano para ecuaciones de tercer grado implica que uno de los posibles valores para $\widetilde x$ es
\begin{eqnarray}
\widetilde x &=& \frac{\chi}{3 a_1^2 r_1 m}-\frac{\sqrt[3]{\psi +\sqrt{\psi^2+4 \left(\varphi -\chi^2\right)^3}}}{3 \sqrt[3]{2}, a_1^2 r_1 m}                             \nonumber \\
    & & +\frac{\sqrt[3]{2} \left(\varphi -\chi^2\right)}{3 a_1^2 r_1 m \sqrt[3]{\psi + \sqrt{\psi^2+4 \left(\varphi -\chi^2\right){}^3}}},                                     \nonumber
\end{eqnarray}
donde 

\begin{eqnarray}
\varphi &=& 3 a_1^2 m r_1\left(c_1 K_1 q_1^2 + m r_1 - 2 a_1 K_1 m r_1 - a_1 K_1 q_1 \mu\right),                                                                                                                              \nonumber \\ 
\chi &=& a_1^2 K_1 m r_1 - 2 a_1 m r_1,                                             \nonumber \\ 
\psi &=& - 18 a_1^3 c_1 K_1 m^2 q_1^2 r_1^2 + 9 a_1^4 c_1 K_1^2 m^2 q_1^2 r_1^2 - 2 a_1^3 m^3 r_1^3 - 6 a_1^5 K_1^2 m^3 r_1^3                                        \nonumber \\
     & & - 6 a_1^4 K_1 m^3 r_1^3 - 2 a_1^6 K_1^3 m^3 r_1^3 - 9 a_1^4 K_1 m^2 q_1 r_1^2 \mu 
         - 9 a_1^5 K_1^2 m^2 q_1 r_1^2 \mu.                                       \nonumber
\end{eqnarray}
El modelo \eqref{4.4} puede presentar soluciones que se aproximen hacia un ciclo límite, y no a un punto de equilibrio asintóticamente estable (se puede revisar el Apéndice \ref{apendicea}). Si el modelo presenta soluciones fluctuantes (soluciones que se aproximan asintóticamente a un ciclo límite), entonces un procedimiento similar permite hallar condiciones para la llegada exitosa de la especie $y$. Se toman en cuenta las respectivas medias temporales a largo plazo para cada una de las especies en cuestión, las cuales se denotan con $\left<x\right>$ y $\left<z\right>$. Se puede asumir que el crecimiento medio de $x$ y de $z$ a largo plazo es cero, y así se obtiene un resultado análogo al de la Observación \ref{obs4.7}.
\begin{obs}
\label{obs4.8}
Si $x$, $z$ coexisten con una dinámica fluctuante, entonces la invasión de $y$ será exitosa ($\dot y > 0$) si y solo si
\begin{eqnarray}
r_2 &>& \alpha_{12} \left<x\right> + q_2 \left<z\right>,                                \nonumber
\end{eqnarray}
o de forma equivalente
\begin{eqnarray}
r_2 &>& \alpha_{12} \left<x\right> + \frac{r_1 q_2}{q_1}\left<\left(1 - \frac{x}{K_1}\right)\left(1 + a_1x\right)\right>.                                                                     \nonumber
\end{eqnarray}
\end{obs}
\begin{ej}
Para ilustrar las Observaciones \ref{obs4.7} y \ref{obs4.8} se considera el sistema \ref{4.1} con los parámetros y las condiciones iniciales
\begin{eqnarray}
\label{4.6}
r_1 \quad = \quad 1,   \quad K_1         &=& 6,     \quad q_1         \quad = \quad 1.5,    \nonumber \\
a_1 \quad = \quad 1,   \quad c_1         &=& 1,     \quad K_2         \quad = \quad 4,      \nonumber \\
q_2 \quad = \quad 1,   \quad a_2         &=& 1,     \quad c_2         \quad = \quad 1,       \\
\mu \quad = \quad 1,  \quad \alpha_{12}       &=& 0.1,   \quad \alpha_{21} \quad = \quad  0.1,   \nonumber  \\ \nonumber \\
 \quad x_0 \quad = \quad 1, \quad y_0  &=& 0.01, \quad z_0         \quad = \quad 1.    \nonumber
\end{eqnarray}

El valor de $r_2$ es tomado como $2$ o como $0.2$, mientras que el valor de $m$ es tomado como $1$ o como $0.01$.
\begin{figure}[h]
    \centering
    \includegraphics[width=242pt]{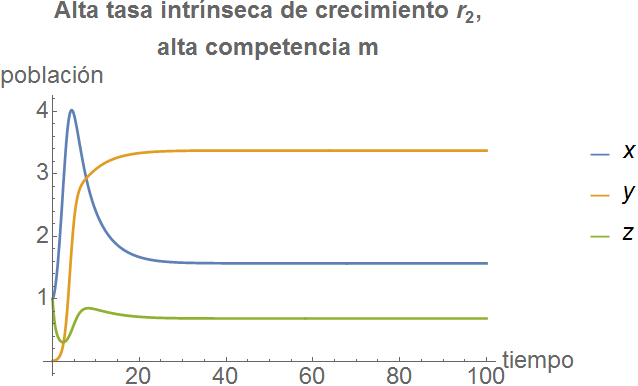}
    \includegraphics[width=242pt]{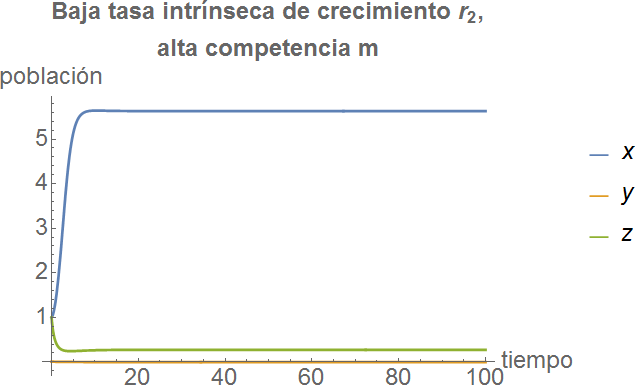}
    \includegraphics[width=242pt]{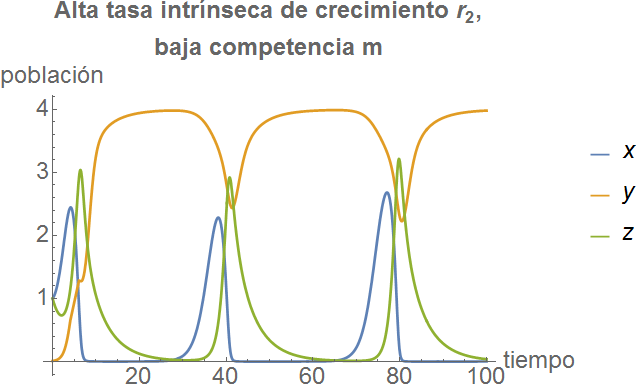}
    \includegraphics[width=242pt]{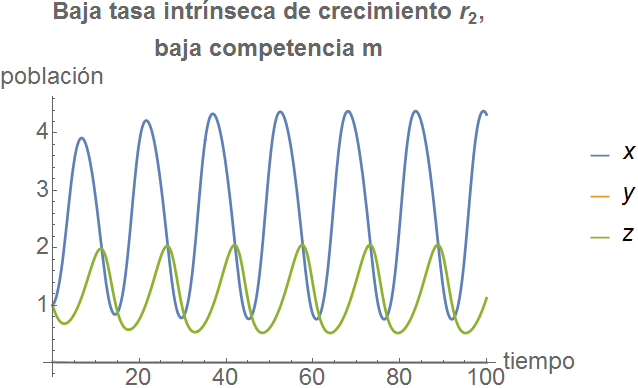}
    \caption{Soluciones del sistema \eqref{4.1} con los parámetros y las condiciones iniciales \eqref{4.6}}
    \label{fig4.1}
\end{figure}

En el primer diagrama de la Figura \ref{fig4.1} se puede apreciar el caso en el que $r_2=2$ y $m=1$. Entonces $\widetilde x = 5.62837$, y el lado derecho de la desigualdad \eqref{4.5} es igual a $0.836537$. Por lo tanto, la desigualdad \eqref{4.5} se cumple, y esto permite que las tres especies alcancen la coexistencia en el entorno. En el segundo caso se considera $r_2=0.02$ y $m=1$. Otra vez se tiene $\widetilde x = 5.62837$, y en este caso la desigualdad \eqref{4.5} no se satisface. A partir de lo anterior, se concluye que la especie denotada con $y$ no cuenta con las características necesarias para competir con las especies residentes. En el tercer y en el cuarto diagrama, se presentan escenarios en los que el coeficiente de competencia del depredador\index{depredador} $m$ es igual a $0.01$. Las curvas no tienden hacia un equilibrio, sino que se aproximan a un ciclo límite. Si $r_2$ es suficientemente grande, entonces la especie $y$ conseguirá ingresar al medio. En caso contrario, la especie $y$ no podrá coexistir con las especies residentes.
\end{ej}
En la Sección \ref{seccion4.3} se muestra como el u.e.e. no juega el mismo papel en los modelos biespecíficos que en los modelos triespecíficos.
\section{Análisis exploratorio}
\label{seccion4.3}
Supóngase que en un escenario biespecífico se da la coexistencia de las especies. Entonces es válido cuestionarse si al llevar ese mismo problema a un nivel triespecífico se podrá concluir lo mismo, o bien si se presentará la exclusión de alguna de las especies. 

Supóngase ahora que se toman en cuenta dos submodelos \eqref{4.4} del modelo original \eqref{4.1}. El primer submodelo está conformado por las especies $x$, $z$
\begin{eqnarray}
\dot x &=& x\left(r_1(1 - x/K_1) - \frac{q_1 z}{1 + a_1 x}\right), \nonumber \\
\dot z &=& z\left(\frac{c_1 q_1 x}{1 + a_1 x} - \mu - m z\right), \nonumber
\end{eqnarray}
y el segundo conformado por las especies $y$, $z$
\begin{eqnarray}
\dot y &=& y\left(r_2(1 - y/K_2) - \frac{q_2 z}{1 + a_2 y}\right), \nonumber \\
\dot z &=& z\left(\frac{c_2 q_2 y}{1 + a_2 y} - \mu - m z\right). \nonumber
\end{eqnarray}

\begin{figure}[h]
    \centering
    \includegraphics[width=179pt]{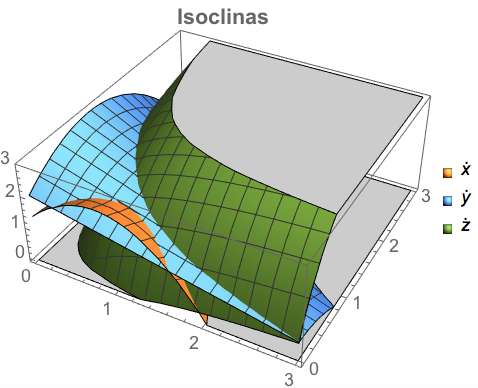}
    \includegraphics[width=179pt]{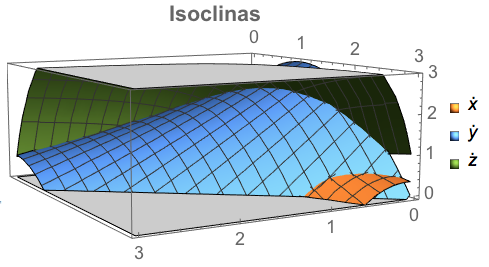}
    \includegraphics[width=179pt]{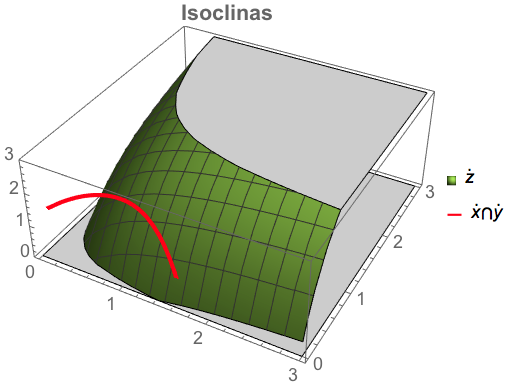}
    \includegraphics[width=179pt]{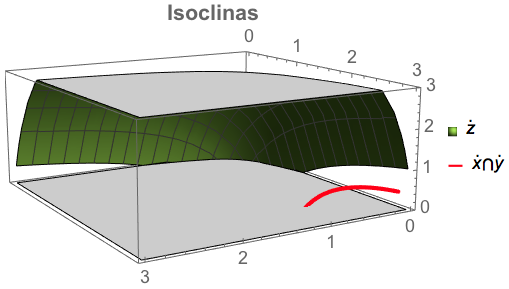}
    \caption{En los diagramas se muestran las superficies isoclinas de \eqref{4.1}. Los parámetros escogidos son $r_1=2$, $q_1=1.5$, $K_1=2$, $a_1=2.3$, $r_2=2$, $q_2=1$, $K_2=3$, $a_2=1.2$, $c_1=2$, $c_2=2$, $\mu =1$, $m=0.3$, $\alpha_{12}=0.5$, $\alpha_{21}=0.5$. Las imágenes en el costado izquierdo muestran la perspectiva frontal del espacio $\mathbb R^3$, mientras que en la derecha aparecen las perspectivas traseras.}
    \label{fig4.2}
\end{figure}

Se puede asumir que las desigualdades 
\begin{eqnarray}
\label{4.7}
0 \quad < \quad x^b \quad < \quad K_1, \nonumber \\
0 \quad < \quad y^b \quad < \quad K_2,
\end{eqnarray}
son ciertas. De acuerdo con la Proposición \ref{propa.1} en el Apéndice \ref{apendicea}, la primera desigualdad \eqref{4.7} implica que las soluciones se mantienen en el interior del cuadrante positivo para el sistema \eqref{4.4}, ya sea con un punto de equilibrio positivo asintóticamente estable, o bien con soluciones fluctuantes. Con base en lo anterior se podría pensar que la primera desigualdad de \eqref{4.7} juega un papel importante en la coexistencia de las especies en \eqref{4.4}, y que por lo tanto las desigualdades \eqref{4.7} podrían tener también un papel similar en el modelo general \eqref{4.1}. Sin embargo, las desigualdades \eqref{4.7} no necesariamente son compatibles con la coexistencia de las especies en el modelo \eqref{4.1}. La Figura \ref{fig4.2} presenta un ejemplo que  muestra como las desigualdades \eqref{4.7} propician la exclusión en \eqref{4.1}. 

\begin{obs}
Las condiciones \eqref{4.7} implican la supervivencia de una de las presas\index{presa} en ausencia de la otra presa\index{presa} y en presencia del depredador\index{depredador}, pero no se puede inferir a partir de eso la coexistencia de las tres especies.
\end{obs} 

En la siguiente sección se analizan las condiciones que deben cumplir los parámetros para que la invasión ocurra de manera exitosa.

Con relación a la cantidad total de puntos de equilibrio en \eqref{4.1}, se tiene lo siguiente.

\begin{obs}
Las desigualdades \eqref{4.7} no implican la existencia de un punto de equilibrio positivo en el sistema \eqref{4.1}. Los parámetros de la especie depredadora\index{depredador} $z$ (que están dados por $c_1$, $c_2$ y $m$) son los que determinan la cantidad de puntos de equilibrio interiores que existen para el sistema \eqref{4.1}. Puede haber ninguno, uno, dos o hasta tres puntos de equilibrio positivos en \eqref{4.1}.
\end{obs}

A continuación los parámetros $K_2$ y $a_2$ son manipulados para conocer las condiciones necesarias para una invasión exitosa.

\begin{figure}[h]
    \centering
    \includegraphics[width=242pt]{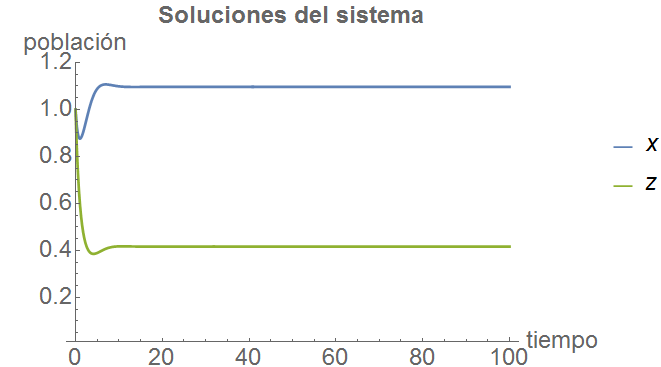}
    \caption{En el modelo \eqref{4.4} con los parámetros y condiciones iniciales \eqref{4.8} se presenta estabilidad asintótica alrededor del punto de equilibrio $E^* = (1.09652,0.41658)$.}
    \label{fig4.3}
\end{figure}

\begin{ej}
Supóngase que una especie y su depredador\index{depredador} se encuentran en coexistencia. Las poblaciones de dichas especies son denotadas con $x$, $z$ respectivamente. Entonces llega al ecosistema\index{ecosistema} una nueva especie con población escasa (denotada con $y$), que compite con $x$ y que es depredada por $z$. Se estudian las soluciones del modelo \eqref{4.1} con los parámetros y condiciones iniciales
\begin{eqnarray}
\label{4.8}
r_1 \quad = \quad 1,   \quad K_1         &=& 1.5,   \quad q_1         \quad = \quad 1,    \nonumber \\ 
a_1 \quad = \quad 0.5, \quad c_1         &=& 2,     \quad r_2         \quad = \quad 1,    \nonumber \\
q_2 \quad = \quad 1,   \quad c_2         &=& 2,     \quad \mu         \quad = \quad 1,              \\
m   \quad = \quad 1,   \quad \alpha_{12} &=& 0.5,   \quad \alpha_{21} \quad = \quad 0.5, \nonumber \\ \nonumber \\
x_0 \quad = \quad 1,   \quad y_0         &=& 0.001, \quad z_0         \quad = \quad 1.    \nonumber
\end{eqnarray}

Si los valores \eqref{4.8} se sustituyen en el submodelo \eqref{4.4}, se obtienen las soluciones que convergen a un estado estacionario en $E^*=(1.09652,0.416458)$, como lo muestra la Figura \ref{fig4.3}. Esto significa que las especies $x$, $z$ alcanzan la coexistencia en ausencia de la especie $y$.

Además de los parámetros exhibidos en \eqref{4.8}, hace falta tomar valores para la capacidad de carga $K_2$ y el tiempo de manejo $a_2$ de la especie depredadora\index{depredador} $z$ respecto a la especie migrante $y$. Los valores de $K_2$ y $a_2$ son tomados dentro de un rango más amplio. En la gráfica aparecen los parámetros $K_2$, $a_2$ en los ejes coordenados, y con un color distinto cada caso según ocurra exclusión del depredador\index{depredador}, exclusión de alguna de las presas\index{presa} o coexistencia de las tres especies.

Para esta simulación se diseñó un programa en \textit{Mathematica}, en el que se calculan las soluciones de \eqref{4.1} con parámetros y condiciones iniciales \eqref{4.8}, y se obtiene la distancia de las soluciones (al tiempo $t=100$) hacia la frontera del espacio $\partial ([0,\infty)^3)$. Con los resultados obtenidos, se construye una gráfica que se muestra en la Figura \ref{fig4.4}. 

\begin{figure}[h]
    \centering
    \includegraphics[width=242pt]{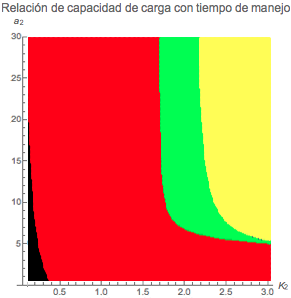}
    \caption{Casos que se presentan de acuerdo a la elección de parámetros $K_2$ y $a_2$. En la región roja es el conjunto de coexistencia, la región verde muestra los casos en los que el depredador\index{depredador} perece, el conjunto negro muestra los casos en el que la especie invasora $y$ es eliminada por las otras especies, y el conjunto amarillo son los casos en el que únicamente la especie $y$ subsiste.}
    \label{fig4.4}
\end{figure}

Después de realizar este experimento surge una pregunta natural: Si las especies endémicas (denotadas con $x$ y $z$) se encuentran de momentáneamente en coexistencia, y la especie inmigrante (denotada con $y$) es introducida al entorno con una tasa intrínseca de crecimiento (denotada con $r_2$) que cumple con la desigualdad de la Observación \ref{obs4.7}, ¿qué condiciones deben cumplir el resto de los parámetros para que las especies endémicas no sean desplazadas del medio? 

\end{ej}
La llegada de la especie invasora puede convertirse en un factor que estabilice el ecosistema\index{ecosistema} pudiendo incluso llegar a salvar de la extinción\index{extinción} a alguna de las especies residentes. Este hecho se muestra en la siguiente simulación. 
\begin{ej}
Considérese el sistema \eqref{4.1} con una selección de parámetros y condiciones iniciales con todos los valores iguales a \eqref{4.8}, con excepción de $a_1$ que cambia a $1.5$.
\begin{eqnarray}
\label{4.9}
r_1 \quad = \quad 1,   \quad K_1         &=& 1.5,   \quad q_1         \quad = \quad 1,    \nonumber \\ 
a_1 \quad = \quad 1.5, \quad c_1         &=& 2,     \quad r_2         \quad = \quad 1,    \nonumber \\
q_2 \quad = \quad 1,   \quad c_2         &=& 2,     \quad \mu         \quad = \quad 1,              \\
m   \quad = \quad 1,   \quad \alpha_{12} &=& 0.5,   \quad \alpha_{21} \quad = \quad 0.5, \nonumber \\ \nonumber \\
x_0 \quad = \quad 1,   \quad y_0         &=& 0.001, \quad z_0         \quad = \quad 1.    \nonumber
\end{eqnarray}

En este caso, las soluciones de \eqref{4.4} con los parámetros \eqref{4.9} no se tienden hacia un estado de coexistencia, como se muestra en la Figura \ref{fig4.5}.

\begin{figure}[h]
    \centering
    \includegraphics[width=242pt]{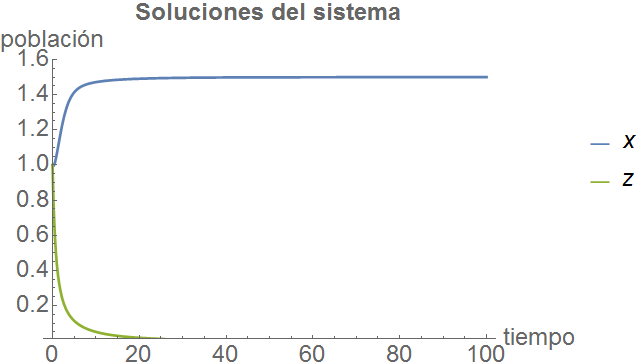}
    \caption{En el modelo \eqref{4.4} con los parámetros y condiciones iniciales \eqref{4.9}, las soluciones se aproximan al punto en la frontera $(1.5,0)$.}
    \label{fig4.5}
\end{figure}

Se procede de forma análoga a la simulación anterior, y se muestran los resultados obtenidos en la Figura \ref{fig4.6}. La pendiente negativa que se forma en la parte derecha de la curva que separa la región verde de la roja es producto de la necesidad que hay de intercambiar capacidad de carga $K_2$ por tiempo de manejo o mecanismo de defensa $a_2$: mientras mayor sea la riqueza o la productividad del ecosistema\index{ecosistema} para la presa\index{presa}, mayor debe ser la eficiencia del depredador\index{depredador} para procesar su alimento, y al revés, si el medio no resulta productivo para la presa\index{presa}, entonces el tiempo de manejo puede ser alto y la coexistencia aún es factible. Nótese que si $K_2$ es menor que $0.15$, entonces el depredador\index{depredador} no sobrevivirá. La llegada de la especie $y$ hace que disminuya la población $x$. La especie $x$ es alimento para la $z$, sin embargo la especie $y$ es también un nuevo alimento para $z$. Otro aspecto que aquí conviene mencionar es que en este escenario se exige una fuerte adaptación de las especies endémicas para que la coexistencia sea alcanzada.
\end{ej}
Estos dos experimentos permiten determinar las regiones de interés en un subespacio de parámetros del sistema \eqref{4.1} y hacer una categorización en el estudio de la colonización de ecosistemas\index{ecosistema}.

\begin{figure}[h]
    \centering
    \includegraphics[width=221pt]{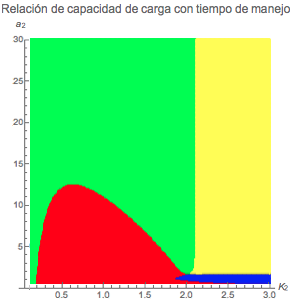}
    \includegraphics[width=221pt]{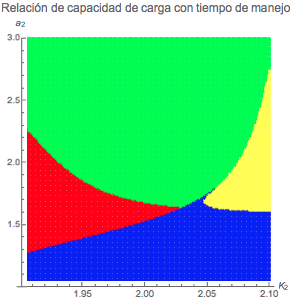}
    \caption{Casos que se presentan de acuerdo a la elección de parámetros $K_2$ y $a_2$. En la región roja es el conjunto de coexistencia, la región verde muestra los casos en los que el depredador\index{depredador} perece, el conjunto azul muestra los casos en el que la presa\index{presa} $x$ es eliminada por las otras especies, y el conjunto amarillo son los casos en el que únicamente la especie $y$ subsiste.}
    \label{fig4.6}
\end{figure}
\section{Dinámica alrededor de los puntos de equilibrio positivos}
\label{seccion4.4}
En la Sección \ref{seccion4.1} se mencionó que el sistema \eqref{4.1} tiene puntos de equilibrio en la frontera en $E_0=(0,0,0)$, $E_1=(K_1,0,0)$ y $E_2=(0,K_2,0)$, los cuales representan la desaparición del depredador\index{depredador} junto con una presa\index{presa}, o junto con las dos presas\index{presa}. En $E_1$ y $E_2$, se aprecia que una presas\index{presa} alcanza la sobrepoblación en ausencia de su depredador\index{depredador} y de su competidor. Además de los casos anteriores, existen otros estados de equilibrio en la frontera, en los que el depredador\index{depredador} y solamente una de las presas\index{presa} sobreviven, o bien en los que ambas presas\index{presa} sobreviven y que el depredador\index{depredador} es eliminado: $E_3=(x_3,y_3,0)$, $E_4=(x_4,0,z_4)$ y $E_5=(0,y_5,z_5)$. Los Teoremas \ref{teo4.2} y \ref{teo4.3} hablan acerca de las propiedades dinámicas de esos puntos de equilibrio.

El tema central en esta sección es el de los puntos de equilibrio positivos del sistema \eqref{4.1}. Como ya se mencionó, estos representan un estado de coexistencia entre las especies de interés. Estos puntos son las soluciones del sistema de tres ecuaciones no lineales y tres incógnitas dado por
\begin{eqnarray}
z &=& \frac{1}{q_1}\left(r_1(1 - x/K_1) - \alpha_{12}y\right)\left(1 + a_1x\right) \quad \colon = \quad f(x,y), \nonumber \\
z &=& \frac{1}{q_2}\left(r_2(1 - y/K_2) - \alpha_{21}x\right)\left(1 + a_2y\right) \quad \colon = \quad g(x,y), \nonumber \\
z &=& \frac{1}{m}\left(\frac{c_1 q_1 x}{1 + a_1 x} + \frac{c_2 q_2 y}{1 + a_2 y} - \mu\right) \quad \quad \quad \colon = \quad h(x,y). \nonumber
\end{eqnarray}
Se iguala $f$ con $g$ para obtener:
\begin{eqnarray}
\label{4.10}
\frac{1}{q_1}\left(r_1(1 - x/K_1) - \alpha_{12}y\right)(1 + a_1 x) &=& \frac{1}{q_2}\left(r_2(1 - y/K_2) - \alpha_{21}x\right)(1 + a_2 y).
\end{eqnarray}
La ecuación anterior determina una hipérbola en $\mathbb R^2$, y si
\begin{eqnarray}
\frac{1}{4}\left(\frac{a_1r_1}{K_1q_1}\left(\frac{r_1}{K_1q_1} - \frac{a_1r_1}{q_1} - \frac{\alpha_{21}}{q_2}\right)^2 + \frac{a_2r_2}{K_2q_2}\left(\frac{r_2}{K_2q_2} - \frac{a_2r_2}{q_2} - \frac{\alpha_{12}}{q_1}\right)^2\right) - \frac{r_1}{q_1} + \frac{r_2}{q_2} \neq 0, \nonumber
\end{eqnarray}
entonces dicha hipérbola es no degenerada. 

Se puede verificar que la pareja ordenada
\begin{eqnarray}
P &=& \left(\frac{r_2(r_1/K_2 - \alpha_{12})}{r_1r_2/(K_1K_2)-\alpha_{12}\alpha_{21}},\frac{r_1(r_2/K_1 - \alpha_{21})}{r_1r_2/(K_1K_2)-\alpha_{12}\alpha_{21}}\right) \nonumber
\end{eqnarray}
pertenece a la hipérbola en \eqref{4.10}, y que satisface
\begin{eqnarray}
f(P) \quad = \quad 0 \quad = \quad g(P). \nonumber
\end{eqnarray}
Además, las coordenadas de $P$ son ambas positivas si y solo si se cumplen las desigualdades 
\begin{eqnarray}
r_1/K_2 \quad > \quad \alpha_{12}, && r_2/K_1 \quad > \quad \alpha_{21}. \nonumber
\end{eqnarray}
De hecho, las coordenadas de $P$ son las primeras dos coordenadas que tiene el punto de equilibrio frontera $E_3=(x_3,y_3,0)$ que fue estudiado en el Teorema \ref{teo4.3}. Por lo tanto, el punto $P$ se puede escribir como
\begin{eqnarray}
P &=& (x_3,y_3). \nonumber
\end{eqnarray}
También considérese
\begin{eqnarray}
x_H &=& \frac{(a_1K_1 - 1)q_2r_1 + K_1q_1\alpha_{21}}{2a_1q_2r_1} \nonumber \\ & & + \frac{\sqrt{4a_1K_1q_2r_1(q_2r_1 - q_1r_2) + ((a_1K_1 - 1)q_2r_1 + K_1q_1\alpha_{21})^2}}{2a_1q_2r_1},                               \nonumber \\
y_H &=& \frac{(a_2K_2 - 1)q_1r_2 + K_2q_2\alpha_{12}}{2a_2q_1r_2} \nonumber \\ & & + \frac{\sqrt{4a_2K_2q_1r_2(q_1r_2 - q_2r_1) + ((a_2K_2 - 1)q_1r_2 + K_2q_2\alpha_{12})^2}}{2a_2q_1r_2}.                               \nonumber
\end{eqnarray}
Se puede verificar que al menos uno de los dos valores $x_H$ o $y_H$ pertenece a $\mathbb R$. Si $x_H$ está en $\mathbb R$, entonces la pareja ordenada 
\begin{eqnarray}
P_1=(x_H,0) \nonumber
\end{eqnarray}
pertenece a la hipérbola de la ecuación \eqref{4.10}. Análogamente si $y_H$ está en $\mathbb R$, entonces 
\begin{eqnarray}
P_2=(0,y_H) \nonumber
\end{eqnarray}
pertenece a la hipérbola \eqref{4.10}. Además, si $x_H$ y $y_H$ están en $\mathbb R$, entonces $P_1$ y $P_2$ están en ramas distintas de la hipérbola \eqref{4.10} (se puede consultar el Apéndice \ref{apendiceb}).

\begin{obs}
\label{obs4.14}
Si $r_1/q_1 \geq r_2/q_2$, entonces $x_H$ está en $\mathbb R$ y $P_1=(x_H,0)$ está en $\mathbb R^2$. Análogamente, si $r_2/q_1 \leq r_1/q_1$, entonces $y_H$ está en $\mathbb R$ y $P_2=(0,y_H)$ está en $\mathbb R^2$.
\end{obs}
Por la Observación \ref{obs4.14}, en adelante se asumirá sin pérdida de generalidad que $r_1/q_1 > r_2/q_2$.

\begin{lema}
\label{lema4.15}
Si $r_2/K_1 > \alpha_{21}$, entonces $0 < x_H < K_1$.
\end{lema}
\begin{proof}
Para demostrar la desigualdad $x_H>0$, considérese el polinomio cuadrático
\begin{eqnarray}
Q(x) &=& \Omega_2x^2 + \Omega_1x + \Omega_0, \nonumber
\end{eqnarray}
con
\begin{eqnarray}
\Omega_2 &=& a_1q_2r_1, \nonumber \\
\Omega_1 &=& (1 - a_1K_1)q_2r_1 - K_1q_1\alpha_{21}, \nonumber \\
\Omega_0 &=& K_1(q_1r_2 - q_2r_1). \nonumber
\end{eqnarray}
Aplíquese la Regla de Descartes\index{Regla de Descartes}, y obsérvese que $Q$ tiene exactamente un cambio de signo, de donde se sigue que $Q$ tiene exactamente una raíz positiva, que es precisamente $x_H$.

Para demostrar la desigualdad $x_H < K_1$ se procede así
\begin{eqnarray}
x_H &=& \frac{(a_1K_1 - 1)q_2r_1 + K_1q_1\alpha_{21}}{2a_1q_2r_1} \nonumber \\ & & + \frac{\sqrt{4a_1K_1q_2r_1(q_2r_1 - q_1r_2) + ((a_1K_1 - 1)q_2r_1 + K_1q_1\alpha_{21})^2}}{2a_1q_2r_1}                              \nonumber \\
&<& \frac{(a_1K_1 - 1)q_2r_1 + q_1r_2}{2a_1q_2r_1} \nonumber \\ & & + \frac{\sqrt{4a_1K_1q_2r_1(q_2r_1 - q_1r_2) + ((a_1K_1 - 1)q_2r_1 + q_1r_2)^2}}{2a_1q_2r_1}                              \nonumber \\
&=& \frac{(a_1K_1 - 1)q_2r_1 + q_1r_2}{2a_1q_2r_1} + \frac{(a_1K_1 + 1)q_2r_1 - q_1r_2}{2a_1q_2r_1} \nonumber \\
&=& K_1. \nonumber
\end{eqnarray}
\end{proof}
\begin{lema} 
\label{lema4.16}
Si $r_2/K_1 > \alpha_{21}$, $r_1/K_2 > \alpha_{12}$, $f(P_1) > h(P_1)$ y $h(P) > 0$, entonces existe $P^*=(x^*,y^*)$ en la hipérbola \eqref{4.10} tal que $f(P^*) = h(P^*) > 0$.
\end{lema}
\begin{proof}
Se puede verificar que los puntos $P$ y $P_1$ están en la misma rama de la hipérbola en \eqref{4.10}. Se define $\varphi \colon [0,1] \longrightarrow [0,\infty)^2$ una función continua e inyectiva que satisface $\varphi (0) = P$, $\varphi(1) = P_1$, y tal que $\varphi(t)$ satisface la ecuación \eqref{4.10} para toda $t \in (0,1)$. Como las funciones $f\circ \varphi$ y $h\circ \varphi$ son continuas, y además $f\circ \varphi(0) = 0 < h\circ \varphi(0)$, $f\circ \varphi(1) > h\circ \varphi(1)$, entonces por el teorema del valor intermedio, existe $t^*\in (0,1)$ tal que
\begin{eqnarray}
f\circ \varphi(t^*) &=& h\circ \varphi(t^*), \nonumber
\end{eqnarray}
de donde se sigue el resultado.
\end{proof}

A partir de los Lemas \ref{lema4.15} y \ref{lema4.16} se pueden establecer condiciones suficientes para la existencia de un punto de equilibrio positivo $E^* = (x^*,y^*,z^*)$.

\begin{teo}[Reyes-García]
\label{teo4.17}
Considérese el sistema \eqref{4.1} con parámetros positivos. Si \newline $r_2/K_1 > \alpha_{21}$, $r_1/K_2 > \alpha_{12}$, y adicionalmente se cumple 
\begin{eqnarray}
c_1p_1(K_1) \quad < \quad \mu \quad < \quad c_1p_1(x_3) + c_2p_2(y_3), \nonumber
\end{eqnarray}
donde
\begin{eqnarray}
p_1(x) &=& \frac{q_1 x}{1 + a_1x}, \nonumber \\
p_2(y) &=& \frac{q_2 y}{1 + a_2y}, \nonumber
\end{eqnarray}
entonces \eqref{4.1} tiene al menos un punto de equilibrio positivo $E^* = (x^*,y^*,z^*)$.
\end{teo}
\begin{proof}
Si $c_1p_1(K_1) < \mu$, entonces $f(x,0) > 0 > h(x,0)$ para toda $x\in[0,K_1)$. Del Lema \ref{lema4.15} se sigue que $f(P_1) > 0 > h(P_1)$.

Además 
\begin{eqnarray}
h(P) &=& m^{-1}\left(c_1p_1(x_3) + c_2p_2(y_3) - \mu\right) > 0, \nonumber
\end{eqnarray}
y por el Lema \ref{lema4.16} se concluye que \eqref{4.1} tiene al menos un punto de equilibrio positivo.
\end{proof}

Nótese que si se considera una selección de parámetros en la cual $\alpha_{12}$, $\alpha_{21} = 0$, entonces se tiene $x_3 = K_1$, $y_3 = K_2$, y así el valor $c_1p_1(K_1)$ es estrictamente menor que $c_1p_1(x_3) + c_2p_2(y_3)$. Por lo tanto, existen valores para $\mu$ en los que se satisface la cadena de desigualdades $c_1p_1(K_1) < \mu < c_1p_1(x_3) + c_2p_2(y_3)$.

El Teorema \ref{teo4.17} provee condiciones suficientes para la existencia de puntos de equilibrio positivos en \eqref{4.1}. En algunos casos se tiene un único de equilibrio positivo, en otros hay dos puntos de equilibrio positivo y en otros hay hasta tres puntos de equilibrio positivos en \eqref{4.1} (véase la Figura \ref{fig4.10}).

En la Sección \ref{seccion4.3} se menciona que una condición necesaria y suficiente para que el modelo de Bazykin \eqref{4.7} tenga puntos de equilibrio positivos es la desigualdad $x^b < K_1$. Lo anterior contrasta (pero no contradice) lo establecido en el Teorema \ref{teo4.17}, ya que la desigualdad  $c_1p_1(K_1) < \mu$ se cumple si y solo si $x^b \notin [0,K_1]$; a partir de lo anterior, se concluye que la presencia de un competidor en el medio es un hecho que favorece la existencia de estados estacionarios positivos, y eso puede dar lugar a una mayor diversidad.

Tomando como base el Teorema \ref{teo4.17}, supóngase que $E^* = (x^*,y^*,z^*)$ es un punto de equilibrio positivo del sistema \eqref{4.1}. A través del cambio de variables 
\begin{eqnarray}
\overline x = x/(4x^*),\quad  
\overline y = y/(4y^*), \quad
\overline z = z/z^*, \quad
\overline r_1 = r_1, \quad 
\overline K_1 = K_1/(4x^*), \quad 
\overline q_1 = q_1 z^*,    \nonumber \\
\overline a_1 = 4a_1 x^*, \quad   
\overline r_2 = r_2, \quad      
\overline K_2 = K_2/(4y^*), \quad
\overline q_2 = q_2 z^*, \quad 
\overline a_2 = 4a_2 y^*, \quad    
\overline c_1 = 4c_1 x^*/z^*, \nonumber \\
\overline c_2 = 4c_2 y^*/z^*, \quad  
\overline \mu = \mu, \quad         
\overline m = m z^*,        \quad
\overline \alpha_{12} = 4 \alpha_{12} y^*, \quad
\overline \alpha_{21} = 4 \alpha_{21} x^*, \nonumber 
\end{eqnarray}
se obtiene el sistema (después de retirar las barras superiores):
\begin{eqnarray}
\label{4.11}
\dot{x} &=& x\left(r_1(1 - x/K_1) - \alpha_{12} y - \frac{q_1 z}{1 + a_1 x}\right),   \nonumber \\
\dot{y} &=& y\left(r_2(1 - y/K_2) - \alpha_{21} x - \frac{q_2 z}{1 + a_2 y}\right),   \\
\dot{z} &=& z\left(\frac{c_1q_1 x}{1 + a_1 x} + \frac{c_2q_2 y}{1 + a_2 y} - \mu - m z \right),                                                          \nonumber 
\end{eqnarray}
que tiene un punto de equilibrio positivo en $E^* = (1/4,1/4,1)$. Además, la matriz jacobiana del sistema \eqref{4.11} en $E^*$ está dada por
\begin{eqnarray}
J(E^*) &=& \left( \begin{array}{ccc}
                A_1                  &  - \alpha_{12}/4      & - B_1  \nonumber \\
                - \alpha_{21}/4      &  A_2                  & - B_2  \nonumber \\
                    C_1B_1           &  C_2B_2               & - m    \nonumber
                \end{array} \right),                                  \nonumber
\end{eqnarray}
con 
\begin{eqnarray}
A_1 = \frac{1}{4}\left(\frac{a_1q_1}{(1 + a_1/4)^2} - \frac{r_1}{K_1} - \alpha_{12}\right), && 
A_2 = \frac{1}{4}\left(\frac{a_2q_2}{(1 + a_2/4)^2} - \frac{r_2}{K_2} - \alpha_{21}\right), 
                                                                                         \nonumber \\
B_1 = \frac{q_1}{4 + a_1}, &&
B_2 = \frac{q_2}{4 + a_2},                                                               \nonumber \\
C_1 = \frac{4c_1}{1 + a_1/4}, &&
C_2 = \frac{4c_2}{1 + a_2/4},                                                            \nonumber
\end{eqnarray}
donde los valores de $q_1$, $q_2$ están dados por
\begin{eqnarray}
\label{4.12}
q_1 &=& \left(r_1\left(1 - \frac{1}{4K_1}\right) - \frac{\alpha_{12}}{4}\right)\left(1 + \frac{a_1}{4}\right), \nonumber \\
q_2 &=& \left(r_2\left(1 - \frac{1}{4K_2}\right) - \frac{\alpha_{21}}{4}\right)\left(1 + \frac{a_2}{4}\right).    
\end{eqnarray}
Para asegurar de que $q_1$ y $q_2$ son ambos positivos, se asume que siempre se cumplen las desigualdades
\begin{eqnarray}
r_1 > \frac{\alpha_{12}K_1}{4K_1 - 1}, & & r_2 > \frac{\alpha_{21}K_2}{4K_2 - 1},     \nonumber \\
K_1 > \frac{1}{4}, & & K_2 > \frac{1}{4}.                                             \nonumber
\end{eqnarray}
El Teorema \ref{teo4.18} da condiciones suficientes para la estabilidad asintótica del sistema \eqref{4.11} en un punto de equilibrio positivo.
\begin{teo}[Reyes-Garcia]
\label{teo4.18}
Supóngase que \eqref{4.11} (con $\alpha_{12}$, $\alpha_{21} \geq 0$, y el resto de los parámetros positivos) tiene un punto de equilibrio positivo en $E^*=(1/4,1/4,1)$. Si los valores 
\begin{eqnarray}
A_1 &=& \frac{1}{4}\left(\frac{a_1q_1}{(1+a_1/4)^2}-\frac{r_1}{K_1}-\alpha_{12}\right), \nonumber \\
A_2 &=& \frac{1}{4}\left(\frac{a_2q_2}{(1+a_2/4)^2}-\frac{r_2}{K_2}-\alpha_{21}\right), \nonumber
\end{eqnarray}
son negativos y se cumplen las desigualdades
\begin{eqnarray}
\label{4.13}
0\leq \alpha_{12} &\leq& \frac{B_1\alpha_{21}}{B_2} \leq -4A_1,           \\
\alpha_{12} &\leq& -\frac{4A_2B_1}{B_2},                        \nonumber
\end{eqnarray}
entonces $E^*$ es un atractor local en \eqref{4.11}.
\end{teo}
\begin{proof}
El polinomio característico de la matriz jacobiana $J(E^*)$ del sistema \eqref{4.11} está dado por
\begin{eqnarray}
P(\lambda) &=& \lambda^3 + \Omega_2\lambda^2 + \Omega_1\lambda + \Omega_0,   \nonumber \\ \nonumber
\end{eqnarray}
donde                                                         
\begin{eqnarray}
\Omega_2 &=& m + \alpha_{12} + \alpha_{21} - A_1 - A_2,                                    \nonumber \\
\Omega_1 &=& B_1^2C_1 + B_2^2C_2 + A_1A_2 - 
A_1\alpha_{21} - A_2\alpha_{12} + m(\alpha_{12} + \alpha_{21} - A_1 -A_2),                 \nonumber \\
\Omega_0 &=& m(A_1A_2 - A_1\alpha_{21} - A_2\alpha_{12}) + (B_1C_1 - B_2C_2)(B_1\alpha_{21} - 
B_2\alpha_{12}) - A_2B_1^2C_1 - A_1B_2^2C_2. \nonumber
\end{eqnarray}
El criterio de Routh-Hurwitz establece que el coeficiente cuadrático $\Omega_2$, el coeficiente independiente $\Omega_0$, así como la expresión $\Omega_2 \Omega_1 - \Omega_0$ deben ser todos positivos para que las raíces de $P$ tengan parte real negativa. A través de un cálculo, se puede comprobar que si las desigualdades \eqref{4.13} se cumplen, entonces las condiciones de Routh-Hurwitz se verifican.
\end{proof}

Nótese que siempre es posible encontrar valores $\alpha_{12}$ y $\alpha_{21}$ que cumplan las desigualdades \eqref{4.13}. En efecto, si los valores $A_i < 0$, $B_i > 0$, $C_i > 0$, entonces la selección de parámetros $\alpha_{12} = \alpha_{21} = 0$ cumple con las desigualdades.

De lo anterior, se concluye que existen puntos de equilibrio positivos en \eqref{4.11} que son asintóticamente estables.

El siguiente resultado dice que el signo de $A_i$ (que aparece en la matriz jacobiana $J(E^*)$ del sistema \eqref{4.11}) depende de la ubicación de la media aritmética de los valores $-1/a_i$ y $K_i$ con respecto al valor de equilibrio $1/4$.
\begin{prop}
Si 
\begin{eqnarray}
\frac{1}{4} &>& \frac{K_i - 1/a_i}{2},  \nonumber
\end{eqnarray}
\end{prop}
entonces $A_i < 0$.
\begin{proof}
Si $1/4 > 2^{-1}(K_i - 1/a_i)$ entonces
\begin{eqnarray}
K_i < \frac{1}{2} + \frac{1}{a_i} &\Rightarrow& \frac{1}{K_i} > \frac{a_i}{1 + a_i/2}     \nonumber \\
&\Rightarrow& \frac{a_ir_i - r_i/K_i - a_ir_i/(2K_i)}{1 + a_i/4} < 0                      \nonumber \\
&\Rightarrow& \frac{a_ir_i(1 - 1/(4K_i)) - r_i/K_i - a_ir_i/(4K_i)}{1 + a_i/4} < 0        \nonumber \\
&\Rightarrow& \frac{a_ir_i(1 - 1/(4K_i))}{1 + a_i/4} - \frac{r_i}{K_i} < 0                \nonumber \\
&\Rightarrow& \frac{a_i(\alpha_{ij}/4 + q_i/(1 + a_i/4))}{1 + a_i/4} - \frac{r_i}{K_i} < 0
\text{, con }q_i\text{ como en }\eqref{4.12}                                                \nonumber \\
&\Rightarrow& \frac{a_i\alpha_{ij}/4}{1 + a_i/4} + \frac{a_iq_i}{(1 + a_i/4)^2} - \frac{r_i}{K_i} < 0                                                                                              \nonumber \\
&\Rightarrow& \frac{a_iq_i}{(1 + a_i/4)^2} - \frac{r_i}{K_i} < 0.                         \nonumber
\end{eqnarray}
\end{proof}
Si las respuestas funcionales de \eqref{4.11} están dadas por 
\begin{eqnarray}
p_1(x) = \frac{q_1x}{1+a_1x}, && p_2(y)=\frac{q_2y}{1+a_2y}, \nonumber 
\end{eqnarray}
entonces los valores de $B_i$ representan, en cada caso, la imagen del valor de equilibrio $1/4$ de las presas\index{presa} bajo las respuestas funcionales $p_i$. Esto es
\begin{eqnarray}
B_i &=& p_i(1/4).     \nonumber
\end{eqnarray}
A continuación se estudian algunas bifurcaciones que presenta el modelo matemático \eqref{4.11}.
\section{Bifurcaciones}
\label{seccion4.5}
Antes de enunciar y demostrar el Teorema \ref{teo4.18}, se mostraron condiciones para  que el sistema \eqref{4.11} tenga un punto de equilibrio positivo en $E^* = (1/4,1/4,1)$. A partir de una sustitución de los valores de las coordenadas de $E^* = (1/4,1/4,1)$ en \eqref{4.11}, se obtiene
\begin{eqnarray}
\left(\frac{1}{4}\right)\left(r_1(1 - 1/(4K_1)) - \alpha_{12}/4 - \frac{q_1}{1 + a_1/4}\right) &=& 0,                                                                                           \nonumber \\
\left(\frac{1}{4}\right)\left(r_2(1 - 1/(4K_2)) - \alpha_{21}/4 - \frac{q_2}{1 + a_2/4}\right) &=& 0,                                                                                           \nonumber \\
\frac{c_1q_1}{4 + a_1} + \frac{c_2q_2}{4 + a_2} - \mu - m &=& 0.                          \nonumber
\end{eqnarray}
Supóngase que algunos parámetros quedan fijos: $K_1 = 1$, $c_2 = 1$, $a_1 = 1$, $a_2 = 1$, $\alpha_{12} = 0$, $\alpha_{21} = 0$. Por lo tanto, las siguientes igualdades se satisfacen
\begin{eqnarray}
\label{4.14}
q_1 &=& \frac{15 r_1}{16},                                                         \nonumber \\
q_2 &=& \frac{5 r_2 (4K_2 - 1)}{16 K_2},                                                     \\
\mu &=& \frac{3c_1r_1}{16} + \frac{r_2(4K_2 - 1)}{16K_2} - m. \nonumber
\end{eqnarray}

Los valores $q_1$, $q_2$, $\mu$ de \eqref{4.14} son positivos, si se cumplen las desigualdades 
\begin{eqnarray}
K_2 &>& 1/4,                                                                           \nonumber \\
0 \quad < \quad m &<& \frac{3c_1r_1}{16} + \frac{r_2(4K_2 - 1)}{16K_2}. \nonumber 
\end{eqnarray}

Se procede a sustituir $q_1$, $q_2$, $\mu$ como en  \eqref{4.14} en el sistema \eqref{4.11}. Entonces se obtiene un nuevo sistema dado por
\begin{eqnarray}
\label{4.15}
\dot x &=& x\left(r_1(1 - x) - \frac{15r_1 z}{16(1 + x)}\right),                                    \nonumber \\
\dot y &=& y\left(r_2(1 - y/K_2) - \frac{5r_2(4K_2 - 1)z}{16K_2(1+y)}\right),                              \\
\dot z &=& z\left(\frac{15c_1r_1 x}{16(1 + x)} + \frac{5r_2(4K_2 - 1)y}{16K_2(1+y)} -
\left(\frac{3c_1r_1}{16} + \frac{r_2(4K_2 - 1)}{16K_2} - m\right) - m z\right),           \nonumber
\end{eqnarray}
que tiene un punto de equilibrio positivo en $E^* = (1/4,1/4,1)$. 

\begin{figure}[h]
    \centering
    \includegraphics[width=432pt]{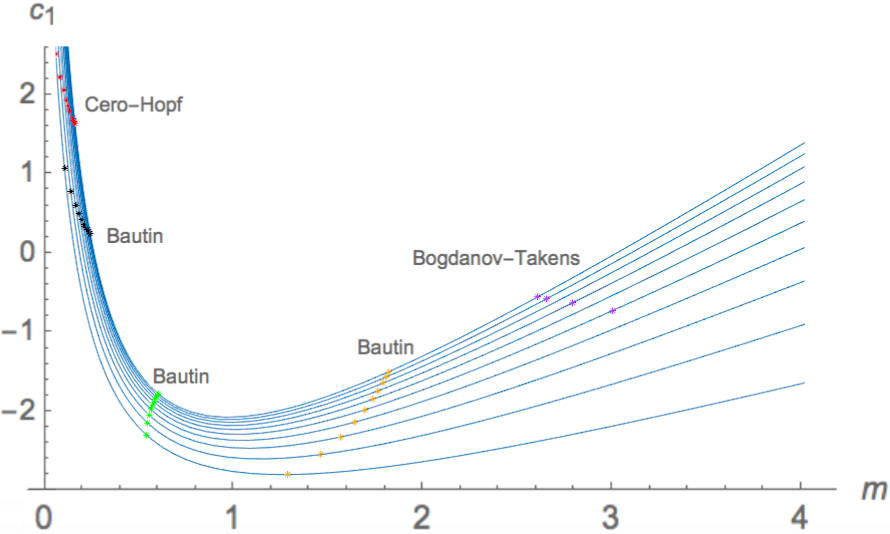}
    \caption{Para $r_1 = 1$, $r_2 = 1.7$ y los niveles $K_2 = 3$, $3.5$, $4$, $4.5$, $5$, $5.5$, $6$, $6.5$, $7$, $7.5$ las curvas formadas por las parejas ordenadas $(m,c_1)$ en las que se presenta un umbral de bifurcación de Hopf (curvas de Hopf) para el sistema \eqref{4.15} se encuentran graficadas de abajo hacia arriba. La variable en el eje horizontal es $m$, y en el eje vertical es $c_1$. Los puntos rojos corresponden a umbrales de una bifurcación de cero-Hopf; los puntos negros, verdes y amarillos son umbrales de bifurcaciones de Bautin; Los puntos morados son bifurcaciones de Bogdanov-Takens. Los umbrales de bifurcación exhibidos en este diagrama, se encuentran por debajo del eje $c_1 = 0$, y por lo tanto están fuera de la región de interés ecológico.}
    \label{fig4.7}
\end{figure}

El Teorema \ref{teo4.20} permite hallar una bifurcación de Hopf en el sistema \eqref{4.15}, ya que provee una configuración de parámetros en la que uno de los valores propios de la matriz jacobiana $J(E^*)$ del sistema es un número real, mientras que los otros dos valores propios son números imaginarios puros.
\begin{teo}[Reyes-García]
\label{teo4.20}
Supóngase que $r_1 > 0$, $r_2 > r_1/2$, $K_2 > 3r_2/(2r_2 - r_1)$, $0 < c_1 < \widetilde c_1$, donde
\begin{eqnarray}
\widetilde c_1 &=& \frac{r_2(2K_2 - 3)\left(4r_1^2K_2^2 + 4r_1r_2K_2(3 - 2K_2) + 5r_2^2(1 - 4K_2)^2\right)}{45r_1^3K_2^3}.
\nonumber
\end{eqnarray}
Si adicionalmente se cumple $m = m_0 < 3c_1r_1/16 + r_2(4K_2 - 1)/(16K_2)$, donde
\begin{eqnarray}
m_0 &=& (80K_2(K_2(r_1 - 2r_2) + 3r_2))^{-1}((r_1^4K_2^4(4 - 45c_1)^2 + 2r_1^2r_2^2K_2^2(45c_1(41 - 88K_2 + 96K_2^2) \nonumber \\ \nonumber \\
& & + 32K_2(1 + 8K_2) - 124) + r_2^4(31 - 8K_2(1 + 8K_2))^2)^{1/2} + r_2^2(-96K_2^2 + 88K_2  - 41) \nonumber \\ \nonumber \\
& & + 8r_1r_2K_2(2K_2 - 3) - r_1^2K_2^2(4 + 45c_1)).                              \nonumber
\end{eqnarray}

entonces el espectro de la matriz jacobiana de \eqref{4.15} evaluada en $E^* = (1/4,1/4,1)$ está formado por un número real $\nu$ y dos números imaginarios conjugados $i \omega$, $-i \omega$.
\end{teo}
Previo a la demostración del Teorema \ref{teo4.20}, se analiza la positividad de la cota superior $\widetilde c_1$ y del umbral de bifurcación $m_0$:

Como $r_1$ es un número positivo, se tiene que
\begin{eqnarray}
K_2 \quad > \quad \frac{3r_2}{2r_2 - r_1} &>& \frac{3}{2}. \nonumber
\end{eqnarray}
Como $K_2 > 3/2$, entonces
\begin{eqnarray}
5r_2^2(1 - 4K_2)^2 > 5r_2^2(3 - 2K_2)^2, &\text{entonces}&
5r_2^2(1 - 4K_2)^2 > r_2^2(3 - 2K_2)^2, \nonumber \\ \nonumber \\
&\text{entonces}& 4r_1^2K_2^2 + 4r_1r_2K_2(3 - 2K_2) + 5r_2^2(1 - 4K_2)^2 \nonumber \\
& & > 4r_1^2K_2^2 + 4r_1r_2K_2(3 - 2K_2) + r_2^2(3 - 2K_2)^2, \nonumber \\ \nonumber \\
&\text{entonces}& 4r_1^2K_2^2 + 4r_1r_2K_2(3 - 2K_2) + 5r_2^2(1 - 4K_2)^2 \nonumber \\
& & > (2r_1K_2 + r_2(3 - 2K_2))^2 \geq 0. \nonumber
\end{eqnarray}
De lo anterior, y del hecho de que las desigualdades $r_2 >0$, $2K_2 - 3 > 0$ siempre se cumplen, se sigue que $\widetilde c_1 > 0$.

Con respecto a $m_0$, se verifica directamente que el argumento de la raíz cuadrada es un número positivo. Por otro lado, el denominador de $m_0$ es positivo, ya que $K_2 > 3r_2/(r_2 - 2r_1) > 0$. 

Obsérvese que $m_0$ es un número positivo: esto sucede porque $m_0$ es la raíz del polinomio de segundo grado
\begin{eqnarray}
a x^2 + b x + c\text{, } \quad \text{donde} \nonumber
\end{eqnarray}
\begin{eqnarray}
a &=& - \frac{r_1}{10} + \frac{(2K_2 - 3)r_2}{10K_2}, \nonumber \\ \nonumber \\
b &=& - \frac{r_1^2}{100} - \frac{9c_1r_1^2}{80} +  \frac{(2 - 3/K_2)r_1r_2}{100} \nonumber \\
& & + \frac{(2K_2 - 3)r_1r_2}{100K_2} - \frac{(1 - 4K_2)^2r_2^2}{80K_2^2} + \frac{(3/K_2 - 2)(2K_2 - 3)r_2^2}{100K_2}, \nonumber \\ \nonumber \\
c &=& - \frac{9c_1r_1^3}{800} + \frac{9c_1(2 - 3/K_2)r_1^2r_2}{800} + \frac{(2K_2 - 3)r_1^2r_2}{1000K_2} \nonumber \\
& & - \frac{9c_1( 2K_2 - 3)r_1^2r_2}{800K_2} - \frac{(2 - 3/K_2)(2K_2 - 3)r_1r_2^2}{1000K_2} + \frac{(2 - 3/K_2)(1 - 4K_2)^2r_2^3}{800K_2^2}. \nonumber \\ \nonumber
\end{eqnarray}

A través de un cálculo con \textit{Mathematica}, se verifica que si $r_1 > 0$, $r_2 > r_1/2$, $K_2 > 3r_2/(2r_2 - r_1)$, $0 < c_1 <\widetilde c_1$, entonces se cumple $a > 0$, $b < 0$, $c > 0$. Por la Regla de Descartes\index{Regla de Descartes} se concluye que $m_0$ es un número positivo.

La tercera ecuación de \eqref{4.14} implica
\begin{eqnarray}
\frac{3c_1r_1}{16} + \frac{r_2(4K_2 - 1)}{16K_2} &>& m \quad > \quad 0. \nonumber
\end{eqnarray}
El umbral de bifurcación $m_0$ debe cumplir también con dicha cadena de desigualdades
\begin{eqnarray}
\frac{3c_1r_1}{16} + \frac{r_2(4K_2 - 1)}{16K_2} &>& m_0 \quad > \quad 0, \nonumber
\end{eqnarray}
Lo anterior ocurre (los cálculos también fueron hechos con \textit{Mathematica}) si y solo si se satisface el conjunto de desigualdades
\begin{eqnarray}
K_2 > \frac{3}{2}\text{, } & & K_2 \neq \frac{3r_2}{2r_2 - r_1}\text{, }     \nonumber \\ \nonumber \\
0 < &c_1& < \frac{r_2(2K_2 - 3)(4r_1^2K_2^2 + 4r_1r_2K_2(3 - 2K_2) + 5r_2^2(1 - 4K_2)^2)}{45r_1^3K_2^3},                                                                  \nonumber
\end{eqnarray}
y esto ocurre siempre que las hipótesis del Teorema \ref{teo4.20} se cumplen.

Por lo tanto, los parámetros $r_1$, $r_2$, $K_2$, $q_1$, $q_2$, $c_1$, $\mu$, $m$, y el umbral de bifurcación $m_0$ son todos números positivos.
\begin{proof}
(del Teorema \ref{teo4.20}) La matriz jacobiana de \eqref{4.15} evaluada en $E^*$ está dada por
\begin{eqnarray}
J(E^*)&=&
\begin{pmatrix}
- r_1/10  & 0                           & - 3r_1/16       \\\\
0         & r_2(2K_2 - 3)/(10K_2) & r_2(1 - 4K_2)/(16K_2) \\\\
3c_1r_1/5 & r_2(4K_2 - 1)/(5K_2)  & -m
\end{pmatrix},
\nonumber
\end{eqnarray}
y el polinomio característico de $J(E^*)$ es
\begin{eqnarray}
P(\lambda) &=& - \lambda^3 + A \lambda^2 + B \lambda + C \text{, } \quad \text{donde}      \nonumber
\end{eqnarray}
\begin{eqnarray}
A &=& \frac{(2 - 3/K_2)r_2 - r_1 - 10m}{10}, \nonumber \\
B &=& \frac{4r_2K_2(2K_2 - 3)(r_1 + 10m) - 5(r_1K_2^2(9c_1r_1 + 8m) + r_2^2(1 - 4K_2)^2)}{400K_2^2}, \nonumber \\
C &=& \frac{r_1r_2(K_2(9c_1r_1 + 8m)(2K_2 -3) - r_2(1 - 4K_2)^2)}{800K_2^2}. \nonumber
\end{eqnarray}
Supóngase que el conjunto de las raíces de $P$ es $\{\nu, i\omega, -i\omega\}$ (con $\nu$, $\omega \in \mathbb R$). Entonces el polinomio característico $P$ se escribe
\begin{eqnarray}
P(\lambda) &=& (\nu - \lambda)(i\omega - \lambda)(- i\omega - \lambda)                     \nonumber \\
           &=& - \lambda^3 + \nu \lambda^2 - \omega^2\lambda + \nu\omega^2.             \nonumber
\end{eqnarray}
Si $B$ es un número negativo, y $A$ tiene el mismo signo que $C$, entonces se escribe
\begin{eqnarray}
A = \nu \text{, } B = -\omega^2 \text{, }C = \nu\omega^2,                               \nonumber
\end{eqnarray}
o de forma equivalente
\begin{eqnarray}
-AB &=& C.                                                                           \nonumber
\end{eqnarray}
Sustituyendo los coeficientes de $P$ en la igualdad anterior, se obtiene
\begin{eqnarray}
-\left(\frac{(2 - 3/K_2)r_2 - r_1 - 10m}{10}\right) \times &&                        \nonumber \\
\left(\frac{4r_2K_2(2K_2 - 3)(r_1 + 10m) - 5(r_1K_2^2(9c_1r_1 + 8m) + r_2^2(1 - 4K_2)^2)}{400K_2^2}\right) &=&                                                                                 \nonumber \\
\frac{r_1r_2(K_2(9c_1r_1 + 8m)(2K_2 -3) - r_2(1 - 4K_2)^2)}{800K_2^2}.              \nonumber
\end{eqnarray}
Se despeja $m$ para obtener
\begin{eqnarray}
m &=& (80K_2(K_2(r_1 - 2r_2) + 3r_2))^{-1}((r_1^4K_2^4(4 - 45c_1)^2 + 2r_1^2r_2^2K_2^2(45c_1(41 - 88K_2 + 96K_2^2) \nonumber \\ \nonumber \\
& & + 32K_2(1 + 8K_2) - 124) + r_2^4(31 - 8K_2(1 + 8K_2))^2)^{1/2} + r_2^2(-96K_2^2 + 88K_2  - 41) \nonumber \\ \nonumber \\
& & + 8r_1r_2K_2(2K_2 - 3) - r_1^2K_2^2(4 + 45c_1)).                              \nonumber
\end{eqnarray}
El lado derecho de la igualdad anterior coincide con $m_0$. Por lo tanto, es claro que si se cumple
\begin{eqnarray}
m &=& m_0, \nonumber
\end{eqnarray}
entonces el espectro de la matriz jacobiana de \eqref{4.15} es de la forma $\{\nu,i\omega,-i\omega\}$. 
\end{proof}

\begin{figure}[h]
\label{fig408}
   \centering
        \subfigure[$m=0.2$]{\includegraphics[width=195pt]{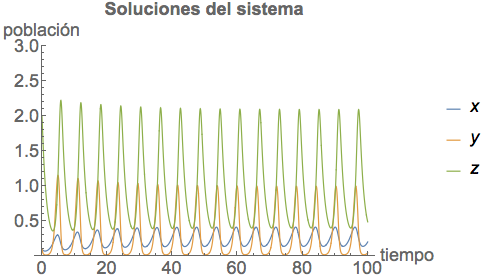}}
        \subfigure[$m=0.25$]{\includegraphics[width=195pt]{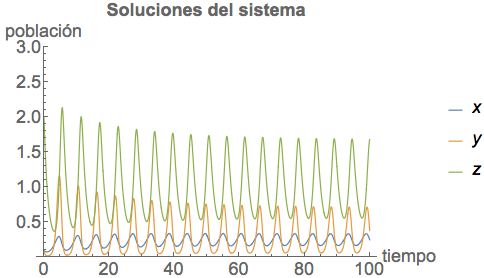}}
        \subfigure[$m=0.3$]{\includegraphics[width=195pt]{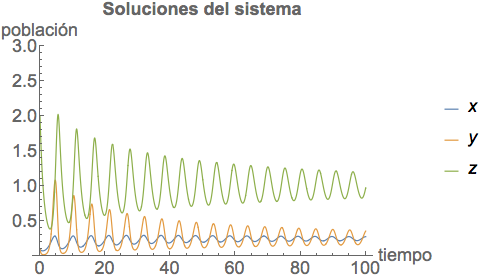}}
        \subfigure[$m=0.35$]{\includegraphics[width=195pt]{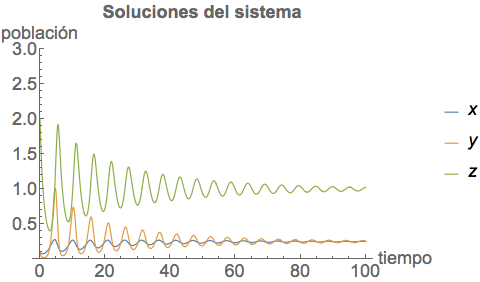}}
        \caption{Soluciones del sistema \eqref{4.15} con condiciones iniciales $x_0=0.1$, $y_0=0.1$, $z_0=2$, y parámetros $r_1=1$, $r_2=3$, $K_2=3$, $c_1=0.5$. En los diagramas (a) y (b), el valor de $m$ es menor que $m_0 = 0.284917$, mientras que en los diagramas (c) y (d) el valor de $m$ es mayor que $m_0$.}
        \label{fig4.8}
\end{figure}

\begin{figure}
\label{fig409}
    \centering
        \subfigure[$m=0.2$]{\includegraphics[width=130pt]{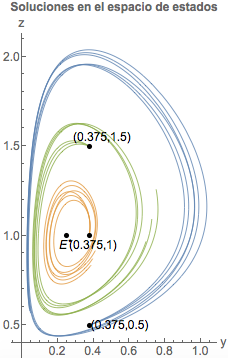}}
        \subfigure[$m=0.25$]{\includegraphics[width=130pt]{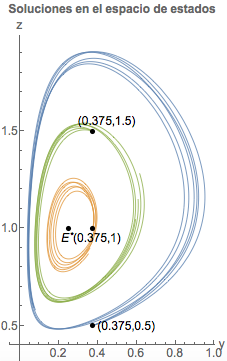}}
        \subfigure[$m=0.3$]{\includegraphics[width=130pt]{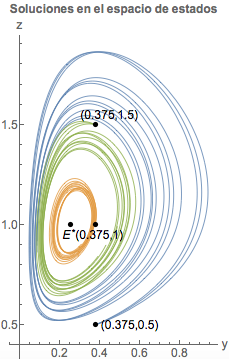}}
        \subfigure[$m=0.35$]{\includegraphics[width=130pt]{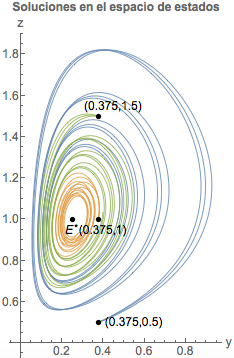}}    
        \caption{Se grafican las coordenadas $y$, $z$ de algunas soluciones particulares de \eqref{4.15}. Los parámetros están dados por $r_1 = 1$, $r_2 = 3$, $K_2 = 3$, $c_1 = 0.5$ y las condiciones $y_0$, $z_0$ se muestran en cada caso, y el valor de $x_0$ es tomado de la lista $0.125$, $0.25$, $0.375$. En los diagramas (a) y (b), el valor de $m$ es menor que $m_0 = 0.284917$, mientras que en los diagramas (c) y (d) el valor de $m$ es mayor que $m_0$.}
    \label{fig4.9}
\end{figure}

A través de simulaciones numéricas, se muestra como la estabilidad asintótica del punto de equilibrio cambia al cruzar el umbral de bifurcación $m_0$. Esto se puede verificar en las Figuras \ref{fig4.8} y \ref{fig4.9}.

Además, en el Apéndice \ref{apendicec} se calcula el primer coeficiente de Lyapunov, y se analiza la condición de transversalidad de la bifurcación de Hopf cuyo umbral es exhibido en el Teorema \ref{teo4.20}.

A continuación se da a conocer una pareja de parámetros $m$ y $c_1$ para los cuales la matriz jacobiana de \eqref{4.15} evaluada en el punto de equilibrio positivo tiene espectro conformado por dos valores propios imaginarios puros conjugados, y un valor propio igual a cero.

\begin{teo}[Reyes-García]
\label{teo4.21}
Supóngase que $r_1>0$, $r_2>0$, $K_2 > 3r_2/(2r_2-r_1)$. Adicionalmente supóngase que 
\begin{eqnarray}
c_1 &=& \frac{1}{45}\left(4 + \frac{r_2(8K_2(8K_2+1)-31)}{r_1K_2(2K_2-3)}\right), \nonumber \\
m &=& \frac{r_2}{5} - \frac{r_1}{10} - \frac{3 r_2}{10 K_2}. \nonumber
\end{eqnarray}
Entonces el espectro de la matriz jacobiana de \eqref{4.15} evaluada en $E^*=(1/4,1/4,1)$ está formado por dos números imaginarios conjugados y un valor propio que es exactamente cero. 
\end{teo}
\begin{proof}
Si los parámetros cumplen con las hipótesis, entonces se puede verificar que los valores $A$, $B$ y $C$ tomados como en la demostración del Teorema \ref{teo4.20} cumplen las condiciones
\begin{eqnarray}
A &=& 0, \nonumber \\
B &>& 0, \nonumber \\
C &=& 0. \nonumber
\end{eqnarray}
Por lo tanto el espectro de la matriz jacobiana de \eqref{4.15} está formado por los valores propios $0$, $i\sqrt B$, $-i\sqrt B$.
\end{proof}

\begin{figure}
    \centering
    \includegraphics[width=241pt]{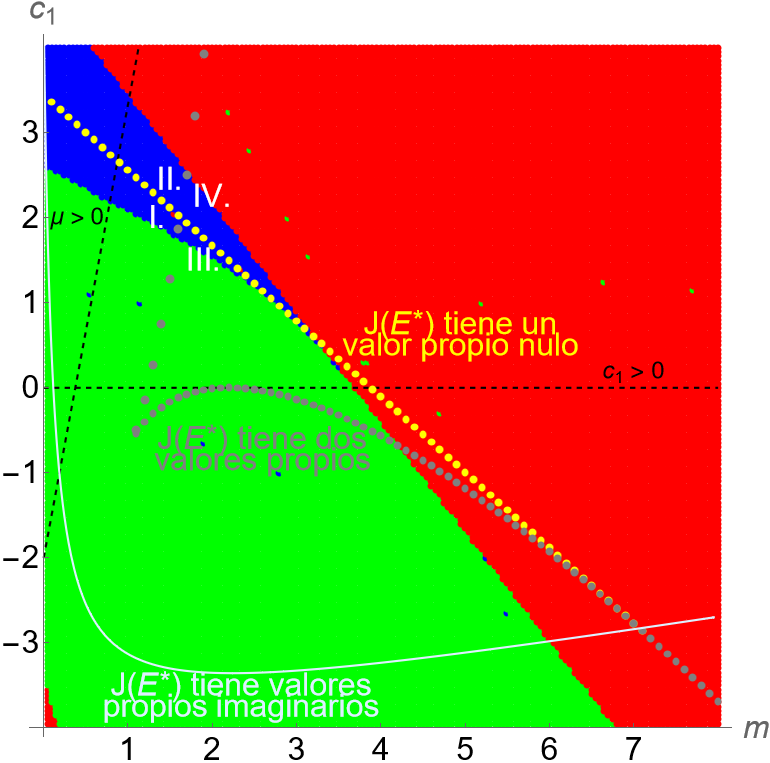}
    \includegraphics[width=241pt]{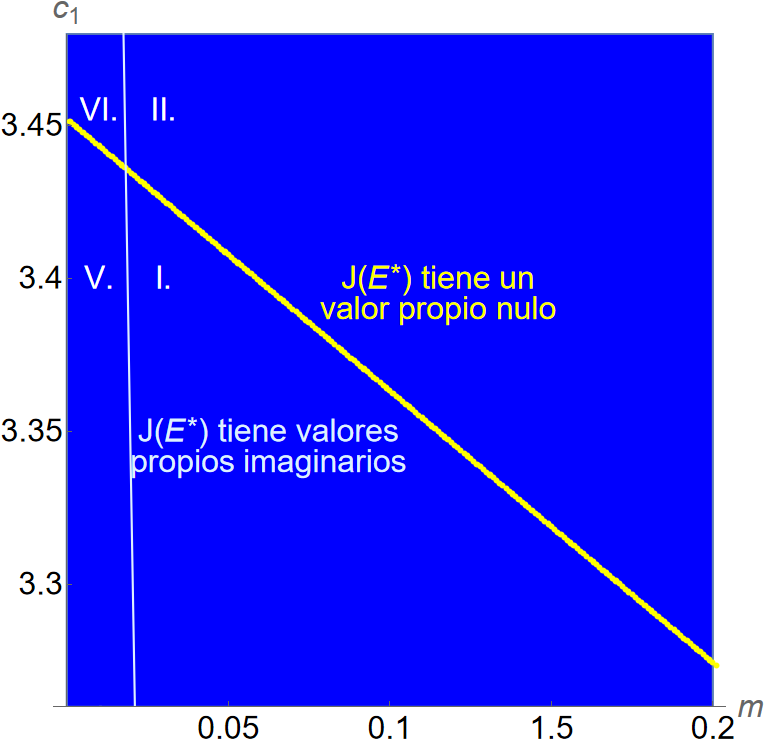}
    \caption{Considérese el sistema \eqref{4.15} con los parámetros $r_1=1$, $r_2=1.7$, $K_2=2.3$, y con $m$, $c_1$ libres (el segundo diagrama es un acercamiento a la intersección de las curvas blanca y amarilla, que es donde está el umbral de la bifurcación de cero-Hopf). En la región marcada con I., $E^*$ es un foco atractor, mientras que en la región III. $E^*$ es nodo atractor. En las regiones II. y IV., $E^*$ es un punto silla, aunque las trayectorias descritas son diferentes. En la región V. $E^*$ es un punto silla y en la región VI. $E^*$ es un foco repulsor. En la región verde, \eqref{4.15} tiene solamente como punto de equilibrio positivo a $E^* = (1/4,1/4,1)$; en la región azul aparecen dos puntos de equilibrio positivo; y en la región roja hay tres puntos de equilibrio positivo. Las líneas discontinuas en negro delimitan las regiones en las que $\mu$ es positivo y $c_1$ es positivo respectivamente, y las líneas discontinuas en color son separatrices.}
    \label{fig4.10}
\end{figure}

\begin{figure}
    \centering
    \includegraphics[width=100pt]{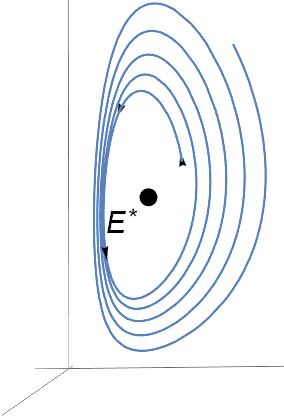}
    \includegraphics[width=110pt]{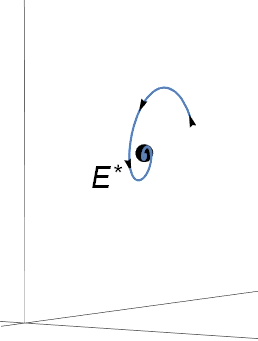}
    \includegraphics[width=75pt]{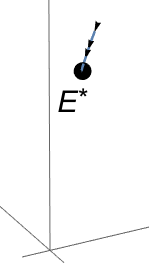}
    \caption{Soluciones del sistema \eqref{4.15} con $r_1=1$, $r_2=1.7$, $K_2=2.3$ con un único punto de equilibrio positivo $E^*$. Puede ocurrir que $E^*$ sea punto silla, foco atractor o nodo atractor.}
    \label{fig4.11}
\end{figure}

\begin{figure}
    \centering
    \includegraphics[width=120pt]{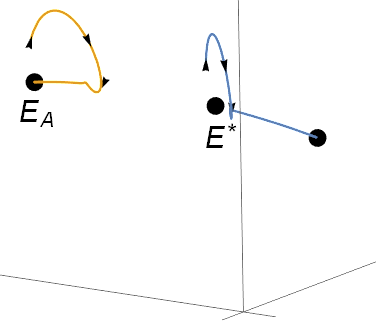}
    \includegraphics[width=120pt]{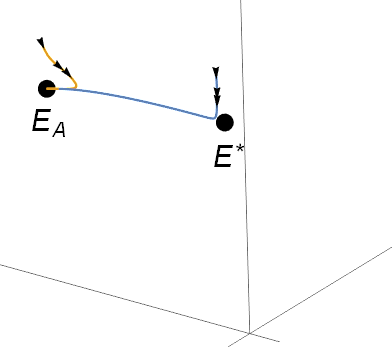}
    \includegraphics[width=50pt]{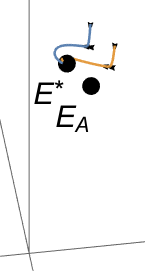}
    \includegraphics[width=70pt]{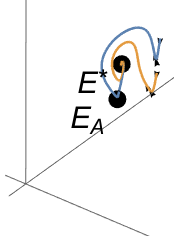}
    \includegraphics[width=70pt]{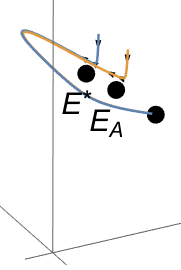}
    \caption{Soluciones del sistema \eqref{4.15} con $r_1=1$, $r_2=1.7$, $K_2=2.3$ con dos puntos de equilibrio positivos $E^*$ y $E_A$. Se presentan casos en que las soluciones convergen a un punto en la frontera. Al igual que en la Figura \ref{fig4.11}, $E^*$ puede ser punto silla, foco atractor o nodo atractor.}
    \label{fig4.12}
\end{figure}

\begin{figure}
    \centering
    \includegraphics[width=80pt]{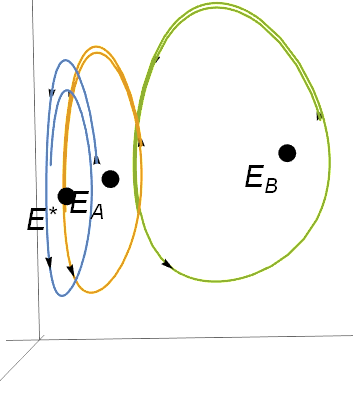}
    \includegraphics[width=110pt]{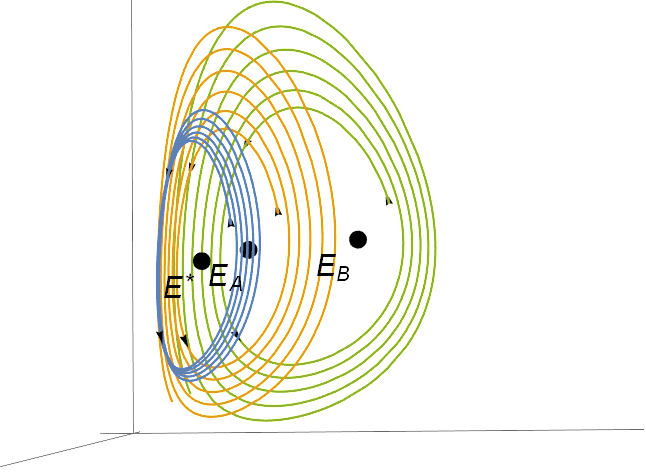}
    \includegraphics[width=70pt]{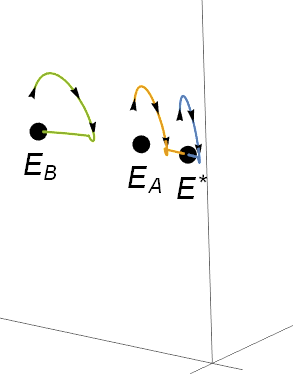}
    \includegraphics[width=90pt]{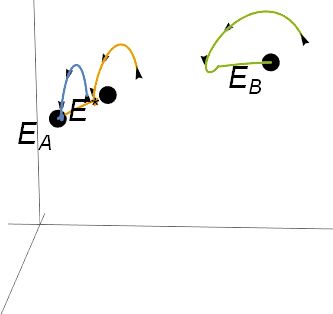}
    \includegraphics[width=80pt]{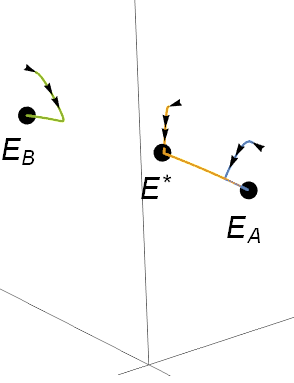}
    \caption{Soluciones del sistema \eqref{4.15} con $r_1=1$, $r_2=1.7$, $K_2=2.3$ con tres puntos de equilibrio positivos $E^*$, $E_A$ y $E_B$. En los primeros dos casos, se puede apreciar la existencia de ciclos límite (parámetros correspondientes a las regiones V. y VII. en la Figura \ref{fig4.10}). En los tres casos restantes, el retrato fase presenta biestabilidad (dos atractores y un punto silla).}
    \label{fig4.13}
\end{figure}

El sistema \eqref{4.15} exhibe una dinámica diversa. El punto de equilibrio $E^* = (1/4,1/4,1)$ puede ser un foco, puede ser un nodo, o también puede ser un punto silla. También es posible identificar las regiones en el espacio de parámetros en los que $E^*$ es asintóticamente estable y en las que no lo es, así como determinar el número de puntos de equilibrio positivos que existen.

Los parámetros del sistema \eqref{4.11} son tomados dentro del conjunto de números no negativos. Sin embargo, es válido preguntarse qué resultados se podrían obtener si se considera $\mu < 0$ ($\mu$ representa la mortalidad de la especie depredadora\index{depredador}), o si se considera $c_1 < 0$ ($c_1$ es la eficiencia en la conversión de biomasa capturada por parte de la especie depredadora\index{depredador}).

En la Figura \ref{fig4.10} se muestran las rectas separatrices que dan lugar a las bifurcaciones presenta el sistema \eqref{4.15}, y también se muestra la región en la que $\mu$ y $c_1$ son mayores que cero (en el sistema \eqref{4.15}, el valor de $\mu$ está dado como en la tercera igualdad de \eqref{4.14}). En las Figuras \ref{fig4.11}, \ref{fig4.12} y \ref{fig4.13} se muestran bosquejos de algunas soluciones de \eqref{4.14} con condiciones iniciales en el octante positivo $(0,\infty)^3$. En las gráficas de la Figura \ref{fig4.11}, los parámetros seleccionados dan lugar a un único punto de equilibrio positivo; en la Figura \ref{fig4.12} se muestran casos en los que hay dos puntos de equilibrio positivo; finalmente, en la Figura \ref{fig4.13} aparecen los casos con tres puntos de equilibrio positivo.

\section{Existencia de órbitas periódicas}
\label{seccion4.6}
En la sección anterior se muestra la existencia de un umbral de bifurcación de cero-Hopf. A continuación se aplica el método del promedio para hallar ciclos límite para el sistema \eqref{4.15} en dicho umbral. También se trabaja con otro sistema que es una variante de \eqref{4.15}. Para conocer más acerca de la teoría del promedio, se pueden consultar \cite{Verhulst1996NonlinearSystems} o el Apéndice \ref{apendiced}. 

En la demostración del Teorema \ref{teo4.20}, aparece la matriz jacobiana de \eqref{4.15} evaluada en el punto de equilibrio $E^*=(1/4,1/4,1)$
\begin{eqnarray}
J(E^*)&=&
\begin{pmatrix}
- r_1/10  & 0                           & - 3r_1/16       \\\\
0         & r_2(2K_2 - 3)/(10K_2) & r_2(1 - 4K_2)/(16K_2) \\\\
3c_1r_1/5 & r_2(4K_2 - 1)/(5K_2)  & -m
\end{pmatrix}.
\nonumber
\end{eqnarray}
Considérese $m$ y $c_1$ como en el Teorema \ref{teo4.21}. La matriz $J(E^*)$ es usada para reescribir \eqref{4.15} con el punto de equilibrio trasladado a $(0,0,0)$ y en forma polinomial
\begin{eqnarray}
\begin{pmatrix}
\dot x \\ \dot y \\ \dot z
\end{pmatrix}
&=& J(E^*)\begin{pmatrix}
x \\ y \\ z \nonumber
\end{pmatrix}
+ O_2(x,y,z),
\end{eqnarray}
donde $O_2$ representa los términos de orden mayor o igual a $2$ en \eqref{4.15}. El cambio a coordenadas cilíndricas $(\rho,\theta,z)$ se realiza a través de las sustituciones 
\begin{eqnarray}
x &=& \rho\cos \theta, \nonumber \\
y &=& \rho\sen \theta. \nonumber
\end{eqnarray}
Se realiza el cálculo de la derivada
\begin{eqnarray}
\dot x &=& -\dot \theta \rho \sen \theta + \dot \rho \cos \theta, \nonumber \\
\dot y &=& \dot \theta \rho \cos \theta + \do \rho \sen \theta. \nonumber
\end{eqnarray}
A través de algunas manipulaciones algebraicas se obtienen las derivadas de $\rho$ y de $\theta$:
\begin{eqnarray}
\dot \rho &=& \dot x\cos\theta + \dot y\sen\theta, \nonumber \\ \nonumber \\
\dot \theta &=& \frac{\dot y\cos\theta - \dot x\sen\theta}{\rho}. \nonumber
\end{eqnarray}
Así, la parte lineal del sistema con el cambio de coordenadas está dada por
\begin{eqnarray}
\dot \rho &=& -\frac{r_1}{10}\rho\cos^2\theta + \left(\frac{r_2(2K_2-3)}{10K_2}\right)\rho\sen^2\theta - \frac{3r_1}{16}z\cos\theta + \frac{r_2(1 - 4K_2)}{16K_2}z\sen\theta, \nonumber \\
\dot \theta &=& \left(\frac{r_2(2K_2 - 3)}{10K_2} + \frac{r_1}{10}\right)\rho \cos\theta \sen\theta + \frac{1}{\rho}\left(\frac{r_2(1 - 4K_2)}{16K_2} - \frac{3r_1}{16}\right)z, \nonumber \\
\dot z &=& \frac{3c_1r_1}{5}\rho\cos\theta + \frac{r_2(4K_2 - 1)}{5K_2}\rho\sen\theta - mz. \nonumber
\end{eqnarray}
Después de hacer el reescalamiento $(\rho,z)\longrightarrow (\varepsilon\rho,\varepsilon z)$, se obtiene el siguiente sistema diferencial
\begin{eqnarray}
\dot \rho &=& \varepsilon\left(-\frac{r_1}{10}\rho\cos^2\theta + \left(\frac{r_2(2K_2-3)}{10K_2}\right)\rho\sen^2\theta - \frac{3r_1}{16}z\cos\theta + \frac{r_2(1 - 4K_2)}{16K_2}z\sen\theta\right), \nonumber \\
\dot z &=& \varepsilon\left(\frac{3c_1r_1}{5}\rho\cos\theta + \frac{r_2(4K_2 - 1)}{5K_2}\rho\sen\theta - mz\right). \nonumber
\end{eqnarray}
Considérese a $\theta$ como una nueva variable independiente. Aquí se usan las ideas y la notación de \cite{Llibre2015Zero-HopfSystem}. Las soluciones del sistema que satisfacen la desigualdad $\dot \theta > 0$ son soluciones para el sistema
\begin{eqnarray}
\frac{d\rho}{d\theta} &=& \varepsilon\left(-\left(\frac{r_1}{10}\right)\rho\cos^2\theta + \left(\frac{r_2(2K_2-3)}{10K_2}\right)\rho\sen^2\theta - \frac{3r_1}{16}z\cos\theta + \frac{r_2(1 - 4K_2)}{16K_2}z\sen\theta\right), \nonumber \\
&: =& \varepsilon F_1(\rho,z), \nonumber \\
\frac{dz}{d\theta} &=& \varepsilon\left(\frac{3c_1r_1}{5}\rho\cos\theta + \frac{r_2(4K_2 - 1)}{5K_2}\rho\sen\theta - mz\right) \nonumber \\
&: =& \varepsilon F_2(\rho,z). \nonumber
\end{eqnarray}
Ahora se calcula el promedio de las funciones $F_1$ y $F_2$ en el costado derecho de las ecuaciones previas, con $\theta$ como el integrador. Esto da lugar a lo siguiente
\begin{eqnarray}
f_1(\rho,z) &:=& \frac{1}{2\pi}\int_0^{2\pi} F_1(\rho,z)d\theta \quad = \quad \frac{\rho}{2}\left(-\left(\frac{r_1}{10}\right) + \left(\frac{r_2(2K_2-3)}{10K_2}\right)\right), \nonumber \\
f_2(\rho,z) &:=& \frac{1}{2\pi}\int_0^{2\pi} F_2(\rho,z)d\theta \quad = \quad
- mz. \nonumber
\end{eqnarray}
Se verifica de forma inmediata que el sistema de ecuaciones
\begin{eqnarray}
f_1(\rho,z) &=& 0, \nonumber \\
f_2(\rho,z) &=& 0, \nonumber
\end{eqnarray}
tiene única solución en $(\rho,z)=(0,0)$. Esto muestra que el retrato fase de \eqref{4.15} exhibe órbitas periódicas en la cercanía de $E^* = (1/4,1/4,1)$. La matriz jacobiana para $f_1$ y $f_2$ con respecto a $\rho$ y a $z$ está dada por
\begin{eqnarray}
\left(
\begin{array}{cc}
\left(-r_1 + r_2(2 - 3/K_2)\right)/20  & 0 \\\\
0  & -m  \end{array}\right). \nonumber
\end{eqnarray}
De lo anterior y de la aplicación de los Teoremas del Apéndice \ref{apendiced}, se concluye lo siguiente.
\begin{teo}[Reyes-García]
\label{teo4.22}
Considérese el sistema \eqref{4.15}, y supóngase que los parámetros $m$, $c_1$ se hallan en las inmediaciones del umbral de bifurcación de cero-Hopf, es decir
\begin{eqnarray}
c_1 &\approx& \frac{1}{45}\left(4 + \frac{r_2(8K_2(8K_2+1)-31)}{r_1K_2(2K_2-3)}\right), \nonumber \\
m &\approx& \frac{r_2}{5} - \frac{r_1}{10} - \frac{3 r_2}{10 K_2}. \nonumber
\end{eqnarray}
Entonces \eqref{4.15} tiene una solución periódica alrededor del punto de equilibrio $E^* = (1/4,1/4,1)$. Más aún, si se cumple la desigualdad
\begin{eqnarray}
r_1 &>& r_2(2 - 3/K_2), \nonumber
\end{eqnarray}
entonces la solución periódica es localmente atractora, y en caso contrario dicha órbita es inestable.
\end{teo}

A continuación se presenta otro modelo que es similar al que se acaba de analizar, y nuevamente se aplica la teoría del promedio para conocer las propiedades de las órbitas periódicas que aparecen en el retrato fase. 

Se realiza una modificación a la respuesta funcional de \eqref{4.15}, y adicionalmente se manipula la tasa intrínseca de crecimiento de las presas\index{presa}. A partir de eso, se obtiene el sistema
\begin{eqnarray}
\label{4.16}
\dot x &=& x\left(r_1(2 - x) - \frac{3 r_1 z}{1+ x + y}\right), \nonumber \\
\dot y &=& y\left(r_2(2 - y) - \frac{3 r_2 z}{1+ x + y}\right),  \\
\dot z &=& z\left(\frac{3 c_1(r_1 x + r_2 y)}{1 + x + y} - (c_1(r_1 + r_2) - m) - m z\right), \nonumber
\end{eqnarray}
Un cambio en la escala del tiempo (véase por ejemplo \cite[Subsección 1.4.1]{Chicone1999OrdinaryApplications}) permite reescribir a \eqref{4.16} como
\begin{eqnarray}
\dot x &=& x \left(r_1 (2-x) (1 + x + y) - 3 r_1 z\right), \nonumber \\
\dot y &=& y \left(r_2 (2-y) (1 + x + y) - 3 r_2 z\right), \nonumber \\
\dot z &=& z \left(3 c_1 \left(r_1 x + r_2 y\right) - (1 + x + y) (c_1(r_1 + r_2) - m + m z)\right), \nonumber
\end{eqnarray}
y que es equivalente a
\begin{eqnarray}
\label{4.17} 
\dot x &=& 2 r_1 x+r_1 x^2+2 r_1 x y-3 r_1 x z-r_1 x^3-r_1 x^2 y, \nonumber \\
\dot y &=& 2 r_2 y+2 r_2 x y+r_2 y^2-3 r_2 y z-r_2 x y^2-r_2 y^3, \\
\dot z &=& \left(m - c_1 \left(r_1+r_2\right)\right)z + \left(m + c_1 \left(2r_1 - r_2\right)\right)xz \nonumber \\
& & + \left(m + c_1 \left(2r_2 - r_1\right)\right)yz - m z^2 - m x z^2 - m y z^2. \nonumber
\end{eqnarray}
Un cálculo muestra que \eqref{4.17} tiene un punto de equilibrio en $(1,1,1)$, y que la matriz jacobiana de \eqref{4.17} evaluada en $(1,1,1)$ está dada por
\begin{eqnarray}
J(1,1,1) &=& \begin{pmatrix}
-2r_1            &  r_1            & -3r_1 \\
r_2              & -2r_2           & -3r_2 \\
c_1(2r_1 - r_2)  & c_1(2r_2 - r_1) & -3m
\end{pmatrix}. \nonumber
\end{eqnarray}
A través de un procedimiento similar al que se realizó en la demostración del Teorema \ref{teo4.21} se verifica que \eqref{4.17} tiene un umbral de bifurcación de cero-Hopf en 
\begin{eqnarray}
(c_1,m) = \left(\frac{2}{3},-\frac{2}{3}(r_1+r_2)\right). \nonumber
\end{eqnarray}
El punto de equilibrio es llevado al origen de coordenadas $(0,0,0)$, y así el sistema se reescribe
\begin{eqnarray}
\begin{pmatrix}
\dot x \\ \dot y \\ \dot z
\end{pmatrix}
&=& J(1,1,1)\begin{pmatrix}
x \\ y \\ z \nonumber
\end{pmatrix}
+ O_2(x,y,z), 
\end{eqnarray}
donde $O_2$ representa los términos de orden mayor o igual a $2$ en \eqref{4.17}. En coordenadas cilíndricas, la parte lineal del sistema queda dada por
\begin{eqnarray}
\dot \rho &=& \rho(-2r_1\cos^2 \theta + (r_1 + r_2)\cos\theta\sen\theta - 2r_2\sen^2\theta) - 3z(r_1\cos\theta + r_2\sen\theta), \nonumber \\
\dot \theta &=& r_2\cos^2\theta - 2(r_1 + r_2)\sen\theta\cos\theta - r_1\sen^2\theta + \frac{3z}{\rho}(-r_2\cos\theta + r_1\sen\theta), \nonumber \\
\dot z &=& c_1\left((2r_1 - r_2)\rho\cos\theta + (2r_2 - r_1)\rho\sen\theta\right) - 3mz. \nonumber
\end{eqnarray}
El rescalamiento $(\rho,z)\longrightarrow (\varepsilon\rho,\varepsilon z)$ da lugar la siguiente pareja de ecuaciones:
\begin{eqnarray}
\dot \rho &=& \varepsilon\left(\rho(-2r_1\cos^2 \theta + (r_1 + r_2)\cos\theta\sen\theta - 2r_2\sen^2\theta) - 3z(r_1\cos\theta + r_2\sen\theta)\right), \nonumber \\
\dot z &=& \varepsilon\left(c_1\left((2r_1 - r_2)\rho\cos\theta + (2r_2 - r_1)\rho\sen\theta\right) - 3mz\right). \nonumber
\end{eqnarray}
Se procede de forma similar a como se hizo anteriormente, usando $\theta$ como nueva variable y considerando las soluciones en el conjunto $\dot \theta > 0$.
\begin{eqnarray}
\frac{\partial \rho}{\partial \theta} &=& \varepsilon\left(\rho(-2r_1\cos^2 \theta + (r_1 + r_2)\cos\theta\sen\theta - 2r_2\sen^2\theta) - 3z(r_1\cos\theta + r_2\sen\theta)\right) \nonumber \\
&:=& \varepsilon F_1(\rho,z), \nonumber \\
\frac{\partial z}{\partial \theta} &=& \varepsilon\left(c_1\left((2r_1 - r_2)\rho\cos\theta + (2r_2 - r_1)\rho\sen\theta\right) - 3mz\right) \nonumber \\
&:=& \varepsilon F_2(\rho,z). \nonumber
\end{eqnarray}
El promedio es calculado
\begin{eqnarray}
f_1(\rho,z) &=& \frac{1}{2\pi}\int_0^{2\pi} F_1(\rho,z)d\theta \quad = \quad -(r_1 + r_2)\rho, \nonumber \\
f_2(\rho, \theta) &=& \frac{1}{2\pi}\int_0^{2\pi} F_2(\rho,z)d\theta \quad = \quad -3mz. \nonumber
\end{eqnarray}
Es claro que el sistema $f_1(\rho,z) = f_2(\rho,z) = 0$ tiene única solución en $(\rho,z) = (0,0)$. Como $m$ tiene signo opuesto a la suma $r_1 + r_2$, entonces se sigue el siguiente resultado.
\begin{teo}[Reyes-García]
\label{teo4.23}
Considérese el sistema \eqref{4.16}, y supóngase que $\displaystyle c_1 \approx \frac{2}{3}$, $\displaystyle m \approx -\frac{2}{3}(r_1 + r_2)$. Entonces \eqref{4.16} tiene una órbita periódica alrededor del punto de equilibrio $(1,1,1)$. Dicha órbita siempre es inestable.
\end{teo}

Existen diversas técnicas que se han implementado para el análisis de la existencia y de las características de las órbitas periódicas en modelos multiespecíficos que presentan bifurcaciones de cero-Hopf. Un ejemplo de lo anterior puede ser hallado en \cite{Al-khedhairi2020Zero-HopfApproach}, donde se realiza un acercamiento a través un reescalamiento usando múltiples variables temporales.
\section{Conclusiones}
Para estudiar la invasión de los ecosistemas\index{ecosistema} a través de la modelación matemática, se propuso el modelo \eqref{4.1} que corresponde a las poblaciones de dos presas\index{presa} y de un depredador\index{depredador} y que es el principal objeto de estudio de este capítulo. Se muestra que el modelo es disipativo, y se muestran sus seis diferentes tipos de puntos de equilibrio: cinco de ellos en la frontera y uno que es positivo.

El análisis que se hace en la Sección \ref{seccion4.2} está basado en el trabajo de Abrams \cite{Abrams1999}, y permite obtener resultado que coincide con el que aparece en el tercer punto del Teorema \ref{teo4.3}.

Los experimentos hechos en la Sección \ref{seccion4.3} muestran que al manipular los parámetros correspondientes a la capacidad de carga y al tiempo de manejo de la especie invasora, es posible alcanzar la coexistencia de las especies. Aunque también hay casos en los que alguna(s) de las especies perece(n).

En comparación con los resultados análogos del Capítulo \ref{capitulo3}, resulta más complicado establecer condiciones para la existencia y la estabilidad de puntos de equilibrio positivos en el sistema \eqref{4.1}. Para analizar la existencia se realiza una construcción geométrica (ver el Apéndice \ref{apendiceb}), y para estudiar la estabilidad se propone una traslación del punto en cuestión a $(1/4,1/4,1)$ para posteriormente usar el teorema de Hartman-Grobman y el criterio de Routh-Hurwitz.

También se mostró que el modelo cuenta con bifurcaciones de tipo Hopf, cero-Hopf, Bautin y Bogdanov-Takens. Los umbrales de Hopf y de cero-Hopf aparecen en los Teoremas \ref{teo4.20} y \ref{teo4.21} respectivamente, y los umbrales de Bautin y de Bogdanov-Takens son mostrados de forma numérica en el diagrama de la Figura \ref{fig4.7}.

Finalmente, usando las ideas expuestas en \cite{Llibre2015Zero-HopfSystem} fue posible encontrar condiciones para la originación de ciclos límites en el modelo.
\chapter{Dinámica Caótica}
\label{capitulo5}
\begin{flushright}
\emph{``Creo que la realidad matemática\\está fuera de nosotros, que nuestra\\función es descubrirla y observarla,\\y que los teoremas que demostramos\\y que describimos grandilocuentemente \\como nuestras creaciones son simples\\notas de nuestras observaciones.''\\$ $\\Godfrey Harold Hardy \cite[p.122-124]{Hardy1967AApology}}
\end{flushright}

Uno de los pioneros en la teoría del caos fue Henri Poincaré\index{Poincaré}, quien demuestra que las ecuaciones de tres cuerpos de Newton no son integrables en 1889, y que enuncia el famoso teorema de la recurrencia en 1890 (que sería demostrado por Constantin Caratheodory\index{Caratheodory} en 1919). El trabajo de Poincaré\index{Poincaré}, fue la base de los resultados del matemático estadounidense George David Birkhoff. Ya en la década de los sesenta del siglo XX, hubo grandes avances en este campo que estuvieron a cargo de autores como Lorenz, Hénon, Smale, Ruelle y Takens (ver \cite[Sección 14.1]{Verhulst1996NonlinearSystems}).

Los sistemas dinámicos no lineales que están conformados por ecuaciones diferenciales ordinarias o ecuaciones diferenciales parciales, han sido y son una parte importante en el desarrollo de múltiples ramas de la física contemporánea (como por ejemplo mecánica celeste, meteorología, entre otras), así como de otras disciplinas como la ecología de poblaciones y la economía (ver \cite[Sección 14.1]{Verhulst1996NonlinearSystems}). El comportamiento de las curvas que representan la solución de un sistema dinámico, y en particular las características de las oscilaciones de estas curvas son un tema de interés. Si es sumamente difícil o imposible predecir el estado de un sistema dinámico al cabo de un cierto número de iteraciones, puede ser que uno se encuentre frente a un sistema con dinámica caótica. 

En términos generales, es complicado decidir si un sistema dinámico presenta caos. Es válido pensar que un flujo caótico es aquel que tiende a ``mezclar'' trayectorias. Además ya se mencionó que el comportamiento de un flujo caótico puede resultar muy difícil o imposible de predecir. Lo anterior conduce a una de las definiciones más aceptadas en la literatura, y que es parte del trabajo de Devaney \cite[Definición 8.5]{Devaney1989AnSystems}.
\begin{df} Se dice que un sistema dinámico
\begin{eqnarray}
\dot x &=& f(x), \nonumber
\end{eqnarray}
es \textbf{caótico} o presenta \textbf{caos}, si su flujo tiene las siguientes tres características:
\begin{enumerate}
\item sensibilidad en condiciones iniciales, 
\item transitividad topológica, 
\item órbitas periódicas que forman un conjunto denso. 
\end{enumerate}
\end{df}

El concepto anterior es presentado y ampliamente analizado en \cite{Banks1992OnChaos}. Además, otras definiciones de caos que también se han utilizado pueden ser halladas en \cite[Sección 7.1]{Meiss2007DifferentialSystems}, \cite[Sección 3.5]{Robinson1998DynamicalChaos} o \cite[Capítulo 23]{Wiggins2003IntroductionChaos}. 

Cabe mencionar que el concepto de caos, a menudo es estudiado a través de la función de la herradura de Smale y del espacio de sucesiones bilaterales 
\begin{eqnarray}
\textbf{s} &=& \hdots s_{-3}s_{-2}s_{-1}s_{0}.s_1s_2s_3\hdots, \nonumber
\end{eqnarray}
donde $s_i = 0$ o $1$ para toda $i\in \mathbb Z$. Es fácil ver que $\textbf s$ corresponde a la escritura en forma binaria de un número real. El lector interesado puede consultar \cite[Capítulo 23]{Wiggins2003IntroductionChaos}.

 El caos ha sido reportado en modelos ecológicos por algunos autores, como por ejemplo \cite{Gakkhar2003} donde se presenta un modelo con dos presas\index{presa} y un depredador\index{depredador}; en \cite{Klebanoff1994ChaosTheory} se analiza la existencia de caos usando la teoría de bifurcaciones; en \cite{Kozlov2013OnApproach} se construye una familia de vectores y se usa aproximación por funciones para estudiar el caos en un modelo de Lotka-Volterra; en \cite{Schaffer1985OrderSystems} se muestra el procedimiento para asociar una función unidimensional a un sistema con caos, para posteriormente vincular estos resultados con otros que están relacionados con el estudio de \textit{Lynx canadensis}.

En este capítulo se retoma el modelo de depredación \eqref{4.1} que fue introducido en el Capítulo \ref{capitulo4}, y se fijan los parámetros para mostrar evidencia numérica de presencia de caos. Concretamente, el modelo que se analiza es el siguiente
\begin{eqnarray}
\label{5.1}
\dot x &=& x\left(r_1\left(1 - x/K_1\right) - \alpha_{12}y - \frac{q_1 z}{1 + a_1x}\right),   \nonumber \\
\dot y &=& y\left(r_2\left(1 - y/K_2\right) - \alpha_{21}x - \frac{q_2 z}{1 + a_2y}\right),    \\
\dot z &=& z\left(\frac{c_1 q_1 x}{1 + a_1x} + \frac{c_2 q_2 y}{1 + a_2y} - \mu - m z\right), \nonumber \\ \nonumber \\
\text{con}& & r_1 = K_1 = r_2 = K_2 = q_2 = \alpha_{12} = \mu = 1, \nonumber \\  
& & c_1 = c_2 = 0.5, \quad q_1 = 10, \quad \alpha_{21} = 1.5, \nonumber \\
& & a_1 = 0.01, \quad a_2 = 0.02, \quad m = 0.1. \nonumber
\end{eqnarray}

Haciendo un cálculo se deduce que los cinco puntos de equilibrio de \eqref{5.1} son
\begin{eqnarray}
E_0 \quad = \quad (0,0,0)\text{, } \quad E_1 &=& (1,0,0)\text{, }\quad E_2 \quad = \quad (0,1,0)\text{, } \nonumber \\ 
E_5 \quad = \quad (0.202418,0,0.0800811)\text{, }&& E^* \quad = \quad (0.11969,0.813857,0.00666116). \nonumber
\end{eqnarray}

La matriz jacobiana de \eqref{5.1} aplicada en el punto de equilibrio en la frontera $E_2 = (0,1,0)$ (véase el Teorema \ref{teo4.2}) está dada por
\begin{eqnarray}
J(E_2) &=& \left( \begin{array}{ccc}
                r_1 - K_2 \alpha_{12} & 0     & 0                                \\\\
                - K_2 r_2             & - r_2 & - K_2 q_2(1 + a_2 K_2)^{-1}                   \\\\
                0                     & 0     & c_2 K_2 q_2(1 + a_2 K_2)^{-1} - \mu \\
                \end{array} \right).                                  \nonumber
\end{eqnarray}

A los puntos de equilibrio que tienen un valor propio con parte real nula se les conoce comúnmente como \textbf{puntos no hiperbólicos}.
Esta característica la tiene el punto de equilibrio en la frontera $E_2 = (0,1,0)$ ya que $r_1/K_2 = \alpha_{12}$. Concretamente, el espectro está dado por
\begin{eqnarray}
\sigma(E_2) &=& \left\{-1,-0.509804,0\right\}, \nonumber
\end{eqnarray}
Lo anterior propicia que algunas trayectorias que pasan suficientemente cerca de $E_2$ permanezcan durante mucho tiempo oscilando alrededor de $E_2$ (lo que podría ser interpretado erróneamente como que $E_2$ es un atractor para dichas trayectorias), para que posteriormente sean repelidos hacia una cuenca de atracción de un conjunto extraño. En el Apéndice \ref{apendicee} se estudia la variedad central\index{variedad central} de $E_2$.

Por otra parte, se puede observar que las desigualdades que aparecen en el punto $1.$ del Teorema \ref{teo4.3} no se cumplen para \eqref{5.1}, y por lo tanto no existe un punto de equilibrio en la frontera de la forma $E_3 = (x_3,y_3,0)$ con $x_3$, $y_3 > 0$. Además de lo anterior, el valor de la tasa de captura del depredador\index{depredador} $z$ hacia la especie $x$ es alta ($q_1 = 10$) con respecto al resto de los parámetros en \eqref{5.1}, y esto implica que $x^b < K_1$, sin embargo al cumplirse la desigualdad $y^b > K_2$ se tiene como consecuencia que no existe un punto de equilibrio en la frontera de la forma $E_4 =(0,y_4,z_4)$ con $y_4$, $z_4 > 0$.

Si el parámetro $q_1$ es manipulado alrededor del valor especificado en \eqref{5.1}, entonces aún es posible hallar dinámica caótica. Ejemplo de lo anterior es tomar $q_1 = 9$, o bien $q_1 = 11$.  Aparentemente se puede conservar el caos en \eqref{5.1} si el parámetro de competencia intraespecífica es llevado a $m=0.8$, o bien si el tiempo de manejo es llevado a $a_1 = 1$. Si los parámetros $a_1$, $a_2$, $m$ se anulan, entonces se obtiene el caso estudiado en \cite{Takeuchi1983}.

\begin{figure}[h]
     \centering
     \includegraphics[width=439pt]{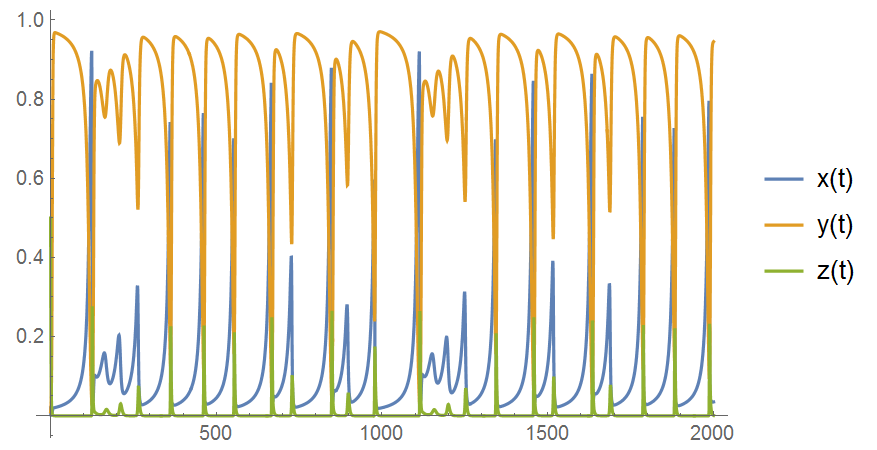}
      \caption{Soluciones de \eqref{5.1} con condiciones iniciales $x_0=0.02$, $y_0=0.01$, $z_0=0.5$ y con parámetro temporal $t$ en el intervalo $[0,2000]$.}
      \label{fig5.1}
\end{figure}

\begin{figure}
    \centering
    \includegraphics[width=300pt]{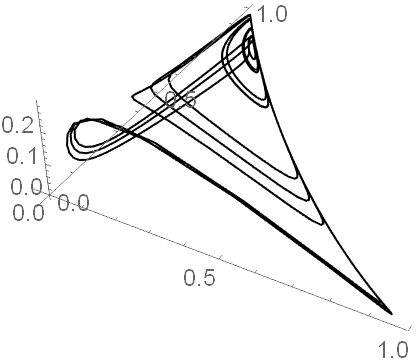}
    \includegraphics[width=300pt]{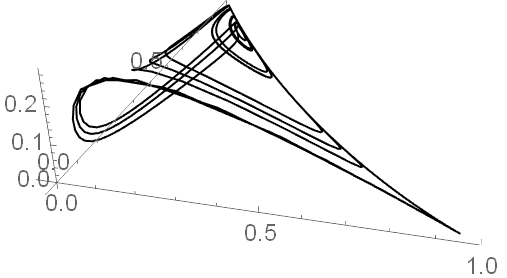}
   \caption{Bosquejo del conjunto atractor extraño de \eqref{5.1} con condiciones iniciales $x_0=0.02$, $y_0=0.01$, $z_0=0.5$, con los parámetros previamente mostrados y con $t\in [9000,10000]$ en el espacio de estados $[0,\infty)^3$.}
   \label{fig5.2}
\end{figure}

En la Figura \ref{fig5.1} se exhibe la evolución de una terna de trayectorias del sistema con las condiciones iniciales que se especifican. La parametrización de las soluciones aparece en la Figura \ref{fig5.2}, así como la aparente convergencia de las curvas hacia un atractor extraño.

A pesar de que no ha sido posible demostrar analíticamente la existencia de caos en \eqref{5.1}, numéricamente se han calculado parámetros como los exponentes de Lyapunov, la dimensión fractal\index{dimensión fractal} y la entropía de correlación.

\begin{figure}
    \centering
    \includegraphics[width=191pt]{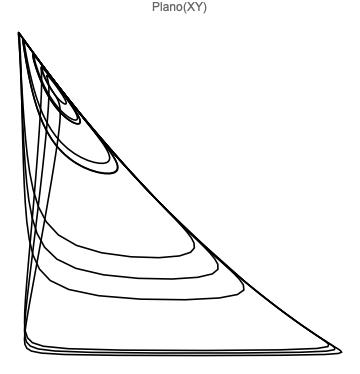}
    \includegraphics[width=191pt]{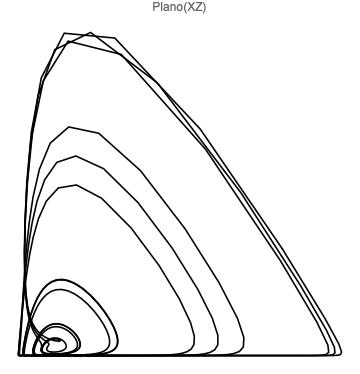}
    \includegraphics[width=191pt]{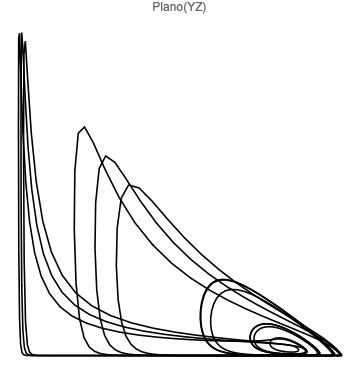}
    \caption{Las proyecciones del conjunto de la Figura \ref{fig5.2} en los planos coordenados.}
    \label{fig5.3}
\end{figure}

\section{Órbitas Heteroclínicas}
Se verifica fácilmente que el sistema \eqref{5.1} cuenta con órbitas heteroclínicas que conectan a cada uno de los puntos de equilibrio $E_1 = (1,0,0)$ $E_2 = (0,1,0)$ con el punto de equilibrio en el origen de coordenadas $E_0 = (0,0,0)$. Específicamente

\begin{eqnarray}
\psi_1(t) &=& \left(\left(1 - c_1e^{-t}\right)^{-1},0,0\right),\quad \text{con }c_1\in(0,1), \nonumber \\
\psi_2(t) &=& \left(0,\left(1 - c_2e^{-t}\right)^{-1},0\right),\quad \text{con }c_2\in(0,1). \nonumber
\end{eqnarray}

Se sabe que $\psi_1$, $\psi_2$ son soluciones de la ecuación del modelo logístico\index{modelo logístico}, y que se pueden calcular usando métodos conocidos.

Aparentemente también existe una órbita de tipo heteroclínico que conecta a los puntos $E_1$ y $E_2$, sin embargo aún no ha sido posible hallar una expresión analítica para dicha órbita.

\begin{figure}[h]
    \centering
    \includegraphics[width=425pt]{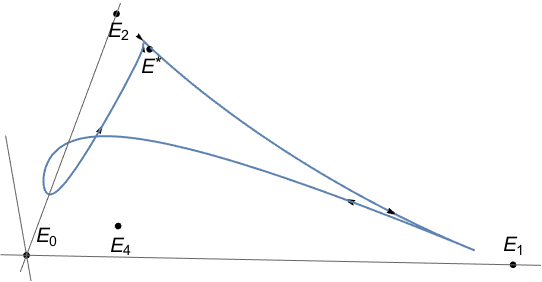}
    \caption{Se muestra un tramo de las trayectorias del sistema \eqref{5.1}. Los puntos de equilibrio tienen coordenadas: $E_0 = (0,0,0)$, $E_1 = (1,0,0)$, $E_2 = (0,1,0)$, $E_4 = (0.202418,0,0.0800811)$, $E^* = (0.11969,0.813857,0.00666116)$. Los valores propios de la matriz jacobiana de \eqref{5.1} aplicada a cada uno de los puntos de equilibrio están dados por: $\sigma(E_0) = \{-1,1,1\}$, $\sigma(E_1) = \{3.90196, -1., -0.5\}$, $\sigma(E_2) = \{-1., -0.50495, 0\}$, $\sigma(E_4) = \{-0.103605 + 0.889713 i, -0.103605 - 0.889713 i, 0.616292\}$, $\sigma(E^*) = \{-0.966449, 0.016224 + 0.158434 i, 0.016224 - 0.158434 i\}$}
    \label{fig5.4}
\end{figure}

En la Figura \ref{fig5.4} se muestra un tramo de las trayectorias junto con la ubicación de los cinco puntos de equilibrio de \eqref{5.1}. La trayectoria pasa en las inmediaciones de los puntos de equilibrio $E_0$, $E_1$ y $E_2$ para diferentes valores de $t$. A partir de eso, se puede intuir que la trayectoria en cuestión está ubicada cerca de órbitas heteroclínicas que conectan a parejas de puntos de equilibrio.

\section{Exponentes de Lyapunov}
Una de las propiedades más notables en la dinámica caótica es la sensibilidad en las condiciones iniciales. En \cite[Sección 5.1]{Kantz2004NonlinearAnalysis}, se describe lo anterior como la imposibilidad de predecir escenarios futuros a pesar de una evolución de tiempo determinista. En \cite[Capítulo 2, Subsección 4.2]{Diks1999Non-linearApplications} se menciona que en estos casos se presenta una separación exponencial entre las trayectorias, y que los exponentes de Lyapunov permiten medir esta divergencia exponencial. En \cite{Sandri1996}, los autores han desarrollado algoritmos en \textit{Mathematica} para el cálculo de dichos exponentes.

El papel que juega el exponente de Lyapunov $\lambda$, está dado por la relación \cite[Sección 5.2]{Kantz2004NonlinearAnalysis}
\begin{eqnarray}
|\varphi(x_0,t) - \varphi(x_0 + h,t)| &\simeq& |h| e^{\lambda t}, \nonumber
\end{eqnarray}
donde $\varphi$ es el flujo del sistema dinámico, y los valores $h$ y $t$ se aproximan a cero y a infinito respectivamente. Con base en eso, se tiene \cite[Capítulo 2, Subsección 4.2]{Diks1999Non-linearApplications}
\begin{eqnarray}
\lambda &=& \lim_{t\to \infty} \left(\lim_{h\to 0}\frac{1}{t}\ln \left|\frac{\varphi(x_0,t) - \varphi(x_0 + h,t)}{h}\right|\right) \nonumber \\
&=& \lim_{t\to \infty} \frac{1}{t}\ln\left|\frac{\partial \varphi (x_0,t)}{\partial x}\right|. \nonumber
\end{eqnarray}

\begin{figure}[h]
    \centering
    \includegraphics[width=242pt]{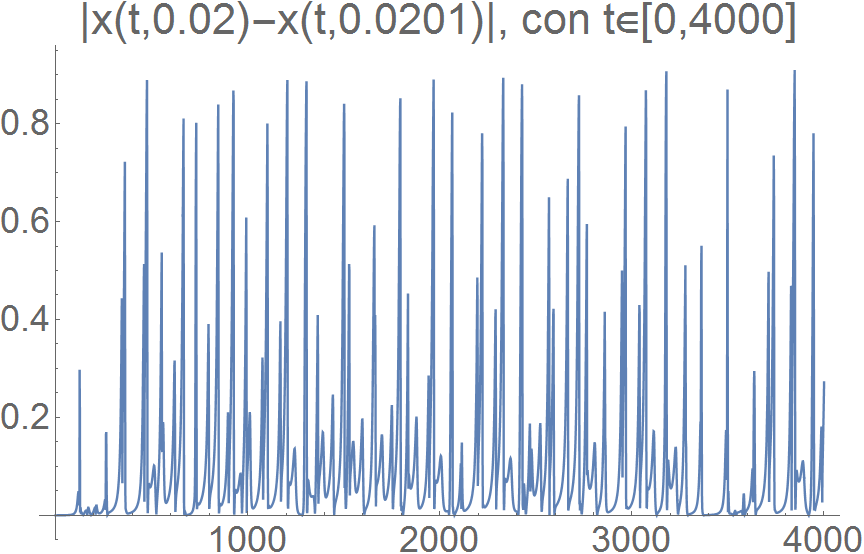}
    \includegraphics[width=242pt]{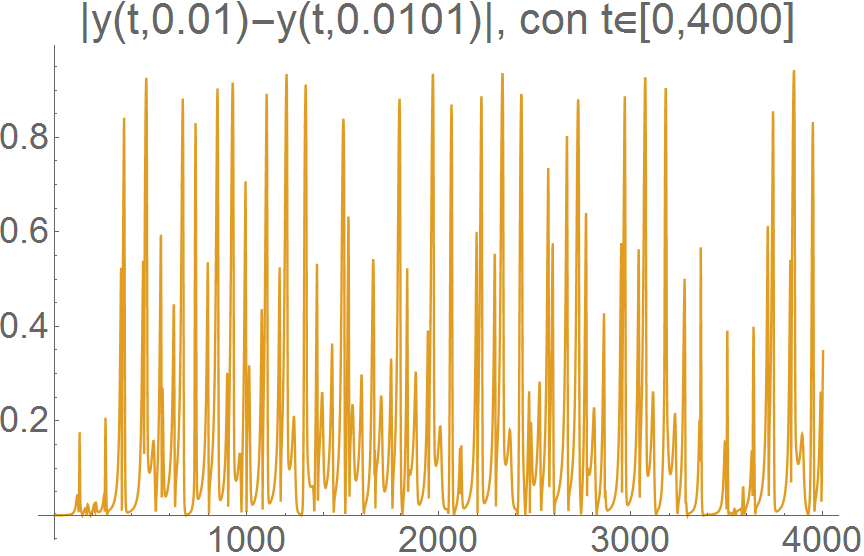}
    \includegraphics[width=242pt]{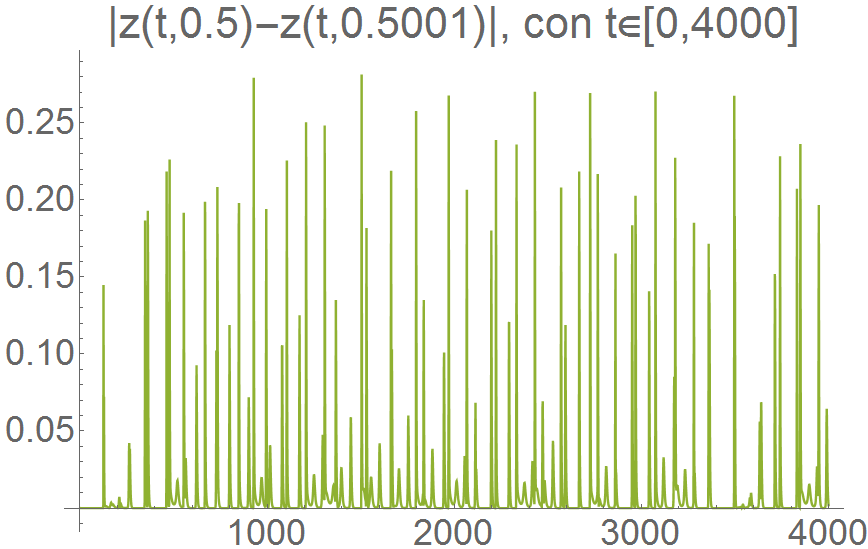}
\caption{Distancias que guardan dos de las soluciones de \eqref{5.1}. Las coordenadas sufren una perturbación de un diezmilésimo de unidad. La gráfica de la coordenada $x$ está en azul, la de $y$ está en amarillo, y la de $z$ está en verde.}
    \label{fig5.5}
\end{figure}

A continuación se realiza el cálculo de los exponentes de Lyapunov de $\eqref{5.1}$. Dadas las condiciones iniciales $x_0 = 0.02$, $y_0 = 0.01$, $z_0=0.5$, se provee un incremento considerablemente pequeño ($h = 0.0001$) y se elabora una gráfica en la que se representa la distancia que guardan las respectivas soluciones, como se muestra en la Figura \ref{fig5.5}.

\begin{figure}
    \centering
    \includegraphics[width=242pt]{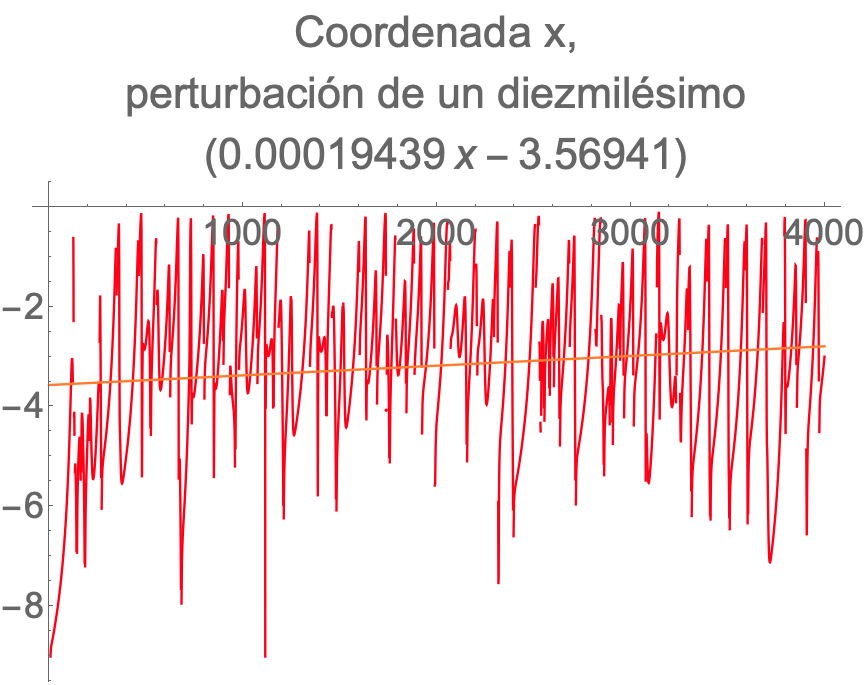}
    \includegraphics[width=242pt]{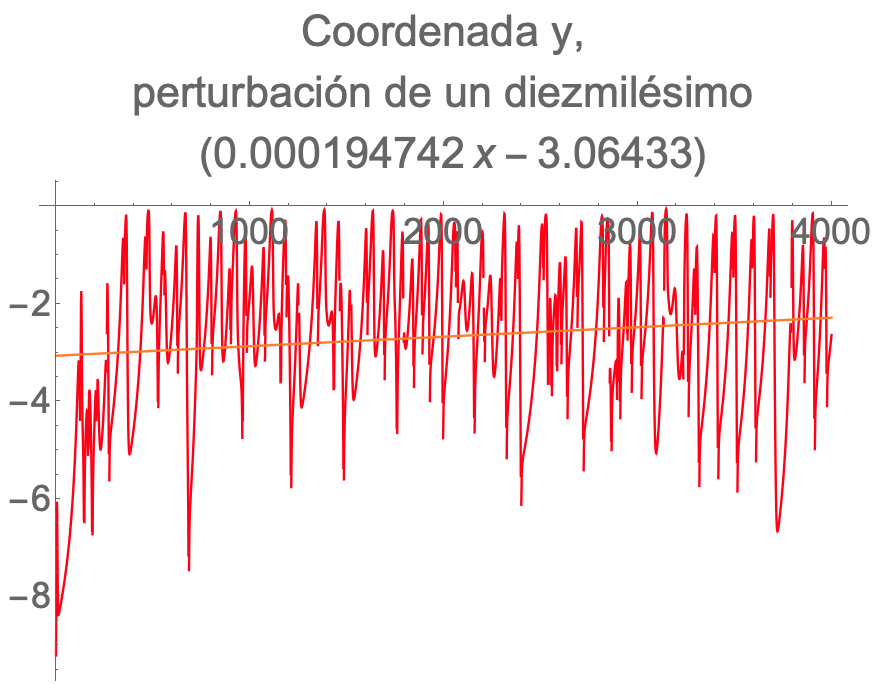}
    \includegraphics[width=242pt]{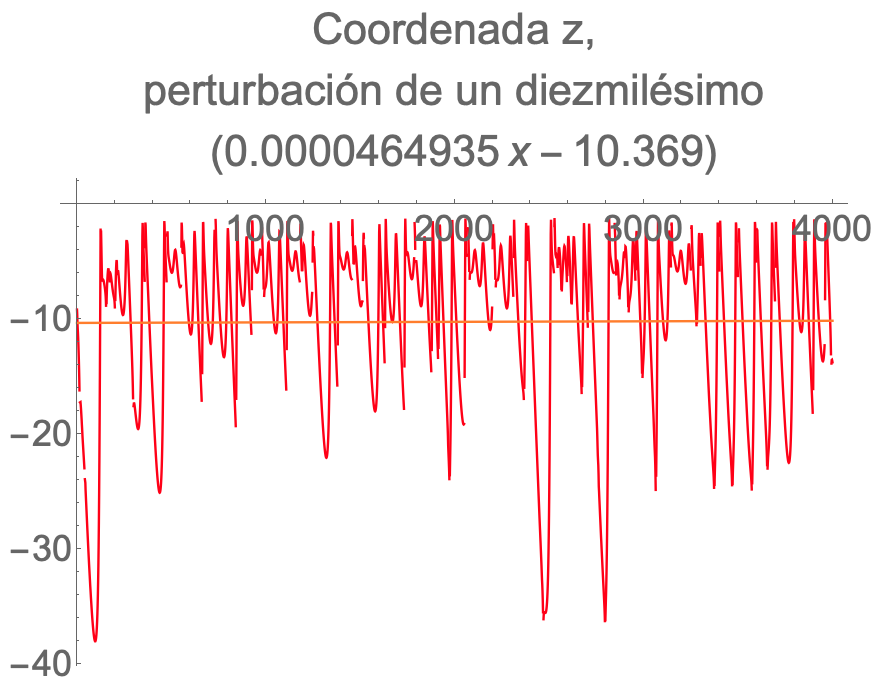}
\caption{Al aplicar logaritmo natural a los valores obtenidos en la Figura \ref{fig5.5}, se obtiene una estimación de los exponentes de Lyapunov, que están dados por la pendiente de la aproximación lineal en cada caso. La primera gráfica corresponde a la coordenada $x$, tiene a los intervalos $(0.000149348,0.000239432)$ y $(-3.67344,-3.46538)$ con un nivel de confianza de $95\%$ y tiene un $P$-value de $3.66381\times 10^{-17}$. La segunda gráfica corresponde a la coordenada $y$, tiene a los intervalos $(0.000151734,0.000237751)$ y $(-3.16366,-2.965)$ con un nivel de confianza de $95\%$ y tiene un $P$-value de $1.01393\times 10^{-18}$. La tercera gráfica corresponde a la coordenada $z$, tiene intervalos de confianza $(-0.00014998,0.000242967)$ y $(-10.8228,-9.91524)$ con un nivel de confianza de $95\%$ y tiene un $P$-value de $0.642712$.}
\label{fig5.6}
\end{figure}

Después de esto, se considera el logaritmo natural de los valores obtenidos en las series de tiempo de la Figura \ref{fig5.5}, y se construye una aproximación lineal. Esta aproximación aparece en la Figura \ref{fig5.6}. El coeficiente lineal en cada una de las ecuaciones obtenidas es una es una estimación del exponente característico de Lyapunov. Los valores de los exponentes de Lyapunov obtenidos son positivos en los casos que aquí fueron analizados.

$ $ \newpage

$ $ \newpage

\section{Dimensión Fractal\index{dimensión fractal}}
En el estudio de la geometría fractal, resulta fundamental el concepto de dimensión. La construcción hecha por Caratheodory\index{Caratheodory} es la base de la noción que Hausdorff\index{Hausdorff} adoptó. En muchas ocasiones, es difícil calcular o estimar la dimensión de Hausdorff\index{Hausdorff} a través de métodos computacionales. En términos coloquiales, la dimensión fractal\index{dimensión fractal} de Hausdorff\index{Hausdorff} nos permite conocer cuanto espacio ocupa un subconjunto de un espacio euclidiano cerca de cada uno de sus puntos \cite[Capítulo 2]{Falconer2003FRACTALApplications}.

Algunos conceptos preliminares de la dimensión fractal\index{dimensión fractal} se presentan a continuación.

\begin{df}
Sea $U\subset \mathbb R^n$ un conjunto abierto no vacío. El diámetro de $U$ está dado por
\begin{eqnarray}
\text{diam} (U) &=& \sup \left\{\|x - y\|\colon x,y \in U\right\} \nonumber
\end{eqnarray}
\end{df}
A partir de lo anterior, se describe lo siguiente
\begin{df}
Sean $F\subset \mathbb R^n$, $s \geq 0$, $\delta > 0$. La $(\delta,s)-$medida de Hausdorff\index{Hausdorff} del conjunto $F$ está dada por
\begin{eqnarray}
H_\delta^s(F) &=& \inf\left\{\sum_{i=i}^\infty \left(\text{diam}(U_i)\right)^s\colon F\subset \bigcup_{i=1}^\infty U_i \text{, }\text{diam}(U_i)\leq \delta\text{ para toda }i\in \mathbb N\right\}. \nonumber
\end{eqnarray}
\end{df}
Obsérvese que si el diámetro $\delta$ tiende a cero, entonces existen menos cubiertas $\{U_i\}_{i=1}^\infty$ de $F$  admisibles. Por lo anterior, el valor de $H_\delta^s(F)$ puede aumentar conforme $\delta$ decrece, y tiene sentido pensar en el valor límite.
\begin{df}
Sean $F\subset \mathbb R^n$, $s \geq 0$. La $s -$ medida de Hausdorff\index{Hausdorff} del conjunto $F$ está dada por
\begin{eqnarray}
H^s(F) &=& \lim_{\delta\to 0}H_\delta^s(F). \nonumber
\end{eqnarray}
\end{df}

Es posible mostrar que $H^s$ es una medida. De hecho, para los casos $s = 1,2,3$, $H^s$ generaliza las nociones conocidas de longitud, área y volumen, respectivamente.

\begin{df}
La dimensión de Hausdorff\index{Hausdorff} de un conjunto está dada por
\begin{eqnarray}
\dim_H(F) &=& \inf\left\{s\geq 0 \colon H^s(F) = 0\right\} \nonumber \\
&=& \sup\left\{s\geq 0 \colon H^s(F) = \infty\right\}. \nonumber
\end{eqnarray}
\end{df}
Para profundizar más en los conceptos expuestos previamente, el lector interesado puede consultar \cite[Capítulo 2]{Falconer2003FRACTALApplications}.

A través de un algoritmo computacional, es posible graficar el número total de píxeles en una celda de dimensión $2$ contra el número de píxeles que cubre a las proyecciones del conjunto en el $Plano(XY)$, en el $Plano(XZ)$ y en el $Plano(YZ)$ que fueron mostradas en la Figura \ref{fig5.3}.

\begin{prop}$ $
\begin{enumerate}
\item La dimensión fractal estimada de la proyección del conjunto en el $Plano(XY)$ es $1.57858$.
\item La dimensión fractal estimada de la proyección del conjunto en el $Plano(XZ)$ es $1.55182$.
\item La dimensión fractal estimada de la proyección del conjunto en el $Plano(YZ)$ es $1.61775$.
\end{enumerate}
\end{prop}
El algoritmo que se usó fue construido por Szabolcs Horvát, Center for Systems Biology Dresden.
\begin{figure}
    \centering
    \includegraphics[width=368pt]{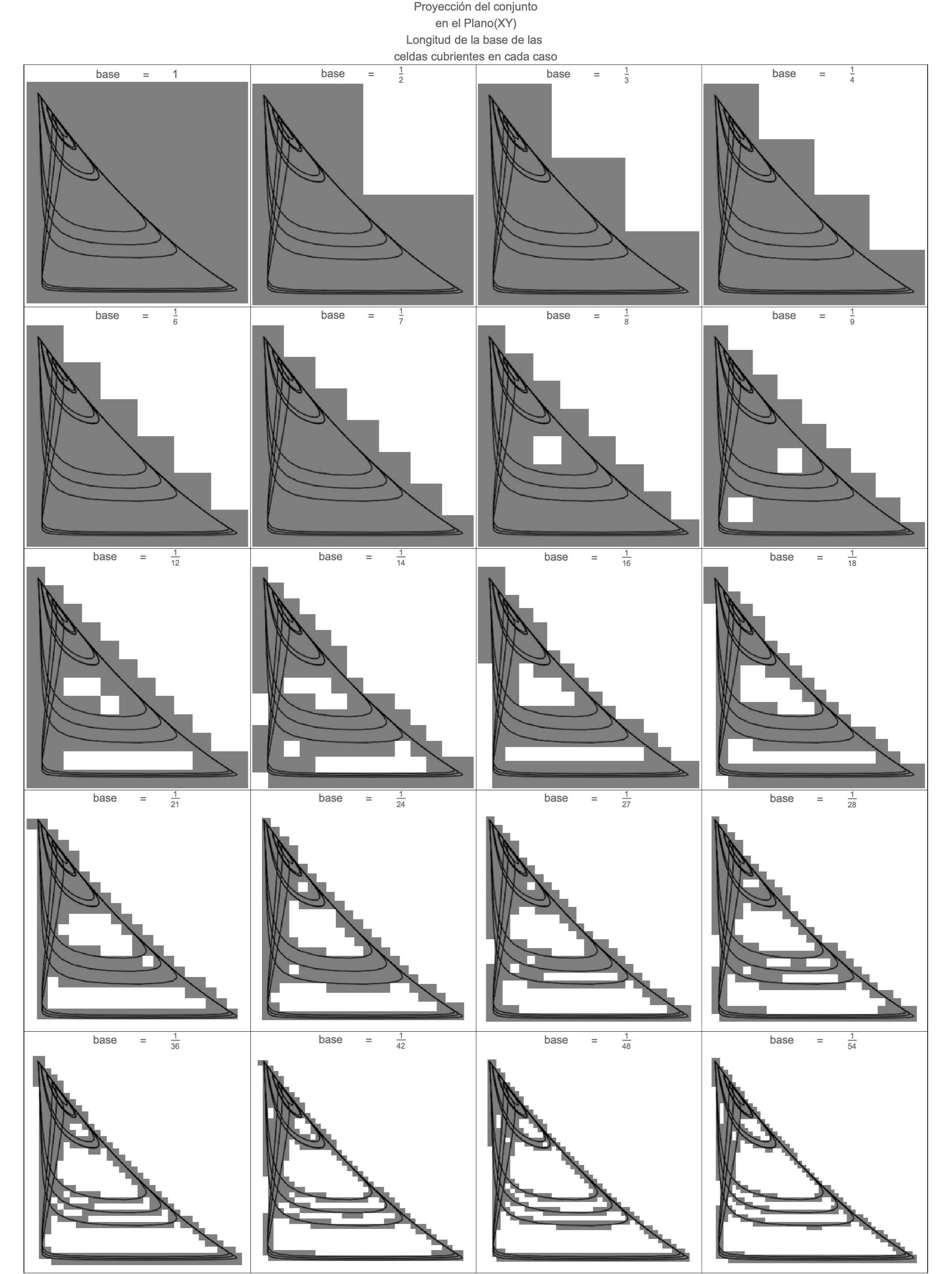}
    \caption{Se consideran 20 tamaños diferentes para las celdas (vecindades canónicas con la norma $\|\cdot\|_1$) que cubren al conjunto extraño. Al disminuir el tamaño de las celdas, aumenta el número de celdas requeridas para cubrir al objeto. }
    \label{fig5.7}
\end{figure}

\begin{figure}
    \centering
    \includegraphics[width=368pt]{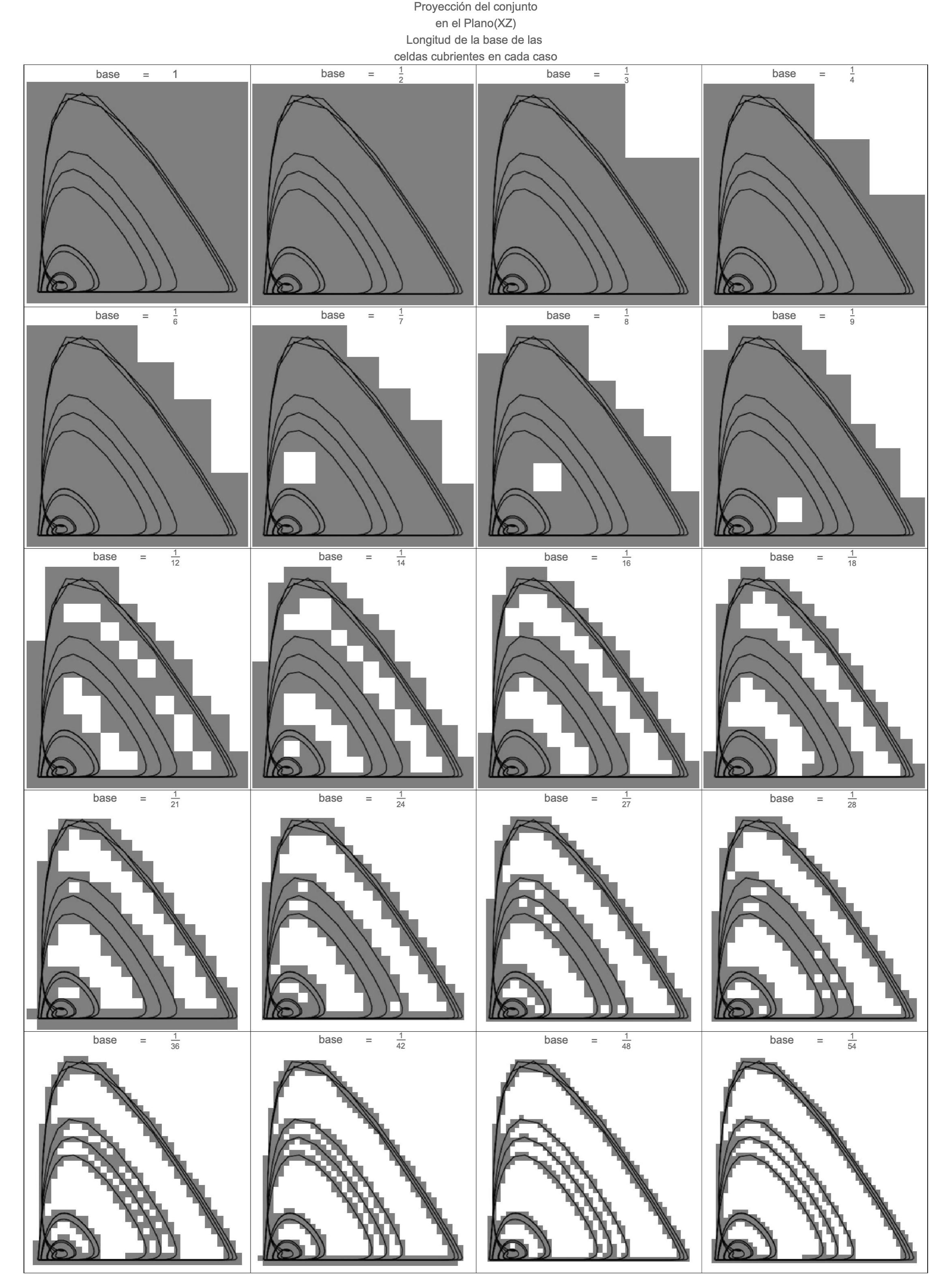}
    \caption{Se consideran 20 tamaños diferentes para las celdas (vecindades canónicas con la norma $\|\cdot\|_1$) que cubren al conjunto extraño. Al disminuir el tamaño de las celdas, aumenta el número de celdas requeridas para cubrir al objeto. }
    \label{fig5.8}
\end{figure}

\begin{figure}
    \centering
    \includegraphics[width=368pt]{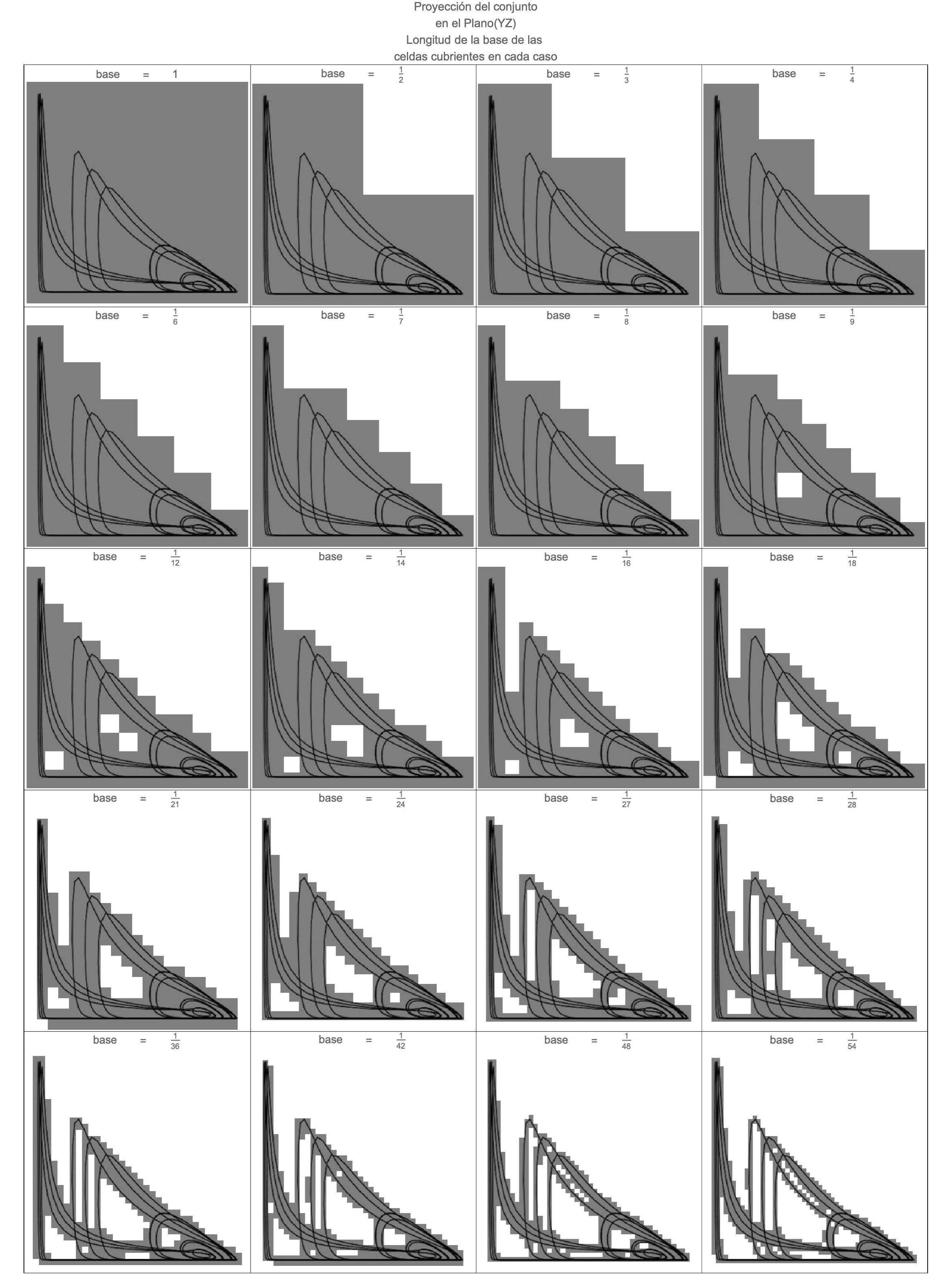}
    \caption{Se consideran 20 tamaños diferentes para las celdas (vecindades canónicas con la norma $\|\cdot\|_1$) que cubren al conjunto extraño. Al disminuir el tamaño de las celdas, aumenta el número de celdas requeridas para cubrir al objeto. }
    \label{fig5.9}
\end{figure}

\begin{figure}
    \centering
    \includegraphics[width=230pt]{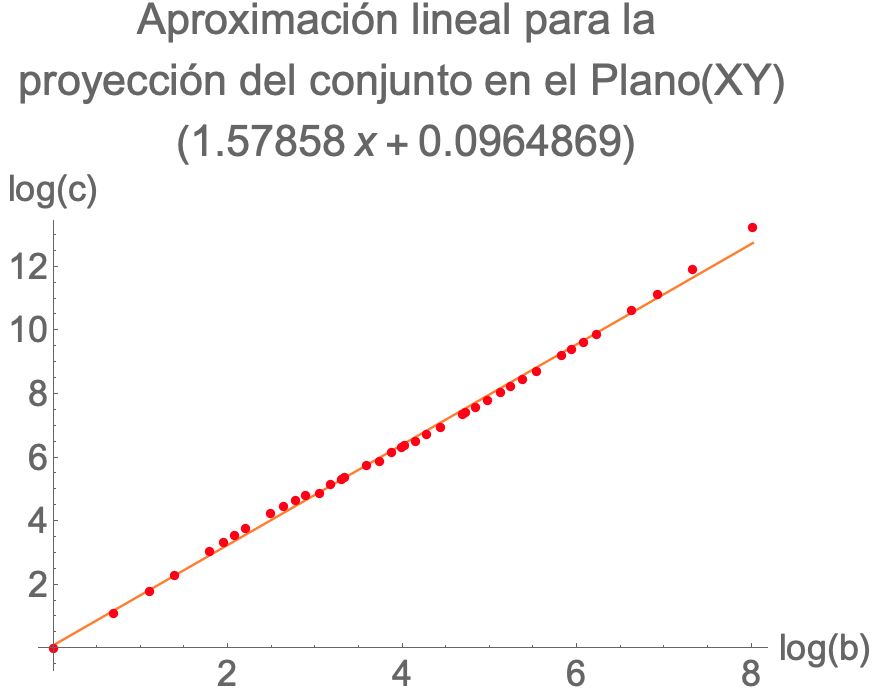}
    \includegraphics[width=230pt]{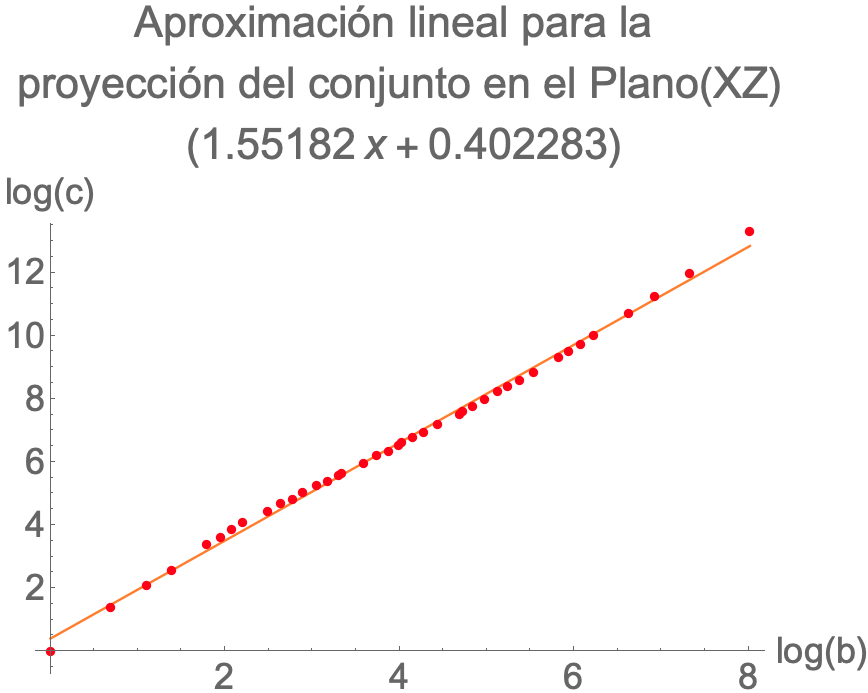}
    \includegraphics[width=230pt]{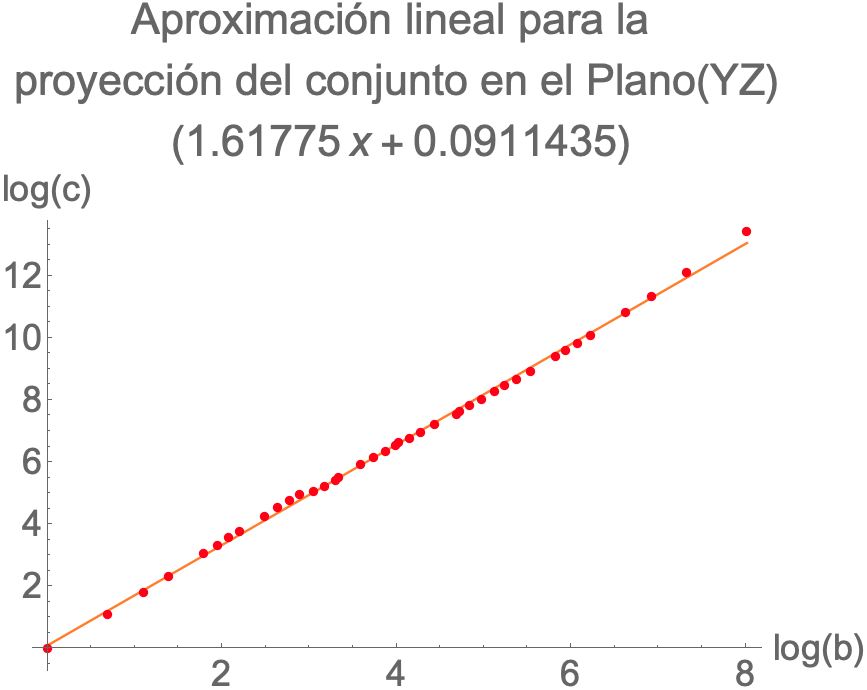}
    \caption{Aproximación lineal para los datos obtenidos en las simulaciones de las Figuras \ref{fig5.7}, \ref{fig5.8} y \ref{fig5.9}. Se toma como abscisa el logaritmo del inverso multiplicativo de la longitud de la base de la celda, y como ordenada el logaritmo del número de celdas usadas para cubrir el conjunto en su totalidad. La pendiente de la recta es $1.57858$ para la proyección en el $Plano(XY)$, $1.55182$ para la proyección en el $Plano(XZ)$, y $1.61775$ para la proyección en el $Plano(YZ)$. En cada caso, la pendiente de la recta es una estimación de la dimensión fractal\index{dimensión fractal} de la proyección correspondiente del conjunto extraño.}
    \label{fig5.10}
\end{figure}

$ $\newpage
\section{Dimensión de Correlación y Entropía de Correlación}
Para el cálculo de la dimensión de correlación y la entropía de correlación, se siguen las ideas de  \cite[Capítulo 2, Secciones 2 y 3]{Diks1999Non-linearApplications}. Como su nombre lo indica, el teorema de reconstrucción (véase por ejemplo \cite{Packard1980GeometrySeries} y \cite{Takens1981DetectingTurbulence}) permite llevar a cabo la reconstrucción de la dinámica asintótica de un sistema a partir de observaciones hechas. En otras palabras, es posible encontrar una correspondencia biunívoca entre un conjunto atractor de un sistema dinámico 
\begin{eqnarray}
\dot x &=& F(x) \nonumber
\end{eqnarray}
con un atractor reconstruido que se obtiene al considerar los vectores de la forma
\begin{eqnarray}
x_n &=& (x(n),x(n+1),\hdots,x(n+m-1))\in \mathbb R^m, \nonumber
\end{eqnarray}
para toda $n\in \mathbb N$.

\begin{df} La integral de correlación $C_m(r)$ de dimensión $m$ y radio $r$ indica la probabilidad con la que dos vectores $\overline x$, $\overline y$ de $\mathbb R^m$ queden ubicados a distancia menor que $r$. Esto es
\begin{eqnarray}
C_m(r) &=& \int \int \Theta (r - \|\overline x - \overline y\|)d\overline yd\overline x, \nonumber
\end{eqnarray}
donde $\Theta$ es la función de Heaviside
\begin{eqnarray}
\Theta(s) &=& \begin{cases}
         0, & \mbox{si $x < 0$,}\\\\
         1, & \mbox{si $x \geq 0$}.\end{cases} \nonumber
\end{eqnarray}
Además, la integral de correlación satisface
\begin{eqnarray}
C_m(r) \approx e^{-m\tau K_2}r^{D_2}\text{, si }m\to \infty\text{, }r\to 0, \nonumber
\end{eqnarray}
donde $\tau$ representa el retardo usado en la reconstrucción (aquí se usa $\tau = 1$), $D_2$ es la dimensión de correlación y $K_2$ la entropía de correlación.
\end{df}

\begin{figure}
    \centering
    \includegraphics[width=255pt]{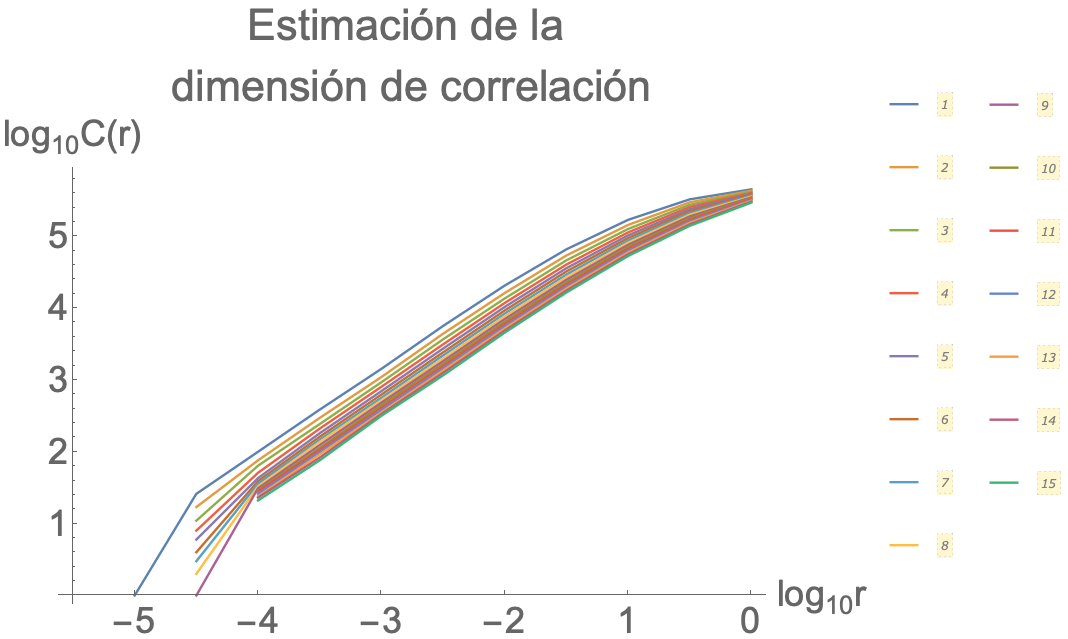}
    \includegraphics[width=255pt]{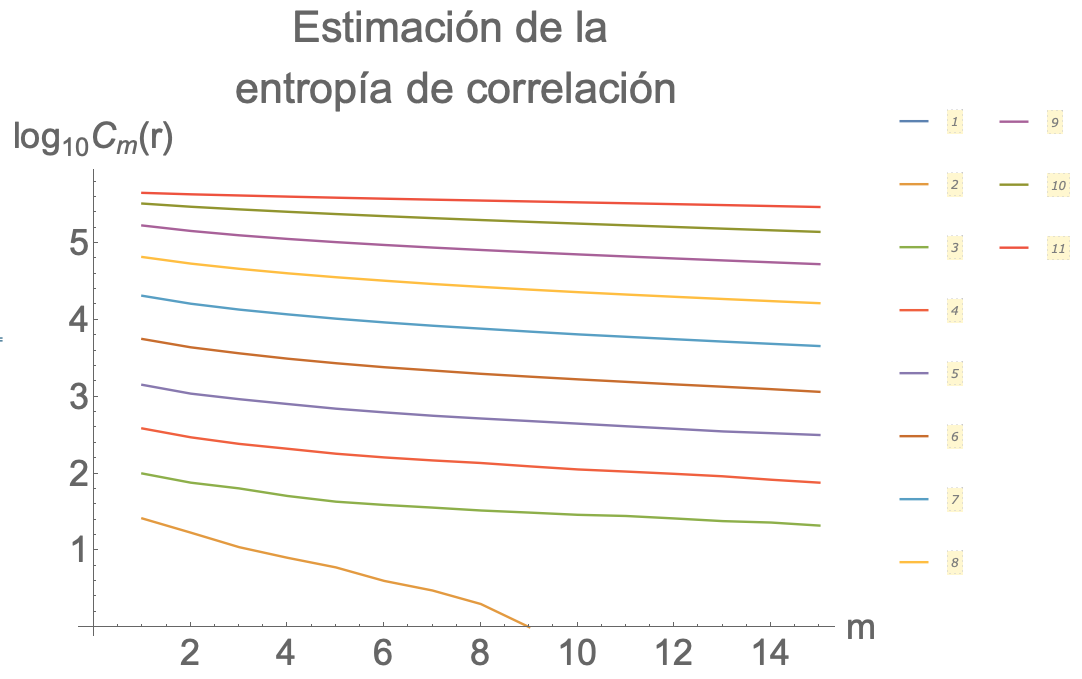}
    \caption{La dimensión de correlación $D_2$ y la entropía de correlación $K_2$ del sistema \eqref{5.1} se pueden aproximar a través las pendientes de las gráficas que se obtuvieron en cada caso y que aquí son mostradas. En la primera gráfica, se considera la dimensión de los vectores $m$ entre $1$ y $15$, mientras que en el segundo caso las dimensiones van de $1$ a $11$.}
    \label{fig5.11}
\end{figure}

Es posible dar una estimación de la integral de correlación $C_m(r)$ a través del método de Grassberger-Procaccia (véase \cite{PhysRevLett.50.346} y \cite{Grassberger1983MeasuringAttractors}), y así conocer un valor aproximado de $D_2$ y $K_2$. Dada la serie de tiempo $\{x(n)\}_{n=1}^L$ de longitud $L$, se consideran $N = L - (m - 1)$  vectores de reconstrucción de la forma $\boldsymbol x_n = (x_n,x_{n+1},\hdots,x_{n-m+1}) \in \mathbb R^m$, y se obtiene
\begin{eqnarray}
\widehat C_m(r) &=& \frac{2}{N(N-1)}\sum_{i=1}^{N-1}\sum_{j=i+1}^N\Theta(r - \|\boldsymbol x_i - \boldsymbol x_j\|). \nonumber
\end{eqnarray}
Usando esto, se estiman los valores de $D_2$ y $K_2$ para el sistema \eqref{5.1} con los parámetros especificados. En las gráficas de la Figura \ref{fig5.11} se muestra que las pendientes de las rectas se aproximan a los valores $D_2 = 1.16$ y $K_2 = -0.08$ para el sistema \eqref{5.1}.

\section{Espectro Potencia}
En esta sección, se analiza el espectro potencia del sistema \eqref{5.1}. En \cite[Sección 11.4]{Wei1989TimeMethods} se establece que una sucesión periódica $\{Z_t\}_{t=1}^\infty$ (con período $n$) se puede escribir como la combinación lineal de las funciones ortogonales seno y coseno de la siguiente forma
\begin{eqnarray}
Z_t &=& \sum_{k = 0}^{n/2} (a_k\cos \omega_k t+ b_k\sen \omega_k t), \nonumber 
\end{eqnarray}
donde
\begin{eqnarray}
a_k &=& \begin{cases}
      \displaystyle \frac{1}{n}\sum_{t=1}^n Z_t\cos \omega_k t, & \displaystyle k = 0\text{, y }k=\frac{n}{2}\text{ si }n\text{ es par,} \\\\
      \displaystyle \frac{2}{n} \sum_{t=1}^n Z_t\cos \omega_k t,& \displaystyle k = 1,2,\hdots,\left[\frac{n-1}{2}\right],
    \end{cases} \nonumber
\end{eqnarray}
\begin{eqnarray}
b_k &=& \frac{2}{n}\sum_{t=1}^n Z_t\sen\omega_k t,\quad k = 1,2,\hdots,\left[\frac{n-1}{2}\right], \nonumber
\end{eqnarray}
con $\omega_k = 2k\pi/n$. De forma similar, se puede escribir también $Z_t$ como la combinación lineal de los exponentes complejos
\begin{eqnarray}
Z_t &=& \begin{cases}
      \displaystyle \sum_{k=-(n-1)/2}^{(n-1)/2} c_ke^{i\omega_k t}, & \displaystyle \text{ si }n\text{ es impar,} \\\\
      \displaystyle \sum_{k=-n/2+1}^{n/2} c_ke^{i\omega_k t}, & \displaystyle \text{ si }n\text{ es par,}
    \end{cases} \nonumber
\end{eqnarray}
donde
\begin{eqnarray}
c_k &=& \frac{1}{n}\sum_{t=1}^n Z_t e^{-i\omega_k t}. \nonumber
\end{eqnarray}
Como las funciones $\sen \omega_k t$, $\cos \omega_k t$ y las funciones exponenciales complejas $e^{i\omega_k t}$ tienen todas período $n$, se concluye que la sucesión $\{Z_t\}_{t=1}^\infty$ tiene también período $n$.

La energía asociada con la sucesión $\displaystyle\{Z_t\}_{t=1}^\infty$ está dada por
\begin{eqnarray}
\sum_{t=1}^n Z_t^2 &=& \begin{cases}
\displaystyle na_0^2 + \frac{n}{2}\sum_{k=1}^{(n-1)/2}(a_k^2 + b_k^2), &  \text{si }n\text{ es impar,} \\\\
\displaystyle na_0^2 + \frac{n}{2}\sum_{k=1}^{(n-1)/2}(a_k^2 + b_k^2) + na_{n/2}^2, & \text{si }n\text{ es par.}
\end{cases}
\nonumber \\ \nonumber \\
&=& \begin{cases}
\displaystyle n\sum_{k=-(n-1)/2}^{(n-1)/2}|c_k|^2, & \text{si }n\text{ es impar,} \\\\
\displaystyle n\sum_{k=-n/2 + 1}^{n/2}|c_k|^2, & \text{si }n\text{ es par.}
\end{cases}
\nonumber
\end{eqnarray}
Lo anterior se sigue de la identidad de Parseval para series de Fourier (ver por ejemplo \cite[Teorema 8.16]{Rudin1980PrincipiosMatematico}). Las ecuaciones previas muestran que el total de energía en una sucesión periódica tomando los valores $t = 1,2,\hdots$ es infinita. Por ende, se considera la energía por unidad de tiempo, que se conoce como la potencia de la sucesión. Esto es
\begin{eqnarray}
\text{Potencia}  &=& \begin{cases}
\displaystyle a_0^2 + \frac{1}{2}\sum_{k=1}^{(n-1)/2}(a_k^2 + b_k^2), &  \text{si }n\text{ es impar,} \\\\
\displaystyle a_0^2 + \frac{1}{2}\sum_{k=1}^{(n-1)/2}(a_k^2 + b_k^2) + a_{n/2}^2, & \text{si }n\text{ es par.}
\end{cases}
\nonumber \\ \nonumber \\
&=& \begin{cases}
\displaystyle \sum_{k=-(n-1)/2}^{(n-1)/2}|c_k|^2, & \text{si }n\text{ es impar,} \\\\
\displaystyle \sum_{k=-n/2 + 1}^{n/2}|c_k|^2, & \text{si }n\text{ es par.}
\end{cases}
\nonumber
\end{eqnarray}
Una de las características principales del espectro potencia es el hecho de que muestra una serie de tiempo con un número de cúspides igual al número de osciladores principales con el que cuenta el sistema subyacente.

\begin{figure}[h]
    \centering
    \includegraphics[width=170pt]{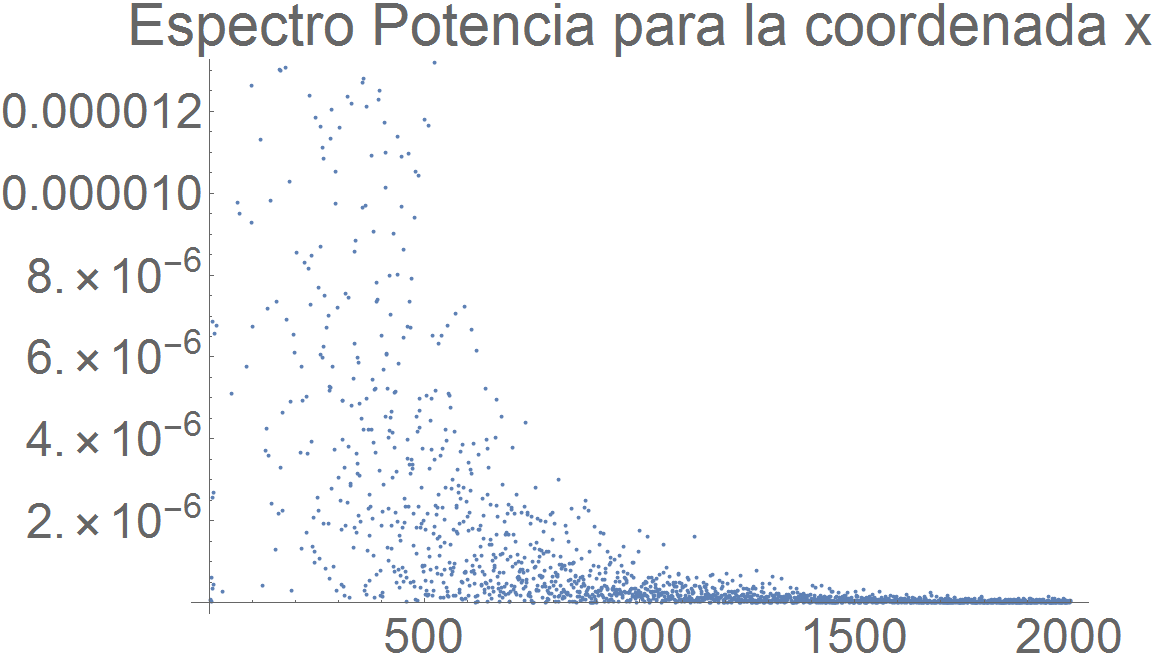}
    \includegraphics[width=170pt]{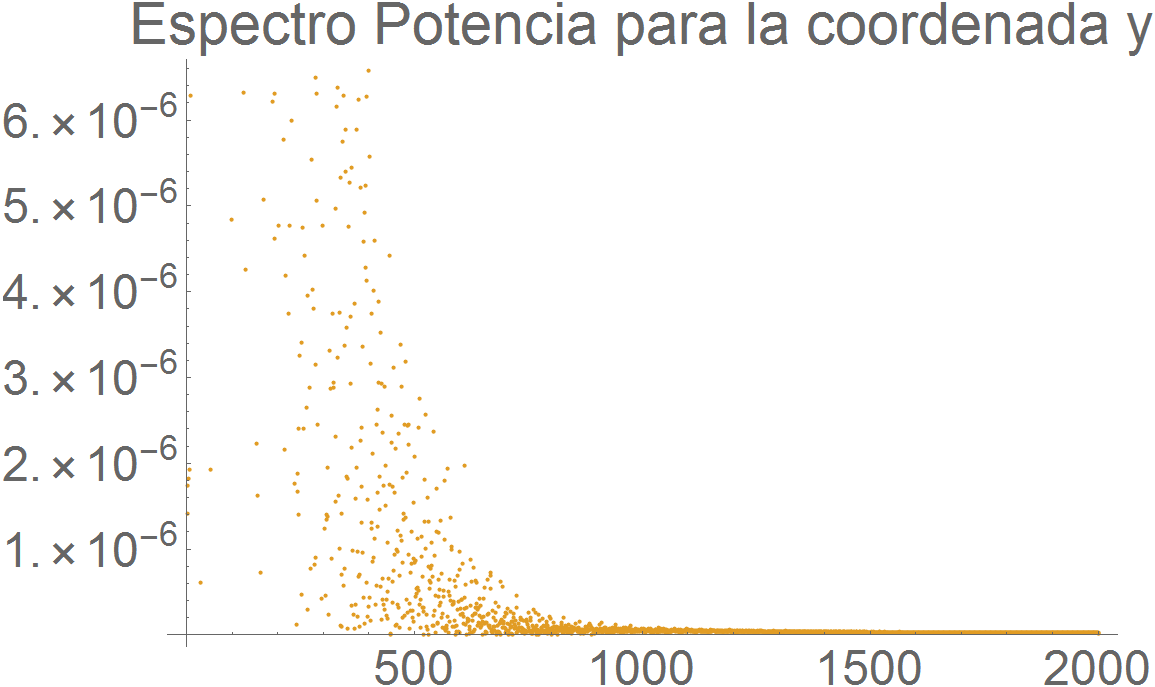}
    \includegraphics[width=170pt]{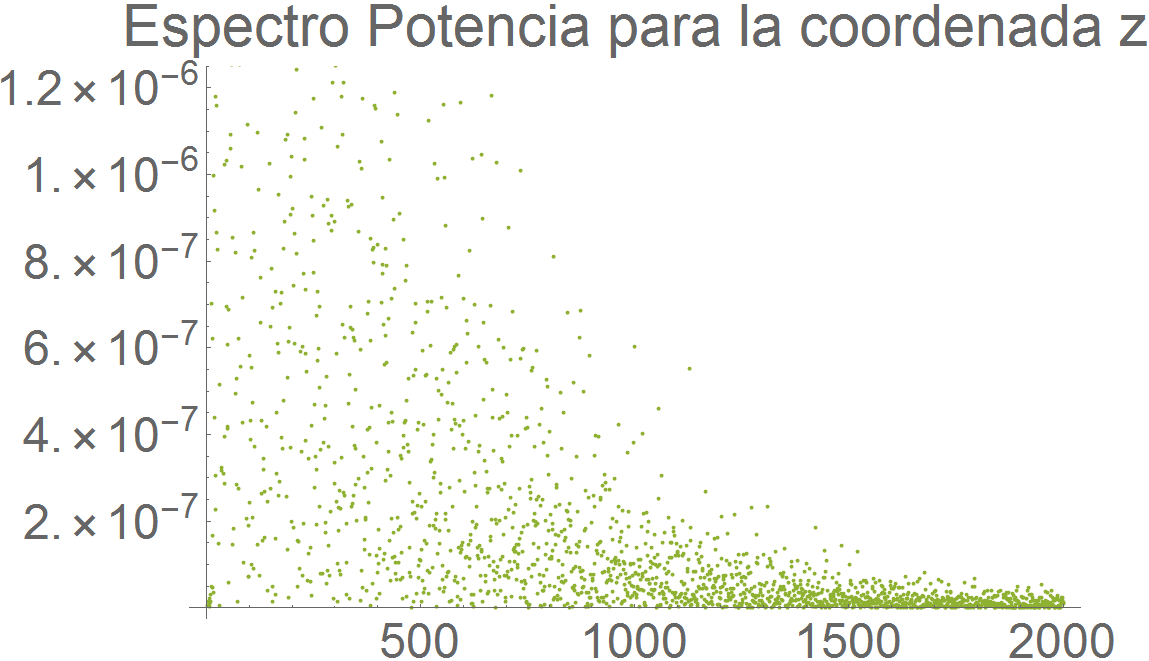}
\caption{En cada caso, se muestra el espectro potencia correspondiente a cada una de las coordenadas de la solución de \eqref{5.1} con los parámetros previamente mencionados. La coordenada $x$ tiene una energía de $143.062$ y una potencia de $0.0357655$; la coordenada $y$ tiene una energía de $2791.01$ y una potencia de $0.697752$; la coordenada $z$ tiene una energía de $2.67374$ y una potencia de $0.000668434$.}
    \label{fig5.12}
\end{figure}

En la Figura \ref{fig5.12} se muestra el cálculo del espectro potencia de \eqref{5.1}. Se puede apreciar la irregularidad de los puntos en la gráfica en cada una de las ilustraciones.

A continuación se presenta una prueba binaria denominada \textbf{Test 0-1}, y que en literatura reciente ha sido implementada para decidir si un sistema determinista tiene dinámica caótica. Se trata de una herramienta en la que se considera solo una parte del espectro potencia.

En \cite{Gottwald2004a} y en \cite{Gottwald2009OnChaos}, se proporciona una descripción del procedimiento que debe seguirse para realizar la prueba. Supóngase que se tienen $\left\{Z_j\right\}_{j=1}^N$, $n\in \{1,2,\hdots,N\}$, $c\in(0,\pi)$. Entonces se definen
\begin{eqnarray}
\alpha(n) \quad = \quad \sum_{j=1}^n Z_j\cos j c, & & \beta(n) \quad = \quad \sum_{j=1}^n Z_j \sen j c. \nonumber
\end{eqnarray}
Así, para $c\in(0,\pi)$, $n\in\{1,2,\hdots,N\}$, el desplazamiento cuadrático medio está dado por
\begin{eqnarray}
M_c(n) &=& \lim_{N\to \infty}\frac{1}{N}\sum_{k=1}^N (\alpha(k + n) - \alpha(k))^2 + (\beta(k + n) - \beta(k))^2. \nonumber
\end{eqnarray}

\begin{figure}
    \centering
    \includegraphics[width=160pt]{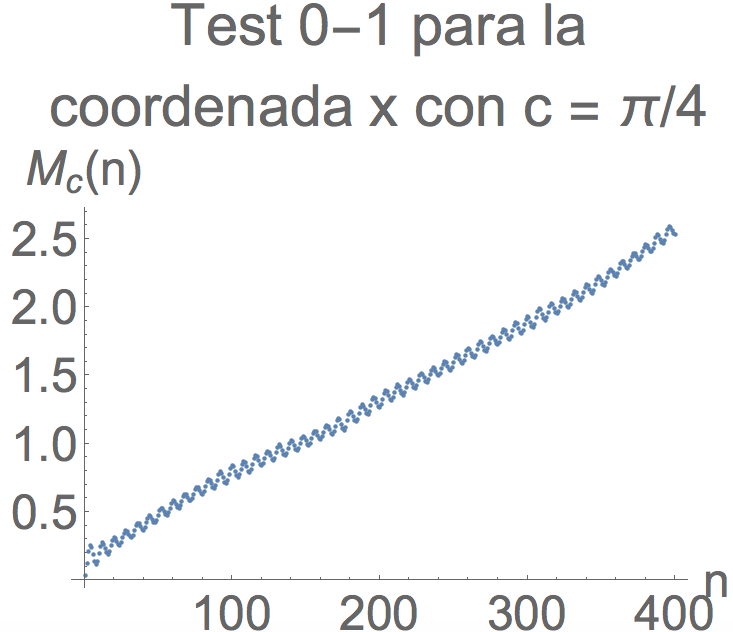}
    \includegraphics[width=160pt]{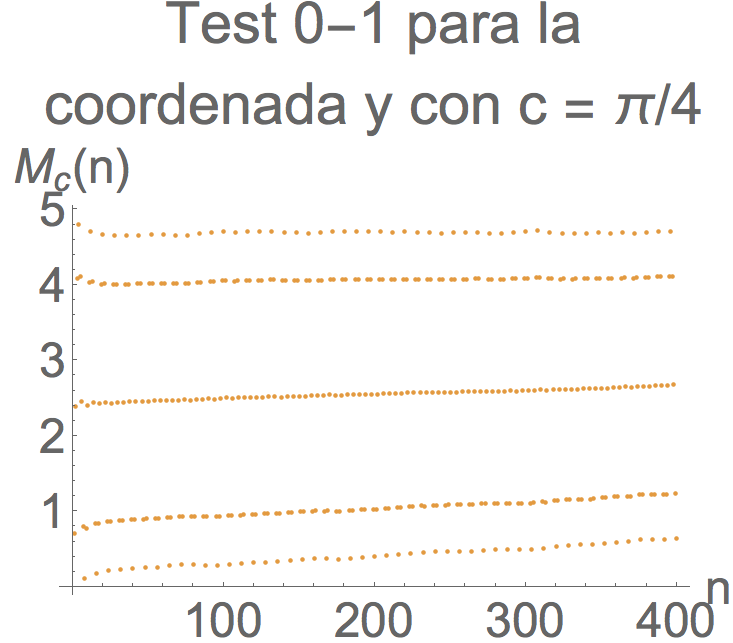}
    \includegraphics[width=160pt]{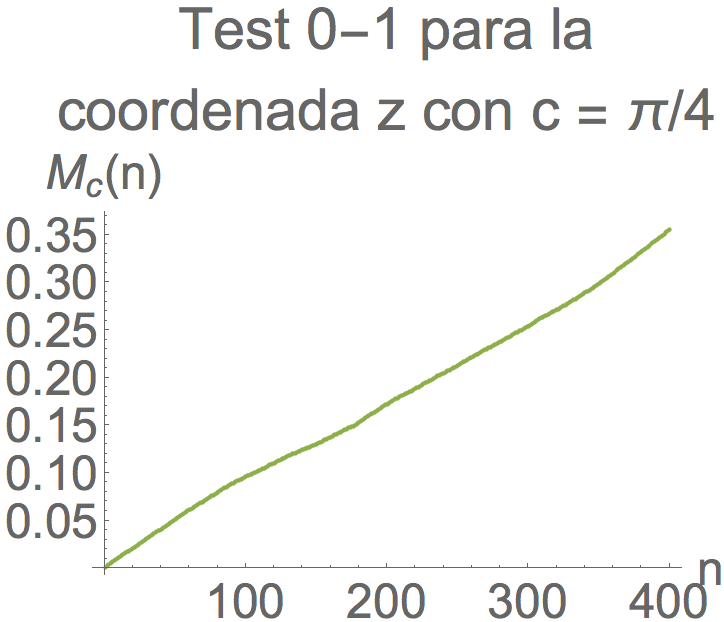}
    \includegraphics[width=160pt]{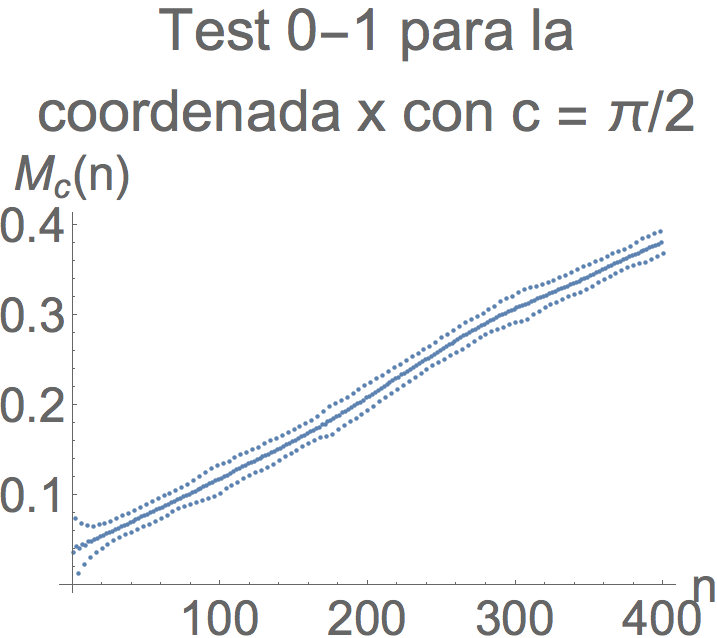}
    \includegraphics[width=160pt]{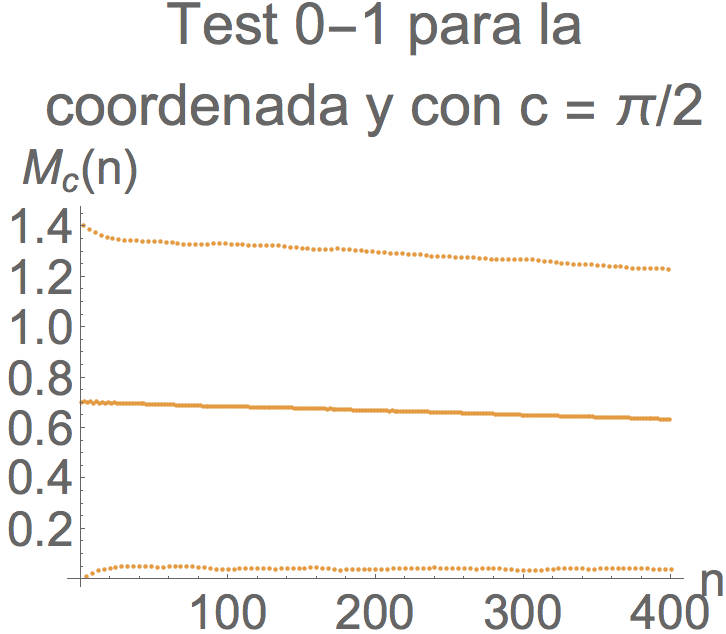}
    \includegraphics[width=160pt]{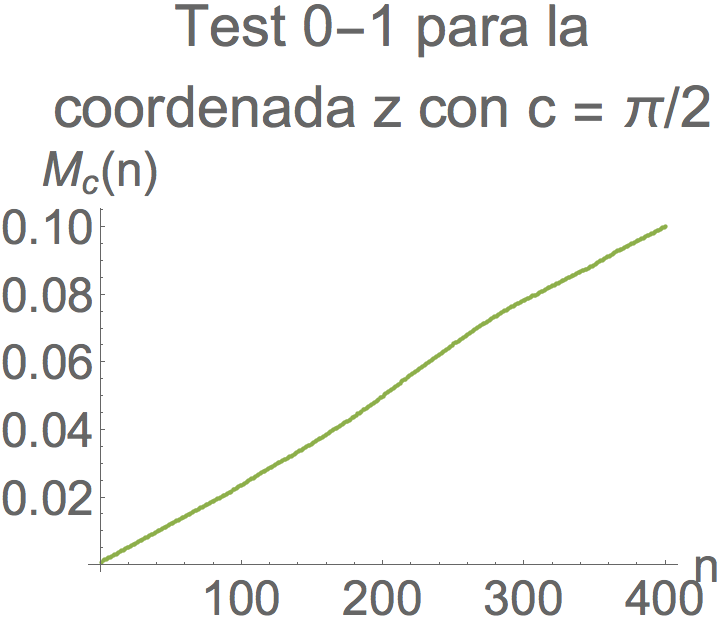}
    \includegraphics[width=160pt]{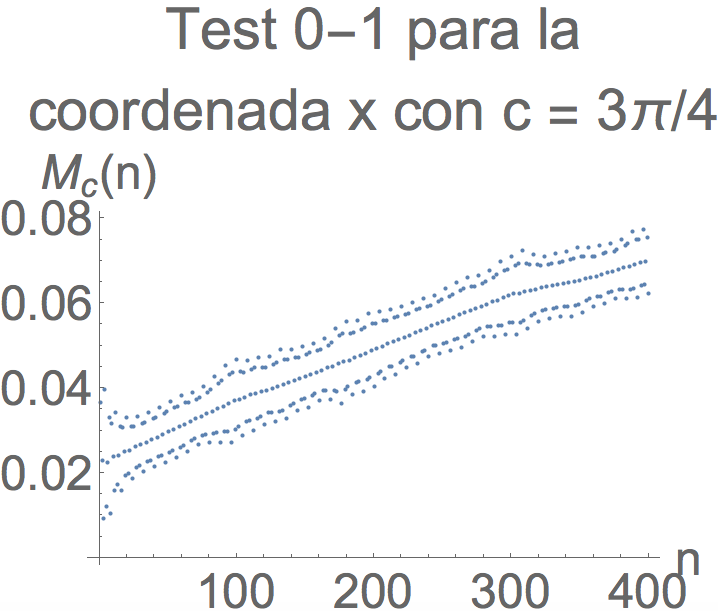}
    \includegraphics[width=160pt]{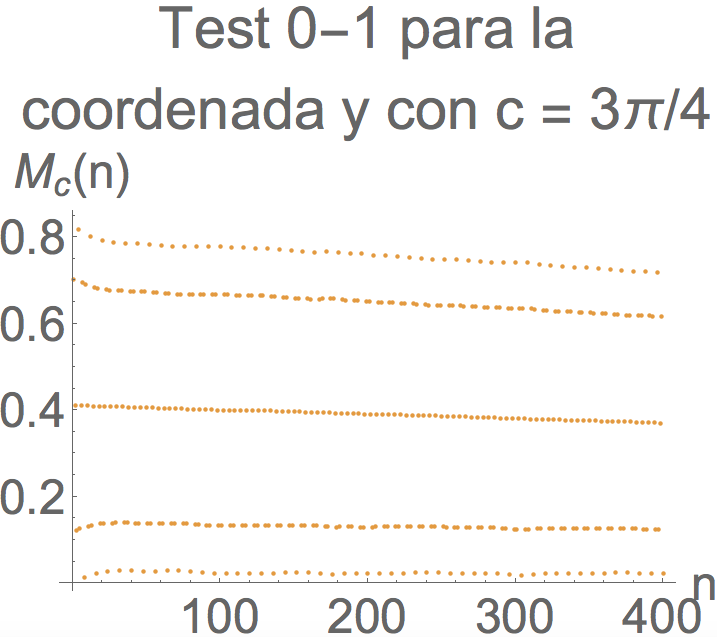}
    \includegraphics[width=160pt]{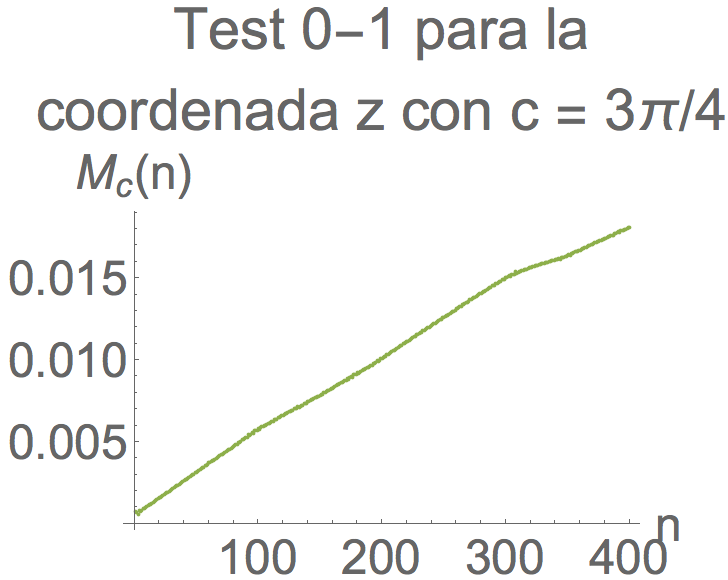}
    \caption{El Test 0-1 fue ejecutado con las tres coordenadas de \eqref{5.1} y los resultados son graficados. Las coordenadas $x$, $z$ parecen tener un crecimiento aproximadamente lineal, mientras que la gráfica de la coordenada $y$ aparentemente está acotada.}
    \label{fig5.13}
\end{figure}

En \cite{Gottwald2004a} se menciona que si $M_c(n)$ está acotada para $n=1,2,3,\hdots$, entonces se tiene una dinámica regular en el sistema original. Si por el contrario $M_c(n)$ no está acotado, entonces la dinámica es caótica.

En la Figura \ref{fig5.13} se muestran los bosquejos de las gráficas del Test 0-1 aplicado a cada una de las entradas de \eqref{5.1}.

\section{Conclusiones}
De acuerdo con la definición de Devaney \cite{Devaney1989AnSystems}, un sistema dinámico en el que hay presencia de caos se caracteriza por la impredecibilidad de las curvas integrales, por su transitividad topológica y por su gran periodicidad. Otra característica de los sistemas con caos, es la existencia de un conjunto extraño al que se aproximan las soluciones. Es bien sabido que el teorema de Poincaré-Bendixson implica que en los sistemas dinámicos planos a tiempo continuo no hay caos, sin embargo esto no ocurre en sistemas con tres o más ecuaciones acopladas, como los que aquí son estudiados. Entre los ejemplos conocidos de sistemas que exhiben caos se encuentran los sistemas de Lorenz \cite{Lorenz1963DeterministicFlow}, de Rössler \cite{Rossler1976AnChaos} y de Chen  \cite{Chen1999YetAttractor}.

Las características de \eqref{5.1} tales como los coeficientes de Lyapunov, la dimensión de correlación y la entropía de correlación, así como la dimensión fractal\index{dimensión fractal} del atractor extraño son analizadas siguiendo los procedimientos esbozados. El espectro potencia es descrito para cada una de las coordenadas de \eqref{5.1}, y posteriormente la prueba Test 0-1 es implementada. 

Algunas de las propiedades que tiene el sistema \eqref{5.1} son la existencia de un punto de equilibrio no hiperbólico, y la ausencia de puntos de equilibrio en la frontera $E_3 = (x_3,y_3,0)$ y $E_4=(0,y_4,z_4)$. 

Si este caso particular es analizado desde la mirada del ecólogo, es posible notar que aquí existe una predilección por parte de la especie depredadora\index{depredador} $z$ por alimentarse de la especie $x$: el valor de $q_1$ es 10 veces mayor que el valor de $q_2$.

\chapter{Conclusiones Generales}
A continuación se hace una última reflexión del trabajo. Todos los modelos matemáticos que aquí fueron presentados son simétricos: en la Parte \ref{partei}, los depredadores\index{depredador} pueden intercambiar posiciones sin afectar el resultado; y en la Parte \ref{parteii} las presas\index{presa} pueden intercambiar lugares sin que esto afecte el resultado. Todos los depredadores\index{depredador} son de tipo especialista, por lo que la ausencia de las presas\index{presa} los conduce hacia la extinción\index{extinción}. Las respuestas funcionales que se eligieron para trabajar fueron de tipo Holling II (que también se conoce como Michaelis-Menten), con excepción del Capítulo \ref{capitulo2}, donde se usa la respuesta funcional de tipo Lotka-Volterra (que resulta ser un caso límite para Holling II). El crecimiento de las presas\index{presa} está dado por el modelo logístico\index{modelo logístico} en los Capítulos \ref{capitulo3}, \ref{capitulo4} y \ref{capitulo5}, mientras que en el Capítulo \ref{capitulo2} se propone trabajar de forma general con una función diferenciable con abscisa al origen positiva y derivada negativa. La competencia entre los depredadores\index{depredador} es considerada en todos los casos, siendo representada con el parámetro $m$, y añadiendo subíndices cuando es conveniente. En la Parte \ref{parteii}, se considera también la competencia interespecífica de las presas\index{presa}, y esta última es denotada con $\alpha_{12}$ y $\alpha_{21}$. En todos los modelos propuestos a los largo de este trabajo, fue posible demostrar disipatividad y llevar a cabo un análisis general de la dinámica de los puntos de equilibrio en la frontera del espacio de estados. También se demuestra persistencia uniforme para el modelo que fue propuesto en el Capítulo \ref{capitulo2}.

La exploración hecha en el Capítulo \ref{capitulo2} busca dar una mayor profundidad al problema que fue propuesto en \cite{Kuang2003BiodiversityCompetition}. Se analizó la incorporación de mecanismos de competencia intra- e interespecífica en los depredadores\index{depredador} del modelo matemático usando una matriz de $n\times n$, donde $n$ es el número de depredadores\index{depredador} involucrados. Los resultados muestran que la única presa\index{presa} que aparece en el modelo cambia su valor de equilibrio (positivo) al variar los coeficientes de competencia de los depredadores\index{depredador}. Esto permite preguntarnos, en qué casos se verá beneficiada (o perjudicada) la especie recurso. La relevancia de la matriz de competencia radica en el hecho de que permite decidir la estabilidad asintótica global del punto de equilibrio positivo. En particular, cuando hay una presa\index{presa} y dos depredadores\index{depredador}, la relación que deben cumplir los parámetros de la matriz se puede expresar en términos de las medias aritmética y geométrica. Además, se muestra que existe un conjunto abierto y conexo en el espacio de matrices de $n\times n$ en el que los teoremas centrales del Capítulo \ref{capitulo2} se cumplen.

En el Capítulo \ref{capitulo3} se retoman las ideas anteriores y se presenta un modelo triespecífico (dos depredadores\index{depredador} y una presa\index{presa}) sin competencia interespecífica. El crecimiento de la presa está dado por el modelo logístico\index{modelo logístico}, y la tasa de captura por parte de los depredadores\index{depredador} es denso-dependiente que se representa con una respuesta funcional de tipo Holling II. A partir de la incorporación la respuesta funcional Holling II surge la posibilidad de que exista más de un punto de equilibrio positivo, y que se presenten casos con multiestabilidad. Entre los resultados de este capítulo, están las condiciones para la unicidad del punto de equilibrio positivo, así como las condiciones para la estabilidad asintótica. Además de lo anterior, fue posible hallar el umbral de una bifurcación de Hopf, y se muestra la existencia de umbrales de Bautin y de cero-Hopf.

En los siguientes capítulos el problema cambia al de un único depredador\index{depredador}. Así, en el Capítulo \ref{capitulo4} aparece un nuevo modelo matemático que en primera instancia es analizado usando las ideas expuestas en \cite{Abrams1999}, y se muestra que la condición para un arribo exitoso de una especie foránea (que se denota con $y$) es la misma que la condición necesaria para que el punto de equilibrio que excluye solamente a dicha especie sea inestable. En la Sección \ref{seccion4.3} se realizaron experimentos numéricos que muestran que la manipulación de los parámetros correspondientes a la capacidad de carga de la especie invasora (denotada con $K_2$) y al tiempo de manejo del depredador\index{depredador} con respecto a la especie invasora (denotada con $a_2$) pueden cambiar la tendencia de las trayectorias que representan a las soluciones del modelo. Desde el punto de vista ecológico, esto implicaría pasar de un escenario de coexistencia, a otro en el que alguna de las especies perezca y viceversa. En el primer experimento se muestra que una capacidad de carga suficientemente grande será vital para la supervivencia de $y$. Por otra parte el segundo experimento muestra que dos especies ($x$, $z$) que aparentemente no pueden coexistir en un escenario biespecífico, podrán alcanzar la coexistencia al incorporarse $y$ al medio. Más adelante, se realiza una construcción geométrica usando las ecuaciones de las curvas isoclinas para demostrar la existencia de puntos de equilibrio positivos. Si el sistema cuenta con un punto de equilibrio, este último se puede trasladar a una posición adecuada para así dar condiciones suficientes para la estabilidad asintótica local. Finalmente se muestra la existencia de puntos en el espacio de parámetros en los que yacen bifurcaciones de Hopf y de cero-Hopf, y se usa la teoría del promedio como se hizo en \cite{Llibre2015Zero-HopfSystem} para hallar ciclos límites que circundan el punto de equilibrio positivo.

En el Capítulo \ref{capitulo5} se expone una selección de parámetros en la que el sistema del capítulo anterior tiene trayectorias que experimentan una dinámica que puede ser clasificada como caótica. Se trata de un escenario en el que el producto de los coeficientes $\alpha_{12}$, $\alpha_{21}$ (que representan la competencia interespecífica de las presas\index{presa}) es mayor que el producto correspondiente a $r_1/K_1$, $r_2/K_2$ (que son los valores de competencia intraespecífica de las presas\index{presa}), lo que implica la ausencia del punto de equilibrio de la forma $(x_3,y_3,0)$, con $x_3$, $y_3 > 0$. Además de lo anterior, no existe un punto de la forma $(0,y_4,z_4)$ con $y_4$, $z_4$ positivos. Las propiedades de las curvas solución y del conjunto atractor son revisadas en este capítulo. En el caso estudiado en el Capítulo \ref{capitulo5}, vale la pena preguntarse si será posible preservar el caos si se manipulan los parámetros $q_i$ (tasa intrínseca de depredación) y $a_i$ (tiempo de manejo) cuidando que la proporción entre ellos no se vea afectada. En los modelos de la Parte \ref{partei} no fue posible hallar caos, por lo tanto otro trabajo futuro es investigar el efecto que tienen los mecanismos de competencia intra- e interespecífica para este propósito.

A lo largo de este trabajo, es notable el papel que juegan los parámetros de competencia intraespecífica de los depredadores\index{depredador} -- denotados con $m$, $m_{ii}$ -- en los modelos propuestos. Si el parámetro correspondiente es suficientemente grande, entonces es probable que se cumpla la existencia de un punto de equilibrio positivo asintóticamente estable. Sin embargo, al hacer demasiado grande el coeficiente de competencia intraespecífica, la coordenada correspondiente a la población de ese depredador\index{depredador} en el punto de equilibrio positivo se aproxima a cero. Esto se puede verificar a través de un análisis de las isoclinas y de la matriz jacobiana de los modelos.
\begin{appendices}
\chapter{El submodelo con una presa\index{presa} y su depredador\index{depredador}}
\label{apendicea}
Un caso que ha sido ampliamente estudiado es el modelo que introdujo Bazykin (véase por ejemplo \cite[p. 67 - 81]{Bazykin1998NonlinearPopulations}). Un análisis profundo de este sistema se hace en \cite{Bazykin1976}. Algunas de las bifurcaciones que el sistema presenta aparecen en \cite{Metzler1985}. En \cite{Zhou2010} los autores estudian un problema en un entorno espacialmente homogéneo. Con esta finalidad, introducen un sistema de ecuaciones con difusión que está basado en el modelo de Bazykin, se hallan soluciones positivas y se estudia el comportamiento asintótico de las mismas.

El modelo de Bazykin, que a continuación se presenta, está dado por
\begin{eqnarray}
\label{a.1}
\dot x &=& x\left(r(1 - x/K) - \frac{q y}{1 + a x}\right), \nonumber \\
\dot y &=&  y\left(\frac{c q x}{1 + a x} - \mu - m y\right).
\end{eqnarray}
Se puede apreciar que se trata de un submodelo de \eqref{3.1} que se obtiene de eliminar a uno de los depredadores\index{depredador} del entorno, o bien de un submodelo de \eqref{4.1} en el que se excluye a una de las presas\index{presa}. El modelo \eqref{a.1} permite estudiar a una presa\index{presa} y un depredor endémicos. Al explorar la dinámica, se puede conocer acerca de las curvas poblacionales de las especies nativas antes de la llegada de una especie exótica al ecosistema.

Considérese las curvas isoclinas de \eqref{a.1}
\begin{eqnarray}
f(x) &=& \frac{r}{q}(1 - x/K)(1 + ax), \nonumber \\
g(x) &=& \frac{1}{m}\left(\frac{c q x}{1 + a x} - \mu\right), \nonumber
\end{eqnarray}
y el u.e.e. de \eqref{a.1}
\begin{eqnarray}
y^b &=& \frac{\mu}{c q - a \mu}. \nonumber
\end{eqnarray}
Se asume que el u.e.e. satisface
\begin{eqnarray}
0 &<& y^b \quad<\quad K. \nonumber
\end{eqnarray}
Para la positividad de $y^b$, es necesario que se satisfaga
\begin{eqnarray}
cq - a\mu &>& 0, \nonumber
\end{eqnarray}
por lo que se asumirá que dicha desigualdad siempre se cumple.

Los siguientes resultados hablan del número de puntos de equilibrio en \eqref{a.1}.
\begin{prop}
\label{propa.1}
Si se cumple $0 < y^b < K$, entonces \eqref{a.1} tiene puntos de equilibrio positivo. Más aún, si se satisface
\begin{eqnarray}
\frac{K - 1/a}{2} \quad \leq \quad y^b \quad < \quad K, \nonumber
\end{eqnarray}
entonces \eqref{a.1} tiene un único punto de equilibrio positivo $E^* = (x^*,y^*)$ que es atractor local. Si adicionalmente se tiene
\begin{eqnarray}
r &<& \mu, \nonumber \\
q &<& \frac{2 - aK(r - \mu)}{c}, \nonumber
\end{eqnarray}
entonces $E^*$ es atractor global.
\end{prop}
\begin{proof}
Se verifica fácilmente que $f$ y $g$ son positivas exactamente en el intervalo $(y^b,K)$. Como $f$ es una función decreciente en $[(K - 1/a)/2, K]$, y $g$ es creciente en $[y^b,K]$, entonces existe un punto $x^*\in (y^b,K)$ que satisface
\begin{eqnarray}
f(x^*) &=& g(x^*). \nonumber
\end{eqnarray}
Se sigue que el punto $E^* = (x^*,f(x^*))$ es un punto de equilibrio positivo para \eqref{a.1}. Además, si $y^b > (K - 1/a)/2$, entonces ese punto es único. La matriz jacobiana del sistema \eqref{a.1} en $E^*$ está dada por
\begin{eqnarray}
J(E^*) &=& \frac{1}{1 + ax^*}\left( \begin{array}{cc}
r K^{-1}(aK - 2ax^* -1)x^* & - qx^* \\\\
cr(1 - x^*/K) & -m y^*(1 + ax^*)
\end{array} \right). \nonumber
\end{eqnarray}
Es fácil verificar que si $x^* > (K - 1/a)/2$, entonces 
$rK^{-1}(aK - 2ax^* - 1)x^* < 0$, y así $\det(J(E^*)) > 0$ y $\text{tr}(J(E^*)) < 0$, y por lo tanto $E^*$ es un atractor local en \eqref{a.1}. Finalmente, se puede apreciar que
\begin{eqnarray}
\nabla \cdot (\dot x, \dot y) &=& r\left(1 - \frac{2x}{K}\right) - \frac{q y}{(1 + ax)^2} + \frac{cqx}{1 + ax} - \mu - 2my \nonumber \\
&<& \frac{cqx}{1 + ax} + r - \frac{2x}{K} - \mu \nonumber \\
&=& \frac{-2ax^2 + (aK(r - \mu) + cq - 2 )x + K(r - \mu)}{K(1 + ax)} \nonumber \\
&<& 0, \nonumber
\end{eqnarray}
si las desigualdades $r < \mu$,
$q < c^{-1}(2 - aK(r - \mu))$ se cumplen. Del teorema de Bendixson-Dulac \cite[Sección 3.9, Teorema 1]{Perko2001}, se sigue que  \eqref{a.1} no tiene órbitas periódicas no triviales en el retrato fase, y por lo tanto $E^*$ es un atractor global.
\end{proof}
En el siguiente resultado, se proveen condiciones necesarias y suficientes para la existencia de tres puntos de equilibrio, lo que permite una bifurcación conocida como tridente.
\begin{teo}
\label{teoa.2}
 Supóngase que los parámetros satisfacen la condición
 \begin{eqnarray}
 cq - a\mu &>& 0. \nonumber
 \end{eqnarray}
 Entonces el sistema \eqref{a.1} tiene tres puntos de equilibrio positivos si y solo si se satisfacen las siguientes condiciones
 \begin{eqnarray}
 K &>& \frac{8 c q + a \mu}{a c q - a^2 \mu}, \nonumber
  \end{eqnarray}
  y $m\in (m_1 - \varepsilon, m_1 + \varepsilon) \cap [m_2,\infty)$, donde
   \begin{eqnarray}
   m_1 &=& (8c(1 + aK)^3r)^{-1}(-2 a^3 c \mu K^3 q - 22 a^2 c \mu K^2 q \nonumber \\
   && +a^2 \mu^2 K (a K + 1)^2 + c^2 K q^2 (a K (a K + 20) - 8) - 20 a c \mu K q), \nonumber
   \end{eqnarray}
   \begin{eqnarray}
   \varepsilon &=& \frac{1}{8} \sqrt{\frac{a K^2 (a \mu K - c K q + \mu) (c q (8 - aK) 
   + a \mu (aK + 1))^3}{c^2 r^2 (a K + 1)^6}}, \nonumber
   \end{eqnarray}
   \begin{eqnarray}
   m_2 &=& (aKr(aK + 2) + r)^{-1}(3Kq(cq - a\mu ) \nonumber
   \end{eqnarray}
\end{teo}
\begin{proof}
Al igualar $f(x)$ y $g(x)$, se obtiene el polinomio cúbico
\begin{eqnarray}
\label{2}
P(x) &=& \Omega_3x^3 + \Omega_2x^2 + \Omega_1x + \Omega_0,                     \nonumber
\end{eqnarray}
con
\begin{eqnarray}
\Omega_3 &=& -a^2 m r, \nonumber \\
\Omega_2 &=& -2 a m r + a^2 K m r, \nonumber \\
\Omega_1 &=& a \mu K q - c K q^2 - m r + 2 a K m r,      \nonumber \\
\Omega_0 &=& \mu K q + K m r. \nonumber               
\end{eqnarray}
Si se calcula la derivada de $P$, y se hace $P'(x) = 0$, entonces se tiene que $P$ alcanza un mínimo local y un máximo local en
\begin{eqnarray}
x_1 \quad = \quad \frac{-\Omega_2 - \sqrt{\Omega_2^2 - 3\Omega_1\Omega_3}}{3\Omega_3}, &&
x_2 \quad = \quad \frac{-\Omega_2 + \sqrt{\Omega_2^2 - 3\Omega_1\Omega_3}}{3\Omega_3}. \nonumber
\end{eqnarray}

Una condición necesaria y suficiente para que el argumento de la raíz cuadrada en las expresiones anteriores no sea negativo ($\Omega_2^2 - 3\Omega_1\Omega_3 \geq 0$) es la siguiente
\begin{eqnarray}
m &\geq& \frac{3 K q(c q - a \mu)}{a K r(a K + 2) + r}. \nonumber
\end{eqnarray}

Si $m \geq (aKr(aK + 2) + r)^{-1}(3Kq(cq - a\mu))$, entonces $x_1$ y $x_2$ están bien definidos. Como $\Omega_0 > 0$ y $\Omega_3 < 0$, entonces se tiene que $P$ tiene tres raíces positivas si y solo si se cumple
\begin{eqnarray}
P(x_1) \quad < \quad 0 &<& P(x_2). \nonumber
\end{eqnarray}

Las desigualdades anteriores se verifican si y solo si
\begin{eqnarray}
 K &>& \frac{8 c q + a \mu}{a c q - a^2 \mu}, \nonumber \\ \nonumber \\
 m_1 - \varepsilon &<& m \quad < \quad  m_1 + \varepsilon. \nonumber
\end{eqnarray}
\end{proof}

En la demostración del Teorema \ref{teoa.2}, se puede apreciar que
\begin{eqnarray}
\Omega_3 \quad < \quad 0 &<& \Omega_0. \nonumber
\end{eqnarray}
Por la Regla de Descartes\index{Regla de Descartes}, se tiene que $P$ tiene al menos una raíz positiva.

Considérese ahora el sistema \eqref{a.1} en el caso límite cuando $m$ es igual a $0$. Este caso fue estudiado por Kuznetsov en \cite[Ejemplo 3.1]{Kuznetsov1998}.
\begin{eqnarray}
\label{a.2}
\dot x &=& x\left(r(1 - x/K) - \frac{q y}{1 + a x}\right),  \\
\dot y &=& y\left(\frac{c q x}{1 + a x} - \mu\right). \nonumber
\end{eqnarray}
El sistema \eqref{a.2} solamente exhibe un punto de equilibrio positivo ubicado, el cual está ubicado en
\begin{eqnarray}
E^* &=& \left(y^b,\frac{r}{q}\left(1 - \frac{y^b}{K}\right)\left(1 + ay^b\right)\right). \nonumber
\end{eqnarray}
\begin{figure}[h]
    \centering
    \includegraphics[width=320pt]{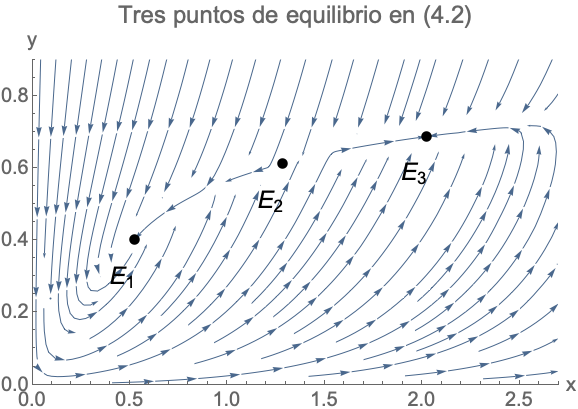}
    \caption{Las soluciones del sistema \eqref{a.1} con parámetros $r = 0.3$, $K = 4.5$, $q = 1.7$, $a = 3$, $c = 2$, $\mu = 0.3$, $m = 0.98$. Se presenta un escenario con multiestabilidad y tres puntos de equilibrio positivos ubicados en $E_1 = (0.526081,0.401793))$, $E_2 = (1.28571,0.612245)$, $E_3 = (2.02154,0.686642)$. Del Teorema \ref{teoa.2} se sigue que si los primeros seis parámetros son así seleccionados, y $m \in (0.96666,1.01401)$, entonces \eqref{a.1} presenta tres equilibrios positivos.}
\end{figure}

Lo anterior es consistente con el Teorema \ref{teoa.2}, ya que $m = 0$, $m_2 > 0$ y entonces $m \notin [m_2,\infty)$. Por otro lado la matriz jacobiana de \eqref{a.2} en $E^*$ tiene la forma
\begin{eqnarray}
J(E^*) &=& \frac{1}{1 + ay^b} \quad \widetilde J(E^*), \nonumber
\end{eqnarray}
donde
\begin{eqnarray}
\widetilde J(E^*) &=& \left( \begin{array}{cc}
r K^{-1}(aK - 2ay^b -1)y^b & - qy^b \\\\
cr(1 - y^b/K) & 0 
\end{array} \right). \nonumber
\end{eqnarray}
\begin{figure}[h]
    \centering
    \includegraphics[width=242pt]{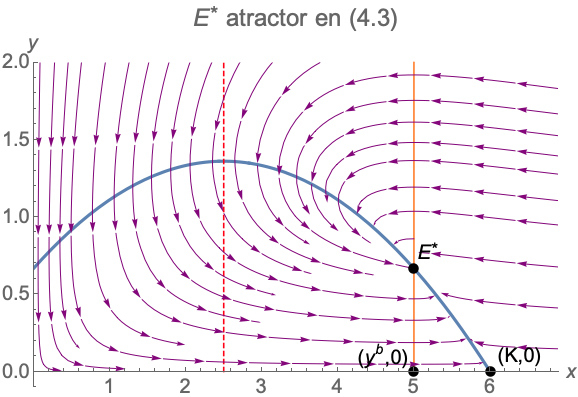}
    \includegraphics[width=242pt]{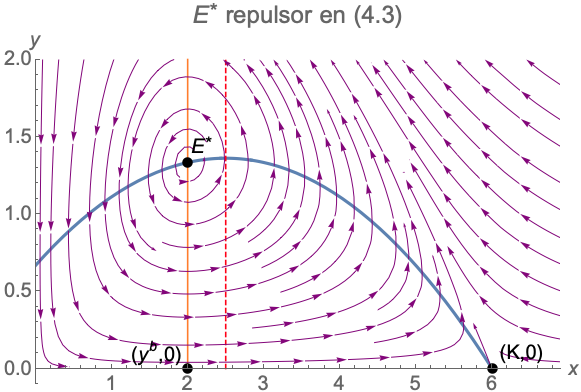}
    \caption{La isoclina correspondiente a $\dot x$ está bosquejada en color azul, la isoclina de $\dot y$ es una recta en amarillo y la línea roja discontinua muestra la latitud del vértice de la curva azul. Ambas intersectan en el punto de equilibrio positivo $E^*$. El punto de equilibrio positivo $E^*$ existe si y solo si $0 < y^b < K$. La gráfica de la izquierda muestra un caso en el que $E^*$ es asintóticamente estable ($E^*$ se ubica a la derecha del vértice de la parábola), mientras que en la imagen derecha $E^*$ es inestable ($E^*$ está a la derecha del vértice).}
    \label{figa.2}
\end{figure}

El polinomio característico de $\widetilde J(E^*)$ está dado por
\begin{eqnarray}
\lambda^2 - \text{tr}(\widetilde J(E^*)) \lambda + \det(\widetilde J(E^*)), \nonumber
\end{eqnarray}
donde
\begin{eqnarray}
\text{tr}(\widetilde J(E^*)) &=& \frac{r}{K}(aK - 2ay^b -1)y^b, \nonumber \\
\det(\widetilde J(E^*)) &=& cqry^b(1 - y^b/K) \quad > \quad 0. \nonumber
\end{eqnarray}

\begin{teo}
\label{teoa.3}
Si $\displaystyle y^b = \frac{K -1/a}{2}$, entonces el sistema \eqref{a.2} tiene un único punto de equilibrio positivo $E^*$ cuyos valores propios son números imaginarios puros conjugados.
\end{teo}
\begin{proof}
Los valores propios de $\widetilde J(E^*)$ son
\begin{eqnarray}
\lambda_{1,2} &=& \frac{\text{tr}(\widetilde J(E^*)) \pm \sqrt{\text{tr}(\widetilde J(E^*))^2 - 4\det(\widetilde J(E^*))}}{2}. \nonumber
\end{eqnarray}

Si $\displaystyle y^b = \frac{K - 1/a}{2}$, entonces
\begin{eqnarray}
\text{tr}(\widetilde J(E^*)) \quad = \quad \frac{r}{K}\left(aK - 2a\frac{(K - 1/a)}{2} - 1\right)\left(\frac{K - 1/a}{2}\right) &=& 0. \nonumber
\end{eqnarray}
De lo anterior, se sigue que
\begin{eqnarray}
\lambda_{1,2} &=& \pm i\sqrt{\det(\widetilde J(E^*))}. \nonumber
\end{eqnarray}
\end{proof}

La existencia de una bifurcación de Hopf cuyo umbral se encuentra ubicado precisamente en el punto $y^b = (K - 1/a)/2$ se puede mostrar a partir del Teorema \ref{teoa.3}. 

Desde el punto de vista geométrico, se puede apreciar que si la recta de la isoclina correspondiente a $\dot y$ se encuentra ubicada del lado derecho del vértice de la parábola de la isoclina de $\dot x$, entonces el punto de equilibrio $E^*$ es localmente asintóticamente estable. Análogamente, si la recta cae del lado izquierdo de la parábola, entonces $E^*$ es inestable. Lo anterior se ilustra en la Figura \ref{figa.2}.
\chapter{Análisis geométrico de las isoclinas}
\label{apendiceb}
En la Sección \ref{seccion4.4}, dos de las superficies isoclinas ($f$, $g$) del sistema \eqref{4.1} se igualan para obtener una ecuación que corresponde a una hipérbola en $\mathbb R^2$:
\begin{eqnarray}
\label{b.1}
\frac{1}{q_1}\left(r_1(1 - x/K_1) - \alpha_{12}y\right)(1 + a_1 x) &=& \frac{1}{q_2}\left(r_2(1 - y/K_2) - \alpha_{21}x\right)(1 + a_2 y),
\end{eqnarray}
donde $r_1$, $q_1$, $a_1$, $\alpha_{12}$, $r_2$, $q_2$, $a_2$, $\alpha_{21}$ son valores positivos fijos. Supóngase que $r_1/K_2 > \alpha_{12}$, $r_2/K_1 > \alpha_{21}$ y supóngase (sin pérdida de generalidad) que $r_1/q_1 > r_2/q_2$. Entonces existen tres casos genéricos distintos para la estructura de la hipérbola:
\begin{enumerate}
    \item La hipérbola \eqref{b.1} intersecta al semieje positivo $X$ en un punto $P_1$, y no intersecta al semieje positivo $Y$. 
    \item La hipérbola \eqref{b.1} intersecta al semieje positivo $X$ en un punto $P_1$, intersecta al semieje positivo $Y$ en dos puntos ($P_0$ y $P_2$), y los puntos $P_0$ y $P_2$ están en la misma rama de la hipérbola.
    \item La hipérbola \eqref{b.1} intersecta al semieje positivo $X$ en un punto $P_1$, intersecta al semieje positivo $Y$ en dos puntos ($P_0$ y $P_2$), y los puntos $P_0$ y $P_1$ están en la misma rama de la hipérbola.
\end{enumerate}
Para hallar el punto $P_1$ que es la intersección de la hipérbola \eqref{b.1} con el semieje positivo $X$, se procede a resolver la ecuación 
\begin{eqnarray}
\label{b.2}
\frac{r_1}{q_1}(1 - x/K_1)(1 + a_1x)&=& \frac{1}{q_2}(r_2 - \alpha_{21}x).
\end{eqnarray}
La solución positiva $x = x_H$ de \eqref{b.2} es la abscisa del punto $P_1$.

De forma análoga, es posible conocer las coordenadas los puntos $P_0$ y $P_2$ donde se intersectan la hipérbola \eqref{b.1} y el semieje positivo $Y$. Se procede a resolver la ecuación
\begin{eqnarray}
\label{b.3}
\frac{r_2}{q_2}(1 - y/K_2)(1 + a_2y) &=& \frac{1}{q_1}(r_1 - \alpha_{12}y).
\end{eqnarray}
Las soluciones reales $y=\widetilde y_H $, $y = y_H$ de \eqref{b.3} son las ordenadas de los puntos $P_0$ y $P_2$ respectivamente.

\begin{figure}[h]
        \centering
        \subfigure[Hipérbola sin intersección en el semieje $Y$.]{
        \includegraphics[width=133pt]{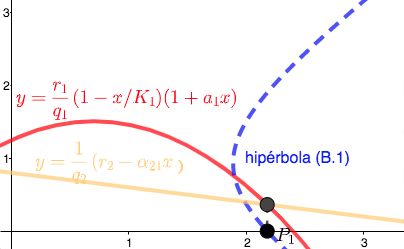}
        \includegraphics[width=133pt]{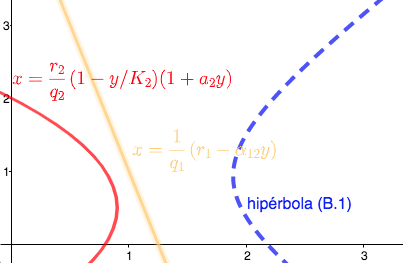}}
        \subfigure[Hipérbola con intersecciones en el semieje $Y$. $P_0$ y $P_1$ pertenecen a la misma rama.]{
        \includegraphics[width=133pt]{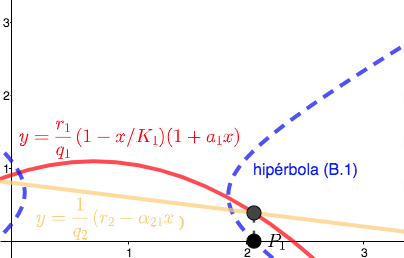}
        \includegraphics[width=133pt]{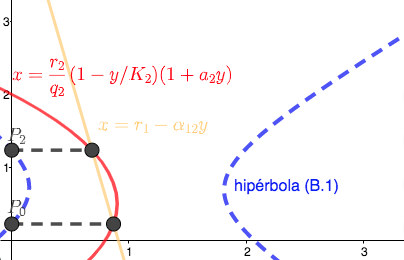}}
        \subfigure[Hipérbola con intersecciones en el semieje $Y$. $P_0$ y $P_2$ pertenecen a la misma rama.]{
        \includegraphics[width=133pt]{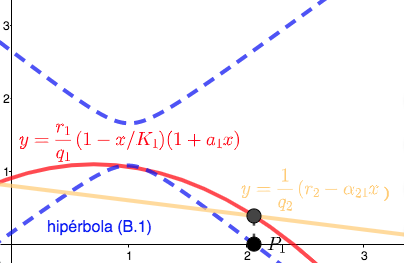}
        \includegraphics[width=133pt]{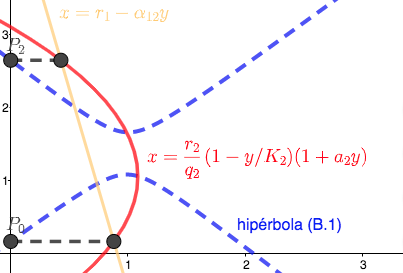}}
        \caption{La hipérbola \eqref{b.1} que en los diagramas aparece en color azul puede estar estructurada de tres maneras diferentes. Las curvas rojas y amarillas se utilizan para hallar la abscisa y la ordenada de las intersecciones de la hipérbola con los ejes coordenados.}
        \label{figb1}
    \end{figure}

En la Figura \ref{figb1} se ilustra en el plano cartesiano cada uno de los casos expuestos.

\chapter[Condiciones de genericidad y transversalidad]{Condiciones de genericidad y transversalidad para la bifurcación de Hopf}
\label{apendicec}
\section*{Traslado del punto de equilibrio y\\ del umbral de bifurcación}
Se procede a calcular el valor del coeficiente principal de Lyapunov (denotado con $\ell_1$) y la primera derivada parcial de la parte real de los valores propios complejos con respecto al parámetro $m$ para el sistema \eqref{4.15} en el punto de equilibrio $E^* = (1/4,1/4,1)$. El sistema está dado por
\begin{eqnarray}
\dot x &=& x\left(r_1(1 - x) - \frac{15r_1 z}{16(1 + x)}\right),                                    \nonumber \\
\dot y &=& y\left(r_2(1 - y/K_2) - \frac{5r_2(4K_2 - 1)z}{16K_2(1+y)}\right),                       \nonumber \\
\dot z &=& z\left(\frac{15c_1r_1 x}{16(1 + x)} + \frac{5r_2(4K_2 - 1)y}{16K_2(1+y)} -
\left(\frac{3c_1r_1}{16} + \frac{r_2(4K_2 - 1)}{16K_2} - m\right) - m z\right).           \nonumber
\end{eqnarray}
Un reescalamiento del tiempo (véase por ejemplo \cite[Subsección 1.4.1]{Chicone1999OrdinaryApplications}) permite reescribir el sistema en forma polinomial,
\begin{eqnarray}
\dot x &=& x\left(16(1 + x)\right)\left(16K_2(1 + y)\right)\left(r_1(1 - x) - \frac{15r_1 z}{16(1 + x)}\right)                                    \nonumber \\
       &=& 256 K_2 r_1 x+256 K_2 r_1 x y-240 K_2 r_1 x z \nonumber \\
       & & -256 K_2 r_1 x^3-240 K_2 r_1 x y z-256 K_2 r_1 x^3 y,     \nonumber \\ \nonumber \\
\dot y &=& y\left(16(1 + x)\right)\left(16K_2(1 + y)\right)\left(r_2(1 - y/K_2) - \frac{5r_2(4K_2 - 1)z}{16K_2(1+y)}\right)                              \nonumber \\
       &=& 256 K_2 r_2 y+256 K_2 r_2 x y+256 (K_2-1) r_2 y^2+80 r_2(1-4K_2)y z \nonumber \\
       & & +256 (K_2-1) r_2 x y^2+80 r_2 (1-4K_2)x y z-256 r_2 y^3-256 r_2 x y^3, \nonumber \\ \nonumber \\
\dot z &=& z\left(16(1 + x)\right)\left(16K_2(1 + y)\right)\times \nonumber \\ & & \left(\frac{15c_1r_1 x}{16(1 + x)} + \frac{5r_2(4K_2 - 1)y}{16K_2(1+y)} -
\left(\frac{3c_1r_1}{16} + \frac{r_2(4K_2 - 1)}{16K_2} - m\right) - m z\right) \nonumber \\ \nonumber \\
      &=& 16(K_2(-3c_1r_1 + 16m -4r_2)+ r_2) z+16(K_2(12c_1r_1 + 16m - 4r_2)+ r_2) x z \nonumber \\
      & & + 16(K_2(16r_2 + 16m -3 c_1r_1)-4r_2) y z +16(K_2(12c_1r_1 + 16m + 16r_2)-4r_2) x y z \nonumber \\
      & & -256 K_2 m x z^2-256 K_2 m y z^2-256 K_2 m z^2-256 K_2 m x y z^2,
           \nonumber
\end{eqnarray}
y se traslada el punto de equilibrio $E^* = (1/4,1/4,1)$ al origen de coordenadas $(0,0,0)$,
\begin{eqnarray}
\dot x &=& -40 K_2 r_1 x-75 K_2 r_1 z-240 K_2 r_1 x^2-32 K_2 r_1 x y-300 K_2 r_1 x z \nonumber \\
       & & -60 K_2 r_1 y z-320 K_2 r_1 x^3-192 K_2 r_1 x^2 y-256 K_2 r_1 x^3 y-240 K_2 r_1 x y z, \nonumber \\ \nonumber \\
\dot y &=& \left(80 K_2 r_2-120 r_2\right)y + \left(25 r_2-100 K_2 r_2\right)z + \left(64 K_2 r_2-96 r_2\right)xy + \left(20 r_2-80 K_2 r_2\right)xz \nonumber \\
       & &+ \left(100 r_2-400 K_2 r_2\right)yz + \left(320 K_2 r_2-560 r_2\right)y^2 + \left(256 K_2 r_2-448 r_2\right)x y^2 \nonumber \\
       & & + \left(80 r_2-320 K_2 r_2\right)xyz -320 r_2 y^3 - 256 r_2 x y^3, \nonumber \\ \nonumber \\ 
\dot z &=& 240 c_1 K_2 r_1 x + \left(320 K_2 r_2-80 r_2\right)y - 400 K_2 m z + \left(192 c_1 K_2 r_1+256K_2 r_2-64 r_2\right)xy \nonumber \\
       & & + \left(240 c_1 K_2 r_1-320 K_2 m\right)xz + \left(-320 K_2 m+320 K_2 r_2-80 r_2\right)yz - 400 K_2 m z^2 \nonumber \\
       & &+ \left(192 c_1 K_2 r_1-256 K_2 m+256 K_2 r_2-64 r_2\right)xyz \nonumber \\
       & & - 320 K_2 m x z^2 - 320 K_2 m y z^2 - 256 K_2 m x y z^2. \nonumber
\end{eqnarray}
Finalmente se traslada el umbral de bifurcación $m_0$ (véase el Teorema \ref{teo4.20}) a cero para obtener
\begin{eqnarray}
\label{c.1}
\dot x &=& F(x,y,z), \nonumber \\
\dot y &=& G(x,y,z),           \\
\dot z &=& H(x,y,z), \nonumber
\end{eqnarray}
donde
\begin{eqnarray}
F(x,y,z) &=& -40 K_2 r_1 x-75 K_2 r_1 z-240 K_2 r_1 x^2-32 K_2 r_1 x y-300 K_2 r_1 x z \nonumber \\
       & & -60 K_2 r_1 y z-320 K_2 r_1 x^3-192 K_2 r_1 x^2 y-256 K_2 r_1 x^3 y-240 K_2 r_1 x y z, \nonumber \\ \nonumber \\
G(x,y,z) &=&\left(80 K_2 r_2-120 r_2\right)y + \left(25 r_2-100 K_2 r_2\right)z + \left(64 K_2 r_2-96 r_2\right)xy \nonumber \\
       & & + \left(20 r_2-80 K_2 r_2\right)xz + \left(100 r_2-400 K_2 r_2\right)yz \nonumber \\
       & & + \left(320 K_2 r_2-560 r_2\right)y^2 + \left(256 K_2 r_2-448 r_2\right)x y^2 \nonumber \\
       & & + \left(80 r_2-320 K_2 r_2\right)xyz -320 r_2 y^3 - 256 r_2 x y^3, \nonumber \\ \nonumber \\ 
H(x,y,z) &=& 240 c_1 K_2 r_1 x + \left(320 K_2 r_2-80 r_2\right)y - 400 K_2 (m + m_0)z \nonumber \\
         & & + \left(192 c_1 K_2 r_1+256K_2 r_2-64 r_2\right)xy + \left(240 c_1 K_2 r_1-320 K_2 (m + m_0) \right)xz \nonumber \\ 
         & & + \left(-320 K_2 (m + m_0) + 320 K_2 r_2-80 r_2\right)yz - 400 K_2 (m + m_0) z^2 \nonumber \\
       & &+ \left(192 c_1 K_2 r_1-256 K_2 (m + m_0)+256 K_2 r_2-64 r_2\right)xyz \nonumber \\
       & & - 320 K_2 (m + m_0) x z^2 - 320 K_2 (m + m _0) y z^2 - 256 K_2 (m + m_0) x y z^2. \nonumber
\end{eqnarray}
y donde
\begin{eqnarray}
m_0 &=& (80K_2(K_2(r_1 - 2r_2) + 3r_2))^{-1}((r_1^4K_2^4(4 - 45c_1)^2 + 2r_1^2r_2^2K_2^2(45c_1(41 - 88K_2 + 96K_2^2) \nonumber \\ \nonumber \\
& & + 32K_2(1 + 8K_2) - 124) + r_2^4(31 - 8K_2(1 + 8K_2))^2)^{1/2} + r_2^2(-96K_2^2 + 88K_2  - 41) \nonumber \\ \nonumber \\
& & + 8r_1r_2K_2(2K_2 - 3) - r_1^2K_2^2(4 + 45c_1)).                              \nonumber
\end{eqnarray}

La matriz jacobiana de \eqref{c.1} evaluada en $(0,0,0)$ está dada por
\begin{eqnarray}
J(0,0,0) &=& \begin{pmatrix} 
F_x & F_y & F_z \\
G_x & G_y & G_z \\
H_x & H_y & H_z   
\end{pmatrix}, \nonumber
\end{eqnarray}
donde
\begin{multicols}{3}
\begin{eqnarray}
F_x &=& -40r_1K_2, \nonumber \\
G_x &=& 0, \nonumber \\
H_x &=& 240 c_1 K_2 r_1, \nonumber
\end{eqnarray}
\begin{eqnarray}
& & \nonumber \\
F_y &=& 0, \nonumber \\
G_y &=& 40r_2(2 K_2-3), \nonumber \\
H_y &=& 80r_2(4 K_2-1), \nonumber
\end{eqnarray}
\begin{eqnarray}
& & \nonumber \\
F_z &=& -75r_1K_2, \nonumber \\
G_z &=& 25r_2(1-4 K_2), \nonumber \\
H_z &=& -400K_2(m + m_0), \nonumber
\end{eqnarray}
\end{multicols}
y tiene valores propios $\lambda_1$, $\lambda_2$, $\lambda_3$ dados por
\begin{eqnarray}
\lambda_1 &=& a - 2\varepsilon,   \nonumber \\
\lambda_2 &=& \varepsilon + i\omega, \nonumber \\
\lambda_3 &=& \varepsilon - i\omega, \nonumber
\end{eqnarray}
donde
\begin{eqnarray}
\varepsilon &=& \frac{a}{3} + \frac{\sqrt[3]{- 2a^3 - 9ab - 27c + 3\sqrt{3}\sqrt{4a^3c - a^2b^2 + 18abc - 4b^3 + 27c^2}}}{6\sqrt[3]{2}}\nonumber \\ \nonumber \\
& & + \frac{a^2 + 3b}{3 \times 2^{2/3}\sqrt[3]{- 2a^3 - 9ab - 27c + 3\sqrt{3}\sqrt{4a^3c - a^2b^2 + 18abc - 4b^3 + 27c^2}}}, \nonumber \\ \nonumber \\ \nonumber \\
\omega &=& - \frac{\sqrt{3}\sqrt[3]{- 2a^3 - 9ab - 27c + 3\sqrt{3}\sqrt{4a^3c - a^2 b^2 + 18abc - 4b^3 + 27c^2}}}{6 \sqrt[3]{2}} \nonumber \\ \nonumber \\
& & + \frac{\sqrt 3 (a^2 + 3b)}{3 \times 2^{2/3} \sqrt[3]{- 2a^3 - 9ab - 27c + 3 \sqrt{3} \sqrt{4a^3c - a^2b^2 + 18abc - 4b^3 + 27c^2}}},\nonumber
\end{eqnarray}
y donde
\begin{eqnarray}
a &=& -400 K_2(m + m_0) - 40K_2r_1 - 120r_2 + 80K_2r_2, \nonumber \\ \nonumber \\
b &=& -18000 c_1 K_2^2 r_1^2-16000 K_2^2 \left(m+m_0\right) r_1+32000 K_2^2 \left(m+m_0\right) r_2 \nonumber \\
& & -48000 K_2 \left(m+m_0\right) r_2+3200 K_2^2 r_1 r_2-4800 K_2 r_1 r_2 \nonumber \\
& &
-320 K_2 r_2 \left(100 K_2 r_2-25 r_2\right) +80 r_2 \left(100 K_2 r_2-25 r_2\right), \nonumber \\ \nonumber \\
c &=& 1440000 c_1 K_2^3 r_1^2 r_2-2160000 c_1 K_2^2 r_1^2 r_2+1280000 K_2^3 \left(m+m_0\right) r_1 r_2 \nonumber \\
& & -1920000 K_2^2 \left(m+m_0\right) r_1 r_2+3200 K_2 r_1 r_2 \left(100 K_2 r_2-25 r_2\right) \nonumber \\
& & -12800 K_2^2 r_1 r_2 \left(100 K_2 r_2-25 r_2\right). \nonumber
\end{eqnarray}
Haciendo uso de los vectores propios de la matriz jacobiana $J(0,0,0)$ se obtiene la matriz cambio de base para $J(0,0,0)$, y a través de ella se llega a la forma normal de Hopf asociada al sistema \eqref{c.1}. Los vectores propios de $J(0,0,0)$ están dados por
\begin{eqnarray}
(u_1,v_1,1), \quad (u_2,v_2,1), \quad (u_3,v_3,1), \nonumber
\end{eqnarray}
donde
\begin{eqnarray}
u_1 &=& \frac{F_z}{\lambda_1 - F_x}, \quad u_2 \quad = \quad \frac{F_z}{\lambda_2 - F_x}, \quad u_3 \quad = \quad \frac{F_z}{\lambda_3 - F_x}, \nonumber \\
v_1 &=& \frac{G_z}{\lambda_1 - G_y}, \quad v_2 \quad = \quad \frac{G_z}{\lambda_2 - G_y}, \quad v_3 \quad = \quad \frac{G_z}{\lambda_3 - G_y}. \nonumber
\end{eqnarray}
Es claro que $\overline{u_2} = u_3$, $\overline{v_2} = v_3$. La matriz cambio de base para la matriz $D$ está dada por
\begin{eqnarray}
Q &=& \begin{pmatrix} 
\text{Re} \hspace{3pt} u_2 & \text{Im} \hspace{3pt} u_2 & u_1 \\
\text{Re} \hspace{3pt} v_2 & \text{Im} \hspace{3pt} v_2 & v_1 \\
1             & 0             & 1   
\end{pmatrix}. \nonumber
\end{eqnarray}
La forma normal de Hopf para \eqref{c.1} es
\begin{eqnarray}
\label{c.2}
\dot x &=& f(x,y,z), \nonumber \\
\dot y &=& g(x,y,z),           \\
\dot z &=& h(x,y,z), \nonumber
\end{eqnarray}
donde
\begin{eqnarray}
\begin{pmatrix} 
f(x,y,z) \\
g(x,y,z) \\
h(x,y,z)  
\end{pmatrix} &=& Q^{-1}\begin{pmatrix} 
F\left(Q\begin{pmatrix} 
x \\
y \\
z  
\end{pmatrix}\right) \\
G\left(Q\begin{pmatrix} 
x \\
y \\
z  
\end{pmatrix}\right) \\
H\left(Q\begin{pmatrix} 
x \\
y \\
z  
\end{pmatrix}\right)  
\end{pmatrix}. \nonumber
\end{eqnarray}
\begin{ej}
Si se considera la selección de parámetros $r_1=1$, $r_2=1.7$, $K_2=2.3$, $c_1=1$, $m=0$, entonces el sistema \eqref{c.2} queda como a continuación se muestra. Los coeficientes lineales (que son los valores propios de la matriz jacobiana de \eqref{c.2} en el origen) aparecen en negritas en cada ecuación diferencial.
\begin{eqnarray}
\dot x &=& -5.43793 x^2 y^2-68.8816 x^2 y z+0.145756 x^3 y+37.1622 x^2 y+468.665 x^2 z^2 \nonumber \\
& & +18.638 x^3 z -421.299 x^2 z+0.115285 x^4-5.88326 x^3+5.29695 x^2 -21.4093 x y^2 z \nonumber \\
& &-0.610502 x y^3+163.15 x y^2+11.3548 x y z^2+2180.13 x y z -684.869 x y +1801.88 x z^3 \nonumber \\
& & -12403.1 x z^2+164.93 x z-785.736 y^2 z^2-34.0897 y^3 z+437.219 y^2 z+0.393711 y^4 \nonumber \\
& & -0.244772 y^3 +149. y^2 -84.6763 y  z^3 +3922.4 y z^2+1302.82 y z \textbf{-687.53 y} \nonumber \\
& & +18125.8 z^4 -3301.85 z^3-13825.1 z^2, \nonumber
\end{eqnarray}
\begin{eqnarray}
\dot y &=& -0.33527 x^2 y^2+20.3026 x^2 y z+0.121061 x^3 y-11.6122 x^2 y-184.208 x^2 z^2-1.08844 x^3 z \nonumber \\
& & -923.669 x^2 z-0.000933574 x^4-5.89212 x^3+150.404 x^2-113.27 x y^2 z-6.3683 x y^3 \nonumber \\
& & +332.407 x y^2+1211.14 x y z^2+3172.82 x y z-1443.88 x y-8379.61 x z^3-22083.9 x z^2 \nonumber \\
& & +913.797 x z\textbf{+687.53 x}-1267.4 y^2 z^2+261.408 y^3 z+1123.45 y^2 z+24.5707 y^4 \nonumber \\
& & -102.819 y^3-86.4572 y^2+24382.3 y z^3-15200.3 y z^2+569.403 y z-113328. z^4 \nonumber \\
& & -38749.4 z^3+6259.31 z^2, \nonumber
\end{eqnarray}
\begin{eqnarray}
\dot z &=& 0.0505015 x^2 y^2-2.1628 x^2 y z-0.0141131 x^3 y-0.0634095 x^2 y-17.0504 x^2 z^2 \nonumber \\
& & -0.240711 x^3 z+3.16536 x^2 z-0.00120377 x^4+0.0536329 x^3-11.2407 x^2+10.6344 x y^2 z \nonumber \\
& & +0.610502 x y^3-5.96823 x y^2-153.839 x y z^2-21.6432 x y z+34.3065 x y-953.383 x z^3 \nonumber \\
& & -1.57877 x z^2-647.665 x z+780.349 y^2 z^2+34.0897 y^3 z-280.038 y^2 z-0.393711 y^4 \nonumber \\
& & +0.244772 y^3+13.5692 y^2+13.3686 y z^3-1801.01 y z^2+202.343 y z-17710.6 z^4 \nonumber \\
& & -8690.48 z^3+1365.39 z^2\textbf{-52.2624 z}. \nonumber 
\end{eqnarray}
\end{ej}
\section*{Cálculo del primer coeficiente de Lyapunov $\ell_1$}
De acuerdo con \cite[Sección 4.A, II.]{Marsden1978}, para encontrar el valor del primer coeficiente de Lyapunov $\ell_1$ se hacen los siguientes cálculos.

\begin{eqnarray}
\ell_1 &=& \frac{3 \pi}{4 \omega}(A_{f_{xxx}}(0,0,0) + A_{f_{xyy}}(0,0,0) + A_{g_{xxy}}(0,0,0) + A_{g_{yyy}}(0,0,0) \nonumber \\
& & - f_{xx}(0,0,0)\times f_{xy}(0,0,0) + g_{yy}(0,0,0) \times g_{xy}(0,0,0) + g_{xx}(0,0,0) \times g_{xy}(0,0,0) \nonumber \\
& & - f_{yy}(0,0,0) \times f_{xy}(0,0,0) + f_{xx}(0,0,0) \times g_{xx}(0,0,0) - f_{yy}(0,0,0) \times g_{yy}(0,0,0)), \nonumber 
\end{eqnarray}
donde
\begin{eqnarray}
A_{f_{xxx}}(0,0,0) &=& f_{xxx}(0,0,0) + 3 (f_{xz}(0,0,0)\times \varphi_{xx}(0,0)), \nonumber \\
A_{f_{xyy}}(0,0,0) &=& f_{xyy}(0,0,0) + 2 (f_{yz}(0,0,0)\times \varphi_{xy}(0,0,0)) + f_{xz}(0,0,0)\times \varphi_{yy}(0,0,0), \nonumber \\
A_{g_{xxy}}(0,0,0) &=& g_{xxy}(0,0,0) + 2 (g_{xz}(0,0,0)\times \varphi_{xy}(0,0,0)) + g_{yz}(0,0,0)\times \varphi_{xx}(0,0,0), \nonumber \\
A_{g_{yyy}}(0,0,0) &=& g_{yyy}(0,0,0) + 3 (g_{yz}(0,0,0)\times \varphi_{yy}(0,0,0)), \nonumber
\end{eqnarray}
donde
\begin{eqnarray}
\begin{pmatrix} 
\varphi_{xx}(0,0,0) \\
\varphi_{xy}(0,0,0) \\
\varphi_{yy}(0,0,0)  
\end{pmatrix}  &=& \nonumber 
\end{eqnarray}
\begin{eqnarray}
& & \Delta^{-1}
\begin{pmatrix} 
2 \omega^2 + (h_z(0,0,0))^2 & -2 \omega h_z(0,0,0) & 2 \omega^2 \\
\omega h_z(0,0,0)           & (h_z(0,0,0))^2       & - \omega h_z(0,0,0) \nonumber \\
2 \omega^2                  & a \omega h(0,0,0)    & 2 \omega^2 + (h_z(0,0,0))^2
\end{pmatrix}
\begin{pmatrix} 
-h_{xx}(0,0,0) \\
-h_{xy}(0,0,0) \\
-h_{yy}(0,0,0)  
\end{pmatrix}, \nonumber
\end{eqnarray}
y donde
\begin{eqnarray}
\Delta &=& h_z(0,0,0)((h_z(0,0,0))^2 + 4\omega^2). \nonumber
\end{eqnarray}
De acuerdo con \cite[Teorema 3.3]{Kuznetsov1998}, la condición de genericidad en una bifurcación de Hopf se cumple si y solo si el primer coeficiente de Lyapunov es distinto de cero ($\ell_1 \neq 0$). Además, la condición de transversalidad consiste en que parte real de la derivada de los valores propios complejos de la matriz jacobiana de \eqref{c.2} sea distinta de cero ($(\text{Re}(\lambda_{2,3})'\neq 0$).

Las siguientes tablas muestran el valor del primer coeficiente de Lyapunov $\ell_1$, así como el valor de la derivada de la parte real de los valores propios complejos $(\text{Re}(\lambda_{2,3}))'$ en algunos casos concretos.
\newpage
Valor de $\ell_1$ para $c_1 = 0.5$, $m=0$

\begin{tabular}{|c||c|c|c|c|c|}
\hline
$r_1 = 1.0$ & $K_2 = 3.0$ & $K_2 = 3.5$ & $K_2 = 4.0$ & $K_2 = 4.5$ & $K_2 = 5.0$ \\
\hline \hline
$r_2 = 1.7$ & 0.691488  &  1.83471 & 2.43701  & 2.91227  & 3.29906   \\
\hline
$r_2 = 2.2$ & -0.703596 & 0.384492 & 0.933515 & 1.3578   & 1.6972    \\
\hline
$r_2 = 2.7$ & -1.38538 & -0.56777 &-0.0419655 & 0.359543 & 0.677529  \\
\hline
$r_2 = 3.2$ & -1.83249 & -1.22852 &-0.716643  &-0.328772 &-0.0235696 \\
\hline
$r_2 = 3.7$ & -2.14332 & -1.70764 &-1.20569   &-0.827395 &-0.531077  \\
\hline \hline
$r_1 = 1.5$ &          &          &          &          &          \\
\hline \hline
$r_2 = 1.7$ & 2.2782   & 4.83453  & 5.7111   & 6.43799  & 7.05387  \\
\hline
$r_2 = 2.2$ & 1.08517  & 2.79753  & 3.45744  & 3.98634  & 4.42229 \\
\hline
$r_2 = 2.7$ & 0.0670446& 1.48879  & 2.07546  & 2.53596  & 2.9091  \\
\hline
$r_2 = 3.2$ & -0.624407& 0.542781 & 1.09631  & 1.52498  & 1.86849 \\
\hline
$r_2 = 3.7$ & -1.12166 & -0.169053& 0.365815 & 0.776237 & 1.10261 \\ 
\hline \hline
$r_1 = 2.0$ &          &          &          &          &         \\
\hline \hline
$r_2 = 1.7$ & 5.8122   & 8.70986  & 10.3108  & 11.7052  & 12.9354 \\ 
\hline
$r_2 = 2.2$ & 2.97409  & 5.12144  & 6.03991  & 6.80496  & 7.45554 \\
\hline
$r_2 = 2.7$ & 1.50701  & 3.38992  & 4.09824  & 4.67166  & 5.14817 \\
\hline
$r_2 = 3.2$ & 0.581613 & 2.21742  & 2.83988  & 3.33404  & 3.73819 \\
\hline
$r_2 = 3.7$ & -0.0810414& 1.32795 & 1.90809  & 2.36237  & 2.72975 \\
\hline \hline
$r_1 = 2.5$ &          &          &          &          &          \\
\hline \hline
$r_2 = 1.7$ & 12.1096  & 15.8126  & 19.1458  & 22.2218  & 25.0772 \\
\hline
$r_2 = 2.2$ & 5.52569  & 8.05191  & 9.51236  & 10.7757  & 11.8839 \\ 
\hline
$r_2 = 2.7$ & 3.09316  & 5.30608  & 6.25287  & 7.04374  & 7.7178  \\ 
\hline
$r_2 = 3.2$ & 1.80565  & 3.79625  & 4.54451  & 5.15445  & 5.66416 \\
\hline
$r_2 = 3.7$ & 0.955262 & 2.73536  & 3.39081  & 3.9156   & 4.34779 \\
\hline \hline
$r_1 = 3.0$ &          &          &          &          &          \\
\hline \hline
$r_2 = 1.7$ & 22.9303  & 29.3943  & 36.9499  & 44.5633  & 52.2538 \\
\hline
$r_2 = 2.2$ & 9.51412  & 12.6763  & 15.2026  & 17.4805  & 19.5514 \\
\hline
$r_2 = 2.7$ & 5.15878  & 7.66484  & 9.04502  & 10.2339  & 11.273  \\ 
\hline
$r_2 = 3.2$ & 3.18015  & 5.43505  & 6.40223  & 7.21168  & 7.90266 \\
\hline
$r_2 = 3.7$ & 2.03058  & 4.09468  & 4.87594  & 5.51601  & 6.05307  \\
\hline
\end{tabular}
\newpage
Valor de $\ell_1$ para $c_1 = 0.75$, $m=0$

\begin{tabular}{|c||c|c|c|c|c|}
\hline
$r_1 = 1.0$ &$K_2 = 3.0$ & $K_2 = 3.5$ & $K_2 = 4.0$ & $K_2 = 4.5$ & $K_2 = 5.0$ \\
\hline \hline
$r_2 = 1.7$ & 0.511728& 2.39184 & 3.10442 & 3.69675 & 4.20016 \\ 
\hline
$r_2 = 2.2$ & -0.361737& 0.776027& 1.36992& 1.84541 & 2.23727 \\
\hline
$r_2 = 2.7$ & -1.12424& -0.263593& 0.287137& 0.71881& 1.06835 \\ 
\hline
$r_2 = 3.2$ & -1.62814& -0.985153& -0.456224& -0.0471141& 0.280496 \\ 
\hline
$r_2 = 3.7$ & -1.98016& -1.50919&-0.993928&-0.598977&-0.285103 \\ 
\hline \hline
$r_1 = 1.5$ &         &         &         &         &         \\
\hline \hline
$r_2 = 1.7$ & 3.70289 & 6.50949 & 7.9919  & 9.34702 & 10.5948 \\
\hline
$r_2 = 2.2$ & 1.72171 & 3.54166 & 4.39771 & 5.12964 & 5.76625 \\
\hline
$r_2 = 2.7$ & 0.516752& 1.99848 & 2.6744  & 3.23093 & 3.70024 \\
\hline
$r_2 = 3.2$ & -0.267761& 0.949439& 1.55251& 2.03719 & 2.43786 \\ 
\hline
$r_2 = 3.7$ & -0.827845& 0.171382& 0.738031& 1.186  & 1.55134 \\
\hline \hline
$r_1 = 2.0$ &         &         &         &         &         \\
\hline \hline
$r_2 = 1.7$ & 10.1951 & 15.1418 & 19.8053 & 24.7975 & 30.1609 \\
\hline
$r_2 = 2.2$ & 4.42867 & 7.01116 & 8.62988 & 10.1248 & 11.5138 \\
\hline
$r_2 = 2.7$ & 2.28869 & 4.31185 & 5.29654 & 6.15416 & 6.91163 \\
\hline
$r_2 = 3.2$ & 1.1291  & 2.8373  & 3.59875 & 4.23864 & 4.78736 \\ 
\hline
$r_2 = 3.7$ & 0.353928& 1.81787 & 2.47883 & 3.0207  & 3.47605 \\
\hline \hline
$r_1 = 2.5$ &         &         &         &         &         \\
\hline \hline
$r_2 = 1.7$ & 25.4393 & 41.4568 & 67.3916 & 115.017 & 231.269 \\
\hline
$r_2 = 2.2$ & 9.14571 & 13.4093 & 17.3023 & 21.3451 & 25.5522 \\
\hline
$r_2 = 2.7$ & 4.66261 & 7.34622 & 9.05996 & 10.6532 & 12.1425  \\
\hline
$r_2 = 3.2$ & 2.71365 & 4.87906 & 5.97481 & 6.94178 & 7.80533 \\ 
\hline
$r_2 = 3.7$ & 1.59523 & 3.46411 & 4.30867 & 5.02941 & 5.65532 \\
\hline \hline
$r_1 = 3.0$ &         &         &         &         &         \\
\hline \hline
$r_2 = 1.7$ & 69.7544 & 225.465 & -595.595& -173.163& -114.766 \\
\hline
$r_2 = 2.2$ & 18.2596 & 27.8921 & 40.3836 & 57.2903 & 81.4792 \\
\hline
$r_2 = 2.7$ & 8.4474  & 12.4386 & 15.9283 & 19.4909 & 23.134 \\
\hline
$r_2 = 3.2$ & 4.83213 & 7.58595 & 9.36954 & 11.0355 & 12.5996 \\
\hline
$r_2 = 3.7$ & 3.04697 & 5.31945 & 6.51032 & 7.57166 & 8.52752 \\
\hline
\end{tabular}
\newpage
Valor de $\ell_1$ para $c_1 = 1$, $m=0$

\begin{tabular}{|c||c|c|c|c|c|}
\hline
$r_1 = 1.0$ &$K_2 = 3.0$ & $K_2 = 3.5$ & $K_2 = 4.0$ & $K_2 = 4.5$ & $K_2 = 5.0$ \\
\hline \hline
$r_2 = 1.7$ & 1.05229    & 3.0864   & 4.0132   & 4.83709  & 5.57824  \\
\hline
$r_2 = 2.2$ & -0.0246356 & 1.19568  & 1.86901  & 2.43196  & 2.91306  \\
\hline
$r_2 = 2.7$ & -0.870904  & 0.047846 & 0.639953 & 1.11873  & 1.51673  \\
\hline
$r_2 = 3.2$ & -1.42963   & -0.739451& -0.184098& 0.255843 & 0.615366  \\
\hline
$r_2 = 3.7$ & -1.82094   & -1.3097  & -0.775157& -0.357456& -0.0199994\\
\hline \hline
$r_1 = 1.5$ & & & & & \\
\hline 
$r_2 = 1.7$ & 5.79276  & 9.86321 & 13.2709  & 17.0631  & 21.3072 \\
\hline 
$r_2 = 2.2$ & 2.48499  & 4.61649 & 5.8931   & 7.08746  & 8.211 \\
\hline
$r_2 = 2.7$ & 0.986855 & 2.60572 & 3.45053  & 4.1896   & 4.84551 \\
\hline
$r_2 = 3.2$ & 0.0855952& 1.3903  & 2.08189  & 2.6636   & 3.16321 \\
\hline
$r_2 = 3.7$ & -0.541735& 0.525146& 1.14617  & 1.65511  & 2.08292 \\
\hline \hline
$r_1 =2.0$ & & & & & \\
\hline
$r_2 = 1.7$ & 20.7782 & 43.1299  & 110.01   & -6708.76 & -150.081 \\ 
\hline  
$r_2 = 2.2$ & 6.75551 & 10.9709  & 14.9977  & 19.651   & 25.0854 \\
\hline
$r_2 = 2.7$ & 3.29151 & 5.7787   & 7.40905  & 8.99409  & 10.5389 \\
\hline
$r_2 = 3.2$ & 1.73917 & 3.65431  & 4.6951   & 5.63805  & 6.50029 \\
\hline
$r_2 = 3.7$ & 0.802485& 2.39044  & 3.20268  & 3.90823  & 4.53065 \\ \hline \hline
$r_1 = 2.5$ & & & & & \\
\hline
$r_2 = 1.7$ & 174.498 & -114.506 & -57.5672 & -43.2723& -36.7703 \\
\hline
$r_2 = 2.2$ & 17.3918 & 33.0661  & 66.5633  & 197.4   & -513.902 \\
\hline
$r_2 = 2.7$ & 7.21365 & 11.7496  & 16.2473  & 21.5868 & 28.0232 \\
\hline
$r_2 = 3.2$ & 3.93499 & 6.71896  & 8.68701  & 10.6604 & 12.6418 \\
\hline
$r_2 = 3.7$ & 2.34882 & 4.50619  & 5.75308  & 6.91535 & 8.00513 \\
\hline \hline
$r_1 = 3.0$ & & & & & \\
\hline
$r_2 = 1.7$ & -80.6253& -45.1195& -35.476  & -31.0165 & -28.4456 \\
\hline
$r_2 = 2.2$ & 59.357  & -3574.28 & -95.7648 & -58.2737 & -45.6068 \\ 
\hline
$r_2 = 2.7$ & 15.5567 & 28.4477  & 52.1839  & 114.76   & 734.458 \\
\hline
$r_2 = 3.2$ & 7.54939 & 12.3265  & 17.1924  & 23.0873  & 30.3691 \\
\hline
$r_2 = 3.7$ & 4.46254 & 7.50202  & 9.78462  & 12.1331  & 14.5521 \\
\hline
\end{tabular}
\newpage
Valor de $(\text{Re}(\lambda_{2,3}))'$ para $c_1 = 0.5$, $m=0$

\begin{tabular}{|c||c|c|c|c|c|}
\hline
$r_1 = 1.0$ & $K_2 = 3.0$ & $K_2 = 3.5$ & $K_2 = 4.0$ & $K_2 = 4.5$ & $K_2 = 5.0$ \\
\hline \hline
$r_2 = 1.7$ & -608.044 & -710.967 & -813.99 & -917.084 & -1020.23   \\
\hline
$r_2 = 2.2$ & -604.657 & -706.444 & -808.301& -910.208 & -1012.15    \\
\hline
$r_2 = 2.7$ & -602.962 & -704.144 & -805.376& -906.644 & -1007.94  \\
\hline
$r_2 = 3.2$ & -602.023 & -702.854 & -803.722& -904.618 & -1005.53 \\
\hline
$r_2 = 3.7$ & -601.459 & -702.071 & -802.712& -903.374 & -1004.05  \\
\hline \hline
$r_1 = 1.5$ &          &          &          &          &          \\
\hline \hline
$r_2 = 1.7$ & -617.135 & -722.791 & -828.593 & -934.496 & -1040.47  \\
\hline
$r_2 = 2.2$ & -610.793 & -714.584 & -818.493 & -922.486 & -1026.54 \\
\hline
$r_2 = 2.7$ & -607.149 & -709.779 & -812.502 & -915.292 & -1018.13  \\
\hline
$r_2 = 3.2$ & -604.979 & -706.878 & -808.85  & -910.875 & -1012.94 \\
\hline
$r_2 = 3.7$ & -603.622 & -705.042 & -806.522 & -908.043 & -1009.6 \\ 
\hline \hline
$r_1 = 2.0$ &          &          &          &          &         \\
\hline \hline
$r_2 = 1.7$ & -625.781 & -733.711 & -841.796 & -949.987 & -1058.26 \\ 
\hline
$r_2 = 2.2$ & -617.983 & -723.875 & -829.916 & -936.059 & -1042.28 \\
\hline
$r_2 = 2.7$ & -612.622 & -716.97  & -821.446 & -926.012 & -1030.64 \\
\hline
$r_2 = 3.2$ & -609.095 & -712.355 & -815.722 & -919.166 & -1022.67 \\
\hline
$r_2 = 3.7$ & -606.75  & -709.249 & -811.837 & -914.49  & -1017.19 \\
\hline \hline
$r_1 = 2.5$ &          &          &          &          &          \\
\hline \hline
$r_2 = 1.7$ & -631.984 & -741.314 & -850.788 & -960.36  & -1070. \\
\hline
$r_2 = 2.2$ & -624.721 & -732.389 & -840.212 & -948.143 & -1056.15 \\ 
\hline
$r_2 = 2.7$ & -618.513 & -724.552 & -830.74  & -937.032 & -1043.4  \\ 
\hline
$r_2 = 3.2$ & -613.906 & -718.635 & -823.498 & -928.454 & -1033.48 \\
\hline
$r_2 = 3.7$ & -610.606 & -714.339 & -818.189 & -922.123 & -1026.12 \\
\hline \hline
$r_1 = 3.0$ &          &          &          &          &          \\
\hline \hline
$r_2 = 1.7$ & -635.347 & -745.252 & -855.282 & -965.395 & -1075.57 \\
\hline
$r_2 = 2.2$ & -630.063 & -738.988 & -848.063 & -957.239 & -1066.49 \\
\hline
$r_2 = 2.7$ & -624.03  & -731.524 & -839.174 & -946.932 & -1054.77  \\ 
\hline
$r_2 = 3.2$ & -618.875 & -725.014 & -831.302 & -937.695 & -1044.16 \\
\hline
$r_2 = 3.7$ & -614.85  & -719.854 & -824.996 & -930.235 & -1035.54  \\
\hline
\end{tabular}
\newpage
Valor de $(\text{Re}(\lambda_{2,3}))'$ para $c_1 = 0.75$, $m=0$

\begin{tabular}{|c||c|c|c|c|c|}
\hline
$r_1 = 1.0$ &$K_2 = 3.0$ & $K_2 = 3.5$ & $K_2 = 4.0$ & $K_2 = 4.5$ & $K_2 = 5.0$ \\
\hline \hline
$r_2 = 1.7$ & -610.314 & -714.063& -817.937& -921.898& -1025.92 \\ 
\hline
$r_2 = 2.2$ & -606.317 & -708.741& -811.258& -913.841& -1016.47 \\
\hline
$r_2 = 2.7$ & -604.145 & -705.798& -807.521& -909.293& -1011.1 \\ 
\hline
$r_2 = 3.2$ & -602.885 & -704.069& -805.307& -906.582& -1007.89 \\ 
\hline
$r_2 = 3.7$ & -602.105 & -702.989& -803.914& -904.869& -1005.85 \\ 
\hline \hline
$r_1 = 1.5$ &         &         &         &         &         \\
\hline \hline
$r_2 = 1.7$ & -619.12 & -725.454& -831.942& -938.535& -1045.2 \\
\hline
$r_2 = 2.2$ & -613.257& -717.92 & -822.722& -927.622& -1032.59 \\
\hline
$r_2 = 2.7$ & -609.298& -712.721& -816.261& -919.884& -1023.57 \\
\hline
$r_2 = 3.2$ & -606.716& -709.277& -811.935& -914.662& -1017.44 \\ 
\hline
$r_2 = 3.7$ & -605.005& -706.969& -809.012& -911.111& -1013.25 \\
\hline \hline
$r_1 = 2.0$ &         &         &         &         &         \\
\hline \hline
$r_2 = 1.7$ & -625.019& -732.775& -840.675& -948.671& -1056.74 \\
\hline
$r_2 = 2.2$ & -619.808& -726.324& -832.994& -939.77 & -1046.62 \\
\hline
$r_2 = 2.7$ & -615.079& -720.282& -825.632& -931.083& -1036.61 \\
\hline
$r_2 = 3.2$ & -611.469& -715.583& -819.828& -924.165& -1028.57 \\ 
\hline
$r_2 = 3.7$ & -608.838& -712.11 & -815.496& -918.963& -1022.49 \\
\hline \hline
$r_1 = 2.5$ &         &         &         &         &         \\
\hline \hline
$r_2 = 1.7$ & -627.315& -735.384& -843.574& -951.845& -1060.17 \\
\hline
$r_2 = 2.2$ & -624.438& -732.074& -839.855& -947.736& -1055.69 \\
\hline
$r_2 = 2.7$ & -620.227& -726.852& -833.632& -940.517& -1047.48  \\
\hline
$r_2 = 3.2$ & -616.293& -721.847& -827.55 & -933.356& -1039.24 \\ 
\hline
$r_2 = 3.7$ & -613.065& -717.669& -822.413& -927.252& -1032.16 \\
\hline \hline
$r_1 = 3.0$ &         &         &         &         &         \\
\hline \hline
$r_2 = 1.7$ & -627    & -734.695& -842.489& -950.348& -1058.25 \\
\hline
$r_2 = 2.2$ & -626.848& -734.906& -843.094& -951.368& -1059.7 \\
\hline
$r_2 = 2.7$ & -624.036& -731.584& -839.28 & -947.076& -1054.94 \\
\hline
$r_2 = 3.2$ & -620.508& -727.206& -834.058& -941.016& -1048.05 \\
\hline
$r_2 = 3.7$ & -617.153& -722.949& -828.896& -934.948& -1041.07 \\
\hline
\end{tabular}
\newpage
Valor de $(\text{Re}(\lambda_{2,3}))'$ para $c_1 = 1$, $m=0$

\begin{tabular}{|c||c|c|c|c|c|}
\hline
$r_1 = 1.0$ &$K_2 = 3.0$ & $K_2 = 3.5$ & $K_2 = 4.0$ & $K_2 = 4.5$ & $K_2 = 5.0$ \\
\hline \hline
$r_2 = 1.7$ & -611.843   & -716.153 & -820.601 & -925.146 & -1029.76  \\
\hline
$r_2 = 2.2$ & -607.641   & -710.574 & -813.617 & -916.737 & -1019.91  \\
\hline
$r_2 = 2.7$ & -605.164   & -707.224 & -809.369 & -911.574 & -1013.82  \\
\hline
$r_2 = 3.2$ & -603.66    & -705.162 & -806.731 & -908.347 & -1010     \\
\hline
$r_2 = 3.7$ & -602.703   & -703.837 & -805.024 & -906.248 & -1007.5   \\
\hline \hline
$r_1 = 1.5$ & & & & & \\
\hline 
$r_2 = 1.7$ & -619.442 & -725.91 & -832.53  & -939.252 & -1046.05 \\
\hline 
$r_2 = 2.2$ & -614.648 & -719.809& -825.12  & -930.533 & -1036.02 \\
\hline
$r_2 = 2.7$ & -610.818 & -714.803& -818.922 & -923.132 & -1027.41 \\
\hline
$r_2 = 3.2$ & -608.081 & -711.164& -814.361 & -917.637 & -1020.97 \\
\hline
$r_2 = 3.7$ & -606.162 & -708.579& -811.093 & -913.674 & -1016.31 \\
\hline \hline
$r_1 =2.0$ & & & & & \\
\hline
$r_2 = 1.7$ & -622.722& -729.824 & -837.055 & -944.374 & -1051.75 \\ 
\hline  
$r_2 = 2.2$ & -619.924& -726.509 & -833.244 & -940.081 & -1046.99 \\
\hline
$r_2 = 2.7$ & -616.257& -721.881 & -827.657 & -933.538 & -1039.49 \\
\hline
$r_2 = 3.2$ & -612.974& -717.633 & -822.437 & -927.339 & -1032.31 \\
\hline
$r_2 = 3.7$ & -610.344& -714.176 & -818.139 & -922.192 & -1026.31 \\ \hline \hline
$r_1 = 2.5$ & & & & & \\
\hline
$r_2 = 1.7$ & -622.541& -729.279 & -836.121 & -943.032& -1049.99 \\
\hline
$r_2 = 2.2$ & -622.518& -729.607 & -836.83  & -944.142& -1051.52 \\
\hline
$r_2 = 2.7$ & -620.208& -726.86  & -833.66  & -940.562& -1047.54 \\
\hline
$r_2 = 3.2$ & -617.273& -723.177 & -829.235 & -935.397& -1041.63 \\
\hline
$r_2 = 3.7$ & -614.472& -719.582 & -824.84  & -930.201& -1035.64 \\
\hline \hline
$r_1 = 3.0$ & & & & & \\
\hline
$r_2 = 1.7$ & -620.612& -726.572 & -832.615 & -938.711 & -1044.84 \\
\hline
$r_2 = 2.2$ & -622.892& -729.836 & -836.893 & -944.025 & -1051.21 \\ 
\hline
$r_2 = 2.7$ & -622.354& -729.425 & -836.633 & -943.93  & -1051.29 \\
\hline
$r_2 = 3.2$ & -620.394& -727.089 & -833.931 & -940.875 & -1047.89 \\
\hline
$r_2 = 3.7$ & -617.962& -724.052 & -830.295 & -936.642 & -1043.06 \\
\hline
\end{tabular}
\chapter{Teoría del Promedio}
\label{apendiced}
En la Sección \ref{seccion4.5} se aplican métodos de teoría del promedio que fueron usados en \cite{Llibre2015Zero-HopfSystem} para hallar y analizar órbitas periódicas en un sistema de ecuaciones diferenciales acopladas. Específicamente, los Teoremas \ref{teo4.22} y \ref{teo4.23} tienen su base en el resultado que a continuación se presenta.

Considérese
\begin{eqnarray}
\label{d.1}
\dot x &=& \varepsilon F(t,x) + \varepsilon^2 G(t,x),
\end{eqnarray}
donde $x\in D$, con $D$ un subconjunto abierto de $\mathbb R^n$, $t\geq 0$. Adicionalmente supóngase que las funciones $F$, $G$ son $T-$periódicas con respecto a $t$. Por otra parte considérese la ecuación promediada en $D$ dada por
\begin{eqnarray}
\label{d.2}
\dot y &=& \varepsilon f(y), 
\end{eqnarray}
donde
\begin{eqnarray}
f(y) &=& \frac{1}{T}\int_0^TF(t,y)dt. \nonumber
\end{eqnarray}
Bajo ciertas condiciones, los puntos de equilibrio de la ecuación promediada corresponden con las soluciones $T-$periódicas de la ecuación \eqref{d.1}.
\begin{teo}
\label{teod.1}
Considérese la ecuación \eqref{d.1} y supóngase que:
\begin{enumerate}
    \item[a.] las funciones $f$, $g$, $\displaystyle \frac{\partial f}{\partial x}$, $\displaystyle \frac{\partial^2 f}{\partial x^2}$, $\displaystyle \frac{\partial g}{\partial x}$ están bien definidas, son continuas y acotadas por una constante independiente de $\varepsilon$ en $[0,\infty)\times D$, $0\leq \varepsilon\leq \varepsilon_0$;
    \item[b.] $F$, $G$ son $T-$periódicas en $t$ (donde $T$ es independiente de $\varepsilon$);
\end{enumerate}
Entonces se cumplen las siguientes afirmaciones:
\begin{enumerate}
    \item Si $p$ es un punto crítico de la ecuación promediada \eqref{d.2} y además
    \begin{eqnarray}
    \left\lvert \frac{\partial f}{\partial y}\right\rvert_{y=p} &\neq& 0, \nonumber
    \end{eqnarray}
    entonces existe una solución $T-$periódica $\varphi$ para la ecuación \eqref{d.1} tal que
    \begin{eqnarray}
    \lim_{\varepsilon\to 0}\varphi(t,\varepsilon) &=& p. \nonumber
    \end{eqnarray}
    \item La estabilidad o inestabilidad de la solución periódica $\varphi$ coincide con la estabilidad de $p$ en \eqref{d.2}.
\end{enumerate}

\end{teo}
 Para ver una demostración, el lector interesado puede consultar \cite[Teoremas 11.5 y 11.6]{Verhulst1996NonlinearSystems}.
 
 El siguiente resultado se encuentra en \cite{Buica2004}, y permite debilitar las hipótesis del Teorema \ref{teod.1} haciendo uso de la teoría del grado de Brouwer. 
\begin{teo} 
\label{teod.2}
Considérese el sistema de ecuaciones diferenciales
\begin{eqnarray}
\label{d.3}
\dot x &=& \varepsilon F(t,x) + \varepsilon^2G(t,x,\varepsilon),
\end{eqnarray}
donde $F\colon \mathbb R\times D \longrightarrow \mathbb R^n$, $G\colon \mathbb R\times D\times (-\varepsilon,\varepsilon)\longrightarrow \mathbb R^n$ son funciones continuas y $T-$periódicas en la primera variable, y $D$ es un subconjunto abierto de $\mathbb R^n$. Si $f\colon D\longrightarrow \mathbb R^n$ está dada por
\begin{eqnarray}
f(z) &=& \int_0^T F(s,z)ds, \nonumber
\end{eqnarray}
y adicionalmente se cumple que
 \begin{enumerate}
     \item $F$ y $G$ son funciones localmente Lipschitz con respecto de $x$.
     \item Para toda $a\in D$ con $f(a) = 0$, existe $V$ una vecindad de $a$ de tal forma que $f(z) \neq 0$ para toda $z\in \overline V\setminus \{a\}$ y $d_B(f,V,0)\neq 0$. 
 \end{enumerate}
 Entonces para $|\varepsilon| > 0$ suficientemente pequeño, existe una solución $T-$periódica $\varphi(\cdot,\varepsilon)$ para el sistema \eqref{d.3} tal que $\varphi(\cdot,\varepsilon)\to a$ si $\varepsilon \to 0$.
\end{teo}
Dado un conjunto abierto y acotado $V\subset \mathbb R^n$ tal que $V\subset D$, y tal que $0\notin f(\partial V,\varepsilon)$ para algún valor de $\varepsilon$, se denota con $d_B(f(\cdot,\varepsilon),V,0)$ el grado de Brouwer de la función $f(\cdot,\varepsilon)$ con respecto al conjunto $V$ y al punto $0$. Lo anterior se define en \cite{Browder1983}.
\chapter[Cálculo de la Variedad Central]{Cálculo de la Variedad Central\index{variedad central} de un Punto No Hiperbólico}
\label{apendicee}
Una característica que exhibe el modelo \eqref{5.1} es que uno de sus puntos de equilibrio no es hiperbólico (tiene un valor propio con parte real nula).

El sistema \eqref{5.1} está dado por

\begin{eqnarray}
\dot x &=& x\left(1 - x - y_1 - \frac{10 y_2}{1 + 0.01 x}\right), \nonumber \\
\dot y_1 &=& y_1\left(1 - y_1 - 1.5x - \frac{y_2}{1 + 0.02 y_1}\right), \nonumber \\
\dot y_2 &=& y_2\left(\frac{5 x}{1 + 0.01 x} + \frac{0.5 y_1}{1 + 0.02 y_1} - 1 - 0.1 y_2\right). \nonumber
\end{eqnarray}
Los cinco puntos de equilibrio de \eqref{5.1} son
\begin{eqnarray}
E_0 \quad = \quad (0,0,0)\text{, } \quad E_1 &=& (1,0,0)\text{, }\quad E_2 \quad = \quad (0,1,0)\text{, } \nonumber \\ 
E_5 \quad = \quad (0.202418,0,0.0800811)\text{, }&& E^* \quad = \quad (0.11969,0.813857,0.00666116). \nonumber
\end{eqnarray}
El espectro de la matriz jacobiana del sistema \eqref{5.1} en  $E_2$ está dado por
\begin{eqnarray}
\sigma(E_2) &=& \left\{-1,-0.509804,0\right\}, \nonumber
\end{eqnarray}
de donde se sigue que $E_2$ es un punto de equilibrio no hiperbólico en \eqref{5.1}. Se realiza el cambio de variable $y_1\to y_1 + 1$ para trasladar el punto $E_2$ al origen $(0,0,0)$ y así obtener el sistema
\begin{eqnarray}
\dot x &=& x\left(- x - y_1 - \frac{10 y_2}{1 + 0.01 x}\right), \nonumber \\
\dot y_1 &=& (1 + y_1)\left(-1.5 x - y_1 - \frac{y_2}{1 + 0.02(1 + y_1)}\right), \nonumber \\
y_2 &=& \left(\frac{5 x}{1 + 0.01 x} + \frac{0.5(1 + y_1)}{1 + 0.02(1 + y_1)} - 1 - 0.1y_2\right). \nonumber
\end{eqnarray}
A continuación se polinomiza el sistema para obtener
\begin{eqnarray}
\dot x &=& - 1.02 x^2 - 1.02 x y_1 - 10.2 x y_2 + O_3(x,y_1,y_2), \nonumber \\
\dot y_1 &=& - 1.53 x - 1.02 y_1 - y_2 - 0.0153 x^2 - 1.04 y_1^2 - 1.5702 x y_1 - 0.01 x y_2 - y_1 y_2 + O_3(x,y_1,y_2), \nonumber \\
\dot y_2 &=& - 0.52 y_2 + 5.0948 x y_2 + 0.48 y_1 y_2 + O_3(x,y_1,y_2), \nonumber  
\end{eqnarray}
y se procede a diagonalizar el sistema (si se desea conocer más al respecto, se puede consultar un procedimiento para la diagonalización de un sistema en el Apéndice \ref{apendicec}).
\begin{eqnarray}
\label{e.1}
\dot x   &=& 0.51 x^2 - 1.02 x y_1 - 8.16 x y_2 + O_3(x,y_1,y_2), \nonumber \\
\dot y_1 &=& - 1.02 y_1 + 0.765 x^2 + 0.0198 x y_1 - 1.04 y_1^2 - 5.1 x y_2 \\
& & + 4.12 y_1 y_2 - 4.284 y_2^2 + O_3(x,y_1,y_2), \nonumber \\
\dot y_2 &=& -0.52 y_2 + 4.3748 x y_2 + 0.48 y_1 y_2 - 1.062 y_2^2 + O_3(x,y_1,y_2). \nonumber
\end{eqnarray}
Ahora se siguen las ideas expuestas en \cite[Secciones 1.3 y 1.4]{Carr1982ApplicationsTheory} para dar una aproximación de la variedad central\index{variedad central} de \eqref{e.1}. Es claro que el sistema \eqref{e.1} se puede reescribir
\begin{eqnarray}
\dot x &=& Ax + f(x,y), \nonumber \\
\dot y &=& By + g(x,y), \nonumber
\end{eqnarray}
donde
\begin{eqnarray}
A &=& 0, \nonumber \\ \nonumber \\
B &=&  \begin{pmatrix}
-1.02 & 0      \\
0     & -0.52 
\end{pmatrix}, \nonumber \\ \nonumber \\
f(x,y) &=& - 1.02xy_1 - 8.16xy_2 + O_3(x,y_1,y_2), \nonumber \\ \nonumber \\
g(x,y) &=& \begin{pmatrix}
0.765x^2 + 0.0198xy_1 - 1.04y_1^2 - 5.1xy_2 + 4.12y_1y_2 - 4.284y_2^2 \\
4.3748xy_2 + 0.48y_1y_2 - 1.062y_2^2
\end{pmatrix} + O_3(x,y_1,y_2). \nonumber
\end{eqnarray}
Supóngase que $y = h(x)$ es la variedad central\index{variedad central} buscada. Entonces se tiene por la Regla de la Cadena
\begin{eqnarray}
\dot y &=& Bh(x) + g(x,h(x))              \nonumber \\
       &=& \left(h(x)\right)'             \nonumber \\
       &=& h'(x)\left[Ax + f(x,y)\right]. \nonumber
\end{eqnarray}
A partir de estas igualdades, se define
\begin{eqnarray}
M_\varphi(x) &: =& \varphi'(x)\left(A x + f(x,\varphi(x)\right) - B\varphi(x) - g(x,\varphi(x)), \nonumber
\end{eqnarray}
y se hace notar que $M_h (x) = 0$ para cualquier $x$. Para aproximar $h$, se asume $\varphi(x) = O_2(x)$, de donde se sigue
\begin{eqnarray}
M_\varphi(x) &=& -\begin{pmatrix}
-1.02 & 0      \\
0     & -0.52 
\end{pmatrix}
\begin{pmatrix}
\varphi_1(x) \\
\varphi_2(x)
\end{pmatrix}
- \begin{pmatrix}
0.765 x^2  \\
 0
\end{pmatrix} 
+ O_3(x). \nonumber
\end{eqnarray}
Se concluye que la aproximación a la variedad central\index{variedad central} $y = h(x)$ del punto de equilibrio $(0,0,0)$ en el sistema \eqref{e.1} está dada por
\begin{eqnarray}
\varphi(x) &=& \begin{pmatrix}
0.75 x^2    \\
0 
\end{pmatrix}. \nonumber
\end{eqnarray}
\begin{teo} El sistema de ecuaciones diferenciales \eqref{5.1} tiene un punto de equilibrio no hiperbólico $E_2 = (0,1,0)$. El espectro del punto $E_2$ es
\begin{eqnarray}
\sigma(E_2) &=& \{-1,-0.509804,0\}, \nonumber
\end{eqnarray}
y la variedad central\index{variedad central} de $E_2$ se puede aproximar con la curva $\varphi\colon \mathbb R\longrightarrow \mathbb R^3$ dada por
\begin{eqnarray}
\varphi(x) &=& (x,1 + 0.75x^2,0). \nonumber
\end{eqnarray}
\end{teo}
\end{appendices}
\bibliographystyle{plain}
\bibliography{bibliografia}
\backmatter
\end{document}